\documentclass[brochure,10pt]{article}
\usepackage{amssymb,amsfonts,amsmath,footnote}
\usepackage{latexsym}
\usepackage{amsmath}
\usepackage{amsthm}
\usepackage{amsfonts}
\usepackage{amssymb}
\usepackage{amsmath,amscd}
\usepackage{graphicx}
\usepackage{graphics}
\usepackage{etex}
\usepackage[matrix,arrow]{xy}
\usepackage{color}
\usepackage{tikz-cd}
\newtheorem{Th}{Theorem}
\newtheorem{Prop}{Proposition}

\newtheorem{Lm}{Lemma}

\newtheorem{Dfi}{Definition}
\newtheorem{Rm}{Remark}
\newtheorem{Con}{Open Problem}
\newcommand{\be}{\begin{equation}}
\newcommand{\ee}{\end{equation}}
\newcommand{\bes}{\begin{equation*}}
\newcommand{\ees}{\end{equation*}}

\newcommand{\R}{\mathbb{R}}
\newcommand{\N}{\mathbb{N}}
\newcommand{\C}{\mathbb{C}}
\newcommand{\Z}{\mathbb{Z}}

\newcommand\res{\mathop{\hbox{\vrule height 7pt width .5pt depth 0pt
\vrule height .5pt width 6pt depth 0pt}}\nolimits}

\newcommand{\reset}{\setcounter{equation}{0}\setcounter{Th}{0}\setcounter{Prop}{0}\setcounter{Co}{0}
\setcounter{Lm}{0}\setcounter{Rm}{0}}
\def\La{\Lambda}
\def\La{\Lambda}
\def\ti{\tilde}
\def\lf{\left}
\def\rg{\right}

\def\al{\alpha}
\def\la{\lambda}

\def\ep{\varepsilon}

\def\ds{\displaystyle}
\def\ov{\overline}
\def\Om{\Omega}
\def\om{\omega}
\def\p{\partial}

\def\r{\mathfrak r}

\def\res{\mathop{\hbox{\vrule height 7pt width .5pt 
depth 0pt\vrule height .5pt width 6pt depth 0pt}}\nolimits}

\begin{document}

\title{Area Variations under Legendrian Constraint}
\author{Tristan Rivi\`ere}

\date{ }
\maketitle

{\bf Abstract : }{\it In any 5-dimensional closed Sasakian manifold we prove that any minmax operation on the area  among Legendrian surfaces is achieved by a continuous conformal Legendrian map from a closed riemann surface $S$ into $N^5$ equipped with an integer multiplicity bounded in $L^\infty$. Moreover this map, equipped with this multiplicity,  satisfies a weak version of the Hamiltonian Minimal Equation.
We conjecture that any solution to this equation is  a smooth branched Legendrian immersion away from isolated Schoen-Wolfson conical singularities with non zero Maslov class.}

\medskip

\noindent{\bf Math. Class. 53D12, 49Q05, 53A10,  58E12, 49Q10}

\section{Introduction}
\subsection{Foreword}
The simplest framework of our investigations in this work could be formulated as follows:  study area variations among maps $u$ from an oriented 2-manifold $\Sigma$ into ${\C}^2$
under the Lagrangian constraint 
\be
\label{1}
u^\ast\om=0\ ,
\ee
 where $\om$ is the standard symplectic form  $\om=dy_1\wedge dy_2+dy_3\wedge dy_4$.  More precisely we ask ourselves how to implement variational operations (minimization under various constraints, minmax) for the 
 standard area Lagrangian
 \[
 A(u):=\int_\Sigma|\p_{x_1}u\wedge\p_{x_2}u|\ dx_1\wedge dx_2
 \]
 among maps satisfying the pointwise constraint\footnote{In case $u$ is a graph $u(x_1,x_2)=(x_1,x_2, v(x_1,x_2))$ the constraint (\ref{1}) is equivalent to the incompressibility condition
 \[
 \mbox{det}(\nabla v)=-1\ .
 \]}
 \[
 du_1\wedge du_2+du_3\wedge du_4=0\ .
 \]
The Lagrangian $A$ is notoriously known in calculus of variations to be delicate to work with (due mostly to the very large invariance group in the domain : the ``gauge group'' of diffeomorphisms of the surface $\Sigma$).
In order to ``break'' this invariance and reduce the gauge group to a finite dimensional one, it is preferable instead to consider the Dirichlet  Energy of the map $u$ for a variable metric $g$ of $\Sigma$
\be
\label{2}
E(u,g):=\frac{1}{2}\int_\Sigma|du|^2_g\ dvol_g\ge A(u)\ .
\ee
with equality if and only if $u$ is conformal with respect to the metric $g$ that is in isothermic or conformal coordinates such that $g=e^{2\la}\ [dx_1^2+dx_2^2]$ there holds
\be
\label{3}
|\p_{x_1}u|^2=|\p_{x_2}u|^2\quad\mbox{ and }\quad \p_{x_1}u\cdot\p_{x_2}u=0\ ,
\ee
where $\cdot$ denotes the standard scalar product in ${\C}^2$. Stationary critical points for fixed $g$ have an holomorphic Hopf differential
\[
H(u):=\p_zu\cdot\p_z u\ dz\otimes dz
\]
while making $g$ vary simultaneously in the Moduli space of constant Gauss curvature metrics is known to imply $H(u)=0$ and equality holds in (\ref{2}). Hence critical points of $E(u,g)$ are critical points of the area.

Assuming $u$ is a conformal critical point of the Dirichlet energy under the Lagrangian constraint, in local conformal coordinates, \underbar{formally} the Euler Lagrange Equation of our constrained variational problem is given by
\be
\label{4}
\lf\{
\begin{array}{l}
\mbox{div}\lf({\mathfrak g}\,\nabla u\rg)=0\\[5mm]
\mbox{div}\lf({\mathfrak g}^{-1}\, \nabla{\mathfrak g}\rg)=0\ ,
\end{array}
\rg.
\ee
where ${\mathfrak g}=e^{i\beta}$ is an $S^1$ valued map\footnote{The function $\beta$ plays the role of a Lagrange multiplier issued from the incompressibility constraint similarly as the role played by the ``pressure'' in the variational formulation of Euler equations.}, which is the Lagrange multiplier associated to the Lagrangian pointwise constraint (\ref{1}). The first equation is a structural equation shared by any conformal parametrization of any Lagrangian surface in ${\C}^2$ while the second one
is the Euler Lagrange equation itself. It is saying that the Lagrange multiplier ${\mathfrak g}$ is an $S^1$ harmonic  map. These equations are obviously calling upon the use of classical elliptic theory for solving regularity or existence questions. The difficulty however is that, while the natural assumption is to consider $u\in W^{1,2}_{loc}(\Sigma,{\C}^2)$ , there is a-priori no information on the function space to which $\mathfrak g$ belongs to in order to start implementing this theory. Hence these questions (regarding existence and regularity) cannot be settled at the strictly ``PDE level'' whose formulation even is problematic. Indeed in order to give a distributional meaning to (\ref{4}) one needs at least an assumption like ${\mathfrak g}$ is in $W^{1,1}$ or $H^{1/2}$. Studying the regularity\footnote{Observe that everywhere discontinuous solutions to the second equation of (\ref{4}) have been constructed in \cite{Alm}.} of solution to the system (\ref{4}) assuming the conformality condition\footnote{A weaker hypothesis would be to assume that the  Hopf differential $H(u)$ is holomorphic, which is equivalent for the critical point to be a stationary critical point that is critical for variation in the domain.  } (\ref{3})  is certainly an interesting elliptic PDE problem by itself but somehow artificial though since these assumptions are not given a-priori by any variational operation we are aiming at doing.
\subsection{Hamiltonian and Lagrangian Stationary Immersions}
A K\"ahler surface is a complex two dimensional manifold $(M^{2},J)$ with a compatible\footnote{A symplectic form is said to be compatible with a complex or with an almost complex structure $J$ if $g(X,Y):=\om(X,JY)$ defines a riemannian metric.} symplectic form $\om$. Lagrangian sub-manifolds are the ones on which the restriction of the symplectic form vanishes. In particular, Lagrangian sub-manifolds are the ones for which the complex structure realizes an isometry between the tangent and the normal spaces at every points.

The study of area variations under Lagrangian constraints has been initiated in a series of work by Yong-Geun Oh  (\cite{Oh2}, \cite{Oh3}). He proposed to minimize the area under the Lagrangian constraint within a Lagrangian homology class with the objective to obtain in this way very special representative of these classes. He called the sub-manifolds critical points of the area within Lagrangian sub-manifolds {\it stationary Lagrangians}. He introduced also a less restrictive class of surfaces\footnote{The class of Hamiltonian stationary surfaces is in general richer and more interesting than the class of Lagrangian stationary surfacesaccording to the author of the present paper.} which are the so called {\it Hamiltonian stationary} sub-manifolds. This later class of surfaces are critical points of the area for local deformations which preserve infinitesimally the Lagrangian condition. These deformations are generated by the {\it Hamiltonian vector-fields} which are of the form
\be
\label{hamilton}
J\,\nabla h\ .
\ee
where $h$ is an arbitrary smooth function. Such surfaces are also called {\it Hamiltonian minimal} or simply {\it $H-$minimal} surfaces. These pioneered works on the subject contain very interesting  and now classical conjectures in the field about possible volume minimizing candidates under hamiltonian deformations within ${\C}^2$ or ${\C}P^2$.

 Y.-G. Oh computed in particular the Euler Lagrange equations associated to this variational problem. In the special case of K\"ahler-Einstein manifolds (satisfying $Ric(g)=\la\, g$), it has been discovered by Dazor \cite{Daz} that the mean curvature of a Lagrangian surface $\Sigma$ is given by
\[
\vec{H}:=J\,\nabla^\Sigma\beta\ ,
\]
where $\beta$ is a locally well defined function on $\Sigma$ and $\nabla^\Sigma$ is the gradient of $\beta$ for the induced metric on $\Sigma$. In the particular case when the manifold is Ricci flat ( which implies that $M^2$ is Calabi Yau), if $\Om$ is a normalized global holomorphic section of the canonical bundle $K:=\wedge^{2,0}M$ the function $\beta$ is given by
\[
e^{i\beta}\iota_\Sigma^\ast\Om=dvol_g
\]
where $\iota_\Sigma$ is the canonical embedding of $\Sigma$ in $M^2$. Hence $\beta $ is a globally well defined function in ${\R}/{2\pi\Z}$ and is called {\it Lagrangian angle}. For general  K\"ahler-Einstein manifolds $\beta$ is still a well defined function in ${\R}/{2\pi\Z}$ in case the induced canonical connection of the restriction to $\Sigma$ of the canonical bundle $\iota_\Sigma^\ast K$ of $M^2$  (which is flat) has no monodromy (see \cite{Wol1}).

In \cite{Oh2} and \cite{Oh3} it is proved that a sub-manifold is {\it hamiltonian stationary} (or $H$-minimal) if the contraction of the mean curvature vector with the K\"ahler form
\[
\om\res\vec{H}=\al_{\vec{H}}
\] 
defines an \underbar{harmonic} one form of the sub-manifold $\Sigma$. This is equivalent to the fact that 
\[
\Delta^\Sigma\beta=0\ ,
\]
which is nothing but the second equation in (\ref{4}).


In the particular case when $\beta$ is a globally univalued function, then $\beta=\beta_0$ and the Euler-Lagrange equation is simply the minimal surface equation $\vec{H}=0$. Such a sub-manifold is called {\it minimal lagrangian}.
If the K\"ahler-Einstein manifolds is Ricci flat  then any {\it minimal lagrangian} sub-manifold is calibrated by $e^{i\beta_0}\,\Om$ and and realizes in this way an absolute minimizer\footnote{The local stability  of arbitrary minimal lagrangian sub-manifolds generalizes to the situation of general K\"ahler-Einstein manifolds of non positive Ricci. This being said we will be interested in the other case
in the present work that is $Ric(g)>0$.}
in its homology class (see \cite{HL}). While it is a very interesting object of research with numerous applications in geometry (and in high energy physics\footnote{Minimal Lagrangian submanifolds of Calabi Yau 3 and 4-folds appear in string theory under the name of  ''super-symmetric cycle'' and they play a key role in the Strominger Yau Zaslow theory of mirror symmetry.}) in general there are no reasons why $d\beta$ should be equal to zero and the general hamiltonian stationary equation in K\"ahler-Einstein Manifold is the 3rd order equation\footnote{The equation is second order with respect to the multivalued Lagrangian angle function $\beta$ but the expression of $d\beta$ involves 2 derivatives with respect to the Lagrangian immersion.}.
It is interesting to observe that, while having $d\beta$ being harmonic we have a ``finite dimensional perturbation'' of $\beta=\beta_0$ but nevertheless H-minimal surfaces have very different features. For instance there exists closed $H-$minimal surfaces in ${\C}^2$ (the Clifford torus $S^1\times S^1\subset {\C}\times {\C}$ is an example of such a surface\footnote{One of the Oh conjecture says that this torus should be minimizing in its Hamiltonian isotopy class : The class of tori which can be obtained through hamiltonian deformations, that is, deformations which infinitesimally are given by (\ref{hamilton}) (see \cite{Oh2}). }) while the existence of closed minimal surfaces is prevented by the maximum principle.

\subsection{The Schoen Wolfson existence and regularity result}

In a breakthrough paper \cite{SW}, Rick Schoen and Jon Wolfson performed the very first variational construction of {\it minimizing stationary lagrangian surfaces} in homology classes. They  proved, in particular, the existence of an area minimizing Lagrangian surface
in any 2 dimensional lagrangian homology class of a K\"ahler surface. The main issue is about proving the regularity of these minimizers. The two authors proved that these area minimizing lagrangian stationary surfaces are branched immersions 
with possible isolated ``conical'' singularities\footnote{The ``conical'' singularities are singularities with non flat {\it tangent cones} and around which $\beta$ is not a univalent ${\R}-$valued map and is contributing to the total Maslov class. They are called nowadays ``Schoen Wolfson cones''.}. At this singularities $a_i$ the $S^1$ valued map $\mathfrak g$ realizes some topological degree $d_i$  and solves in conformal coordinates the equation
\[
\mbox{curl}\lf[{\mathfrak g}^{-1}\,\nabla {\mathfrak g}   \rg]=\sum_{i\in I} d_i\, \delta_{a_i}\ .
\]
In fact these conical singularities, whose description is fully given in \cite{SW}, are proved to exist even for minimizers inside some spherical Lagrangian homology classes (see \cite{MiWo} and \cite{Wol}) in some K\"ahler-Einstein surfaces. This last fact is a surprise and is one of the reasons why the variational analysis of Lagrangian surfaces, even in the lowest co-dimension 2, is particularly involved.

The starting point for performing such an analysis is the search of a monotonicity formula. In classical minimal surface theory it says that the quantity of area of a closed minimal surface inside a ball of radius $r$
divided by $r^2$ is increasing. At this point Schoen and Wolfson made the following crucial observation that such a monotonicity formula is not true in ${\C}^2$ for  $H-$minimal surfaces in general (even minimizing !).
They give an explicit simple counter-example in \cite{SW0} (section 4). The existence of such a counter-example could be explained by the fact that the generator of dilation\footnote{The existence of monotonicity formula, which is a conservation law, is intimately linked to the existence of infinitesimal isotropy of the space (see \cite{All}) and technically in general is obtained by considering the variation of the surface with respect to the generator of dilations ``exploring'' this isotropy.} in ${\C}^2$ is not hamiltonian :
\[
y_1\,\p_{y_1}+y_2\,\p_{y_2}+y_3\,\p_{y_3}+y_4\,\p_{y_4}\ne i\nabla h
\]
The major observation made by Schoen and Wolfson is that by adding the Legendrian coordinate (if it locally exists) this problem can be overcome. If $u$ is Lagrangian in ${\C}^2$ there holds
\[
d(u_1\,d{u_2}-u_2\, d{u_1}+u_3\,d{u_4}-u_4\,d{u_3})=0
\]
and one introduces the local Legendrian coordinate
\[
d\varphi_u:=u_1\,d{u_2}-u_2\, d{u_1}+u_3\,d{u_4}-u_4\,d{u_3}\ .
\]
In ${\R}^5$ one introduce the form
\[
\al:=-d\varphi+y_1\,d{y_2}-y_2\, d{y_1}+y_3\,d{y_4}-y_4\,d{y_3}
\]
If $u$ is Lagrangian the map $v:=(\varphi_u,u_1,u_2,u_3,u_4)$ is cancelling the form $\al$ which is non integrable $\al\wedge d\al\wedge d\al\ne 0$. The map $v$ is called {\it Legendrian}. Infinitesimal variations 
among Legendrian maps in ${\R}^5$ are given by Hamiltonian vector fields in the contact space $({\R}^5,\al)$ of the form
\be
\label{ham}
\vec{X}:=J_H\, \nabla^H h-2\, h\,\p_{\varphi}\ .
\ee
where $H$ is the 4-plane given by the kernel of $\al$, $\nabla^H$ is the projection on $H$ of the gradient in ${\R}^5$ with respect to the metric such that  $\pi_\ast$, the differential of the canonical projection $\pi(\varphi,y_1,y_2,y_3,y_4) =(y_1,y_2,y_3,y_4)$ from $H$ into ${\C}^2$, realizes an isometry and the length of $\p_\varphi$ is 1. ${\R}^5$ equipped with this metric defines the Heisenberg group ${\mathbb H}^2$. Finally $J_H$ is the complex structure on $H$ such that it projects by $\pi_\ast$ onto the standard complex structure of ${\C}^2$. The authors in this work introduce the notion of $H-$minimal Legendrian for the immersion which are critical point of the area with respect to all infinitesimal variations generated by Hamiltonian vector fields in ${\R}^5$ of the form (\ref{ham}). One then observe
\[
u=\pi\circ v\mbox{ is H-Minimal in }{\C}^2\quad\Longleftrightarrow\quad v \mbox{ is H-Minimal Legendrian in }({\R}^5,\al)\ .
\]
Testing $h=-\varphi$ produces the generator of dilations in ${\C}^2$ plus dilations in the 5th Legendrian direction but with a different scale.
\[
-\,i\, \nabla^H \varphi+2\, \varphi\,\p_{\varphi}=y_1\,\p_{y_2}-y_2\, \p_{y_1}+y_3\,\p_{y_4}-y_4\,\p_{y_3}+2\, \varphi\,\p_{\varphi}\ .
\]
This fact is an encouragement for looking for a monotonicity formula for H-Minimal Legendrian surfaces. However, as a matter of fact,  it cannot be deduced from a simple insertion of this generator of the dilation in the stationarity condition. The derivation of a monotonicity formula for  H-Minimal Legendrian immersions in $({\R}^5,\al)$
proposed in \cite{SW}, which is maybe one of the main achievement of this work, is quite indirect and involved (going through the resolution of a wave equation, the use of Bessel functions...etc). Another drawback is that the formula is semi-explicit.

Recently the author in \cite{Riv} obtained a new monotonicity formula by inserting in the Legendrian stationarity condition an explicit hamiltonian : $h:=\chi(\r)\,\arctan(2\,\varphi/|y|^2)$, where $\r$ is the 
Folland-Kor\'anyi gauge given by $\r^4=|y|^4+2\,\varphi^2$ and $\chi$ is some well chosen cut-off function. Because of its explicit nature the formula is more effective and is the starting point of the present work and the proof of the main result which was out of reach without it.
\subsection{Main result}

Let $(N^5,g)$ be a 5 dimensional oriented Riemannian manifold and $\al$ a non degenerate $1-$form on $N^5$  (i.e. $(N^5,\al)$ is contact) satisfying
\[
\al\wedge d\al\wedge d\al>0\ ,
\]
The triple $(N^5,g,\al)$ is called a Sasakian structure if the cone $(N^5\times {\R}_+, k:=dt^2+t^2 g)$ with the  non degenerate symplectic form $\Om:=2^{-1}\, d(t^2\,\al)$ is K\"ahler. Let ${ J}$ be the compatible complex structure ($\Om(\cdot,{ J}\cdot )=g$). The tangent vector field along $N^5$ given by
$\vec{R}:=J(t\p_t)$ is unit on $N^5$ and is orthogonal to the ``horizontal hyperplanes'' given by $H=\mbox{Ker}(\al)$. It is called {\it Reeb Vector-field} of the distribution $H$. This distribution of plane is invariant under the action of $J$ in $N^5\times {\R}_+$ and we shall denote $J_H$ its restriction on $H$. Such a structure is called Sasakian structure (see \cite{spar}). Classical examples are the Heisenberg group ${\mathbb H}^2$, the unit sphere $S^5$ for the Reeb vector-field tangent to the fibers of the Hopf fibration into ${\C}P^2$ or the Stiefel manifold $V_2({\R}^4)\simeq S^3\times S^2$ of orthonormal 2-frames in ${\R}^4$ for the Reeb vector-field given by fibers of the tautological projection onto the Grassman manifold $G_2({\R}^4)$ of 2-planes in ${\R}^4$.


A  Sasakian structure being given, we introduce the Sobolev space of Legendrian $W^{2,4}$ immersions  of a closed oriented surface $\Sigma$ in $N^5$
\[
{\mathfrak E}_{\Sigma,Leg}^{2,4}(\Sigma,N^5):=\lf\{\vec{\La}\in W^{2,4}(\Sigma,N^5)\quad ;\quad |d\vec{\La}\,\dot{\wedge}\,d\vec{\La}|>0\ \mbox{ and }\ \vec{\La}^\ast\al=0\quad\mbox{ on }\Sigma\rg\}\ .
\]
We prove in section II that ${\mathfrak M}:={\mathfrak E}_{\Sigma,Leg}^{2,4}(\Sigma,N^5)$ has the structure of Banach manifold and possesses a compatible Finsler structure for which the associated {\it Palais distance} is complete.

An admissible family ${\mathcal A}$ in ${\mathfrak M}$ is a set of subsets of ${\mathfrak M}$ which is invariant under isotopies in ${\mathfrak M}$. The min-max value or the ``width'' associated to such a family is the number given by
\[
\beta:=\inf_{A\in{\mathcal A}}\ \sup_{\vec{\La}\in A} \int_{\Sigma} dvol_{\vec{\La}}\ .
\]
Our main result in the present work says that any such min-max is achieved by a continuous conformal Legendrian map from a closed riemann surface $S$ into $N^5$ equipped with an integer multiplicity bounded in $L^\infty$ and satisfying a weak version of the $H-$minimal equation. Such a surface is called {\it Legendrian Hamiltonian stationary parametrized integral varifold}. Precisely we have:
\begin{Th}
\label{th-1}
Let ${\mathcal A}$ be an admissible family for the space of $W^{2,4}$ Legendrian immersions of a closed oriented surface $\Sigma$ into a closed Sasakian manifold $(N^5,g,\al)$. Assume that the associated width is strictly positive
\[
\beta:=\inf_{A\in{\mathcal A}}\ \sup_{\vec{\La}\in A} \int_{\Sigma} dvol_{\vec{\La}}\ >0\ .
\]
Then there exists a closed Riemann surface\footnote{We equip the surface with a compatible  metric.} $(S,h)$ with genus$(S)\,\le$ genus$(\Sigma)$ a map $\vec{\La}\in C^0(S,N^5)\cap W^{1,2}(S,N^5)$ and $N\in L^\infty(S,{\N}^\ast)$ such that
\begin{itemize}
\item[ i)]
\[
\vec{\La}^\ast\al=0\ ,
\]
\item[ii)]
\[
\p\vec{\La}\,\dot{\otimes}\,\p\vec{\La}:=\p_z\vec{\La}\cdot\p_z\vec{\La}\ dz\otimes dz=0\quad\quad\mbox{ a. e. on }S\ ,
\]
\item[iii)]
\[
\int_{S} N\, dvol_{\vec{\La}}= \frac{1}{2}\,\int_{S} N\, |\nabla\vec{\La}|^2_h\ dvol_h=\beta\ ,
\]
\item[iv)]
\be
\label{ham-stat}
\begin{array}{l}
\ds\forall \,h\in C^3(S)\quad,\quad\forall f\in C^1(S)\quad\mbox{ for a.e. }\la\in {\R}\quad\mbox{ s.t.} \quad \vec{\La}(f^{-1}\{\la\})\cap\mbox{Supp}(h)=\emptyset  \\[5mm]
\ds\int_{f^{-1}((\la,+\infty))}N\, \nabla\vec{\La}\cdot \nabla\lf[\vec{X}_h\circ\vec{\La}\rg]\ dvol_h=0\ ,
\end{array}
\ee
where $\vec{X}_h$ is the Hamiltonian vector-field associated to $h$ and given by
\be
\label{ham-11}
\vec{X}_h:=J_H(\nabla^H h)- 2\, h\, \vec{R}\ ,
\ee
where $\nabla^H$ is the projection onto $H$ of the gradient in $(N^5,g)$. 

\end{itemize}
\end{Th}
In \cite{Riv-1} the corresponding result  to theorem~\ref{th-1} was proven without the pointwise Legendrian constraint. Solutions to (\ref{ham-stat}) for arbitrary $\vec{X}$ (non necessarily Hamiltonian) were called {\it parametrised stationary integral varifolds}. In an analogous way we shall be calling any triple $(S,\vec{\La}, N)$ solving (\ref{ham-stat}) for any Hamiltonian vector-field of the form (\ref{ham-11})  {\it Legendrian parametrised Hamiltonian stationary integral varifolds}. 

We stress that the stationarity identity (\ref{ham-stat}) holds for a.e. $\la\in {\R}$ foor which the image by $\vec{\La}$ of the level set $f^{-1}\{\la\}$ avoids the support of the testing
vector-field $\vec{X}$. This  is the called ``localisation hypothesis'' (It has been first introduced in \cite{Riv-00}). Combining the nullity of the variation of the area for perturbation $\vec{X}$ in the target with the ``localisation hypothesis'' is a substitute to the Euler Lagrange equation (\ref{4}), which has no satisfying weak formulation as we discussed in the foreword.

The strategy adopted to establish theorem~\ref{th-1} is the penalisation approach introduced by the author in \cite{Riv-1} under the name ``viscosity method''. This method has then been implemented for free boundary surfaces in \cite{Piga} and we will adopt some of the improvements of the method brought in this work. Precisely we study the smoothed minmax operation
\[
\beta(\ep):=\inf_{A\in{\mathcal A}}\ \sup_{\vec{\La}\in A} E_\ep(\vec{\La})\ ,
\]
with
\[
E_\ep(\vec{\La}):=\int_{\Sigma} dvol_{\vec{\La}}+\ep^4\,\int_\Sigma (1+|\vec{\mathbb I}_{\vec{\La}}|^2)^2\ dvol_{\vec{\La}}\ ,
\]
where $\vec{\mathbb I}_{\vec{\La}}$ is the second fundamental form of $\vec{\La}$.  In section II We construct almost critical points $\vec{\La}_k$ of  $E_{\ep_k}$  for a sequence $\ep_k\rightarrow 0$ such that
\[
E_{\ep_k}(\vec{\La}_k)-\beta(\ep_k)\longrightarrow 0\quad\mbox{ and }\quad \ep^4_k\,\int_\Sigma (1+|\vec{\mathbb I}_{\vec{\La}}|^2)^2\ dvol_{\vec{\La}_k}=o\lf(\frac{1}{\log\ep_k^{-1}}\rg)\ .
\]
The main difficulty to overcome in sections III and IV  consists in passing to the limit for some well chosen subsequence of $\vec{\La}_k$ in order to obtain the map $\vec{\La}$ satisfying the conclusions of theorem~\ref{th-1}.  While we broadly follow  the same scheme as the one introduced in \cite{Riv-1}, the Legendrian constraint however introduces new serious challenges for the proof of theorem~\ref{th-1}.

\medskip

In \cite{PiRi1} it is proved that any {\it parametrised stationary varifolds} is a smooth branched minimal immersion equipped with a smooth multiplicity. We make the following conjecture.

\begin{Con}
\label{con-reg}
Prove that 2-dimensional Legendrian parametrised  Hamiltonian stationary integral varifolds are  branched Legendrian immersions and smooth away from isolated Schoen Wolfson conical singularities. 
\end{Con}

\medskip

We first present the proof of theorem~\ref{th-1} for $N^5=V_2({\R}^4)$ which has motivated the present works for reasons explained in \cite{Riv-1}  in relation with the Willmore conjecture.
The case $N^5=V_2({\R}^4)$ has the advantage of clarity in the presentation due to the symmetric nature of the Stiefel space of orthonormal 2 frames in ${\R}^4$. This being said the global symmetry of $V_2({\R}^4)$ has no incidence at all on the analysis presented
in sections II, IV and V. What is crucial is the common nature of the asymptotic of any arbitrary Sasakian manifold while dilating at any arbitrary point. It always converge to to the Heisenberg Group ${\mathbb H}^2$.
Hence the proof of theorem~\ref{th-1} for general target is formally identical and requires some obvious adaptation of the computation of sections III, IV and V. This is explained in the last section VI.

\section{Notations and Preliminaries on $V_2({\R}^4)$}
\reset
\subsection{Contact structure;  complex structure on  horizontal planes}
We denote by $V_2({\R}^4)$ the Stiefel Manifold of orthonormal 2 frames in ${\R}^4$ :
\[
V_2({\R}^4):=\lf\{ (\vec{a},\vec{b})\in S^3\times S^3\quad;\quad \vec{a}\cdot\vec{b}=0  \rg\}
\]
We consider on $V_2({\R}^4)$ the following 2-form
\[
\al:=\vec{a}\cdot d\vec{b}-\vec{b}\cdot d\vec{a}=\sum_{i=1}^4 a^i\,db^i-b^i\,da^i\quad.
\]
where $\vec{a}=(a^1,\cdots,a^4)$ and $\vec{b}=(b^1,\cdots,b^4)$. The Stiefel manifold $V_2({\R}^4)$ defines an $S^1-$bundle over $G_2({\R}^4)\simeq S^2_+\times S^2_-$. The projection ``Hopf map'' $\Pi$ which to $(\vec{a},\vec{b})$ assigns
$((\vec{a}\wedge \vec{b})^+,(\vec{a}\wedge \vec{b})^-)$ where $(\vec{a}\wedge \vec{b})^\pm$ are respectively the dual and anti-dual parts of $\sqrt{2}\,\vec{a}\wedge\vec{b}$ that is
 where
\[
(\vec{a}\wedge \vec{b})^\pm=\frac{1}{\sqrt{2}}\,\lf[\vec{a}\wedge \vec{b}\pm\ast (\vec{a}\wedge \vec{b}) \rg]
\]
and $\ast$ is the Hodge operator on 2 vectors in ${\R}^4$ such that 
\[
\vec{e}\wedge \vec{f}\wedge\ast (\vec{e}\wedge \vec{f})=|\vec{e}\wedge\vec{f}|^2\ \vec{\ep}_1\wedge\vec{\ep}_2\wedge\vec{\ep}_3\wedge\vec{\ep}_4\ ,
\]
where $(\vec{\ep}_1\cdots \vec{\ep}_4)$ is the canonical basis of ${\R}^4$.

So we have
where we recall that $\Pi$ is the tautological projection from $V_2({\R}^4)$ into $G_2({\R}^4)\simeq S^2_+\times S^2_-$ given by
\[
\Pi\ :\ (\vec{a},\vec{b})\ \longrightarrow\ ((\vec{a}\wedge\vec{b})^+,(\vec{a}\wedge\vec{b})^-)= \lf(\frac{\vec{a}\wedge\vec{b}+\ast\, \vec{a}\wedge\vec{b}}{\sqrt{2}}\ ,\  \frac{\vec{a}\wedge\vec{b}-\ast\, \vec{a}\wedge\vec{b}}{\sqrt{2}}    \rg)
\]
We shall denote  $\pi_{+}$ and $\pi_-$  the canonical projections
\[
\pi_{+}\ :\ S^2_+\times S^2_-\ \longrightarrow\ S^2_+\quad\mbox{ and }\quad\pi_{-}\ :\ S^2_+\times S^2_-\ \longrightarrow\ S^2_-\ .
\]
Using the quaternionic representatives of $\vec{a}$ and $\vec{b}$ that we denote respectively
\[
{\mathbf a}:=a^1+a^2\,{\mathbf i}+a^3\,{\mathbf j}+a^4\,{\mathbf k}\quad\mbox{ and }\quad{\mathbf b}:=b^1+b^2\,{\mathbf i}+b^3\,{\mathbf j}+b^4\,{\mathbf k}
\]
the projection $\Pi$ identifies to the following map
\[
\Pi\ :\ \vec{a}\wedge\vec{b}\in G_2({\R}^4)\longleftrightarrow\lf(\lf(
\begin{array}{l}
<{\mathbf a},{\mathbf i}{\mathbf b}>\\[3mm]
<{\mathbf a},{\mathbf j}{\mathbf b}>\\[3mm]
<{\mathbf a},{\mathbf k}{\mathbf b}>
\end{array}
\rg),\lf(
\begin{array}{l}
<{\mathbf a},{\mathbf b}{\mathbf i}>\\[3mm]
<{\mathbf a},{\mathbf b}{\mathbf j}>\\[3mm]
<{\mathbf a},{\mathbf b}{\mathbf k}>
\end{array}
\rg)\rg)=({\bf a}{\bf b}^\ast,{\bf b}^\ast {\bf a})\in S^2_+\times S^2_-
\]
where $<\cdot,\cdot>$ is the scalar product on ${\mathbf H}\simeq {\R}^4$. This diffeomorphism will naturally induce the complex structure of $\C P^{1}\times \ov{\C P^1}$ on $G_2({\R}^4)$ as we will see where ${\C} P^1$ is referring to the  reversed complex structure  (i.e reversed orientation) with respect to the canonical one $J_+$ on $S^2$ : $J_- X=-J_+ X$.

\medskip

Let ${\vec \La}=(\vec{a},\vec{b})$ we have
\[
T_{\vec{\La}}V_2({\R}^4)=\lf\{ (\vec{V},\vec{W}) \ ;\ \vec{V}\cdot\vec{a}=0\ ,\ \vec{W}\cdot\vec{b}=0\ \mbox{ and }\ \vec{a}\cdot\vec{W}+\vec{V}\cdot\vec{b}=0  \rg\}\ .
\]
Let $(\vec{c},\vec{d})\in V_2({\R}^4)$ such that $\ast( \vec{a}\wedge\vec{b})=\vec{c}\wedge\vec{d}$. We denote $\vec{G}=\vec{a}\wedge\vec{b}$ and note that
\[
T_{\vec{G}}G_2({\R}^4)=\lf\{\vec{a}\wedge \vec{W}+\vec{V}\wedge\vec{b}\quad \mbox{ where }\quad \vec{V}\,,\,\vec{W}\in \mbox{Span}(\vec{a},\vec{b})^\perp=\mbox{Span}(\vec{c},\vec{d})\rg\}\ ,
\]
and
\[
\begin{array}{rl}
\Pi_\ast T_{\vec{\La}}V_2({\R}^4)&\ds=\lf\{((\vec{a}\wedge \vec{W})^++(\vec{V}\wedge\vec{b})^+,(\vec{a}\wedge \vec{W})^-+(\vec{V}\wedge\vec{b})^-)\quad \mbox{ where }\quad \vec{V}\,,\,\vec{W}\in \mbox{Span}(\vec{a},\vec{b})^\perp\rg\}\\[5mm]
 &\ds=TS^2_+\oplus TS^2_-\ .
 \end{array}
\]
We denote respectively $\vec{G}^+:=\pi_+(\Pi(\vec{\La}))=(\vec{a}\wedge \vec{b})^+$ and $\vec{G}^-:=\pi_-(\Pi(\vec{\La}))=(\vec{a}\wedge \vec{b})^-$. Denote by $\pi_\perp$ the projection onto  the span of $\{\vec{c},\vec{d}\}$. First we introduce
We equip the unit spheres of unit self-dual (resp. anti-self-dual) 2-vectors of ${\R}^4$ with the following complex structure
\[
J_{(\vec{a}\wedge\vec{b})^\pm}((\vec{a}\wedge\vec{c})^\pm)=\pm(\vec{a}\wedge\vec{d})^\pm\ .
\]
Then we equip the product $S^2_+\times S^2_-$ with the product complex structure
\[
J_{\vec{G}^{\,\pm}}(\vec{a}\wedge\vec{c}\pm\vec{d}\wedge \vec{b})=\pm[\vec{a}\wedge\vec{d}\mp\vec{c}\wedge \vec{b}]
\]
For $\vec{X}\in \mbox{Span}(\vec{c},\vec{d})$ we have $(\vec{a}\wedge \vec{X})^\pm= \sqrt{2^{-1}}\ (\vec{a}\wedge \vec{X}\pm\vec{X}^\perp\wedge\vec{b}) $ where $\vec{X}^\perp:=\star \lf(\vec{a}\wedge \vec{b}\wedge X\rg)$ or in other words $\vec{a}\wedge \vec{b}\wedge \vec{X}\wedge \vec{X}^\perp=|X|^2\ \star 1$.  With this notation we have
\[
\forall\, \vec{X}\in \,\mbox{Span}(\vec{c},\vec{d})\quad\quad J_{{\vec{G}^{\,\pm}}}\lf((\vec{a}\wedge\vec{X})^\pm\rg)=\pm(\vec{a}\wedge\vec{X}^{\perp})^\pm\ .
\]
Let $X,Y\in T_{\vec{\La}}V_2({\R}^4)$, we have\footnote{The - sign comes from the fact that the relation between the K\"ahler metric and its associated K\"ahler form is given by $g(u,v)=\om(u,Jv)$ and hence $-g(u,Jv)=\om(u,v)$.}
\[
\ds-\lf((\vec{a}\wedge\vec{b})^+\rg)^\ast\om_{S^2}\,(X,Y)=\lf(d_X(\vec{a}\wedge\vec{b})^+,J_{\vec{G}^+}(d_Y(\vec{a}\wedge\vec{b})^+)\rg)
\]
where $d_X(\vec{a}\wedge\vec{b})^+$ is denoting the push forward by $\Pi$ of the tangent vector $X$ at $(\vec{a},\vec{b})$. Let $P^+(\vec{a},\vec{b})=\vec{a}$ and 
$P^- (\vec{a},\vec{b})=\vec{b}$   the push forward of $X$ by $P^+$ (resp. by $P^-$) will simply be denoted $d_X\vec{a}$ (resp. $d_X\vec{b}$) 
\[
\begin{array}{l}
\ds d_X(\vec{a}\wedge\vec{b})^+=\frac{1}{\sqrt{2}}(d_X(\vec{a}\wedge\vec{b})+\ast d_X(\vec{a}\wedge\vec{b}))=d_X\vec{a}\cdot\vec{c}\ (\vec{c}\wedge\vec{b})^++d_X\vec{a}\cdot\vec{d}\ (\vec{d}\wedge\vec{b})^+\\[5mm]
\ds +d_X\vec{b}\cdot\vec{c}\ (\vec{a}\wedge\vec{c})^++d_X\vec{b}\cdot\vec{d}\ (\vec{a}\wedge\vec{d})^+\\[5mm]
\ds=[d_X\vec{a}\cdot\vec{c}-  d_X\vec{b}\cdot\vec{d} ]\ (\vec{c}\wedge\vec{b})^++[d_X\vec{a}\cdot\vec{d}+d_X\vec{b}\cdot\vec{c}] \ (\vec{d}\wedge\vec{b})^+
\end{array}
\]
Hence
\[
\begin{array}{l}
\ds J_{\vec{G}^+}(d_Y(\vec{a}\wedge\vec{b})^+)=[d_Y\vec{a}\cdot\vec{c}-d_Y\vec{b}\cdot\vec{d} ]\ J_{\vec{G}^+}(\vec{c}\wedge\vec{b})^++[d_Y\vec{a}\cdot\vec{d}+d_Y\vec{b}\cdot\vec{c}]\ J_{\vec{G}^+}(\vec{d}\wedge\vec{b})^+\\[5mm]
\ds =[d_Y\vec{a}\cdot\vec{c}-d_Y\vec{b}\cdot\vec{d} ]\ (\vec{d}\wedge\vec{b})^+-[d_Y\vec{a}\cdot\vec{d}+d_Y\vec{b}\cdot\vec{c}]\ (\vec{c}\wedge\vec{b})^+\ .
\end{array}
\]
We deduce
\[
\begin{array}{l}
\ds\lf(d_X(\vec{a}\wedge\vec{b})^+,J_{\vec{G}^+}(d_Y(\vec{a}\wedge\vec{b})^+)\rg)\\[5mm]
=\lf([d_X\vec{a}\cdot\vec{c}-  d_X\vec{b}\cdot\vec{d} ]\ (\vec{c}\wedge\vec{b})^++[d_X\vec{a}\cdot\vec{d}+d_X\vec{b}\cdot\vec{c}]\ (\vec{d}\wedge\vec{b})^+,\rg.\\[5mm]
\ds \lf .[d_Y\vec{a}\cdot\vec{c}-d_Y\vec{b}\cdot\vec{d} ]\ (\vec{d}\wedge\vec{b})^+-[d_Y\vec{a}\cdot\vec{d}+d_Y\vec{b}\cdot\vec{c}]\ (\vec{c}\wedge\vec{b})^+\  \rg)\\[5mm]
\ds=-\,[d_X\vec{a}\cdot\vec{c}-  d_X\vec{b}\cdot\vec{d} ]\ [d_Y\vec{a}\cdot\vec{d}+d_Y\vec{b}\cdot\vec{c}]+[d_X\vec{a}\cdot\vec{d}+d_X\vec{b}\cdot\vec{c}]\,[d_Y\vec{a}\cdot\vec{c}-d_Y\vec{b}\cdot\vec{d} ]\\[5mm]
=-\,d_X\vec{a}\cdot\vec{c}\ d_Y\vec{a}\cdot\vec{d}+d_X\vec{a}\cdot\vec{d}\,d_Y\vec{a}\cdot\vec{c}-\,d_X\vec{b}\cdot\vec{c}\ d_Y\vec{b}\cdot\vec{d}+d_X\vec{b}\cdot\vec{d}\,d_Y\vec{b}\cdot\vec{c}
-d_X\vec{a}\cdot d_{Y}\vec{b}+d_{Y}\vec{a}\cdot d_X\vec{b}\\[5mm]
\ds=-\,\star[\vec{a}\wedge\vec{b}\wedge d_X\vec{a}\wedge d_{Y}\vec{a}]-\,\star[\vec{a}\wedge\vec{b}\wedge d_X\vec{b}\wedge d_{Y}\vec{b}]-d_X\vec{a}\cdot d_{Y}\vec{b}+d_{Y}\vec{a}\cdot d_X\vec{b}
\end{array}
\]
We have
\[
d\vec{a}\wedge d\vec{a}= 2\,(d\vec{a}\cdot \vec{b})\wedge (d\vec{a}\cdot \vec{c})\ \ \vec{b}\wedge\vec{c}+2\,(d\vec{a}\cdot \vec{b})\wedge (d\vec{a}\cdot \vec{d})\ \ \vec{b}\wedge\vec{d}+2\,(d\vec{a}\cdot \vec{c})\wedge (d\vec{a}\cdot \vec{d})\ \ \vec{c}\wedge\vec{d}
\]
Hence
\[
\frac{1}{2}\,\vec{a}\wedge\vec{b}\wedge d\vec{a}\wedge d\vec{a}=(d\vec{a}\cdot \vec{c})\wedge (d\vec{a}\cdot \vec{d})\ \vec{a}\wedge\vec{b}\wedge\vec{c}\wedge\vec{d}\ .
\]
Similarly
\[
\frac{1}{2}\,\vec{a}\wedge\vec{b}\wedge d\vec{b}\wedge d\vec{b}=(d\vec{b}\cdot \vec{c})\wedge (d\vec{b}\cdot \vec{d})\ \vec{a}\wedge\vec{b}\wedge\vec{c}\wedge\vec{d}\ .
\]
Denote 
\be
\label{dcd}
d\vec{c}\,\dot{\wedge}\, d\vec{d}:=\sum_{l=1}^4dc^l\wedge dd^l\ ,
\ee
there holds
\[
d\vec{c}\,\dot{\wedge}\, d\vec{d}=(d\vec{c}\cdot\vec{a})\,{\wedge}\, (\vec{a}\cdot d\vec{d})+(d\vec{c}\cdot\vec{b})\,{\wedge}\, (\vec{b}\cdot d\vec{d})=(\vec{c}\cdot d\vec{a})\,{\wedge}\, (d\vec{a}\cdot \vec{d})+(\vec{c}\cdot d\vec{b})\,{\wedge}\, (d\vec{b}\cdot \vec{d})\ .
\]
Hence
\[
d\vec{c}\,\dot{\wedge}\, d\vec{d}=\frac{1}{2}\,\star_{\wedge{\R}^4}\ \vec{a}\wedge\vec{b}\wedge \lf[d\vec{a}\wedge d\vec{a}+d\vec{b}\wedge d\vec{b}\rg]
\]
Thus
\[
-\lf((\vec{a}\wedge\vec{b})^+\rg)^\ast\om^+_{S^2}\,(X,Y)=-d\vec{c}\,\dot{\wedge}\, d\vec{d}(X,Y)-d\vec{a}\dot{\wedge} d\vec{b}(X,Y)
\]
Similarly one computes
\[
\ds-\lf((\vec{a}\wedge\vec{b})^-\rg)^\ast\om^-_{S^2}\,(X,Y)=\lf(d_X(\vec{a}\wedge\vec{b})^-,J_{\vec{G}^-}(d_Y(\vec{a}\wedge\vec{b})^-)\rg)
\]
but
\[
\begin{array}{l}
\ds d_X(\vec{a}\wedge\vec{b})^-=\frac{1}{\sqrt{2}}(d_X(\vec{a}\wedge\vec{b})-\ast d_X(\vec{a}\wedge\vec{b}))=d_X\vec{a}\cdot\vec{c}\ (\vec{c}\wedge\vec{b})^-+d_X\vec{a}\cdot\vec{d}\ (\vec{d}\wedge\vec{b})^-\\[5mm]
\ds +d_X\vec{b}\cdot\vec{c}\ (\vec{a}\wedge\vec{c})^-+d_X\vec{b}\cdot\vec{d}\ (\vec{a}\wedge\vec{d})^-\\[5mm]
\ds=[d_X\vec{a}\cdot\vec{c}+  d_X\vec{b}\cdot\vec{d} ]\ (\vec{c}\wedge\vec{b})^-+[d_X\vec{a}\cdot\vec{d}-d_X\vec{b}\cdot\vec{c}] \ (\vec{d}\wedge\vec{b})^-
\end{array}
\]
Hence
\[
\begin{array}{l}
\ds J_{\vec{G}^-}(d_Y(\vec{a}\wedge\vec{b})^-)=[d_Y\vec{a}\cdot\vec{c}+  d_Y\vec{b}\cdot\vec{d} ]\ J_{\vec{G}^-}(\vec{c}\wedge\vec{b})^-+[d_Y\vec{a}\cdot\vec{d}-d_Y\vec{b}\cdot\vec{c}] \ J_{\vec{G}^-}(\vec{d}\wedge\vec{b})^-\\[5mm]
\ds =-\,[d_Y\vec{a}\cdot\vec{c}+  d_Y\vec{b}\cdot\vec{d} ]\ (\vec{d}\wedge\vec{b})^-+[d_Y\vec{a}\cdot\vec{d}-d_Y\vec{b}\cdot\vec{c}] \ (\vec{c}\wedge\vec{b})^-\ .
\end{array}
\]
We deduce
\[
\begin{array}{l}
\ds\lf(d_X(\vec{a}\wedge\vec{b})^-,J_{\vec{G}^-}(d_Y(\vec{a}\wedge\vec{b})^-)\rg)=[d_X\vec{a}\cdot\vec{c}+  d_X\vec{b}\cdot\vec{d} ]\,[d_Y\vec{a}\cdot\vec{d}-d_Y\vec{b}\cdot\vec{c}]\\[5mm]
\ds-[d_X\vec{a}\cdot\vec{d}-d_X\vec{b}\cdot\vec{c}]\,[d_Y\vec{a}\cdot\vec{c}+  d_Y\vec{b}\cdot\vec{d} ]\\[5mm]
\ds=d_X\vec{a}\cdot\vec{c}\ d_Y\vec{a}\cdot\vec{d}-d_X\vec{a}\cdot\vec{d}\,d_Y\vec{a}\cdot\vec{c}+\,d_X\vec{b}\cdot\vec{c}\ d_Y\vec{b}\cdot\vec{d}-d_X\vec{b}\cdot\vec{d}\,d_Y\vec{b}\cdot\vec{c}
-d_X\vec{a}\cdot d_{Y}\vec{b}+d_{Y}\vec{a}\cdot d_X\vec{b}\\[5mm]
\ds=\star[\vec{a}\wedge\vec{b}\wedge d_X\vec{a}\wedge d_{Y}\vec{a}]+\,\star[\vec{a}\wedge\vec{b}\wedge d_X\vec{b}\wedge d_{Y}\vec{b}]-d_X\vec{a}\cdot d_{Y}\vec{b}+d_{Y}\vec{a}\cdot d_X\vec{b}\ .
\end{array}
\]
Hence
\[
-\lf((\vec{a}\wedge\vec{b})^-\rg)^\ast\om^-_{S^2}\,(X,Y)=d\vec{c}\,\dot{\wedge}\, d\vec{d}(X,Y)-d\vec{a}\dot{\wedge} d\vec{b}(X,Y)\ .
\]
We deduce the following lemma
\begin{Lm}
\label{lm-structure} {\bf[The contact Form]}
Let $(\vec{a},\vec{b})\in V_{2}({\R}^4)$ and $(\vec{c},\vec{d})\in V_{2}({\R}^4)$ such that
\[
\ast \vec{a}\wedge\vec{b}=\vec{c}\wedge \vec{b}
\]
Then the following identities hold
\[
\lf((\vec{a}\wedge\vec{b})^\pm\rg)^\ast\om^\pm_{S^2}=\pm d\vec{c}\,\dot{\wedge}\, d\vec{d}+ d\vec{a}\dot{\wedge} d\vec{b}
\]
and
\[
d\vec{c}\,\dot{\wedge}\, d\vec{d}=\frac{1}{2}\,\star_{\wedge{\R}^4}\ \vec{a}\wedge\vec{b}\wedge \lf[d\vec{a}\wedge d\vec{a}+d\vec{b}\wedge d\vec{b}\rg]\ .
\]
In particular\footnote{Observe that $\pi_+^\ast\om_{S^2}+\pi_-^\ast\om_{S^2}=\pi_+^\ast\om_{{\C}P^1}-\pi_-^\ast\om_{\ov{{\C}P^1}}$.}
\[
d\al=2\,d\vec{a}\dot{\wedge}d\vec{b}= \lf((\vec{a}\wedge\vec{b})^+\rg)^\ast\om_{S^2}+\lf((\vec{a}\wedge\vec{b})^-\rg)^\ast\om_{S^2}= \Pi^\ast \lf(\pi_+^\ast\om^+_{S^2}+\pi_-^\ast\om^-_{S^2}\rg)=\Pi^\ast\om_{{\C}P^1\times \ov{{\C}P^1}}\ .
\]
and
\[
\al\wedge d\al\wedge d\al=2\, \al\wedge \lf((\vec{a}\wedge\vec{b})^+\rg)^\ast\om^+_{S^2}\wedge \lf((\vec{a}\wedge\vec{b})^-\rg)^\ast\om^-_{S^2}=2\,\al\wedge \Pi^\ast \pi_+^\ast\om^+_{S^2}\wedge\Pi^\ast\pi_-^\ast\om^-_{S^2}
\]
\end{Lm}
We introduce  the Reeb Vector Field
\[
\vec{R}:=\frac{\p}{\p\theta}(\vec{a}^{\,\theta},\vec{b}^{\,\theta})_{\theta=0}:=\lf.\frac{d}{d\theta}(\cos\theta\,\vec{a}+\sin\theta\,\vec{b}\ ,\ -\sin\theta\,\vec{a}+\cos\theta\,\vec{b})\rg|_{\theta=0}=\lf(\vec{b},-\vec{a}\rg)
\]
then
\[
\al\wedge d\al\wedge d\al\res \vec{R}=-\,2\, \Pi^\ast\lf[ \pi_+^\ast\om_{S^2}\wedge\pi_-^\ast\om_{S^2}\rg]=-\,2\, \Pi^\ast\, \Om_{ S^2_+\times S^2_-}
\]
where $\Om_{ S^2_+\times S^2_-}$ is the volume form on $S^2_+\times S^2_-$. Denote by $H_{(\vec{a},\vec{b})}:=\mbox{Ker}(\Pi_{\ast})^\perp=\mbox{Span}(\vec{R})^\perp$. We have
\[
H_{(\vec{a},\vec{b})}=\lf\{ (\vec{V},\vec{W})\in {\R}^4\times {\R}^4\ ;\ \vec{V},\vec{W}\in\mbox{Span}(\vec{a},\vec{b})^\perp  \rg\}
\]
and
\[
\Pi_\ast(\vec{V},\vec{W})=((\vec{V}\wedge\vec{b})^++(\vec{a}\wedge\vec{W})^+,(\vec{V}\wedge\vec{b})^-+(\vec{a}\wedge\vec{W})^-)
\]
Thus $\Pi_\ast$ realises an isometry between $H_{(\vec{a},\vec{b})}$ and $T_{\Pi(\vec{a},\vec{b})} S^2_+\times S^2_-$. On $H$ we define the complex structure $J_H$ such that
\[
\forall \,\vec{X}\in H\quad \Pi_\ast J_H \vec{X}:=J_{{\C}P^1\times \ov{{\C}P^1}}\,\Pi_\ast \vec{X}\ .
\]
This gives in particular
\[
\forall \vec{X}\,,\vec{Y}\in H\quad \vec{Y}\cdot J_H \vec{X}=\om_{{\C}P^1\times \ov{{\C}P^1}}(\vec{X},\vec{Y})\ .
\]
Finally we are proving a last preliminary result
\begin{Lm}
\label{lm-II.2} {\bf [The transverse Complex structure]}
Under the previous notations,  $(\vec{V},\vec{W})\in H_{(\vec{a},\vec{b})}$ being given there holds
\[
J_H(\vec{V},\vec{W})=(-\vec{W},\vec{V})\ .
\]
\end{Lm}
\noindent{\bf Proof of Lemma~\ref{lm-II.2}  }It remains to prove the last assertion on $J_H$. Let $(\vec{V},\vec{W})\in H_{(\vec{a},\vec{b})}$ we have
\[
\begin{array}{l}
\ds\Pi_\ast J_H(\vec{V},\vec{W})=J_{\vec{a}\wedge\vec{b}}\Pi_\ast(\vec{V},\vec{W})=J_{\vec{a}\wedge\vec{b}}((\vec{V}\wedge\vec{b})^++(\vec{a}\wedge\vec{W})^+,(\vec{V}\wedge\vec{b})^-+(\vec{a}\wedge\vec{W})^-)\\[5mm]
\ds=((\vec{V}^\perp\wedge\vec{b})^++(\vec{a}\wedge\vec{W}^\perp)^+,-(\vec{V}^\perp\wedge\vec{b})^--(\vec{a}\wedge\vec{W}^\perp)^-)
\end{array}
\]
We have respectively
\[
\begin{array}{l}
\ds(\vec{V}^\perp\wedge\vec{b})^+=\frac{1}{2}\,[\vec{V}^\perp\wedge\vec{b}+\ast (\vec{V}^\perp\wedge\vec{b})]=\frac{1}{2}\,[\vec{V}^\perp\wedge\vec{b}+\vec{a}\wedge\vec{V}]=(\vec{a}\wedge\vec{V})^+\\[5mm]
\ds(\vec{a}\wedge\vec{W}^\perp)^+=\frac{1}{2}\,[ \vec{a}\wedge\vec{W}^\perp+\ast(\vec{a}\wedge\vec{W}^\perp)] =\frac{1}{2}\,[\vec{a}\wedge\vec{W}^\perp+\vec{b}\wedge\vec{W}]=(\vec{b}\wedge\vec{W})^+\\[5mm]
\ds-(\vec{V}^\perp\wedge\vec{b})^-=-\frac{1}{2}\,[ \vec{V}^\perp\wedge\vec{b}-\ast(\vec{V}^\perp\wedge\vec{b})]=-\frac{1}{2}\,[ \vec{V}^\perp\wedge\vec{b}-\vec{a}\wedge\vec{V}]=(\vec{a}\wedge\vec{V})^-\\[5mm]
\ds-(\vec{a}\wedge\vec{W}^\perp)^-=-\frac{1}{2}\,[ \vec{a}\wedge\vec{W}^\perp-\ast(\vec{a}\wedge\vec{W}^\perp)]=-\frac{1}{2}\,[ \vec{a}\wedge\vec{W}^\perp-\vec{b}\wedge\vec{W}]=(\vec{b}\wedge\vec{W})^-
\end{array}
\]
Thus
\be
\label{almost-comp}
\Pi_\ast J_H(\vec{V},\vec{W})=\lf((\vec{a}\wedge\vec{V})^++(\vec{b}\wedge\vec{W})^+,(\vec{a}\wedge\vec{V})^-+(\vec{b}\wedge\vec{W})^-\rg)=\Pi_\ast(-\vec{W},\vec{V})
\ee

\subsection{The space  of $W^{2,4}$ Legendrian immersions of a surface }
Let $\Sigma$ be a closed oriented $2-$dimensional manifold. In \cite{Riv-1} section II the author introduces the space of Sobolev Immersions 
\[
{\mathfrak E}_{\Sigma}^{2,4}:=W^{2,4}_{imm}(\Sigma,V_2({\R}^4)):=\lf\{\vec{\La}\in W^{2,4}(\Sigma,V_2({\R}^4))\ \mbox{ s. t. rank}(d\vec{\La})(x)=2\quad\forall\, x\in\Sigma\rg\}\ .
\]
which is proven  to be a $C^2$ submanifold of the Banach vector space $W^{2,4}(\Sigma,{\R}^4\times{\R}^4)$ ; moreover  \cite{Riv-1} introduces  a Finsler Structure over ${\mathfrak E}_{\Sigma}^{2,4}$ for which the associated Palais distance defines a complete metric space structure on ${\mathfrak E}_{\Sigma,2}$. We  now prove the following theorem.

\begin{Th}
\label{th-subman}
The space of $W^{2,4}$ Legendrian immersions of a closed oriented surface $\Sigma$ into $V_2({\R}^4)$  (i.e. elements from ${\mathfrak E}_{\Sigma}^{2,4}$ satisfying $\vec{\La}^\ast\al=0$) realises a $C^2$ sub-manifold of $W^{2,4}_{imm}(\Sigma,V_2({\R}^4))$.
\end{Th}
\noindent{\bf Proof of theorem~\ref{th-subman}.}
We consider the map
\[
{\mathfrak H}\ :\ \vec{\La}\in {\mathfrak E}_{\Sigma}^{2,4}\quad \longrightarrow\quad \vec{\La}^{\,\ast}\al\in \Gamma_{W^{1,4}}(T^\ast \Sigma)\ .
\]
We claim that ${\mathfrak H}$ realizes a $C^2$ submersion from the Banach manifold ${\mathfrak E}_{\Sigma}^{2,4}$ onto the Banach vector space of $W^{1,4}$ one forms on $\Sigma$. Moreover there holds
\[
\begin{array}{c}
\ds\forall\ \vec{\La}=(\vec{\Phi},\vec{n})\in {\mathfrak E}_{\Sigma}^{2,4}\quad \forall \ (\vec{V},\vec{W})\in W^{2,4}(\Sigma,{\R}^4\times {\R}^4)\ ;\quad\vec{V}\cdot\vec{\Phi}=0\ ,\ \vec{W}\cdot\vec{n}=0\ \mbox{ and }\ \vec{\Phi}\cdot\vec{W}+\vec{n}\cdot\vec{V}=0 \\[5mm]
\ds\quad d{\mathfrak H}_{\vec{\La}}\cdot(\vec{V},\vec{W})=\vec{V}\cdot d\vec{n}+\vec{\Phi}\cdot d\vec{W}-\vec{n}\cdot d\vec{V}-\vec{W}\cdot d\vec{\Phi}\ .
\end{array}
\]
Indeed, first, ${\mathfrak H}$ is the restriction to the $C^2$ sub-manifold ${\mathfrak E}_{\Sigma}^{2,4}$ of the Banach vector space $W^{2,4}(\Sigma,{\R}^4\times{\R}^4)$ of the (obviously $C^\infty$) bilinear form
\[
B(\vec{a},\vec{b}):=\vec{a}\cdot d\vec{b}-\vec{b}\cdot d\vec{a}\ .
\]
It remains to prove in a second step that ${\mathfrak H}$ realises a submersion that is at every point $\vec{\La}=(\vec{\Phi},\vec{n})\in {\mathfrak E}_{\Sigma}^{2,4}$ the linear map $d{\mathfrak H}_{\vec{\La}}$ is surjective from $T_{\vec{\La}}{\mathfrak E}_{\Sigma}^{2,4}$ onto $\Gamma_{W^{1,4}}(T^\ast \Sigma)$.

\medskip

Instead of proving it directly, we introduce the map
\[
{\mathfrak K}\ :\ \vec{\La}\in {\mathfrak E}_{\Sigma}^{2,4}\quad \longrightarrow\quad \lf(\vec{\La}^{\,\ast}d\al\,, \,\Pi_0(\vec{\La}^\ast\al)\,, \,\int_{\Gamma_1}\vec{\La}^\ast\al,\cdots  ,\int_{\Gamma_{2\,g(\Sigma)}}\vec{\La}^\ast\al\rg)\ ,
\]
where $\Pi_0(\vec{\La}^\ast\al)$ is the exact part in the Hodge decomposition of $\vec{\La}^\ast \al$ with respect to the smooth reference metric $g_0$
\[
\Pi_0(\vec{\La}^\ast\al):=d^{\ast_{g_0}}\Delta_{g_0}^{-1}\lf[\vec{\La}^\ast \al\rg]
\]
and $\Gamma_j$ is a fixed basis of closed loops generating $H_1(\Sigma,{\Z})$. We have that 
\[
{\mathfrak K}\ :\ {\mathfrak E}_{\Sigma}^{2,4}\quad \longrightarrow\quad W^{1,4}(\wedge^2\Sigma)\times W^{2,4}(\Sigma)\times{\R}^{2\,g(\Sigma)}\ .
\]
Since ${\mathfrak H}$ is $C^1$we deduce that ${\mathfrak K}$ is also $C^1$ and one has
\[
\begin{array}{c}
\ds\forall\ \vec{\La}=(\vec{\Phi},\vec{n})\in {\mathfrak E}_{\Sigma}^{2,4}\quad \forall \ (\vec{V},\vec{W})\in W^{2,4}(\Sigma,{\R}^4\times {\R}^4)\ ;\quad\vec{V}\cdot\vec{\Phi}=0\ ,\ \vec{W}\cdot\vec{n}=0\ \mbox{ and }\ \vec{\Phi}\cdot\vec{W}+\vec{n}\cdot\vec{V}=0 \\[5mm]
\ds\quad d{\mathfrak K}_{\vec{\La}}\cdot(\vec{V},\vec{W})=\lf(2\,d\vec{V}\dot{\wedge}d\vec{n}+2\,d\vec{\Phi}\dot{\wedge} d\vec{W},\Pi_0(d{\mathfrak H}_{\vec{\La}}\cdot(\vec{V},\vec{W})),2\,\int_{\Gamma_1} \vec{V}\cdot d\vec{n}-\vec{W}\cdot d\vec{\Phi}\cdots 2\,\int_{\Gamma_{2\,g(\Sigma)}} \vec{V}\cdot d\vec{n}-\vec{W}\cdot d\vec{\Phi}\rg) \ .
\end{array}
\]
Let $(\Om,\phi,s_1\cdots s_{2\,g(\Sigma)})\in W^{1,4}_0(\wedge^2\Sigma)\times W^{2,4}(\Sigma)\times{\R}^{2\,g(\Sigma)}$ where $W^{1,4}_0(\wedge^2\Sigma)$ is the space of $W^{2,4}-$forms on $\Sigma$ with integral 0. Assuming $\vec{\La}^\ast\al=0$, we are looking for $(\vec{V},\vec{W})\in W^{2,4}(\vec{\La}^{-1}TV_2({\R}^4))$ such that
\[
 d{\mathfrak K}_{\vec{\La}}\cdot(\vec{V},\vec{W})=(\Om,\phi,s_1\cdots s_{2\,g(\Sigma)})\ .
\]
For any $\vec{\La}=(\vec{\Phi},\vec{n})\in W^{2,4}_{imm}(\Sigma,V_2({\R}^4))$ we introduce the corresponding Gauss Map
\[
\vec{T}:=\frac{\p_{x_1}\vec{\La}\wedge\p_{x_2}\vec{\La}}{|\p_{x_1}\vec{\La}\wedge\p_{x_2}\vec{\La}|}\in W^{1,4}(\Sigma,G_2(TV_2({\R}^4)))\hookrightarrow C^{0,1/2}(\Sigma,G_2(TV_2({\R}^4)))
\]
{\bf Claim :} {\it There exists a radius $r>0$ depending only on the $W^{1,4}$ norm of $\vec{T}$ such that on any ball\footnote{We have fixed a reference smooth metric $g_0$ on $\Sigma$ with respect to which the spaces $W^{k,p}(\Sigma)$ are defined and we take the geodesic balls with respect to this metric.} $B_r(x_0)$ for $x_0$ in $\Sigma$ the image of $B_r(x_0)$ by $\vec{\La}$ is injective and described by a local graph in $V_2({\R}^4)$ : that is there exists a local chart $(y_1,y_2,y_3,y_4,y_5)$ in an open set ${\mathcal O}$ containing $\vec{\La}(x)$ such that
\[
\vec{\La}(B_r(x))=\lf\{(y_1,y_2, f(y_1,y_2), g(y_1,y_2),h(y_1,y_2))\in\,y(\om)\ ;\ (y_1,y_2)\in U\rg\}
\]
where $U$ is an open set in ${\R}^2$. Moreover $f,g,h\in W^{2,4}(U,{\R})$.} 

\medskip

\noindent{\bf Proof of the claim :} We can produce this chart by taking smooth coordinates $(y_1,y_2)$ in the tangent plate $\vec{\La}_\ast T_{x_0}\Sigma$ with $y_1(\vec{\La}(x))=0$ and $y_2(\vec{\La}(x))=0$ that we complete in order to have a local chart inn the neighborhood of $\vec{\La}(x)$.
Using the local inversion theorem we have that the map which assigns to any point $x\in B_r(x_0)$ the two first coordinates $\ti{\La}(x):=({\La}_1(x),\La_2(x))$ of $\vec{\La}(x)$ realizes a $C^1$ diffeomorphism between $B_r(x)$ and $U$, moreover since $\mbox{det}(d\ti{\La})\in C^0(B_r(x))$ and since $\mbox{det}(d\ti{\La})(x_0)>0$, for $r>0$ small enough we have the existence of $C>1$ such that
\[
\forall x\in B_r(x_0)\quad C^{-1}\le \mbox{det}(d\ti{\La})(x)\le C
\]
Moreover using the chain rule in Sobolev space we have
\[
\p^2_{y_iy_j}\ti{\La}^{-1}=\p_{y_i}\lf[ (d\ti{\La})^{-1}_j\circ\ti{\La}^{-1}(y)\rg]=\p_{x_l}(d\ti{\La})^{-1}_j\ \p_{y_j}(\ti{\La}^{-1})_l(y)\ .
\]
Hence
\[
\int_{B_r(x)}|\p^2_{y_iy_j}\ti{\La}^{-1}|^4(\ti{\La}(x))\ \mbox{det}(d\ti{\La})(x)\ dx^2\le C\, \int_{B_r(x)}|\nabla_x(d\ti{\La})^{-1}|^4(x)\ \|d(\ti{\La})^{-1}\|^4_\infty dx^2<+\infty
\] 
We deduce using the coarea formula that  $\ti{\La}^{-1}\in W^{2,4}(U,B_r(x))$ and we take respectively
\[
f(y_1,y_2):=\La_3\circ\ti{\La}^{-1}(y_1,y_2)\ ,\ g(y_1,y_2):=\La_4\circ\ti{\La}^{-1}(y_1,y_2)\ \mbox{ and }\ h(y_1,y_2):=\La_5\circ\ti{\La}^{-1}(y_1,y_2)\ .
\]
Using Chain rules (successively the classical one of $C^1$ functions and then the one in Sobolev Spaces) we have for any $m=1,2,3$ and $i,j=1,$
\[
\p^2_{y_iy_j}[\La_m\circ\ti{\La}^{-1}]=\p^2_{x_lx_k}\La_m\circ\ti{\La}^{-1}\ \p_{y_i}\ti{\La}_l^{-1}\ \p_{y_j}\ti{\La}_k^{-1}+\p_{x_l}\La_m\circ\ti{\La}^{-1}\ \p^2_{y_iy_j}\ti{\La}_l^{-1}
\]
and using again the co-area formula we have that $f,g,h$ are all three in $W^{2,4}(U)$. This ends the proof of the claim.

\medskip

We fix a finite covering of $\Sigma$ by balls $(B_{r}(x_i))_{i=1\cdots Q}$ and we consider an adapted partition of unity made of $C^\infty$ non negative functions $\chi_i$ with support in $B_r(x_i)$
\[
\sum_{i=1}^Q\chi_i\equiv 1
\]
There exists $\beta\in W^{2,4}(\wedge^1\Sigma)$ such that $d\beta=\Om$. Let $(h_k^0)_{k=1\cdots 2\, g(\Sigma)}$ be a basis of harmonic 1 forms for $g_0$ generating $H^1(\Sigma)$. These forms are known to b smooth. We modify possibly $\beta$ by adding the linear combination of $h_k^0$ such that
\be
\label{norm-harm}
\forall k=1\cdots 2\, g(\Sigma)\quad\quad\int_{\Gamma_k}\beta= s_k\ .
\ee
We write $\Om_i:=d(\chi_i\,\beta)\in W^{1,4}_0(\wedge^1B_r(x_i))$ and
\[
\Om_i= \om_i(x)\ dx_1\wedge dx_2
\]
Hence
\[
(\ti{\La}^{-1})^\ast\Om_i=\om_i(\ti{\La}^{-1}(y))\ \mbox{det}(\nabla_y\ti{\La}^{-1})\ dy_1\wedge dy_2=\ti{\om}_i(y)\ dy_1\wedge dy_2
\]
where we have $\ti{\om}_i\in W^{1,4}_0(U)$ (using again chain rule in Sobolev spaces as well as co-area formula). Let $\gamma_i=(\ti{\La}^{-1})^\ast (\chi_i\,\beta)\in W^{1,4}_0(\wedge^1U)$. We have
\[
(\ti{\La}^{-1})^\ast\Om_i= d\gamma_i\ .
\]
In particular $\int_U(\ti{\La}^{-1})^\ast\Om_i=0$. Let $\xi_i\in W^{3,4}(U)$ solving
\[
\lf\{
\begin{array}{l}
\ds \Delta \xi_i=\ti{\om}_i(y)\quad\mbox{ in }{\mathcal D}'(U)\\[3mm]
\ds \frac{\p\xi_i}{\p\nu}=0\quad \mbox{ on }\p U\ .
\end{array}
\rg.
\]
We have in particular $d(\ast d\xi_i)=(\ti{\La}^{-1})^\ast\Om_i$ and $\ast d\xi_i$ can  continuously be extended by $0$ outside $U$. Observe that on $U$ we have the existence of $\phi_i$ such that
\[
\chi_i\circ(\ti{\La}^{-1})\ (\ti{\La}^{-1})^\ast\beta=\ast d\xi_i+ d\phi_i
\]
On $U\setminus \mbox{supp}(\chi_i\circ(\ti{\La}^{-1}))$ both $\xi_i$ and $\phi_i$ are harmonic. We take a $C^\infty$ regularisation $\ov{\phi}_i$ of $\phi_i$ such that $\ov{\phi}_i-\phi_i$ is equal to zero
in an open neighbourhood of $\p U$ in $U$ and we have that $\sigma_i:=\ast d\xi_i+ d\ov{\phi}_i$ is $W^{2,4}$ in $U$ with support strictly inside $U$ and satisfies
\[
d\sigma_i=\Om_i\ .
\]
Recall that on $V_2({\R}^4)$ we have
\[
\al\wedge d\al\wedge d\al>0
\]
We have seen that the restriction to the Horizontal plane $H$ of $d\al$ is non degenerate and that
\[
\mbox{Ker}(d\al)(\vec{a},\vec{b})=\mbox{Span}\lf\{ (\vec{b},-\vec{a})\rg\}=\mbox{Span}\lf\{\vec{R}\rg\}
\]
while in the opposite way
\[
\mbox{Ker}(\al)(\vec{a},\vec{b})=H_{(\vec{a},\vec{b})}\quad\mbox{ and }\quad\al(\vec{R})=-2\ .
\]
Hence for any $(\vec{a},\vec{b})\in V_2({\R}^4)$
\[
\Xi_{\vec{a},\vec{b}}\ : \ \vec{X}=(\vec{V},\vec{W})\ \longrightarrow d\al\res\vec{X}+\al(\vec{X})\, \al
\]
realises an isomorphism from $T_{(\vec{a},\vec{b})}V_2({\R}^4)$ into $T^\ast_{(\vec{a},\vec{b})}V_2({\R}^4)$. 
The explicit form of $\Xi^{-1}$
is the following : Let ${\gamma}\in \Gamma(T^\ast V_2({\R}^4))$  we take
\[
\vec{X}=\vec{X}^H+\vec{X}^V
\]
where
\[
\lf\{
\begin{array}{l}
\ds\vec{X}^V(\vec{a},\vec{b}):= 2^{-1}\,<\al,\gamma>\ \vec{R}(\vec{a},\vec{b})\\[5mm]
\ds\vec{X}^{\,H}(\vec{a},\vec{b})\in H_{(\vec{a},\vec{b})} \quad\mbox{s. t. }\quad\forall\,\vec{Y}\in H_{(\vec{a},\vec{b})}\quad \vec{Y}\cdot\vec{X}^H:=d\al(\vec{X}^H,J_H\vec{Y})=\gamma(\vec{Y})\ .
\end{array}
\rg.
\]
where we have used that on $H$ the restriction of the two form $d\al$ is non degenerate $d\al\wedge d\al>0$ on $H$. The explicit forms above imply $\Xi_{(\vec{a},\vec{b})}$ as a map from $\Gamma(TV_2({\R}^4))$ into  $\Gamma(T^\ast V_2({\R}^4))$  as well as $\Xi^{-1}_{(\vec{a},\vec{b})}$ are obviously smooth. 

\medskip

We denote
\[
\vec{X}_i\circ y^{-1}=(\vec{V}_i\circ y^{-1},\vec{W}_i\circ y^{-1})=\Xi^{-1}_{y^{-1}(y_1,y_2,y_3,y_4,y_5)}( y^{\ast}\,\sigma_i)
\]
It satisfies
\[
y^{\ast}\sigma_i=2 \,\lf[\vec{V}_i\cdot d\vec{b}-\vec{W}_i\cdot d\vec{a}\rg]+ (\vec{a}\cdot\vec{W}_i-\vec{b}\cdot\vec{V}_i)\ [\vec{a}\cdot d\vec{b}-\vec{b}\cdot d\vec{a}]
\]
The restriction of $\vec{X}_i$ to $y\circ\vec{\La}(B_r(x_i))$ satisfies
\[
\vec{X}_i\circ y^{-1}=(\vec{V}_i\circ y^{-1},\vec{W}_i\circ y^{-1})=\Xi^{-1}_{y^{-1}(y_1,y_2,f(y_1,y_2),g(y_1,y_2),h(y_1,y_2))}( y^{\ast}\,\sigma_i)
\]
Since $y$ and $y^{-1}$ are smooth, since $f,g,h$ are in $W^{2,4}(U)$ and since $\sigma_i\in W^{2,4}(\wedge^1U)$, the chain rule for Sobolev maps combined with the area formula as above
gives that the extensions by zero outside $B_r(x_i)$ of $(\vec{V}_i\circ \vec{\La},\vec{W}_i\circ \vec{\La})$ belong to $W^{2,4}( \vec{\La}^{-1}TV_2({\R}^4))$.
By restricting this identity on $\vec{\La}(B_r(x_i))$, assuming $\vec{\La}$ is Legendrian we obtain
\[
\Om_i=d\lf[\vec{\La}^\ast y^{\ast}\sigma_i\rg]=2 \,\lf[d[\vec{V}_i\circ{\vec{\La}}]\dot{\wedge} d\vec{n}-d[\vec{W}_i\circ{\vec{\La}}]\dot{\wedge} d\vec{\Phi}\rg]
\]
Summing over $i$ gives $$(\vec{V},\vec{W})=\sum_i(\vec{V}_i\circ \vec{\La},\vec{W}_i\circ \vec{\La})\ \mbox{ in }\ W^{2,4}( \vec{\La}^{-1}TV_2({\R}^4))$$ and we have
\[
\Om=2 \,\lf[d\vec{V}\dot{\wedge} d\vec{n}-d\vec{W}\dot{\wedge} d\vec{\Phi}\rg]\ .
\]
Observe that
\[
\sigma_i=\chi_i\circ(\ti{\La}^{-1})\ (\ti{\La}^{-1})^\ast\beta+d(\phi_i-\ov{\phi}_i)
\]
where $\phi_i-\ov{\phi}_i$ is supported i the open set $U_i$. Hence each $(\phi_i-\ov{\phi}_i)\circ \ti{\La}$ extends by 0 outside $B_r(x_i)$ as a $W^{1,4}$ function on the whole $\Sigma$ and we have
\[
2\,\lf[\vec{V}\cdot d\vec{n}-\vec{W}\cdot d\vec{\Phi}\rg]=\beta +d\lf[\sum_{i}\phi_i-\ov{\phi}_i\rg]
\]
Recall that
\[
d{\mathfrak H}_{\vec{\La}}\cdot(\vec{V},\vec{W})=2\, [\vec{V}\cdot d\vec{n}-\vec{W}\cdot d\vec{\Phi}]+d\lf[\vec{\Phi}\cdot \vec{W}-\vec{n}\cdot \vec{V}\rg]
\]
Observe moreover that for any $f\in W^{2,4}(\Sigma)$ there holds
\[
d{\mathfrak H}_{\vec{\La}}\cdot(- f\,\vec{n}, f\,\vec{\Phi})=2\, df\ .
\]
Consider then $f$ such that 
\[
2\,df=d^{\ast_{g_0}}\Delta^{-1}_{g_0}\lf[d{\mathfrak H}_{\vec{\La}}\cdot(\vec{V},\vec{W})\rg]-d\phi\in W^{1,4}(\wedge^1\Sigma)\ .
\]
This gives
\[
d^{\ast_{g_0}}\Delta^{-1}_{g_0}\lf[d{\mathfrak H}_{\vec{\La}}\cdot(\vec{V}+f\vec{n},\vec{W}-f\vec{\Phi})\rg]=d\phi\ ,
\]
where we have used the fact that
\[
d^{\ast_{g_0}}\Delta^{-1}_{g_0}\lf[df\rg]=df\ .
\]
Observe moreover
\[
2 \,\lf[d[\vec{V}+f\vec{n}]\dot{\wedge} d\vec{n}-d[\vec{W}-f\vec{\Phi}]\dot{\wedge} d\vec{\Phi}\rg]=2 \,\lf[d\vec{V}\dot{\wedge} d\vec{n}-d\vec{W}\dot{\wedge} d\vec{\Phi}\rg]=\Om
\]
and
\[
\forall k=1\cdots 2g(\Sigma)\quad\quad 2\,\int_{\Gamma_k}(\vec{V}+f\vec{n})\cdot d\vec{n}-(\vec{W}-f\vec{\Phi})\cdot d\vec{\Phi}=s_k\ .
\]
Hence we have proved that if $\vec{\La}$ is a $W^{2,4}$ Legendrian immersion into $V_2({\R}^4)$ we have 
\[
\begin{array}{c}
\ds\forall\, (\Om,\phi,s_1,\cdots, s_{2\,g(\Sigma)})\in W^{1,4}(\wedge^2\Sigma)\times W^{2,4}(\Sigma)\times {\R}^{2\,g(\Sigma)}\quad\\[5mm] \ds\exists\, (\vec{V},\vec{W})\in \Gamma_{W^{2,4}}(\vec{\La}^{-1}T V_2({\R}^4))
\quad\quad d{\mathfrak K}_{\vec{\La}}(\vec{V},\vec{W})=(\Om,\phi,s_1,\cdots, s_{2\,g(\Sigma)})\ .
\end{array}
\]
This implies that ${\mathfrak K}$ realises s submersion at every point  $\vec{\La}\in {\mathfrak E}^{2,4}_\Sigma$ satisfying $\vec{\La}^\ast\al=0$. 

We claim that at a point $\vec{\La}\in {\mathfrak E}^{2,4}_\Sigma$ satisfying $\vec{\La}^\ast\al=0$ the Kernel of $d{\mathfrak K}_{\vec{\La}}$ is closed in $\Gamma_{W^{2,4}}(\vec{\La}^{-1}T V_2({\R}^4))$ and splits, that is there exists a closed supplement of Ker$\,d{\mathfrak K}_{\vec{\La}}$. We have
\[
(\vec{V},\vec{W})\in \mbox{Ker}\,d{\mathfrak K}_{\vec{\La}}\quad\Longleftrightarrow\quad \vec{V}\cdot d\vec{n}-\vec{W}\cdot d\vec{\Phi}=0\ \mbox{ on }\Sigma\quad\ .
\]
The Splitting of Ker$\,d{\mathfrak K}_{\vec{\La}}$ can be seen by taking the intersection of its orthogonal complement in $L^2(\Sigma,{\R}^8)$ with $\Gamma_{W^{2,4}}(\vec{\La}^{-1}T V_2({\R}^4))$

\medskip

Observe that 
\[
{\mathfrak K}(\vec{\La})=0\quad\quad\Longleftrightarrow\quad\quad \vec{\La}^\ast\al=0\ .
\]
Hence we have proved that the space of $W^{2,4}$ Legendrian immersions of a surface $\Sigma$ is a $C^2$ sub-manifold of ${\mathfrak E}^{2,4}_\Sigma$ and this concludes the proof of theorem~\ref{th-subman}.\hfill $\Box$

\medskip

We equip the tangent bundle of ${\mathfrak E}^{2,4}_{\Sigma,Leg}$ with the Finsler structure introduced in \cite{Riv-1} (Proposition II.2).
\[
\forall\ \vec{w}\in T_{\vec{\La}}{\mathfrak E}^{2,4}_{\Sigma,Leg}\quad\quad\|\vec{w}\|_{\vec{\La}}:=\lf[\int_\Sigma\lf[|\nabla^2\vec{w}|^2+|\nabla\vec{w}|^2+|\vec{w}|^2\rg]^2\ dvol_{g_{\vec{\La}}}\rg]^{1/4}+\||\nabla\vec{w}|_{g_{\La}}\|_{L^\infty(\Sigma)}\ .
\] 
The restriction to ${\mathfrak E}^{2,4}_{\Sigma,Leg}$ of the same Finsler structure is inducing a Palais distance $d_{Leg}$ which is clearly larger than the restriction to ${\mathfrak E}^{2,4}_{\Sigma,Leg}$ of the Palais distance $d$ on ${\mathfrak E}^{2,4}_{\Sigma}$. Hence since this last one is complete on ${\mathfrak E}^{2,4}_{\Sigma}$ (see  \cite{Riv-1}  proposition II.3) and since the membership to ${\mathfrak E}^{2,4}_{\Sigma,Leg}$ passes to the limit under distance control, $({\mathfrak E}^{2,4}_{\Sigma,Leg},d_{Leg})$ is also complete (i.e. Cauchy sequences for $d_{Leg}$ are obviously Cauchy for $d$).

\medskip

We recall the definition of {\it homotopic admissible families} or simply {\it admissible family}.

\begin{Dfi}
\label{df-admi}
Let ${\mathfrak M}$ be a  Banach manifold. A subset ${\mathcal A}$ of the power set of ${\mathfrak M}$ (i.e. $A$ is a set of subsets of ${\mathfrak M}$) is called admissible family if it is stable under isotopy that is
: forall $H\in C^0([0,1]\times {\mathfrak M},{\mathfrak M})$ such that $H_t$ is an homeomorphism of ${\mathfrak M}$ for any $t\in [0,1]$ and for any $A\in {\mathcal A}$ we have $H_1(A)\in{\mathcal A}$.
\end{Dfi}

We recall the proposition which was central in the viscosity method introduced in \cite{Riv-1} and more precisely the version given in \cite{Piga} (which permits to avoid the use of the Palais-Smale condition which should hold for our problem but whose proof might be a bit tedious in the present context).

\begin{Prop}
\label{pr-stmono}
Let $E_\ep$ be a family of $C^1$ functionals on a complete Finsler manifold $\mathfrak M$, with $E_\ep$ differentiable with respect to $\ep$ and $\ep\rightarrow E_\ep(x)$, $\ep\rightarrow \frac{dE_\ep(x)}{d\ep}$ both increasing in $\ep$ for any $x\in {\mathfrak M}$. Assume also that
\be
\label{cond-der}
\|dE_{\ep_i}(x_i)-dE_\ep(x_i)\|_{T^\ast{\mathfrak M}}\longrightarrow 0\ ,
\ee
whenever $0<\ep\le\p_i$,  $\ep_i\rightarrow \ep$ and $\limsup_{i\rightarrow +\infty} E_\ep(x_i)<+\infty$.

Then for any admissible family ${\mathcal A}$ we define the corresponding minmax value (or width) for $E_\ep$
\[
\beta(\ep):=\inf_{A\in {\mathcal A}}\ \sup_{x\in A}E_\ep(x)
\]
there exists a sequence $\ep_k\rightarrow 0$ and $x_k\in{\mathfrak M}$ such that
\[
E_{\ep_k}(x_k)-\beta(\ep_k)\rightarrow 0,\quad\|dE_{\ep_k}(x_k)\|_{T^\ast{\mathfrak M}}<f(\ep_k)\ ,\quad \ep_k\,\log\ep_k^{-1}\, \lf.\frac{dE_\ep(x_k)}{d\ep}\rg|_{\ep=\ep_k}\rightarrow 0
\]
where $f$ is any function fixed in advance.
\end{Prop}

\subsection{Local Geometry of Legendrian Immersions}
\subsubsection{Infinitesimal variations of Legendrian Immersions}
Let $\vec{\La}=(\vec{\Phi},\vec{n})$ be a $C^2$ immersion of a closed oriented surface $\Sigma$ into $V_2({\R}^4)$. We shall restrict to Legendrian immersions which are the immersions satisfying 
\[
\forall\, x\in \Sigma\quad\quad\vec{\La}_\ast T_x\Sigma\subset H_{\vec{\La}(x)}\ .
\]
Let $\Psi_t$ be the flow of a vector-field $X$ in $V_2({\R}^4)$ such that $\Psi_t$ is preserving the kernel of $\al$ for every $t$ (i.e. is a contacto-morphism). We have then for all $t$
\[
\quad\mbox{Ker}(\Psi_t^\ast\al)=\mbox{Ker}(\al)\ .
\]
Two non zero one forms have the same kernel if and only if they are parallel to each other. Hence, since $|\al|=\sqrt{2}$ at every point one has 
\[
\Psi_t^\ast\al= f(t)\ \al\quad\mbox{ where }\quad 2\,f(t):=<\al,\Psi_t^\ast\al>\ .
\] 
Hence
\[
\frac{d \Psi_t^\ast\al}{dt}={\mathcal L}_X\al=\dot{f}\,\al\ .
\]
Thus
\[
d(\al(X))+d\al\res X=\dot{f}\,\al
\]
We decompose $\vec{X}=\vec{X}^V+\vec{X}^H$ where $\vec{X}^H$ is the orthogonal projection of $\vec{X}$ onto $H$. This gives in particular
\[
d(\al(\vec{X}^V))+(\Pi^\ast \om_{{\C}P^1\times \ov{{\C}P^1}})\res \vec{X}^H=\dot{f}\,\al
\]
Denote $h:=\al(\vec{X}^V)$ that is
\[
\vec{X}^V=\frac{h}{2}\ \vec{R}\ ,
\]
and let $\vec{Y}=\vec{Y}^V+\vec{Y}^H$ we have
\[
d_{\vec{Y}}h+\vec{Y}^H\cdot J_H\vec{X}^H=\dot{f}\ \al(\vec{Y}^V)\ .
\]
This gives
\[
\vec{X}^H=J_H\nabla^H h\ ,
\]
where $\nabla^H h$ is the horizontal projection of the gradient of $h$. To conclude, infinitesimal variations preserving the Legendrian constraint are given by vector fields of the form
\[
\vec{X}_h:=J_H\nabla^H h+\frac{h}{2}\ \vec{R}\ ,
\]
where $h$ is an arbitrary function (i.e. Hamiltonian).

\medskip

The projection of a Legendrian immersion is satisfying $\vec{\La}^\ast \al=0$ hence
\[
0=d\vec{\La}^\ast \al=\vec{\La}^\ast \Pi^\ast\om_{{\C}P^1\times \ov{{\C}P^1}}=(\Pi\circ\vec{\La})^\ast \om_{{\C}P^1\times \ov{{\C}P^1}}
\]
hence $\vec{G}:=\Pi\circ\vec{\La}$ is Lagrangian into ${\C}P^1\times \ov{{\C}P^1}$. If $\vec{e}_1,\vec{e}_2$ is an orthonormal basis of $\vec{G}_\ast T\Sigma$ we have since the restriction of $\om_{{\C}P^1\times \ov{{\C}P^1}}$ to $\vec{G}_\ast T\Sigma=0$ we have
\[
\vec{e}_2\cdot J\vec{e}_1=-\vec{e}_1\cdot J\vec{e}_2=\om(\vec{e}_1,\vec{e}_2)=0
\]
hence $J\vec{e}_1,J,\vec{e}_2$ realises an orthonormal basis of the normal bundle to the immersion $\vec{G}$. By uniqueness of the horizontal lifts and the fact that $\Pi_\ast$ realises an isometry
between $H_(\vec{a},\vec{b})$ and $T_{\Pi(\vec{a},\vec{b})}{\C}P^1\times \ov{{\C}P^1}$ the same holds for the horizontal lifts as well and $(J_H\vec{e}_1^{\,H},J_H\vec{e}_2^{\,H})$ realises an orthonormal basis of the normal plane to $\vec{\La}_\ast T\Sigma$ in $H$.

Let $\vec{c}\in \mbox{Span}\{\vec{a},\vec{b}\}^\perp$. We have
\[
\Pi_\ast(\vec{c},0)=\frac{d}{dt}\lf[ \Pi\circ(\vec{a}+t\,\vec{c},\vec{b})  \rg]=\frac{d}{dt}\lf( ((\vec{a}+t\vec{c})\wedge\vec{b})^+, ((\vec{a}+t\vec{c})\wedge\vec{b})^- \rg)=\lf((\vec{c}\wedge\vec{b})^+,(\vec{c}\wedge\vec{b})^-\rg)\ .
\]
Hence
\[
\Pi_\ast J_H(\vec{c},0)=\lf(J_{\vec{G}^+}((\vec{c}\wedge\vec{b})^+),J_{\vec{G}^-}((\vec{c}\wedge\vec{b})^-)\rg)=\lf((\vec{d}\wedge\vec{b})^+,-(\vec{d}\wedge\vec{b})^-\rg)=\lf( (\vec{a}\wedge\vec{c})^+,(\vec{a}\wedge\vec{c})^- \rg)=\Pi_\ast(0,\vec{c})
\]
where
\[
\vec{d}:=\star\lf[\vec{a}\wedge\vec{b}\wedge\vec{c}\rg]\ .
\]
Since $J_H^2=-I_H$
\[
 J_H(0,\vec{c})=-(\vec{c},0)
\]


\subsubsection{The second fundamental Form of a Legendrian immersion}

We assume that $\vec{\La}=(\vec{\Phi},\vec{n})$ is conformal Legendrian from the unit disc $D^2$ into $V_2({\R}^4)$. Hence the unit Gauss Map is a 3-vector given by
\[
\vec{N}:=\ \frac{e^{-2\la}}{\sqrt{2}}\,\vec{R}(\vec{\La})\wedge J_H\,\frac{\p \vec{\La}}{\p x_1}\wedge\ J_H\,\frac{\p \vec{\La}}{\p x_2}\ ,
\]
where
\[
e^{2\la}=\lf|\frac{\p \vec{\La}}{\p x_1}\rg|^2=\lf|\frac{\p \vec{\La}}{\p x_2}\rg|^2\ .
\]
We see $V_2({\R}^4)$ canonically isometrically embedded in $S^3\times S^3\hookrightarrow {\R}^4\times {\R}^4$. The second fundamental form of the immersion $\vec{\La}$ is given by
\[
\vec{\mathbb I}_{ij}\ dx_i\otimes dx_j:=\pi_{\vec{N}}\p^2_{x_ix_j}\vec{\La}\ dx_i\otimes dx_j
\]
where $\pi_{\vec{N}}$ is the orthogonal projection onto the normal space to $\vec{\La}_\ast T\Sigma$ in $TV_2({\R}^4)$. Assuming as just above that $\vec{\La}$ is in local conformal coordinates, this normal space is generated by the
orthonormal basis
\[
\sqrt{2}^{-1}\,\vec{R}(\vec{\La})\ ,\ e^{-\la}\,J_H\p_{x_1}\vec{\La}\ ,\ e^{-\la}\,J_H\p_{x_2}\vec{\La}\ .
\] 
Hence we have in these conformal coordinates
\[
\vec{\mathbb I}_{ij}= e^{-2\la}\ J_H\p_{x_1}\vec{\La}\cdot\p^2_{x_ix_j}\vec{\La}\ J_H\p_{x_1}\vec{\La}+e^{-2\la}\ J_H\p_{x_2}\vec{\La}\cdot\p^2_{x_ix_j}\vec{\La}\ J_H\p_{x_2}\vec{\La}+
\frac{1}{2}\, \vec{R}\cdot\p^2_{x_ix_j}\vec{\La}\ \vec{R}
\]
Observe that the Legendrian condition is implying
\be
\label{leg-cond}
\p_{x_k}\p_{x_l} \vec{\La}\cdot\vec{R}=\p_{x_k}\p_{x_l} \vec{\Phi}\cdot\vec{n}-\p_{x_k}\p_{x_l} \vec{n}\cdot\vec{\Phi}=\p_{x_k}\lf[\p_{x_l} \vec{\Phi}\cdot\vec{n}-\p_{x_l} \vec{n}\cdot\vec{\Phi}\rg]=0\ .
\ee
Hence
\[
\vec{\mathbb I}_{ij}= e^{-2\la}\ J_H\p_{x_1}\vec{\La}\cdot\p^2_{x_ix_j}\vec{\La}\ J_H\p_{x_1}\vec{\La}+e^{-2\la}\ J_H\p_{x_2}\vec{\La}\cdot\p^2_{x_ix_j}\vec{\La}\ J_H\p_{x_2}\vec{\La}
\]
In particular
\[
|\vec{\mathbb I}_{ij}|^2=\lf|\pi_T\lf[ J_H\nabla_{x_i}\p_{x_j}\vec{\La}\rg]\rg|^2=\lf|\pi_{\vec{N}}\p^2_{x_ix_j}\vec{\La}\rg|^2\ ,
\]
where $\pi_T$ is the orthogonal projection onto $\vec{\La}_\ast T\Sigma$,  $\nabla$ is the covariant derivative on $V_2({\R}^4)$ and we denote $\nabla_{x_k}:=\nabla_{\p_{x_k}\vec{\La}}$. We have
in particular
\[
|\vec{\mathbb I}|^2_g=e^{-4\la}\, \sum_{i,j=1}^2\lf|\pi_{\vec{N}}\p^2_{x_ix_j}\vec{\La}\rg|^2\ .
\]
\[
\begin{array}{l}
\ds\nabla_{x_i}\vec{N}=\ \frac{1}{\sqrt{2}}\,\nabla_{x_i}[\vec{R}(\vec{\La})]\wedge J_H\,e^{-\la}\frac{\p \vec{\La}}{\p x_1}\wedge\ J_H\,e^{-\la}\frac{\p \vec{\La}}{\p x_2}\\[5mm]
\ds\quad +\frac{1}{\sqrt{2}}\,\vec{R}(\vec{\La})\wedge \nabla_{x_i}\lf[J_H\,e^{-\la}\frac{\p \vec{\La}}{\p x_1}\rg]\wedge\ J_H\,e^{-\la}\frac{\p \vec{\La}}{\p x_2}\\[5mm]
\ds\quad +\frac{1}{\sqrt{2}}\,\vec{R}(\vec{\La})\wedge J_H\,e^{-\la}\frac{\p \vec{\La}}{\p x_1}\wedge\ \nabla_{x_i}\lf[J_H\,e^{-\la}\frac{\p \vec{\La}}{\p x_2}\rg]\\[5mm]
\end{array}
\]
We have first for $i,k=1,2$
\[
\nabla_{x_i}[\vec{R}(\vec{\La})]\cdot \vec{R}(\vec{\La})=0\quad,\quad \nabla_{x_i}[\vec{R}(\vec{\La})]\cdot \p_{x_k}\vec{\La}=\p_{x_k}\vec{\Phi}\cdot\p_{x_i}\vec{n}-\p_{x_i}\vec{\Phi}\cdot\p_{x_k}\vec{n}=0\ ,
\]
hence
\[
\nabla_{x_i}[\vec{R}(\vec{\La})]\wedge J_H\,e^{-\la}\frac{\p \vec{\La}}{\p x_1}\wedge\ J_H\,e^{-\la}\frac{\p \vec{\La}}{\p x_2}=0\ .
\]
We have also for $i,k=1,2$
\[
\nabla_{x_i}\lf[J_H\,e^{-\la}\frac{\p \vec{\La}}{\p x_k}\rg]\cdot \p_{x_l}\vec{\La}=-J_H\,e^{-\la}\frac{\p \vec{\La}}{\p x_k}\cdot \p^2_{x_lx_i}\vec{\La}
\]
Thus
\[
\begin{array}{l}
\ds\nabla_{x_i}\vec{N}=-\sum_{j=1}^2\frac{e^{-2\la}}{\sqrt{2}}\,J_H\,e^{-\la}\frac{\p \vec{\La}}{\p x_1}\cdot \p^2_{x_jx_i}\vec{\La}\ \vec{R}(\vec{\La})\wedge\p_{x_l}\vec{\La}\wedge\ J_H\,e^{-\la}\frac{\p \vec{\La}}{\p x_2}\\[5mm]
\ds-\sum_{j=1}^2\frac{e^{-2\la}}{\sqrt{2}}\,J_H\,e^{-\la}\frac{\p \vec{\La}}{\p x_2}\cdot \p^2_{x_jx_i}\vec{\La}\ \vec{R}(\vec{\La})\wedge\ J_H\,e^{-\la}\frac{\p \vec{\La}}{\p x_1}\wedge \p_{x_l}\vec{\La}
\end{array}
\]
This is implying in particular
\[
\sum_{i=1}^2|\nabla_{x_i}\vec{N}|^2= e^{-2\la}\ \sum_{i,j=1}^2|\pi_{\vec{N}}\p^2_{x_ix_j}\vec{\La}|^2\ .
\]
We recall that the  mean curvature vector is given by
\[
\vec{H}_{\vec{\La}}:=\frac{1}{2} \mbox{Tr}_g{\mathbb I}=\frac{e^{-2\la}}{2}\pi_{\vec{N}}\lf[\p^2_{x_1^2}\vec{\La}+\p^2_{x_2^2}\vec{\La}\rg]\ .
\]
Hence
\[
2\,\vec{H}_{\vec{\La}}=e^{-4\la}\ J_H\p_{x_1}\vec{\La}\cdot\Delta\vec{\La}\ J_H\p_{x_1}\vec{\La}+e^{-4\la}\ J_H\p_{x_2}\vec{\La}\cdot\Delta\vec{\La}\ J_H\p_{x_2}\vec{\La}\]
Let
\[
\gamma:=e^{-2\la}\ J_H\p_{x_1}\vec{\La}\cdot\Delta\vec{\La}\ dx_1+e^{-2\la}\ J_H\p_{x_2}\vec{\La}\cdot\Delta\vec{\La}\ dx_2\ .
\]
We have
\[
\begin{array}{l}
\ds d\gamma=-\p_{x_2}\lf[e^{-2\la}\ J_H\p_{x_1}\vec{\La}\cdot\Delta\vec{\La}\rg]+\p_{x_1}\lf[e^{-2\la}\ J_H\p_{x_2}\vec{\La}\cdot\Delta\vec{\La}\rg]\ dx_1\wedge dx_2\\[5mm]
\ds=\p_{x_2}\lf[e^{-2\la}\ \nabla(J_H\p_{x_1}\vec{\La})\cdot\nabla\vec{\La}\rg]-\p_{x_1}\lf[e^{-2\la}\ \nabla(J_H\p_{x_2}\vec{\La})\cdot\nabla\vec{\La}\rg]\ dx_1\wedge dx_2
\end{array}
\]
Observe that
\[
\begin{array}{l}
\ds\nabla(J_H\p_{x_i}\vec{\La})\cdot\nabla\vec{\La}=-\,\star \nabla\lf(\vec{\Phi}\wedge\vec{n}\wedge \p_{x_i}\vec{\Phi}\rg)\wedge\nabla\vec{\Phi}-\,\star\, \nabla\lf(\vec{\Phi}\wedge\vec{n}\wedge \p_{x_i}\vec{n}\rg)\wedge\nabla\vec{n}\\[5mm]
-\,\star \lf(\vec{\Phi}\wedge\vec{n}\wedge \p_{x_i}\nabla\vec{\Phi}\rg)\wedge\nabla\vec{\Phi}-\,\star\, \lf(\vec{\Phi}\wedge\vec{n}\wedge \p_{x_i}\nabla\vec{n}\rg)\wedge\nabla\vec{n}\\[5mm]
\ds\quad=J_H\nabla^H\p_{x_i}\vec{\La}\cdot\nabla\vec{\La}=-\nabla\p_{x_i}\vec{\La}\cdot J_H\nabla\vec{\La}
\end{array}
\]
Hence
\be
\label{mean-1}
\begin{array}{l}
\ds d\gamma=-\p_{x_2}\lf[e^{-2\la}\ \nabla\p_{x_1}\vec{\La}\cdot J_H\nabla\vec{\La}\rg]+\p_{x_1}\lf[e^{-2\la}\ \nabla\p_{x_2}\vec{\La}\cdot J_H\nabla\vec{\La}\rg]\ dx_1\wedge dx_2\\[5mm]
\ds \quad =\lf[2\, \p_{x_2}\la \,e^{-2\la}\ \nabla\p_{x_1}\vec{\La}\cdot J_H\nabla\vec{\La}-2\,\p_{x_1}\la\,e^{-2\la}\ \nabla\p_{x_2}\vec{\La}\cdot J_H\nabla\vec{\La}\rg]\ dx_1\wedge dx_2\\[5mm]
\ds\quad\quad-\,e^{-2\la}\ \lf[\nabla\p_{x_1}\vec{\La}\cdot \p_{x_2}(J_H\nabla\vec{\La})-\nabla\p_{x_2}\vec{\La}\cdot \p_{x_1}(J_H\nabla\vec{\La})\rg]\ dx_1\wedge dx_2
\end{array}
\ee
We have
\[
\begin{array}{l}
\ds\nabla_{x_k}\p_{x_l} \vec{\La}=e^{-2\la}\, \lf[\p_{x_k}\p_{x_l} \vec{\La}\cdot\p_{x_k}\vec{\La}\ \p_{x_k}\vec{\La}+  \p_{x_k}\p_{x_l} \vec{\La}\cdot\p_{x_{k+1}}\vec{\La}\ \p_{x_{k+1}}\vec{\La}\rg]\\[5mm]
\ds +e^{-2\la}\, \lf[\p_{x_k}\p_{x_l} \vec{\La}\cdot J_H(\vec{\La})\,\p_{x_k}\vec{\La}\ J_H(\vec{\La})\,\p_{x_k}\vec{\La}+  \p_{x_k}\p_{x_l} \vec{\La}\cdot\,J_H(\vec{\La})\p_{x_{k+1}}\vec{\La}\ J_H(\vec{\La})\,\p_{x_{k+1}}\vec{\La}\rg]\\[5mm]
\ds+2^{-1}\ \p_{x_k}\p_{x_l} \vec{\La}\cdot\vec{R}\ \vec{R}
\end{array}
\]
(\ref{leg-cond}) gives
\[
\lf\{
\begin{array}{l}
\ds\nabla_{x_k}\p_{x_k} \vec{\La}=\p_{x_k}\la\ \p_{x_k}\vec{\La}-\p_{x_{k+1}}\la\ \p_{x_{k+1}}\vec{\La}
\ds+e^{-2\la}\, \lf[\p_{x_k}\p_{x_k} \vec{\La}\cdot J_H\,\nabla \vec{\La}\rg] \,J_H\nabla\, \vec{\La}\\[5mm]
\ds\nabla_{x_k}\p_{x_{k+1}} \vec{\La}=\p_{x_{k+1}}\la\ \p_{x_k}\vec{\La}+\p_{x_{k}}\la\ \p_{x_{k+1}}\vec{\La}
\ds+e^{-2\la}\, \lf[\p_{x_k}\p_{x_{k+1}} \vec{\La}\cdot J_H\,\nabla \vec{\La}\rg] \,J_H\nabla\, \vec{\La}
\end{array}
\rg.
\]
which implies in particular
\[
\begin{array}{l}
\ds\nabla\p_{x_1}\vec{\La}=\nabla\la\ \p_{x_1}\vec{\La}+\nabla^\perp\la\ \p_{x_2}\vec{\La}+e^{-2\la}\ \nabla\p_{x_1}\vec{\La}\cdot J_H(\vec{\La})\p_{x_1}\vec{\La}\ J_H(\vec{\La})\p_{x_1}\vec{\La}\\[5mm]
\ds +e^{-2\la}\ \nabla\p_{x_1}\vec{\La}\cdot J_H(\vec{\La})\p_{x_2}\vec{\La}\ J_H\p_{x_2}\vec{\La}
\end{array}
\]
Observe also that
\[
\begin{array}{l}
J_H\p_{x_j}\vec{\La}\cdot \p_{x_k}\lf(J_H\p_{x_l}\vec{\La}\rg)=
\p_{x_j}\vec{\La}\cdot \p^2_{x_kx_l}\vec{\La}
\end{array}
\]
We deduce
\[
\nabla\p_{x_1}\vec{\La}\cdot \p_{x_2}(J_H\nabla\vec{\La})-\nabla\p_{x_2}\vec{\La}\cdot \p_{x_1}(J_H\nabla\vec{\La})=2\, \nabla\p_{x_1}\vec{\La}\cdot J_H\nabla\p_{x_2}\vec{\La}
\]
Combining the previous gives then
\[
\begin{array}{l}
\ds\nabla\p_{x_1}\vec{\La}\cdot \p_{x_2}(J_H\nabla\vec{\La})-\nabla\p_{x_2}\vec{\La}\cdot \p_{x_1}(J_H\nabla\vec{\La})=2\,\nabla\la\ \p_{x_1}\vec{\La}\cdot J_H\nabla\p_{x_2}\vec{\La}\\[5mm]
\ds+2\,\nabla^\perp\la\ \p_{x_2}\vec{\La}\cdot J_H\nabla\p_{x_2}\vec{\La}+2\,e^{-2\la}\ \nabla\p_{x_1}\vec{\La}\cdot J_H\p_{x_1}\vec{\La}\ \p_{x_1}\vec{\La}\cdot \nabla\p_{x_2}\vec{\La}\\[5mm]
\ds+2\,e^{-2\la}\ \nabla\p_{x_1}\vec{\La}\cdot J_H\p_{x_2}\vec{\La}\ \p_{x_2}\vec{\La}\cdot \nabla\p_{x_2}\vec{\La}
\end{array}
\]
Observe that
\[
e^{-2\la}\, \p_{x_2}\vec{\La}\cdot \nabla\p_{x_2}\vec{\La}=\nabla\la\quad\mbox{ and }\quad e^{-2\la}\,\p_{x_1}\vec{\La}\cdot \nabla\p_{x_2}\vec{\La}=-\nabla^\perp\la
\]
Hence
\[
\begin{array}{l}
\ds\nabla\p_{x_1}\vec{\La}\cdot \p_{x_2}(J_H(\vec{\La})\nabla\vec{\La})-\nabla\p_{x_2}\vec{\La}\cdot \p_{x_1}(J_H(\vec{\La})\nabla\vec{\La})\\[5mm]
\ds=2\,\nabla\la\,\lf[ -  J_H(\vec{\La})\p_{x_1}\vec{\La}\cdot \nabla\p_{x_2}\vec{\La}+J_H(\vec{\La})\p_{x_2}\vec{\La}\cdot \nabla^\perp\p_{x_2}\vec{\La}+J_H(\vec{\La})\p_{x_1}\vec{\La}\cdot \nabla^\perp\p_{x_1}\vec{\La} +
J_H(\vec{\La})\p_{x_2}\vec{\La}\cdot \nabla\p_{x_1}\vec{\La}\rg]

\end{array}
\]
Observe that
\[
J_H(\vec{\La})\p_{x_2}\vec{\La}\cdot \nabla\p_{x_1}\vec{\La}=\nabla\lf[ J_H(\vec{\La})\p_{x_2}\vec{\La}\cdot \p_{x_1}\vec{\La}\rg]-\nabla(J_H(\vec{\La})\p_{x_2}\vec{\La})\cdot\p_{x_1}\vec{\La}=J_H(\vec{\La})\p_{x_1}\vec{\La}\cdot \nabla\p_{x_2}\vec{\La}\ .
\]
Thus
\be
\label{mean-2}
\begin{array}{l}
\ds\nabla\p_{x_1}\vec{\La}\cdot \p_{x_2}(J_H(\vec{\La})\nabla\vec{\La})-\nabla\p_{x_2}\vec{\La}\cdot \p_{x_1}(J_H(\vec{\La})\nabla\vec{\La})\\[5mm]
\ds=2\,\nabla\la\,\lf[ J_H(\vec{\La})\p_{x_2}\vec{\La}\cdot \nabla^\perp\p_{x_2}\vec{\La}+J_H(\vec{\La})\p_{x_1}\vec{\La}\cdot \nabla^\perp\p_{x_1}\vec{\La} \rg] \\[5mm]
\ds=-\,2\,\p_{x_1}\la \lf[J_H(\vec{\La})\p_{x_2}\vec{\La}\cdot \p^2_{x_2x_2}\vec{\La}+J_H(\vec{\La})\p_{x_1}\vec{\La}\cdot \p^2_{x_1x_2}\vec{\La}\rg]+2\,\p_{x_2}\la \lf[J_H(\vec{\La})\p_{x_2}\vec{\La}\cdot \p^2_{x_1x_2}\vec{\La}+J_H(\vec{\La})\p_{x_1}\vec{\La}\cdot \p^2_{x_1x_1}\vec{\La}\rg]\\[5mm]
=-\,2\,\p_{x_1}\la \lf[J_H(\vec{\La})\nabla\vec{\La}\cdot \nabla\p_{x_2}\vec{\La}\rg]+2\,\p_{x_2}\la \lf[J_H(\vec{\La})\nabla\vec{\La}\cdot \nabla\p_{x_1}\vec{\La}\rg]\ .
\end{array}
\ee
Combining (\ref{mean-1}) and (\ref{mean-2}) we obtain that $\gamma$ is closed and we have established the following lemma which was already known since the work of Dazor \cite{Daz}.
\begin{Lm}
\label{lm-mean}
Let $\vec{\La}$ be a Legendrian immersion into $V_2({\R}^4)$ then the second fundamental form is horizontal and
\be
\label{mean-3}
d\lf(J_H(\vec{\La})\vec{H}_{\vec{\La}}\rg)^\sharp=0\ ,
\ee
where for any vector $\vec{X}\in H$ we assign the one form $\vec{X}^{\,\sharp}$ on $H$ such that 
\[
\forall\, \vec{Y}\in H\quad\quad\vec{X}^{\,\sharp}(\vec{Y}):=\vec{X}\cdot\vec{Y}\ .
\]
\end{Lm}
The  co-homology class in $\Sigma$ defined by $\lf(J_H(\vec{\La})\vec{H}_{\vec{\La}}\rg)^\sharp$ is invariant under the action of Hamiltonian isotopies (this is one of the main results in \cite{Oh3}).

\medskip

Let $d\beta=2^{-1}\,\gamma$ ($\beta$ is locally defined). We have
\[
\lf\{
\begin{array}{l}
\ds\p_{x_1}\beta=2^{-1}\,e^{-2\la}\ J_H\p_{x_1}\vec{\La}\cdot\Delta\vec{\La}\\[5mm]
\ds \p_{x_2}\beta=2^{-1}\,e^{-2\la}\ J_H\p_{x_2}\vec{\La}\cdot\Delta\vec{\La}
\end{array}
\rg.
\]
This gives
\[
J_H\lf(\p_{x_1}\beta\,\p_{x_1}\vec{\La}+\p_{x_2}\beta\,\p_{x_2}\vec{\La}  \rg)=2^{-1}\,\pi_H\lf[\Delta\vec{\La}\rg]\ .
\]
We have
\[
\nabla^\Sigma\beta=<d\beta\cdot d\vec{\La}>_g=e^{-2\la}\lf(\p_{x_1}\beta\,\p_{x_1}\vec{\La}+\p_{x_2}\beta\,\p_{x_2}\vec{\La}  \rg)=-\,2^{-1}\,J_H\,\Delta_g\vec{\La}
\]
Hence
\[
\int_{\Sigma}<d\vec{w},d\vec{\La}>_g\ dvol_g=-\int_{\Sigma}\vec{w}\cdot\Delta_g\vec{\La}\ dvol_g=2\, \int_\Sigma\, \vec{w}\cdot J_H\nabla^\Sigma\beta\ dvol_g
\]
Assume $\vec{w}$ is of the form $\vec{w}:=J^H\nabla^Hh+2^{-1}\,h\,\vec{R}$ this then gives
\be
\label{forme-aire}
\int_{\Sigma}<d\vec{w},d\vec{\La}>_g\ dvol_g=2\, \int_\Sigma\, <dh,d\beta>_g\ dvol_g\ .
\ee


\subsubsection{The variations of the Gauss unit multi-vector of a Legendrian immersion}

We consider a smooth perturbation of a conformal Legendrian immersion and we denote by $\vec{w}:={d\vec{\La}}/{dt}(0)$. We see $\vec{\La}_t$ as an immersion into $V_2({\R}^4)\hookrightarrow S^3\times S^3\hookrightarrow {\R}^4\times {\R}^4$. Hence $\vec{N}\in \wedge^3{\R}^8$. We have
\[
\frac{d \vec{R}(\vec{\La}_t)}{dt}=\frac{d}{dt}\lf(\begin{array}{c} \vec{n}_t\\[3mm]-\vec{\Phi}_t\end{array}\rg)= I\,\vec{w}
\]
where
\[
I:=\lf(\begin{array}{cc}0 &\,1 \\[3mm]-1 &0 \end{array}\rg)
\]
Let
\[
\vec{T}:=\frac{\p_{x_1}\vec{\La}\wedge\p_{x_2}\vec{\La}}{|\p_{x_1}\vec{\La}\wedge\p_{x_2}\vec{\La}|}\ .
\]
The link between $\vec{T}$ and $\vec{N}$ is made as follows. We have that
\[
\vec{\La}\cdot J_H\p_{x_k}\vec{\La}=\lf(\begin{array}{c}
\vec{\Phi}\\[3mm]
\vec{n}
\end{array}\rg)\cdot J_H\p_{x_k}\vec{\La}=\star \,\vec{\Phi}\wedge\vec{n}\wedge \p_{x_k}\vec{\Phi}\wedge\vec{\Phi}+\star\, \vec{\Phi}\wedge\vec{n}\wedge \p_{x_k}\vec{n}\wedge\vec{n}=0
\]
we have similarly
\[
\lf(\begin{array}{c}
\vec{\Phi}\\[3mm]
-\vec{n}\end{array}\rg)\cdot J_H\p_{x_k}\vec{\La}=0\quad\lf(\begin{array}{c}
\vec{\Phi}\\[3mm]
-\vec{n}\end{array}\rg)\cdot \vec{R}(\vec{\La})=0\quad\cdots
\]
Hence
\[
\sqrt{2}^{-1}\,\vec{\La}=\sqrt{2}^{-1}\,\lf(\begin{array}{c}
\vec{\Phi}\\[3mm]
\vec{n}
\end{array}\rg)\ ,\ \sqrt{2}^{-1}\, \vec{U}:= \sqrt{2}^{-1}\, \lf(\begin{array}{c}
\vec{\Phi}\\[3mm]
-\vec{n}\end{array}\rg)\ ,\   \sqrt{2}^{-1}\, \vec{V}:=\sqrt{2}^{-1}\, \lf(\begin{array}{c}
\vec{n}\\[3mm]
\vec{\Phi}\end{array}\rg)\ ,  \sqrt{2}^{-1}\, \vec{R}:=\sqrt{2}^{-1}\, \lf(\begin{array}{c}
\vec{n}\\[3mm]
-\vec{\Phi}\end{array}\rg)\ \]
together with
\[
 e^{-\la}\,\p_{x_1}\vec{\La}\quad ,\quad e^{-\la}\,\p_{x_2}\vec{\La}\quad , \quad e^{-\la}\,J_H\p_{x_1}\vec{\La}\quad \mbox{and}\quad  e^{-\la}\,J_H\p_{x_1}\vec{\La}\
\]
realises an orthonormal basis of ${\R}^8$. We deduce that
\[
\vec{N}:=\pm \,\sqrt{2}^{-3}\,\star\lf[\vec{\La}\wedge\vec{U}\wedge\vec{V}\wedge\vec{T}\rg]\ .
\]
In particular
\[
|d\vec{N}|^2_g\le 2^{-3}\, [|d\vec{T}|^2_g+2]\ .
\]
We have
\[
\begin{array}{l}
\ds\frac{d}{dt} |\p_{x_1}\vec{\La}_t\wedge\p_{x_2}\vec{\La}_t|^2=
\ds +2\, \lf[\p_{x_1}\vec{w}\wedge\p_{x_2}\vec{\La}\rg]\cdot\lf[\p_{x_1}\vec{\La}\wedge\p_{x_2}\vec{\La}\rg]+2\, \lf[\p_{x_1}\vec{\La}\wedge\p_{x_2}\vec{w}\rg]\cdot\lf[\p_{x_1}\vec{\La}\wedge\p_{x_2}\vec{\La}\rg]\\[5mm]
\ds=2\, e^{2\la}\, \lf[\p_{x_1}\vec{w}\cdot\p_{x_1}\vec{\La}+\p_{x_2}\vec{w}\cdot\p_{x_2}\vec{\La}\rg]=2\, e^{2\la}\,  \nabla\vec{w}\cdot\nabla\vec{\La}
\end{array}
\]
This is implying
\[
\begin{array}{l}
\ds\frac{d\vec{T}}{dt}={e^{-2\la}}\ \lf[\p_{x_1}\vec{w}\wedge\p_{x_2}\vec{\La}+\p_{x_1}\vec{\La}\wedge\p_{x_2}\vec{w}  \rg]-\frac{1}{2} \lf|\p_{x_1}\vec{\La}\wedge\p_{x_2}\vec{\La}\rg|^{-3}\ \frac{d}{dt} |\p_{x_1}\vec{\La}_t\wedge\p_{x_2}\vec{\La}_t|^2\ \lf[\p_{x_1}\vec{\La}\wedge\p_{x_2}\vec{\La}\rg]\\[5mm]
\ds\quad={e^{-2\la}}\ \lf[\p_{x_1}\vec{w}\wedge\p_{x_2}\vec{\La}+\p_{x_1}\vec{\La}\wedge\p_{x_2}\vec{w}  \rg]- e^{-2\la}\,\nabla\vec{w}\cdot\nabla\vec{\La}\ \vec{T}
\end{array}
\]
Now we compute the derivative with respect to $t$ of the metric of $\vec{\La}_t$.
\be
\label{metric}
\frac{dg_{ij}}{dt}=\p_{x_i}\vec{w}\cdot\p_{x_j}\vec{\La}+\p_{x_i}\vec{\La}\cdot\p_{x_j}\vec{w}\ ,
\ee
which gives
\[
\frac{dg^{ij}}{dt}=-\, e^{-4\la}\, \lf[\p_{x_i}\vec{w}\cdot\p_{x_j}\vec{\La}+\p_{x_i}\vec{\La}\cdot\p_{x_j}\vec{w}\rg]\ .
\]
Hence
\[
\begin{array}{l}
\ds\frac{d}{dt}\lf[|d\vec{T}|^2_{g}  \rg]=-\, e^{-4\la}\, \lf[\p_{x_i}\vec{w}\cdot\p_{x_j}\vec{\La}+\p_{x_i}\vec{\La}\cdot\p_{x_j}\vec{w}\rg]\ \p_{x_i}\vec{T}\,\p_{x_j}\vec{T}+2\,e^{-2\la}\,\nabla\vec{T}\cdot\nabla\frac{d\vec{T}}{dt}\\[5mm]
\ds \quad=-2\,\lf<d\vec{T}\,\dot{\otimes}\,d\vec{T}, d\vec{w}\,\dot{\otimes}\,d\vec{\La}\rg>_g+2\, e^{-2\la}\, \nabla\vec{T}\cdot\nabla\lf<d\vec{w}\wedge d\vec{\La}\rg>_g-2\, e^{-2\la}\, \nabla\vec{T}\cdot\nabla\lf[<d\vec{w}\cdot d\vec{\La}>_g\ \vec{T}\rg]\\[5mm]
\ds=-2\,\lf<d\vec{T}\,\dot{\otimes}\,d\vec{T}, d\vec{w}\,\dot{\otimes}\,d\vec{\La}\rg>_g+2\,d\vec{T}\cdot d\lf<d\vec{w}\wedge d\vec{\La}\rg>_g- 2\, |d\vec{T}|^2_g\ <d\vec{w}\cdot d\vec{\La}>_g\ .
\end{array}
\]
We deduce from (\ref{metric}) that
\[
\frac{d}{dt}\lf[dvol_g\rg]=\frac{d}{dt}\lf[\sqrt{g_{11}\,g_{22}-g_{12}^2}\rg] dx_1\wedge dx_2=<d\vec{w}\cdot d\vec{\La}>_g\ dvol_g\ .
\]
Hence combining the previous gives
\[
\begin{array}{l}
\ds\frac{d}{dt}\int_\Sigma (1+|d\vec{T}|^2_{g})^2\ dvol_g=4\, \int_\Sigma (1+|d\vec{T}|^2_{g})\ \lf[-\,\lf<d\vec{T}\,\dot{\otimes}\,d\vec{T}, d\vec{w}\,\dot{\otimes}\,d\vec{\La}\rg>_g+\,d\vec{T}\cdot d\lf<d\vec{w}\wedge d\vec{\La}\rg>_g\rg]\ dvol_g\\[5mm]
\ds\quad+\int_\Sigma\lf [(1+|d\vec{T}|^2_{g})^2- 4\, |d\vec{T}|^2_g\, (1+|d\vec{T}|^2_{g})\rg] <d\vec{w}\cdot d\vec{\La}>_g\ dvol_g\

\end{array}
\]
\subsubsection{Almost critical point of $E_\ep$ with entropy estimate.}
We are considering Legendrian immersions of a closed oriented surface $\Sigma$. We denote simply by $g$ the metric induced by the immersion. We  shall be using the following definition
\begin{Dfi}
\label{df-admi-pt}
We call ``admissible sequence of almost critical point'' of 
\[
E_\ep(\vec{\La}):=\int_\Sigma dvol_g+\ep^4 \int_\Sigma (1+|d\vec{T}|^2_{g})^2\ dvol_g
\]
any sequence of elements $\vec{\La}_k\in {\mathfrak E}_{\Sigma,Leg}^{2,4}$ and parameter $\ep_k\rightarrow 0$ satisfying respectively
\be
\label{energybound}
\limsup_{k\rightarrow +\infty}E_{\ep_k}(\vec{\La}_k)<+\infty\ ,
\ee
moreover
\be
\label{alm-crit}
\|dE_{\ep_k}(\vec{\La}_k)\|_{\vec{\La}_k}\le e^{-\ep_k^{-2}}\ ,
\ee
and satisfying the entropy condition
\be
\label{entropy}
\ep_k^4\,\int_\Sigma (1+|d\vec{T}_k|^2_{g_k})^2\ dvol_{g_k}=o\lf(\frac{1}{\log \ep_k^{-1}}\rg)\ .
\ee
\end{Dfi}
Hence in particular for any subdomain $\Om$ of $\Sigma$ and any smooth function $h$ in $V_2({\R}^4)$ such that supp$(h)\subset V_2({\R}^4)\setminus \vec{\La}_k(\p\Om)$ we have 
\be
\label{euler-lagrange}
\begin{array}{l}
\ds e^{-\ep_k^{-2}}\, O(\|w\|_{\vec{\La}_k})=\int_{\Om}<d\vec{w}\cdot d\vec{\La}_k>_{g_k}\ dvol_{g_k}\\[5mm]
\ds\quad+4\,\ep^4_k\, \int_\Om (1+|d\vec{T}_k|^2_{g_k})\ \lf[-\,\lf<d\vec{T}_k\,\dot{\otimes}\,d\vec{T}_k, d\vec{w}\,\dot{\otimes}\,d\vec{\La}_k\rg>_g+\,d\vec{T}_k\cdot d\lf<d\vec{w}\wedge d\vec{\La}_k\rg>_{g_k}\rg]\ dvol_{g_k}\\[5mm]
\ds\quad+\ep^4_k\ \int_\Om\lf [(1+|d\vec{T}_k|^2_{g_k})^2- 4\, |d\vec{T}_k|^2_{g_k}\, (1+|d\vec{T}_k|^2_{g_k})\rg] <d\vec{w}\cdot d\vec{\La}_k>_{g_k}\ dvol_{g_k}\
\end{array}
\ee
where $\vec{T}_k$ and $g_k$ are denoting respectively the tangent 2-vector to $\vec{\La}_k$ and the induced metric by $\vec{\La}_k$ moreover
\[
\vec{w}_h:=\vec{X}_h(\vec{\La}_k)=\lf[J_H\,\nabla^Hh+\frac{h}{2}\,\vec{R}\rg]\circ\vec{\La}_k\ .
\]
Observe that since $|d\vec{\La}_k|^2_{g_k}=2$ at every point we have - interpreting simply $\vec{N}_k$ as a section of $\vec{\La}^{-1} T{\R}^8$ -
\be
\label{pointwize-w}
|d\vec{w}|_g\le \ |\nabla\vec{X}_h|\quad\mbox{ and }\quad |\nabla^g d\vec{w}|_g\le C\ |\nabla^2\vec{X}_h|+C\ |\nabla\vec{X}_h|\ |d\vec{N}|_g\le C\ \lf[|\nabla^2\vec{X}_h|+|\nabla\vec{X}_h|\ [1+|d\vec{T}|_g]\rg]
\ee
Hence the entropy condition is implying
\be
\label{smoother}
\begin{array}{l}
\ds \lim_{k\rightarrow +\infty}4\,\ep^4_k\, \int_\Om (1+|d\vec{T}_k|^2_{g_k})\ \lf[-\,\lf<d\vec{T}_k\,\dot{\otimes}\,d\vec{T}_k, d\vec{w}\,\dot{\otimes}\,d\vec{\La}_k\rg>_g+\,d\vec{T}_k\cdot d\lf<d\vec{w}\wedge d\vec{\La}_k\rg>_g\rg]\ dvol_g\\[5mm]
\ds\quad+\ep^4_k\ \int_\Om\lf [(1+|d\vec{T}_k|^2_{g_k})^2- 4\, |d\vec{T}_k|^2_{g_k}\, (1+|d\vec{T}_k|^2_{g_k})\rg] <d\vec{w}\cdot d\vec{\La}_k>_{g_k}\ dvol_{g_k}=0\ .
\end{array}
\ee
Observe that 
\be
\label{horiz}
\p_{x_i}\vec{T}\res\vec{R}=\p_{x_i}\lf[ \frac{\p_{x_1}\vec{\La}\wedge\p_{x_2}\vec{\La}}{|\p_{x_1}\vec{\La}\wedge\p_{x_2}\vec{\La}|}  \rg]\res\vec{R}=0
\ee
since $\p^2\vec{\La}$ and $\p\vec{\La}$ are horizontal. Hence
\be
\label{ptw-horiz}
\lf|d\vec{T}_k\cdot d\lf<d\vec{w}\wedge d\vec{\La}_k\rg>_g\rg|\le \lf|d\vec{T}_k\rg|_g^2 |\nabla^H\vec{X}_h|+\lf|d\vec{T}_k\rg|_g\ |\nabla^H\nabla\vec{X}_h|
\ee
where $\nabla^H\vec{Y}=\nabla\vec{Y}-2^{-1}\,\nabla\vec{Y}\cdot\vec{R}\ \vec{R}$. We have in particular $\nabla^H\vec{R}=\nabla\vec{R}$ and then
\[
\nabla^H\nabla\vec{X}=\nabla^H\nabla^H\vec{X}+2^{-1}\nabla\vec{X}\cdot\vec{R}\ \nabla\vec{R}\ .
\]
We have at the point $\vec{p}=(\vec{a},\vec{b})\in V_2({\R}^4)$
\[
\begin{array}{l}
\ds\nabla\vec{R}(\vec{p})=\nabla^{{\R}^8}\vec{R}-{2}^{-1}\nabla^{\R^8}\vec{R}\cdot\vec{p}\, \vec{p}-{2}^{-1}\nabla^{\R^8}\vec{R}\cdot\vec{U}\, \vec{U}-{2}^{-1}\nabla^{\R^8}\vec{R}\cdot\vec{V}\, \vec{V}\\[5mm]
\ds=\lf(\begin{array}{c}   d\vec{b}\\[3mm] -\,d\vec{a} \end{array} \rg)-2^{-1}\,[ d\vec{b}\cdot\vec{a}-d\vec{a}\cdot\vec{b} ]\,\lf(\begin{array}{c}   \vec{a}\\[3mm] \vec{b} \end{array} \rg)-2^{-1}\,[ d\vec{b}\cdot\vec{a}+d\vec{a}\cdot\vec{b} ]\,\lf(\begin{array}{c}   \vec{a}\\[3mm] -\vec{b} \end{array} \rg)\ .
\end{array}
\]
In particular if $\vec{Z}$ is horizontal, we have
\be
\label{cova-reeb}
\nabla_{\vec{Z}}\vec{R}=\lf(\begin{array}{c}   d_{\vec{Z}}\vec{b}\\[3mm] -\,d_{\vec{Z}}\vec{a} \end{array} \rg)\ .
\ee
This permits to refine a bit (\ref{ptw-horiz}) in order to deduce
\be
\label{ptw-horiz-ref}
\lf|d\vec{T}_k\cdot d\lf<d\vec{w}\wedge d\vec{\La}_k\rg>_g\rg|\le \lf|d\vec{T}_k\rg|_g^2 |\nabla^H\vec{X}_h|+\lf|d\vec{T}_k\rg|_g\ [|\nabla^H\nabla^H\vec{X}_h|+|\nabla\vec{X}_h\cdot\vec{R}|]
\ee
It remains to estimate $ e^{-\ep_k^{-2}}\, O(\|w\|_{\vec{\La}_k})$. One has
\[
\|w\|_{\vec{\La}_k}=\lf[\int_\Sigma\lf[|\nabla^2\vec{w}|^2_g+|\nabla\vec{w}|^2_g+|\vec{w}|^2\rg]^2\ dvol_{g}\rg]^{1/4}+\||\nabla\vec{w}|_{g}\|_{L^\infty(\Sigma)}
\]
Using (\ref{pointwize-w}) we obtain
\be
\label{sss}
\|w\|_{\vec{\La}_k}\le \lf[\int_\Sigma\lf[|\nabla^2 \vec{X}_h|^2\circ{\vec{\La}_k}+|\nabla \vec{X}_h|^2\circ{\vec{\La}_k}\ |1+|d\vec{T}_k||^2+|\vec{X}_h|^2\rg]^2\ dvol_{g}\rg]^{1/4}+\|\nabla X_h\|_\infty\ .
\ee
Combining (\ref{forme-aire}), (\ref{euler-lagrange}), (\ref{smoother}) and (\ref{sss}) 

\begin{Lm}
\label{lm-passage}
For any choice of smooth function $h\in C^3(V_2({\R}^4))$ and any choice of admissible sequence of almost critical points of $E_\ep$ there holds
\be
\label{passage}
\lim_{k\rightarrow +\infty}\int_\Sigma\, <dh,d\beta_k>_{g_k}\ dvol_{g_k}=0\ .
\ee
\end{Lm}
Let $G_2^{Leg}(V_2({\R}^4))$ be the Grassmann bundle over $V_2({\R}^4)$ of 2-planes which are Legendrians. Denote by $\pi$ the projection from  $G_2^{Leg}(V_2({\R}^4))$ onto $V_2({\R}^4)$ associated to this bundle. A point in $G_2^{Leg}(V_2({\R}^4))$ is denoted $({\mathcal P},\vec{p})$ where ${\mathcal P}$ is a 2 dimensional subspace of $H_{\vec{p}}$ the horizontal 2-plane at $\vec{p}$. 
Considering now an admissible sequence of almost critical points of $E_\ep$ we shall denote ${\mathbf v}_{k}$ the varifold given by 
\[
\forall \ \Xi\in C^0(G_2^{Leg}(V_2({\R}^4)))\quad {\mathbf v}_{k}(\Xi):=\int_\Sigma\Xi((\vec{\La}_k)_\ast T_x\Sigma,\vec{\La}_k(x))\ dvol_{g_k}
\]
Modulo extraction of a subsequence the non negative measure ${\mathbf v}_k$ is converging in Radon measure to a limiting measure ${\mathbf v}_\infty$.  Given a $C^1$ vector-field $\vec{X}$ in $V_2({\R}^4)$ and a 2 plane
$\mathcal P$ in a fiber of the tangent bundle $TV_2({\R}^4)$ we define the divergence of $\vec{X}$ along $\mathcal P$ (see \cite{Sim} page 42) to be the number
\[
\mbox{div}_{\mathcal P}\vec{X}:=\sum_{j=1}\vec{e}_j\cdot \nabla_{\vec{e}_j}\vec{X}
\]
where $(\vec{e}_j)_{j=1,2}$ is an arbitrary orthonormal basis of ${\mathcal P}$ (one easily verify that the definition is independent of the choice of such a basis). If $\vec{X}$ is $C^1$ the map
\[
(\mathcal P,\vec{p})\ \longrightarrow\ \mbox{div}_{\mathcal P}\vec{X}
\]
is in particular a continuous function on the Grassmann bundle $G_2(V_2({\R}^4))$. Observe that
\[
\int_{\Sigma}<d\vec{w}\cdot d\vec{\La}_k>_{g_k}\ dvol_{g_k}=\int_{({\mathcal P},\vec{p})\in G_2^{Leg}(V_2({\R}^4))} \mbox{div}_{\mathcal P}\vec{X}_h\ d{\mathbf v}_k({\mathcal P},\vec{p})
\]
Hence we come naturally to the following definition
\begin{Dfi}
\label{df-leg-sta-var} A Radon measure $\mathbf v$ on the Grassmann bundle $G^{Leg}_2(V_2({\R}^4))$ of Legendrian 2-planes in $TV_2({\R}^4)$ is called Legendrian stationary or contact stationary if for any $C^2$ function $h$ on $V_2({\R}^4)$
there holds
\[
\int_{({\mathcal P},\vec{p})\in G_2^{Leg}(V_2({\R}^4))} \mbox{div}_{\mathcal P}\vec{X}_h\ d{\mathbf v}({\mathcal P},\vec{p})=0
\]
where
\[
\vec{X}_h:=J_H\,\nabla^Hh+\frac{h}{2}\,\vec{R}\ .
\]
\end{Dfi}
The weak limit ${\mathbf v}_\infty$ of ${\mathbf v}_k$ is clearly legendrian stationary thanks to (\ref{euler-lagrange}),  (\ref{smoother}) and the weak convergence in Radon measures. The goal of the sections below is to describe
${\mathbf v}_\infty$ in more details.


\section{The almost monotonicity formula at the limit}
\reset
In this section $\vec{\La}_k$ denotes an admissible sequence of almost critical points of $E_\ep$.
\subsection{The ``truncated'' almost monotonicity formula at the level $k$}
In this subsection we shall omit to write the subscript $k$. We assume $\vec{\La}$ is passing through the point $\vec{p}:=(\vec{\ep}_1,\vec{\ep}_2)$ where $\vec{\ep}_i$ is the canonical basis of ${\R}^4$. In the neighbourhood of $\vec{p}$ we introduce the {\it Legendrian coordinate}
\[
\varphi:=\vec{a}\cdot\vec{\ep}_2-\vec{b}\cdot\vec{\ep}_1=a_2-b_1\ ,
\]

\[
\rho^2:=|(\vec{a},\vec{b})-\vec{p}|^2
\]
In the neighbourhood of $\vec{p}$ we write $\vec{a}=\vec{\ep}_1+\vec{v}$ and $\vec{b}=\vec{\ep}_2+\vec{w}$. Using $|\vec{a}|^2=1$, $|\vec{b}|^2=1$ and $\vec{a}\cdot\vec{b}=0$ gives 
\be
\label{tan}
\lf\{
\begin{array}{l}
\ds 2\,v^1=-|\vec{v}|^2\\[3mm]
\ds2\,w_2=-|\vec{w}|^2\\[3mm]
\ds w^1+v^2=-v^3\,w^3-v^4\,w^4-v^1\,w^1-v^2\,w^2
\end{array}
\rg.
\ee
This gives
\[
\lf\{
\begin{array}{l}
v^1=O(\rho^2)\\[3mm]
w^2=O(\rho^2)\\[3mm]
w^1(1+O((v^2)^2))+v^2(1+O((w^1)^2))=O(\rho^2)\\[3mm]
v^2-w^1=\varphi
\end{array}
\rg.
\]
which finally implies
\[
\lf\{
\begin{array}{l}
v^1=O(\rho^2)\\[3mm]
w^2=O(\rho^2)\\[3mm]
w^1+v^2=O(\rho^2)\\[3mm]
v^2-w^1=\varphi
\end{array}
\rg.
\]
We deduce
\be
\label{liouv}
\al=\vec{a}\cdot d\vec{b}-\vec{b}\cdot d\vec{a}=-\,d\varphi+v^1\,dw^1-w^1\,dv^1+v^2\,dw^2-w^2\,dv^2+v^3\,dw^3-w^3\,dv^3+v^4\,dw^4-w^4\,dv^4\ .
\ee
We introduce  the multivalued function $\beta$ \underbar{defined on $\Sigma$} such that
\[
-d\beta:=(J_H(\vec{\La})\vec{H}_{\vec{\La}})^\sharp\ .
\]
With these notations, in conformal coordinates we have
\[
J_H\,\vec{H}_{\vec{\La}}=\frac{e^{-2\la}}{2} \sum_{k=1}^2J_H\,\nabla_{x_k}\p_{x_k}\vec{\La}=<d\beta, d\vec{\La}>_g=-{e^{-2\la}}\,\sum_{k=1}^2\,\p_{x_k}\beta\,\p_{x_k}\vec{\La}\ .
\]
We take the scalar product of this identity with $\vec{\La}-\vec{p}$. This gives
\be
\label{equation-courb}
\sum_{k=1}^2[\vec{\La}-\vec{p}]\cdot J_H\,\nabla_{x_k}\p_{x_k}\vec{\La}=-\sum_{k=1}^2\,\p_{x_k}\beta\,\p_{x_k}|\vec{\La}-\vec{p}|^2
\ee
Observe that
\[
\begin{array}{l}
\ds[\vec{\La}-\vec{p}]\cdot J_H\,\nabla_{x_k}\p_{x_k}\vec{\La}=\p^2_{x_k^2}\vec{\Phi}\cdot(\vec{n}-\vec{\ep}_2)-\p^2_{x_k^2}\vec{n}\cdot(\vec{\Phi}-\vec{\ep}_1)=-\p_{x_k}\lf[\p_{x_k}\vec{\Phi}\cdot\vec{\ep}_2-\p_{x_k}\vec{n}\cdot\vec{\ep}_1\rg] =-\p^2_{x_k^2}\varphi\ .
\end{array}
\]
Hence
\be
\label{first-equ}
\Delta\varphi=\sum_{k=1}^2\,\p_{x_k}\beta\,\p_{x_k}|\vec{\La}-\vec{p}|^2\ .
\ee
Multiplying now by $J_H(\vec{\ep}_2,-\vec{\ep}_1)^H$ where $(\vec{\ep}_2,-\vec{\ep}_1)^H$ is the orthogonal projection of $(\vec{\ep}_2,-\vec{\ep}_1)$ onto $H$
\[
\sum_{k=1}^2J(\vec{\ep}_2,-\vec{\ep}_1)\cdot\nabla_{x_k}\p_{x_k}\vec{\La}=2\,\sum_{k=1}^2\,\p_{x_k}\beta\, (\vec{\ep}_2,-\vec{\ep}_1)\cdot \p_{x_k}\vec{\La}=2\,\sum_{k=1}^2\,\p_{x_k}\beta\,  \p_{x_k}{\varphi}
\]
Recall that $\nabla_{x_k}\p_{x_k}\vec{\La}$ is horizontal and that
\[
\begin{array}{l}
\ds\sum_{k=1}^2\nabla_{x_k}\p_{x_k}\vec{\La}=\Delta\vec{\La}-2^{-1}\,\vec{\La}\cdot\Delta\vec{\La}\ \vec{\La}-2^{-1}\,\vec{U}\cdot\Delta\vec{\La}\ \vec{U}-2^{-1}\,\vec{V}\cdot\Delta\vec{\La}\ \vec{V}\\[5mm]
\ds=\Delta\vec{\La}+2^{-1}\,\lf[|\nabla\,\vec{\La}|^2\,\vec{\La}+\nabla\vec{U}\cdot\nabla\vec{\La}\ \vec{U}+\nabla\vec{V}\cdot\nabla\vec{\La}\ \vec{V}\rg]
\end{array}
\]
We have
\[
\begin{array}{l}
\ds 2^{-1}\,J(\vec{\ep}_2,-\vec{\ep}_1)\cdot\lf[ \,|\nabla\,\vec{\La}|^2\,\vec{\La}+\nabla\vec{U}\cdot\nabla\vec{\La}\ \vec{U}+\nabla\vec{V}\cdot\nabla\vec{\La}\ \vec{V} \rg]\\[5mm]
\ds \ =2\, e^{2\la}+  (v^1+w^2)\ e^{2\la}+ 2^{-1}\,\nabla\vec{U}\cdot\nabla\vec{\La}\ (v^1-w^2)+2^{-1}\,\nabla\vec{V}\cdot\nabla\vec{\La}\ (w^1+v^2)=2\, e^{2\la} \,\lf(1+O(\rho^2)+O(\varphi)\rg)\ .
\end{array}
\]
hence
\[
\Delta[v^1+w^2]+2\, e^{2\la} \,\lf(1+O(\rho^2)+O(\varphi)\rg)=2\,\sum_{k=1}^2\,\p_{x_k}\beta\,  \p_{x_k}{\varphi}
\]
We have
\[
2\, (v^1+w^2)=-\sum_{i=1}^4(v^i)^2+(w^i)^2=-|\vec{\La}-\vec{p}|^2
\]
This gives finally
\be
\label{systeme}
\lf\{
\begin{array}{l}
\ds 1+O(\rho^2)+O(\varphi)-\,4^{-1}\Delta_g\rho^2=<d\beta,d\varphi>_g\\[5mm]
\ds \Delta_g\varphi= <d\beta,d\rho^2>_g
\end{array}
\rg.
\ee
We have
\[
\begin{array}{l}
\ds|\vec{\La}-\vec{p}|^2=2^{-1}\,|(\vec{\La}-\vec{p})\cdot\vec{\La}|^2+2^{-1}\,|(\vec{\La}-\vec{p})\cdot\vec{R}(\vec{\La})|^2+2^{-1}\,|(\vec{\La}-\vec{p})\cdot\vec{U}(\vec{\La})|^2+2^{-1}\,|(\vec{\La}-\vec{p})\cdot\vec{V}(\vec{\La})|^2\\[5mm]
\ds+e^{-2\la}\,\lf|(\vec{\La}-\vec{p})\cdot\nabla\vec{\La}\rg|^2+e^{-2\la}\,\lf|(\vec{\La}-\vec{p})\cdot J_H\nabla\vec{\La}\rg|^2
\end{array}
\]
We have respectively
\[
|(\vec{\La}-\vec{p})\cdot\vec{\La}|^2=|v^1+w^2|^2=O(\rho^4)\ ,\ |(\vec{\La}-\vec{p})\cdot\vec{R}(\vec{\La})|^2=|w^1-v^2|^2=\varphi^2\ ,
\]
and
\[
|(\vec{\La}-\vec{p})\cdot\vec{U}(\vec{\La})|^2=|v^1-w^2|^2=O(\rho^4) ,\ |(\vec{\La}-\vec{p})\cdot\vec{V}(\vec{\La})|^2=|w^1+v^2|^2=O(\rho^4)+O(\varphi^2)\ .
\]
Moreover
\[
e^{-2\la}\,\lf|(\vec{\La}-\vec{p})\cdot\nabla\vec{\La}\rg|^2=4^{-1}\, \lf|d\rho^2\rg|^2_g\ ,\ e^{-2\la}\,\lf|(\vec{\La}-\vec{p})\cdot J_H\nabla\vec{\La}\rg|^2=|d(v^2-w^1)|^2_g=|d\varphi|^2_g
\]
Combining the previous is implying
\be
\label{structure}
1=|d\rho|^2_g+\rho^{-2}|d\varphi|_g^2+O(\rho^2)+\rho^{-2}\,O(\varphi^2)\ .
\ee
Away from $\rho=0$ we introduce $$\sigma:=\frac{2\,\varphi}{\rho^2}\ .$$ We now follow the main computations of \cite{Riv}. First of all we have (see (III.10) in \cite{Riv})
\[
\begin{array}{l}
\ds <d\sigma,d\beta>_g=\lf<[2\,\rho^{-2}\,d\varphi-2\,\varphi\,\rho^{-4}\,d\rho^2]\,,\,d\beta\rg>_g\\[5mm]
\ds=4\rho^{-2}\ (1+O(\rho^2)+O(\varphi)+\rho^{-2}\,O(\varphi^2))+4^{-1}\,\rho^{-4}\,\Delta_g\r\\[5mm]
\ds =4\,\frac{\sqrt{1+\sigma^2}}{\r^2}\ \lf(1+O(\rho^2)+O(\varphi)+\rho^{-2}\,O(\varphi^2)\rg)-4^{-1}\,\frac{1+\sigma^2}{\r^4}\, \Delta_g\r^4
\end{array}
\]
where $\r$ is the {\it Folland-Kor\'anyi gauge} given by $\r:=(\rho^4+4\,\varphi^2)^{1/4}$. This gives in particular
\be
\label{pre-monoton}
\lf<\frac{d\sigma}{1+\sigma^2},d\beta\rg>_g=\frac{4}{\r^2}\frac{1}{\sqrt{1+\sigma^2}}\,\lf(1+O(\rho^2)+O(\varphi)+\rho^{-2}\,O(\varphi^2)\rg)-4\,|d\log\r|^2_g-\Delta_g\log\r\ .
\ee
Observe that from (\ref{liouv})
\[
\nabla^H\varphi=\sum_{i=1}^4 v^i\, \nabla^Hw^i-w^i\, \nabla^Hv^i\ .
\]
We have in ${\R}^8$ at a point $\vec{Q}=(\vec{\ep}_1+\vec{v},\vec{\ep}_2+\vec{w})\in V_2({\R}^8)$
\[
\begin{array}{l}
\ds\lf|\sum_{i=1}^4 v^i\, \nabla^Hw^i-w^i\, \nabla^Hv^i\rg|^2=\lf|\sum_{i=1}^4 v^i\, \nabla w^i-w^i\, \nabla v^i\rg|^2-2^{-1}\,\lf|\lf(\sum_{i=1}^4 v^i\, \nabla w^i-w^i\, \nabla v^i\rg)\cdot\vec{R}\rg|^2\\[5mm]
\ds -2^{-1}\,\lf|\lf(\sum_{i=1}^4 v^i\, \nabla w^i-w^i\, \nabla v^i\rg)\cdot\vec{U}\rg|^2 -2^{-1}\,\lf|\lf(\sum_{i=1}^4 v^i\, \nabla w^i-w^i\, \nabla v^i\rg)\cdot\vec{V}\rg|^2 -2^{-1}\,\lf|\lf(\sum_{i=1}^4 v^i\, \nabla w^i-w^i\, \nabla v^i\rg)\cdot\vec{Q}\rg|^2\\[5mm]
\ds=\rho^2-2^{-1}|-v^1\, (1+v^1)-(w^1)^2-(v^2)^2-w^2\, (1+w^2) -(v^3)^2-(v^4)^2 -(w^3)^2-(w^4)^2 |^2\\[5mm]
\ds-\, 2^{-1}\,|v^2+w^1+2\,\sum_{i=1}^4 v^i\,w^i|^2\\[5mm]
\ds\ -\,2^{-1}|v^1\, (1+v^1)+(v^2)^2+(v^3)^2+(v^4)^2-(w^1)^2-w^2\,(1+w^2)-(w^3)^2-(w^4)^2|^2-\, 2^{-1}\,|v^2-w^1|^2\\[5mm]
\ds=\rho^2+O(\rho^4)+O(\varphi^2)=\rho^2+O(\r^4)
\end{array}
\]
We have also
\[
\begin{array}{l}
\ds|\rho\,\nabla^H\rho|^2=\lf|\sum_{i=1}^4 v^i\, \nabla^Hv^i+w^i\, \nabla^Hw^i\rg|^2=\lf|\sum_{i=1}^4 v^i\, \nabla v^i+w^i\, \nabla w^i\rg|^2-2^{-1}\,\lf|\lf(\sum_{i=1}^4 v^i\, \nabla v^i+w^i\, \nabla w^i\rg)\cdot\vec{R}\rg|^2\\[5mm]
\ds -2^{-1}\,\lf|\lf(\sum_{i=1}^4 v^i\, \nabla v^i+w^i\, \nabla w^i\rg)\cdot\vec{U}\rg|^2 -2^{-1}\,\lf|\lf(\sum_{i=1}^4 v^i\, \nabla v^i+w^i\, \nabla w^i\rg)\cdot\vec{V}\rg|^2 -2^{-1}\,\lf|\lf(\sum_{i=1}^4 v^i\, \nabla v^i+w^i\, \nabla w^i\rg)\cdot\vec{Q}\rg|^2\\[5mm]
\ds=\rho^2-2^{-1}\,|v^2-w^1|^2-2^{-1}\,\lf|v^1-w^2+\sum_{i=1}^4(v^i)^2-(w^i)^2\rg|^2-2^{-1}\,\lf|v^2+w^1+2\,\sum_{i=1}^2v^i\,w^i\rg|^2-2^{-1}\,|v^1+w^2+\rho^2|^2\\[5mm]
\ds=\rho^2+O(\varphi^2)+O(\rho^4)=\rho^2+O(\r^4)\ .
\end{array}
\]
and similarly we prove
\[
\nabla^H\varphi\cdot\nabla^H\rho=O(\r^4)\ .
\]
Following again \cite{Riv} (III.13) we obtain
\be
\label{horizont-gradient}
|\nabla^H\r|^2=\frac{\rho^6}{\r^6}|\nabla^H\rho|^2+\frac{|\nabla^H\varphi^2|}{\r^6}+4\,\frac{\varphi\,\rho^3}{\r^6}\, \nabla^H\varphi\cdot\nabla^H\rho=\frac{\rho^2}{\r^2}+O(\r^2)=\frac{1}{\sqrt{1+\sigma^2}}+O(\r^2)\ .
\ee
We denote by $(\nabla^\Sigma\r)^\perp:=\nabla^H\r-\nabla^\Sigma\r$ where $\nabla^\Sigma\r$ is denoting the gradient of $\r\circ\vec{\La}$ that is the orthogonal projection of $\nabla^H\r$ onto the tangent plane $\vec{\La}_\ast T\Sigma$. Hence one has in particular
\[
|d\log\r|^2_g=|\nabla^\Sigma\r|^2\quad\mbox{ and }\quad|(\nabla^\Sigma\r)^\perp|^2=|\nabla^H\r|^2-|\nabla^\Sigma\r|^2=\frac{1}{\sqrt{1+\sigma^2}}+O(\r^2)-|\nabla^\Sigma\r|^2\ .
\]
Combining this fact with (\ref{pre-monoton}) and (\ref{horizont-gradient}) and arguing as in \cite{Riv} we finally obtain away from $\rho=0$
\[
\begin{array}{l}
\ds <d\arctan\sigma,d\beta>_g+\Delta_g\log\r=4\,\frac{|(\nabla^\Sigma\r)^\perp|^2}{\r^2}+\frac{4}{\r^2}\frac{1}{\sqrt{1+\sigma^2}}\,\lf(O(\r^2)+\rho^{-2}\,O(\varphi^2)\rg)+O(1)
\end{array}
\]
Observe that
\[
\frac{1}{\rho^2}\,\frac{1}{\sqrt{1+\sigma^2}}=\frac{1}{\sqrt{\rho^4+4\varphi^2}}=\frac{1}{\r^2}\ .
\]
Hence
\[
\frac{4}{\r^2}\frac{1}{\sqrt{1+\sigma^2}}\,\lf(O(\r)^2+\rho^{-2}\,O(\varphi^2)\rg)=O(1)\ ,
\]
and we get away from $\rho=0$
\be
\label{identite-monotone}
<d\arctan\sigma,d\beta>_g+\Delta_g\log\r=4\,\frac{|(\nabla^\Sigma\r)^\perp|^2}{\r^2}+O(1)
\ee
As in \cite{Riv} again we compute away from $\rho=0$
\be
\label{autre-identite}
\begin{array}{l}
\ds \r^3\,<d\r,d\beta>_g=\rho^3\, <d\rho,d\beta>_g+<d\varphi^2, d\beta>_g\\[5mm]
\ds=2^{-1}\,\rho^2\, \Delta_g\varphi+2\, \varphi\ [1+O(\r^2)-\,4^{-1}\Delta_g\rho^2]\\[5mm]
\ds=2\varphi \,[1+O(\r^2)]+2^{-1} \mbox{div}^\Sigma(\rho^2\,\nabla^\Sigma\varphi-\varphi\nabla^\Sigma\rho^2)\\[5mm]
\ds=\r^2\,\frac{\sigma}{\sqrt{1+\sigma^2}}+O(\r^4)+4^{-1}\,\mbox{div}^\Sigma\lf(\r^4\frac{\nabla^\Sigma\sigma}{1+\sigma^2}\rg)
\end{array}
\ee
Let $0<r<1$ arbitrary and let $\chi$ be a cut-off function such that
\[
\chi(t)=\lf\{
\begin{array}{l}
1\quad\mbox{ for }t<1\\[3mm]
0\quad\mbox{ for }t>2
\end{array}
\rg.\quad,\quad \chi'\le 0\quad\mbox{on }{\R}_+\ .
\]
Let $1>r>\eta>0$ and consider the  hamiltonian function given on a neighbourhood of $\vec{p}$ in $V_2({\R}^4)$ by
\[
h_{r,\eta}:=[\chi(\r/r)-\chi(\r/\eta)]\, \arctan\sigma\ .
\]
We claim  that away from $\r=0$ the function $\arctan\sigma$ defines a $C^l$ function for any $l\in {\N}$. Indeed  for $\rho=0$ and $\varphi>0$ it extends continuously by $\pi/2$ while for $\rho=0$ and $\varphi<0$ it extends continuously by $-\pi/2$.
We have moreover
\be
\label{id-sigma}
\nabla\arctan\sigma=\frac{2}{1+4\,\frac{\varphi^2}{\rho^4}}\,\lf[\frac{\nabla\varphi}{\rho^2}-2\,\frac{\varphi}{\rho^3}\ \nabla\rho   \rg]=2\,\frac{\rho^2}{\r^4}\,\nabla\varphi- \frac{2}{\r^4}\, \varphi\,\nabla\rho^2\ .
\ee
Hence, away from $\r=0$, $\nabla\arctan\sigma$ extends continuously by $0$ on $\rho=0$. $\r$ is obviously smooth away from $\r=0$ as well as $\rho^2\,\nabla\varphi-\varphi\,\nabla\rho^2$ which is proving the claim.
We multiply (\ref{identite-monotone}) by $[\chi(\r/r)-\chi(\r/\eta)]$ and (\ref{autre-identite}) by $\r^{-3}\,\arctan\sigma\, [r^{-1}\chi'(\r/r)-\eta^{-1}\chi'(\r/\eta)]$, this gives
\[
\begin{array}{l}
\ds<d h_{r,\eta},d\beta>_g+[\chi(\r/r)-\chi(\r/\eta)]\,\Delta_g\log\r    =   4\,[\chi(\r/r)-\chi(\r/\eta)]\,\frac{|(\nabla^\Sigma\r)^\perp|^2}{\r^2}+O(1)\,[\chi(\r/r)-\chi(\r/\eta)]\\[5mm]
\ds+[r^{-1}\chi'(\r/r)-\eta^{-1}\chi'(\r/\eta)]\,\r^{-1}\,\arctan\sigma\,\frac{\sigma}{\sqrt{1+\sigma^2}}+O(\r)\, [r^{-1}\chi'(\r/r)-\eta^{-1}\chi'(\r/\eta)]\\[5mm]
\ds+4^{-1}\,\r^{-3} \ [r^{-1}\chi'(\r/r)-\eta^{-1}\chi'(\r/\eta)]\ \arctan\sigma\,\mbox{div}^\Sigma\lf(\r^4\frac{\nabla^\Sigma\sigma}{1+\sigma^2}\rg)
\end{array}
\]
Integrating over $\Sigma$ gives
\[
\begin{array}{l}
\ds\int_\Sigma<d h_{r,\eta},d\beta>_g\ dvol_g-\int_\Sigma\lf[\frac{\r}{r}\,\chi'\lf(\frac{\r}{r}\rg)-\frac{\r}{\eta}\,\chi'\lf(\frac{\r}{\eta}\rg)\rg]\frac{|\nabla^\Sigma\r|^2}{\r^2}\ dvol_g\\[5mm]
\ds+\int_\Sigma\ O(1)\ [\chi(\r/r)-\chi(\r/\eta)]+[O(\r/r)\,\chi'(\r/r)-O(\r/\eta)\chi'(\r/\eta)]\ dvol_g \\[5mm]
\ds+4^{-1}\,\int_\Sigma \nabla^\Sigma\lf[\r^{-3} \ [r^{-1}\chi'(\r/r)-\eta^{-1}\chi'(\r/\eta)]\ \arctan\sigma\rg]\,\r^4\frac{\nabla^\Sigma\sigma}{1+\sigma^2}\ dvol_g \\[5mm]
\ds-\int_\Sigma\ \lf[\frac{\r}{r}\,\chi'\lf(\frac{\r}{r}\rg)-\frac{\r}{\eta}\,\chi'\lf(\frac{\r}{\eta}\rg)\rg]\,\frac{1}{\r^{2}}\,\frac{\sigma\,\arctan\sigma}{\sqrt{1+\sigma^2}}\ dvol_g =4\,\int_\Sigma[\chi(\r/r)-\chi(\r/\eta)]\,\frac{|(\nabla^\Sigma\r)^\perp|^2}{\r^2}\ dvol_g
\end{array}
\]
From which we deduce the ``truncated\footnote{``Truncated'' is referring to the fact that we are not making $\eta$ go to zero at this stage and that we have to ``wait'' for  $k$ going to infinity first.} '' almost monotonicity formula satisfied by $\vec{\La}_k$ :
\be
\label{truncated}
\begin{array}{l}
\ds\int_\Sigma<d h_{r,\eta},d\beta>_g\ dvol_g-\int_\Sigma\frac{\r}{r}\,\chi'\lf(\frac{\r}{r}\rg)\,\frac{1}{\r^2}\,\lf[{|\nabla^\Sigma\r|^2} +\frac{\sigma\,\arctan\sigma}{\sqrt{1+\sigma^2}} \rg]\ dvol_g\\[5mm]
\ds+\int_\Sigma\ \lf[O(1)\chi\lf(\frac{\r}{r}\rg)+O\lf(\frac{\r}{r}\rg)\,\chi'\lf(\frac{\r}{r}\rg)\rg]\ dvol_g +\frac{1}{4}\,\int_\Sigma \frac{\r^2}{r^{2}}\,\chi''\lf(\frac{\r}{r}\rg) \,\arctan\sigma\,\frac{\nabla^\Sigma\r}{\r}\cdot\frac{\nabla^\Sigma\sigma}{1+\sigma^2}\ dvol_g \\[5mm]
\ds-\frac{3}{4}\,\int_\Sigma  \ \frac{\r}{r}\,\chi'\lf(\frac{\r}{r}\rg) \,\arctan\sigma\,\frac{\nabla^\Sigma\r}{\r}\cdot\frac{\nabla^\Sigma\sigma}{1+\sigma^2}\ dvol_g +\frac{1}{4}\,\int_\Sigma \frac{\r}{r}\chi'\lf(\frac{\r}{r}\rg)\,\frac{|\nabla^\Sigma\sigma|^2}{(1+\sigma^2)^2}\ dvol_g \\[5mm]
\ds=4\,\int_\Sigma[\chi(\r/r)-\chi(\r/\eta)]\,\frac{|(\nabla^\Sigma\r)^\perp|^2}{\r^2}\ dvol_g-\int_\Sigma\frac{\r}{\eta}\,\chi'\lf(\frac{\r}{\eta}\rg)\,\frac{1}{\r^2}\,\lf[{|\nabla^\Sigma\r|^2} +\frac{\sigma\,\arctan\sigma}{\sqrt{1+\sigma^2}} \rg]\ dvol_g\\[5mm]
\ds+4\,\int_\Sigma\frac{\r}{\eta}\, \chi'\lf(\frac{\r}{\eta}\rg)\,\frac{|\nabla^\Sigma\sigma|^2}{(1+\sigma^2)^2}\ dvol_g-\frac{3}{4}\int_\Sigma\,\frac{\r}{\eta}\,\chi'\lf(\frac{\r}{\eta}\rg)\, \,\arctan\sigma\,\frac{\nabla^\Sigma\r}{\r}\cdot\frac{\nabla^\Sigma\sigma}{1+\sigma^2}\ dvol_g  \\[5mm]
\ds+\frac{1}{4}\int_\Sigma\frac{\r^2}{\eta^{2}}\,\chi''\lf(\frac{\r}{\eta}\rg)\, \,\arctan\sigma\,\frac{\nabla^\Sigma\r}{\r}\cdot\frac{\nabla^\Sigma\sigma}{1+\sigma^2}\ dvol_g +\int_\Sigma\ \lf[O(1)\chi\lf(\frac{\r}{\eta}\rg)+O\lf(\frac{\r}{\eta}\rg)\,\chi'\lf(\frac{\r}{\eta}\rg)\rg]\ dvol_g
\end{array}
\ee
\subsection{Passing to the limit $k\rightarrow +\infty$ in the ``truncated'' almost monotonicity formula  (\ref{truncated}).}
We have seen (\ref{passage}) that
\[
\lim_{k\rightarrow +\infty}\int_\Sigma<d h_{r,\eta},d\beta_k>_{g_k}\ dvol_{g_k}=0
\]
More precisely, from (\ref{euler-lagrange}), (\ref{entropy}), (\ref{pointwize-w}), (\ref{ptw-horiz-ref}) and from (\ref{forme-aire}) we obtain
\[
\begin{array}{l}
\ds\lf|\int_\Sigma<d h_{r,\eta},d\beta_k>_{g_k}\ dvol_{g_k}\rg|\le C\ \ep_k^4\,\int_{ \Sigma}(1+|d\vec{T}_k|^2_{g_k})^2 \ |\nabla^H \vec{X}_{h_{r,\eta}}|\circ\vec{\La}_k\ dvol_{g_k}\\[5mm]
\ds\quad\quad+\, C\ \ep_k^4\,\int_{ \Sigma}(1+|d\vec{T}_k|^2_{g_k})\  |d\vec{T}_k|_{g_k}\ \lf[|\nabla^H\nabla^H \vec{X}_{h_{r,\eta}}|\circ\vec{\La}_k+|\nabla\vec{X}_{h_{r,\eta}}\cdot\vec{R}|\circ\vec{\La}_k\rg]\ dvol_{g_k}
\end{array}
\]
Recall
\[
\vec{X}_h=J_H\nabla^Hh+\frac{h}{2}\,\vec{R}\ .
\]
We have in particular
\[
J_H\nabla^H\lf[\chi(\r/\eta)\,\arctan\sigma\rg]=\frac{J_H\nabla^H\r}{\eta}\,\chi'(\r/\eta)\,\arctan\sigma+\chi(\r/\eta)\,\frac{J_H\nabla^H\sigma}{1+\sigma^2}
\]
We recall that a subsequence has been extracted in such a way that the varifold ${\mathbf v}_k$ defined by $\vec{\La}_k$ is converging towards a limiting varifold ${\mathbf v}_\infty$. 
We also denote by $\mu_k$ (resp. $\mu_\infty$ ) the so called {\it weight} of ${\mathbf v}_k$ (resp. ${\mathbf v}_\infty$) . It is given by
\[
\forall \ A\subset V_2({\R}^4)\quad\mbox{ Borel }\quad\mu_k(A):={\mathbf v_k}(\pi^{-1}(A))
\]
where $\pi$ denotes the bundle projection of the Grassmann bundle $G_2^{Leg}(V_2({\R}^4)$ over $V_2({\R}^4)$. 
The Radon measure convergence
in $G_2^{Leg}(V_2({\R}^4))$ is implying the following passages to the limit for any smooth functions $F$ and $G$ on $V_2({\R}^4)$ : we have respectively
\[
\lim_{k\rightarrow +\infty} \int_\Sigma F(\vec{\La}_k(x))\ dvol_{g_k}=\int_{V_2({\R}^4)} F(\vec{p})\ d\mu_\infty(\vec{p})\ ,
\]
we have also
\[
\begin{array}{l}
\ds\lim_{k\rightarrow +\infty} \int_\Sigma\ G(\vec{\La}_k)\  |\nabla^\Sigma(F(\vec{\La}_k(x)))|^2\ dvol_{g_k}=\lim_{k\rightarrow +\infty} \int_{G_2^{Leg}(V_2({\R}^4))}\ G(\vec{p})\  |\nabla^{\mathcal P} F(\vec{p}) |^2\ d{\mathbf v}_k({\mathcal P},\vec{p})\\[5mm]
\ds \quad=\int_{G_2^{Leg}(V_2({\R}^4))}\ G(\vec{p})\  |\nabla^{\mathcal P} F(\vec{p}) |^2\ d{\mathbf v}_\infty({\mathcal P},\vec{p})
\end{array}
\]
where $\nabla^{\mathcal P} F$ is denoting the orthogonal projection of the gradient of $F$ on the Legendrian 2-plane ${\mathcal P}$ in $T_{\vec{p}}V_2({\R}^4)$. Passing to the limit in (\ref{truncated}) gives
\be
\label{truncated-limit}
\begin{array}{l}
\ds-\int_{G_2^{Leg}(V_2({\R}^4))}\ \frac{\r(\vec{p})}{r}\,\chi'\lf(\frac{\r(\vec{p})}{r}\rg)\,\frac{1}{\r^2(\vec{p})}\ \lf[{|\nabla^{\mathcal P}\r|^2}+ \frac{\sigma(\vec{p})\,\arctan\sigma(\vec{p})}{\sqrt{1+\sigma^2(\vec{p})}} \rg]\ d{\mathbf v}_\infty({\mathcal P},\vec{p}) \\[5mm]
\ds+\int_{G_2^{Leg}(V_2({\R}^4))} \lf[O(1)\,\lf(\chi\lf(\frac{\r(\vec{p})}{r}\rg)-\chi\lf(\frac{\r(\vec{p})}{\eta}\rg)\rg)+O\lf(\frac{\r(\vec{p})}{r}\rg)\,\chi'\lf(\frac{\r(\vec{p})}{r}\rg)\rg]\ d{\mathbf v}_\infty({\mathcal P},\vec{p})\\[5mm]
\ds+\frac{1}{4}\,\int_{G_2^{Leg}(V_2({\R}^4))}\frac{\r^2(\vec{p})}{r^{2}}\,\chi''\lf(\frac{\r(\vec{p})}{r}\rg) \,\arctan\sigma(\vec{p})\ \frac{\nabla^{\mathcal P}\r}{\r}\cdot\frac{\nabla^{\mathcal P}\sigma}{1+\sigma^2}\ d{\mathbf v}_\infty({\mathcal P},\vec{p}) \\[5mm]
\ds-\frac{3}{4}\,\int_{G_2^{Leg}(V_2({\R}^4))}  \ \frac{\r(\vec{p})}{r}\,\chi'\lf(\frac{\r(\vec{p})}{r}\rg) \,\arctan\sigma(\vec{p})\,\frac{\nabla^{\mathcal P}\r}{\r}\cdot\frac{\nabla^{\mathcal P}\sigma}{1+\sigma^2(\vec{p})}\ \ d{\mathbf v}_\infty({\mathcal P},\vec{p}) \\[5mm]
\ds+\frac{1}{4}\,\int_{G_2^{Leg}(V_2({\R}^4))} \frac{\r(\vec{p})}{r}\chi'\lf(\frac{\r(\vec{p})}{r}\rg)\,\frac{|\ds\nabla^{\mathcal P}\sigma|^2}{(1+\sigma^2(\vec{p}))^2}\ \ d{\mathbf v}_\infty({\mathcal P},\vec{p}) \\[5mm]
\ds=4\,\int_{G_2^{Leg}(V_2({\R}^4))} \lf[\chi\lf(\frac{\r(\vec{p})}{r}\rg)-\chi\lf(\frac{\r(\vec{p})}{\eta}\rg)\rg]\,\frac{|(\nabla^{\mathcal P}\r)^\perp|^2}{\r(\vec{p})^2}\ d{\mathbf v}_\infty({\mathcal P},\vec{p})\\[5mm]
\ds-\int_{G_2^{Leg}(V_2({\R}^4))} \frac{\r(\vec{p})}{\eta}\,\chi'\lf(\frac{\r(\vec{p})}{\eta}\rg)\,\frac{1}{\r^2(\vec{p})}\,\lf[{|\nabla^{\mathcal P}\r|^2} +\frac{\sigma\,\arctan\sigma}{\sqrt{1+\sigma^2}}(\vec{p}) \rg]\ \ d{\mathbf v}_\infty({\mathcal P},\vec{p})\\[5mm]
\ds+\,4\,\int_{G_2^{Leg}(V_2({\R}^4))}\frac{\r(\vec{p})}{\eta}\, \chi'\lf(\frac{\r(\vec{p})}{\eta}\rg)\,\frac{|\nabla^{\mathcal P}\sigma|^2}{(1+\sigma^2(\vec{p}))^2}\ d{\mathbf v}_\infty({\mathcal P},\vec{p})\\[5mm]
\ds-\frac{3}{4}\int_{G_2^{Leg}(V_2({\R}^4))}\,\frac{\r(\vec{p})}{\eta}\,\chi'\lf(\frac{\r(\vec{p})}{\eta}\rg)\, \,\arctan\sigma(\vec{p})\ \frac{\nabla^{\mathcal P}\r}{\r(\vec{p})}\cdot\frac{\nabla^{\mathcal P}\sigma}{1+\sigma^2(\vec{p})}\ d{\mathbf v}_\infty({\mathcal P},\vec{p})  \\[5mm]
\ds+\frac{1}{4}\int_{G_2^{Leg}(V_2({\R}^4))}\frac{\r^2(\vec{p})}{\eta^{2}}\,\chi''\lf(\frac{\r(\vec{p})}{\eta}\rg)\, \,\arctan\sigma(\vec{p})\ \frac{\nabla^{\mathcal P}\r}{\r(\vec{p})}\cdot\frac{\nabla^{\mathcal P}\sigma}{1+\sigma^2(\vec{p})} \ d{\mathbf v}_\infty({\mathcal P},\vec{p})\\[5mm]
\ds+\int_{G_2^{Leg}(V_2({\R}^4))}\ \lf[O\lf(\frac{\r(\vec{p})}{\eta}\rg)\,\chi'\lf(\frac{\r(\vec{p})}{\eta}\rg)\rg]\ \ d{\mathbf v}_\infty({\mathcal P},\vec{p})
\end{array}
\ee
We first deduce the following lemma
\begin{Lm}
\label{lm-dens-finie}
Assume
\be
\label{dens-finie}
\liminf_{\eta\rightarrow 0}\frac{1}{\eta^2}\int_{ \eta<\r<2\eta}d\mu_\infty(\vec{p})<+\infty\ ,
\ee
then
\be
\label{bound-arctan-sigma}
\lim_{\eta\rightarrow 0}\int_{ \eta<\r<1}\frac{|\nabla^{\mathcal P}\sigma|^2}{(1+\sigma^2(\vec{p}))^2}\ d{\mathbf v}_\infty({\mathcal P},\vec{p})<+\infty\ .
\ee
\end{Lm}
\noindent{\bf Proof of lemma~\ref{lm-dens-finie}}. We have from (\ref{id-sigma}) away from $\r=0$ and for any $({\mathcal P},\vec{p})\in G_2^{Leg}(V_2({\R}^4))$
\be
\label{df-1}
\lf|\frac{\nabla^{\mathcal P}\sigma}{1+\sigma^2(\vec{p})}\rg|=\frac{1}{\r^4(\vec{p})}\, \lf|\rho^2(\vec{p})\,\nabla^{\mathcal P}\varphi-\varphi(\vec{p})\,\nabla^{\mathcal P}\rho^2\rg|\le \frac{|\nabla^{\mathcal P}\varphi|}{\r^2(\vec{p})}+\frac{|\nabla^{\mathcal P}\rho|}{\r(\vec{p})}\le \frac{2}{\r(\vec{p})}\ .
\ee
Observe that we have respectively
\be
\label{df-2}
\chi\lf(\frac{\r(\vec{p})}{r}\rg)-\chi\lf(\frac{\r(\vec{p})}{\eta}\rg)\ge 0\ ,
\ee
moreover
\be
\label{df-3}
-\chi'\lf(\frac{\r(\vec{p})}{\eta}\rg)\ge 0\ ,
\ee
and
\be
\label{df-4}
\frac{\sigma\,\arctan\sigma}{\sqrt{1+\sigma^2}}(\vec{p})\ge 0\ .
\ee
Using (\ref{df-1}) we have respectively
\be
\label{df-5}
\begin{array}{l}
\ds\,4\,\lf|\int_{G_2^{Leg}(V_2({\R}^4))}\frac{\r(\vec{p})}{\eta}\, \chi'\lf(\frac{\r(\vec{p})}{\eta}\rg)\,\frac{|\nabla^{\mathcal P}\sigma|^2}{(1+\sigma^2(\vec{p}))^2}\ d{\mathbf v}_\infty({\mathcal P},\vec{p})\rg|\\[5mm]
\ds\quad\quad\le\, \frac{C}{\eta^2}\int_{ \eta<\r<2\eta}d\mu_\infty(\vec{p})
\end{array}
\ee
moreover
\be
\label{df-6}
\begin{array}{l}
\ds\lf|\frac{3}{4}\int_{G_2^{Leg}(V_2({\R}^4))}\,\frac{\r(\vec{p})}{\eta}\,\chi'\lf(\frac{\r(\vec{p})}{\eta}\rg)\, \,\arctan\sigma(\vec{p})\ \frac{\nabla^{\mathcal P}\r}{\r(\vec{p})}\cdot\frac{\nabla^{\mathcal P}\sigma}{1+\sigma^2(\vec{p})}\ d{\mathbf v}_\infty({\mathcal P},\vec{p})  \rg|\ \\[5mm]
\ds\quad\quad\le\, \frac{C}{\eta^2}\int_{ \eta<\r<2\eta}d\mu_\infty(\vec{p})
\end{array}
\ee
and
\be
\label{df-7}
\begin{array}{l}
\ds\frac{1}{4}\int_{G_2^{Leg}(V_2({\R}^4))}\frac{\r^2(\vec{p})}{\eta^{2}}\,\chi''\lf(\frac{\r(\vec{p})}{\eta}\rg)\, \,\arctan\sigma(\vec{p})\ \frac{\nabla^{\mathcal P}\r}{\r}\cdot\frac{\nabla^{\mathcal P}\sigma}{1+\sigma^2(\vec{p})} \ d{\mathbf v}_\infty({\mathcal P},\vec{p})\\[5mm]
\ds\quad\quad\le\, \frac{C}{\eta^2}\int_{ \eta<\r<2\eta}d\mu_\infty(\vec{p})
\end{array}
\ee
Combining (\ref{truncated-limit}), (\ref{df-5}), (\ref{df-6}) and (\ref{df-7}) we finally obtain
\be
\label{df-8}
\begin{array}{l}
\ds\limsup_{\eta\rightarrow 0}\int_{G_2^{Leg}(V_2({\R}^4))} 4\, \lf[\chi\lf(\frac{\r(\vec{p})}{r}\rg)-\chi\lf(\frac{\r(\vec{p})}{\eta}\rg)\rg]\,\frac{|(\nabla^{\mathcal P}\r)^\perp|^2}{\r(\vec{p})^2}\ d{\mathbf v}_\infty({\mathcal P},\vec{p})\\[5mm]
\ds-\int_{G_2^{Leg}(V_2({\R}^4))} \frac{\r(\vec{p})}{\eta}\,\chi'\lf(\frac{\r(\vec{p})}{\eta}\rg)\,\frac{1}{\r^2(\vec{p})}\,\lf[{|\nabla^{\mathcal P}\r|^2} +\frac{\sigma\,\arctan\sigma}{\sqrt{1+\sigma^2}}(\vec{p}) \rg]\ d{\mathbf v}_\infty({\mathcal P},\vec{p})<+\infty\ .
\end{array}
\ee
We then deduce for any $r<1$
\be
\label{df-9}
\lim_{\eta\rightarrow 0}\int_{\r>\eta} \chi\lf(\frac{\r(\vec{p})}{r}\rg)\,\frac{|(\nabla^{\mathcal P}\r)^\perp|^2}{\r(\vec{p})^2}\ d{\mathbf v}_\infty({\mathcal P},\vec{p})<+\infty\ .
\ee
At $\r=0$ we have that $\nabla v_i$ and $\nabla w_i$ for $i=3,4$ are orthogonal to $(\vec{\ep}_1,\vec{\ep}_2)$ and hence horizontal. We have moreover thanks to (\ref{almost-comp})
\[
J_H\nabla^Hv_i=J_H(\vec{\ep}_i,0)=(0,\vec{\ep}_i)=\nabla^Hw_i
\]
Then, in a neighbourhood of $(\vec{\ep}_1,\vec{\ep}_2)$ at the point $(\vec{\ep}_1+\vec{v},\vec{\ep}_2+\vec{w})$ we have
\[
J_H(\nabla^H \varphi)=J_H\lf(\sum_{i=1}^4 v_i\,\nabla^H w_i-w_i\,\nabla^Hv_i\rg)=-\rho\,\nabla^H\rho+O(\r^2)\ .
\]
This gives
\be
\label{id-sig-folland}
\begin{array}{l}
\ds\r^3\,J_H(\nabla^H\r)=\rho^3\,J_H(\nabla^H\rho)+2\,\varphi\,J_H(\nabla^H\varphi)=\rho^2\,\nabla^H \varphi-2\,\varphi\,\rho\,\nabla^H\rho+O(\r^4)\\[5mm]
\ds\quad=-\rho^4\,\nabla^H\lf(\frac{\varphi}{\rho^2}\rg)+O(\r^4)\ .
\end{array}
\ee
For every ${\mathcal P}\in G_2^{Leg}(V_2({\R}^4))$ we have that $J_H({\mathcal P})$ is also horizontal, orthogonal to $\mathcal P$ and
\[
H_{(\vec{a},\vec{b})}={\mathcal P}\oplus J_H({\mathcal P})\ .
\]
Hence we deduce from (\ref{id-sig-folland})
\[
-\rho^4\,\nabla^{\mathcal P}\lf(\frac{\varphi}{\rho^2}\rg)+O(\r^4)=\r^3\,J_H((\nabla^{\mathcal P}\r)^\perp)
\]
This gives
\be
\label{gg}
\frac{(\nabla^{\mathcal P}\r)^\perp}{\r}=\frac{1}{2}\,\frac{\rho^4}{\r^4}\, J_H(\nabla^{\mathcal P}\sigma)+O(1)=\frac{1}{2}\,J_H\lf(\frac{\nabla^{\mathcal P}\sigma}{1+\sigma^2}   \rg)+O(1)\ ,
\ee
from which we deduce
\[
\int_{\r>\eta} \chi\lf(\frac{\r(\vec{p})}{r}\rg)\,\frac{|(\nabla^{\mathcal P}\r)^\perp|^2}{\r(\vec{p})^2}\ d{\mathbf v}_\infty({\mathcal P},\vec{p})=\frac{1}{4}\,\int_{\r>\eta} \chi\lf(\frac{\r(\vec{p})}{r}\rg)\,\lf|\frac{\nabla^{\mathcal P}\sigma}{1+\sigma^2(\vec{p})}   \rg|^2
\ d{\mathbf v}_\infty({\mathcal P},\vec{p})+O(\mu_\infty(\{\r<2r\}))
\]
Combining this last identity with (\ref{df-9}) we get (\ref{bound-arctan-sigma}) and lemma~\ref{lm-dens-finie} is proved. \hfill$\Box$

\medskip

We deduce the following lemma

\begin{Lm}
\label{lm-density}
Assume
\be
\label{dens-finie-2}
\limsup_{\eta\rightarrow 0}\frac{1}{\eta^2}\int_{ \eta<\r<2\eta}d\mu_\infty(\vec{p})<+\infty\ ,
\ee
then for any function $\phi\in C^\infty_0({\R}_+,{\R}_+)$ such that supp$(\phi)\subset {\R}_+\setminus \{0\}$  and $\int_0^{+\infty}\phi(t)=1$
\be
\label{dens-finie-3}
\lim_{\eta\rightarrow 0}\frac{1}{\eta^2}\int_{V_2({\R}^4)}\,\frac{\eta}{\r(\vec{p})}\,\phi\lf(\frac{\r(\vec{p})}{\eta}\rg)\,\frac{1+\sigma\,\arctan\sigma}{\sqrt{1+\sigma^2}}(\vec{p}) \ d\mu_\infty=\theta_0
\ee
exists, is finite and is independent of $\phi$ .
\end{Lm}
\begin{Rm}
\label{rm-II.1}
Observe that
\[
\frac{d}{ds}\lf[   \frac{1+s\,\arctan s}{\sqrt{1+s^2}} \rg]=\frac{\arctan s}{(1+s^2)^{3/2}} \ .
\]
Hence the function $(1+s\,\arctan s)/\sqrt{1+s^2}$ is decreasing on ${\R}_-$ and increasing on ${\R}_+$ and we deduce
\[
\forall\ s\in{\R}\quad\quad\quad 1\le\frac{1+s\,\arctan s}{\sqrt{1+s^2}} \le \frac{\pi}{2}\ .
\]
\end{Rm}
\begin{Rm}
\label{rm-II.2}
Observe that if ${\mathbf v}_\infty$ is the varifold associated to a smooth Legendrian immersion $\Sigma$ we have $\sigma\rightarrow 0$ as $\r\rightarrow 0$ (see \cite{Riv} section III) hence
\[
\begin{array}{l}
\ds\theta_0=\lim_{\eta_\rightarrow 0}\frac{1}{\eta^2}\ \int_{V_2({\R}^4)}\frac{\eta}{\r(\vec{p})}\,\phi\lf(\frac{\r(\vec{p})}{\eta}\rg)\,\ d\mu_\infty=\lim_{\eta_\rightarrow 0}\frac{1}{\eta^2}\ \int_{V_2({\R}^4)}\frac{\eta}{\r(\vec{p})}\,\phi\lf(\frac{\r(\vec{p})}{\eta}\rg)\,|\nabla^\Sigma\r|\ d\mu_\infty\\[5mm]
\ds\quad=\lim_{\eta_\rightarrow 0}\ \int_{0}^\infty\,\phi\lf(\frac{s}{\eta}\rg)\,\frac{{\mathcal H}^1(\{\r=s\})}{s}\ \frac{ds}{\eta}=2\,\pi \, \mbox{Card}(\{\r^{-1}(\{0\})\})\ \int_0^\infty \phi(t)\ dt=2\,\pi \, \mbox{Card}(\{\r^{-1}(\{0\})\})\ .
\end{array}
\]
\end{Rm}
\noindent{\bf Proof of lemma~\ref{lm-density}.}
Let $d>1$ such that Supp$(\phi)\in [d^{-1},d]$. We introduce $\chi(t):=\int_t^{+\infty}\phi(s)\ ds$ and we apply (\ref{truncated-limit}). Because of (\ref{bound-arctan-sigma}) we have respectively
\[
\begin{array}{l}
\ds\lf|  \int_{G_2^{Leg}(V_2({\R}^4))}\frac{\r(\vec{p})}{\eta}\, \chi'\lf(\frac{\r(\vec{p})}{\eta}\rg)\,\frac{|\nabla^{\mathcal P}\sigma|^2}{(1+\sigma^2(\vec{p}))^2}\ d{\mathbf v}_\infty({\mathcal P},\vec{p}) \rg|\\[5mm]
\ds\quad\le C\, \int_{d^{-1}\eta<\r<d\,\eta}\frac{|\nabla^{\mathcal P}\sigma|^2}{(1+\sigma^2(\vec{p}))^2}\ d{\mathbf v}_\infty({\mathcal P},\vec{p}) \ \longrightarrow 0
\end{array}
\]
moreover
\[
\begin{array}{l}
\ds\lf|\int_{G_2^{Leg}(V_2({\R}^4))}\,\frac{\r(\vec{p})}{\eta}\,\chi'\lf(\frac{\r(\vec{p})}{\eta}\rg)\, \,\arctan\sigma(\vec{p})\ \frac{\nabla^{\mathcal P}\r}{\r(\vec{p})}\cdot\frac{\nabla^{\mathcal P}\sigma}{1+\sigma^2(\vec{p})}\ d{\mathbf v}_\infty({\mathcal P},\vec{p})\rg|  \\[5mm]
\ds\quad\le C\, \lf[\eta^{-2}\,\int_{d^{-1}\eta<\r<d\,\eta}\ d\mu_\infty\rg]^{1/2}\ \lf[ \int_{d^{-1}\eta<\r<d\,\eta}\frac{|\nabla^{\mathcal P}\sigma|^2}{(1+\sigma^2(\vec{p}))^2}\ d{\mathbf v}_\infty({\mathcal P},\vec{p})   \rg]^{1/2}\ \longrightarrow 0
\end{array}
\]
and
\[
\begin{array}{l}
\ds\lf|\int_{G_2^{Leg}(V_2({\R}^4))}\frac{\r^2(\vec{p})}{\eta^{2}}\,\chi''\lf(\frac{\r(\vec{p})}{\eta}\rg)\, \,\arctan\sigma(\vec{p})\ \frac{\nabla^{\mathcal P}\r}{\r(\vec{p})}\cdot\frac{\nabla^{\mathcal P}\sigma}{1+\sigma^2(\vec{p})} \ d{\mathbf v}_\infty({\mathcal P},\vec{p})\rg|
\\[5mm]
\ds\quad\le C\, \lf[\eta^{-2}\,\int_{d^{-1}\eta<\r<d\,\eta}\ d\mu_\infty\rg]^{1/2}\ \lf[ \int_{d^{-1}\eta<\r<d\,\eta}\frac{|\nabla^{\mathcal P}\sigma|^2}{(1+\sigma^2(\vec{p}))^2}\ d{\mathbf v}_\infty({\mathcal P},\vec{p})   \rg]^{1/2}\ \longrightarrow 0
\end{array}
\]
We have also
\[
\lf|\int_{G_2^{Leg}(V_2({\R}^4))}\ O\lf(\frac{\r(\vec{p})}{\eta}\rg)\,\chi'\lf(\frac{\r(\vec{p})}{\eta}\rg)\ \ d{\mathbf v}_\infty({\mathcal P},\vec{p})\rg|\le C\,\eta^2\ \lf[  \eta^{-2}\,\int_{d^{-1}\eta<\r<d\,\eta}\ d\mu_\infty\rg]\longrightarrow 0
\]
and
\[
\begin{array}{l}
\ds\int_{G_2^{Leg}(V_2({\R}^4))} \lf[O(1)\,\lf(\chi\lf(\frac{\r(\vec{p})}{r}\rg)-\chi\lf(\frac{\r(\vec{p})}{\eta}\rg)\rg)\rg] \ d{\mathbf v}_\infty({\mathcal P},\vec{p})\\[5mm]
\ds=\int_{G_2^{Leg}(V_2({\R}^4))\setminus\pi^{-1}\{\r=0\}}\,O(1)\,\chi\lf(\frac{\r(\vec{p})}{r}\rg)\ d{\mathbf v}_\infty({\mathcal P},\vec{p})\\[5mm]
\ds-\sum_{j=0}^\infty\int_{G_2^{Leg}(V_2({\R}^4))}O(1)\,\lf[ \chi\lf(\frac{\r(\vec{p})}{2^{-j}\eta}\rg)-    \chi\lf(\frac{\r(\vec{p})}{2^{-j-1}\eta}\rg)\rg] \ d{\mathbf v}_\infty({\mathcal P},\vec{p})
\end{array}
\]
We have
\[
\begin{array}{l}
\ds\lf|\sum_{j=0}^\infty\int_{G_2^{Leg}(V_2({\R}^4))}O(1)\,\lf[ \chi\lf(\frac{\r(\vec{p})}{2^{-j}\eta}\rg)-    \chi\lf(\frac{\r(\vec{p})}{2^{-j-1}\eta}\rg)\rg] \ d{\mathbf v}_\infty({\mathcal P},\vec{p})\rg|\\[5mm]
\ds\le C\,\eta^2\, \sum_{j=0}^\infty 2^{-2j}\ (2^{-j}\eta)^{-2}\,\int_{d^{-1}2^{-j-1}\eta<\r<d\,2^{-j}\eta}\ d\mu_\infty\ \longrightarrow 0
\end{array}
\]
Finally we observe that 
\[
4\,\int_{G_2^{Leg}(V_2({\R}^4))} \lf[\chi\lf(\frac{\r(\vec{p})}{r}\rg)-\chi\lf(\frac{\r(\vec{p})}{\eta}\rg)\rg]\,\frac{|(\nabla^{\mathcal P}\r)^\perp|^2}{\r(\vec{p})^2}\ d{\mathbf v}_\infty({\mathcal P},\vec{p})
\]
is a monotone function of $\eta$. Hence
\[
\begin{array}{l}
\ds\lim_{\eta\rightarrow 0}4\,\int_{G_2^{Leg}(V_2({\R}^4))} \lf[\chi\lf(\frac{\r(\vec{p})}{r}\rg)-\chi\lf(\frac{\r(\vec{p})}{\eta}\rg)\rg]\,\frac{|(\nabla^{\mathcal P}\r)^\perp|^2}{\r(\vec{p})^2}\ d{\mathbf v}_\infty({\mathcal P},\vec{p})\\[5mm]
\ds=4\,\int_{G_2^{Leg}(V_2({\R}^4))\setminus\pi^{-1}(\{0\})} \chi\lf(\frac{\r(\vec{p})}{r}\rg)\,\frac{|(\nabla^{\mathcal P}\r)^\perp|^2}{\r(\vec{p})^2}\ d{\mathbf v}_\infty({\mathcal P},\vec{p})\ .
\end{array}
\]
 Combining all the previous we obtain that
\[
-\int_{G_2^{Leg}(V_2({\R}^4))} \frac{\r(\vec{p})}{\eta}\,\chi'\lf(\frac{\r(\vec{p})}{\eta}\rg)\,\frac{1}{\r^2(\vec{p})}\,\lf[{|\nabla^{\mathcal P}\r|^2} +\frac{\sigma\,\arctan\sigma}{\sqrt{1+\sigma^2}}(\vec{p}) \rg]\ \ d{\mathbf v}_\infty({\mathcal P},\vec{p})
\]
is converging as $\eta$ converges to zero. Hence
\be
\label{dens-finie-3-b}
\lim_{\eta\rightarrow 0}\int_{G_2^{Leg}(V_2({\R}^4))} \frac{\r(\vec{p})}{\eta}\,\phi\lf(\frac{\r(\vec{p})}{\eta}\rg)\,\frac{1}{\r^2(\vec{p})}\,\lf[{|\nabla^{\mathcal P}\r|^2} +\frac{\sigma\,\arctan\sigma}{\sqrt{1+\sigma^2}}(\vec{p}) \rg]\ d{\mathbf v}_\infty({\mathcal P},\vec{p})=\theta_0(\phi)
\ee
exists and is finite. Recall that for any ${\mathcal P}\in G_2^{Leg}(V_2({\R}^4))$
\[
|\nabla^H\r|^2=|\nabla^{\mathcal P}\r|^2+|(\nabla^{\mathcal P}\r)^\perp|^2
\]
From (\ref{gg}) and using also (\ref{horizont-gradient}) we have
\[
\frac{|\nabla^{\mathcal P}\r|^2}{\r^2}=\frac{|\nabla^H\r|^2}{\r^2}-\frac{1}{4}\,\frac{|\nabla^{\mathcal P}\sigma|^2}{(1+\sigma^2)^2}+O(1)=\frac{1}{\r^2}\,\frac{1}{\sqrt{1+\sigma^2}}-\frac{1}{4}\,\frac{|\nabla^{\mathcal P}\sigma|^2}{(1+\sigma^2)^2}+O(1)
\]
We deduce
\[
\lim_{\eta\rightarrow 0}\int_{G_2^{Leg}(V_2({\R}^4))} \frac{\r(\vec{p})}{\eta}\,\phi\lf(\frac{\r(\vec{p})}{\eta}\rg)\,\frac{1}{\r^2(\vec{p})}\,\lf[\frac{1+\sigma\,\arctan\sigma}{\sqrt{1+\sigma^2}}(\vec{p}) \rg]\ d{\mathbf v}_\infty({\mathcal P},\vec{p})=\theta_0(\phi)
\]
By taking $h_{r,\eta}:=[\chi(\r/r)-\int^{+\infty}_{\r/\eta}\phi(s)\ ds] \arctan\sigma$ where $\chi$ is fixed and $\phi$ is arbitrary satisfying the assumption of the lemma
we observe that
\[
\begin{array}{l}
\ds\theta_0(\phi)=-4\,\int_{G_2^{Leg}(V_2({\R}^4))\setminus\pi^{-1}(\{0\})} \chi\lf(\frac{\r(\vec{p})}{r}\rg)\,\frac{|(\nabla^{\mathcal P}\r)^\perp|^2}{\r(\vec{p})^2}\ d{\mathbf v}_\infty({\mathcal P},\vec{p})\\[5mm]
\ds-\int_{G_2^{Leg}(V_2({\R}^4))}\ \frac{\r(\vec{p})}{r}\,\chi'\lf(\frac{\r(\vec{p})}{r}\rg)\,\frac{1}{\r^2(\vec{p})}\ \lf[{|\nabla^{\mathcal P}\r|^2}+ \frac{\sigma(\vec{p})\,\arctan\sigma(\vec{p})}{\sqrt{1+\sigma^2(\vec{p})}} \rg]\ d{\mathbf v}_\infty({\mathcal P},\vec{p}) \\[5mm]
\ds+\int_{G_2^{Leg}(V_2({\R}^4))\setminus\pi^{-1}(\{0\})} \lf[O(1)\,\lf(\chi\lf(\frac{\r(\vec{p})}{r}\rg)\rg)+O\lf(\frac{\r(\vec{p})}{r}\rg)\,\chi'\lf(\frac{\r(\vec{p})}{r}\rg)\rg]\ d{\mathbf v}_\infty({\mathcal P},\vec{p})\\[5mm]
\ds+\frac{1}{4}\,\int_{G_2^{Leg}(V_2({\R}^4))}\frac{\r^2(\vec{p})}{r^{2}}\,\chi''\lf(\frac{\r(\vec{p})}{r}\rg) \,\arctan\sigma(\vec{p})\ \frac{\nabla^{\mathcal P}\r}{\r}\cdot\frac{\nabla^{\mathcal P}\sigma}{1+\sigma^2}\ d{\mathbf v}_\infty({\mathcal P},\vec{p}) \\[5mm]
\ds-\frac{3}{4}\,\int_{G_2^{Leg}(V_2({\R}^4))}  \ \frac{\r(\vec{p})}{r}\,\chi'\lf(\frac{\r(\vec{p})}{r}\rg) \,\arctan\sigma(\vec{p})\,\frac{\nabla^{\mathcal P}\r}{\r}\cdot\frac{\nabla^{\mathcal P}\sigma}{1+\sigma^2(\vec{p})}\ \ d{\mathbf v}_\infty({\mathcal P},\vec{p}) \\[5mm]
\ds+\frac{1}{4}\,\int_{G_2^{Leg}(V_2({\R}^4))} \frac{\r(\vec{p})}{r}\chi'\lf(\frac{\r(\vec{p})}{r}\rg)\,\frac{|\ds\nabla^{\mathcal P}\sigma|^2}{(1+\sigma^2(\vec{p}))^2}\ \ d{\mathbf v}_\infty({\mathcal P},\vec{p})
\end{array}
\]
Obviously $\theta_0(\phi)$ is independent of $\phi$ and this concludes the proof of lemma~\ref{lm-density}.\hfill$\Box$

\medskip

We have the following lemma
\begin{Lm}
\label{lm-quasi-mono}
There exists a universal constant $C>1$ such that, for any $0<2s<r<1$, assuming we have
\be
\label{dens-finie-2-2}
\limsup_{\eta\rightarrow 0}\frac{1}{\eta^2}\int_{ \eta<\r<2\eta}d\mu_\infty(\vec{p})<+\infty\ ,
\ee
then
\be
\label{mono}
C^{-1} \lf[\theta_0+\int_{0<\r<s/2}\lf|\frac{\nabla^{\mathcal P}\sigma}{1+\sigma^2(\vec{p})}   \rg|^2
\ d{\mathbf v}_\infty({\mathcal P},\vec{p})\rg]\le \frac{1}{s^2}\ \int_{0<\r<s}\ d\mu_\infty\le C\  \frac{1}{r^2}\ \int_{\r<r}\ d\mu_\infty\ .
\ee
\end{Lm}
We replace $r$ by $s$ in (\ref{truncated-limit}) and we make $\eta$ converge to zero. Thanks to the previous lemma we obtain
\be
\label{truncated-limit-2}
\begin{array}{l}
\ds-\int_{G_2^{Leg}(V_2({\R}^4))}\ \frac{\r(\vec{p})}{s}\,\chi'\lf(\frac{\r(\vec{p})}{s}\rg)\,\frac{1}{\r^2(\vec{p})}\ \lf[{|\nabla^{\mathcal P}\r|^2}+ \frac{\sigma(\vec{p})\,\arctan\sigma(\vec{p})}{\sqrt{1+\sigma^2(\vec{p})}} \rg]\ d{\mathbf v}_\infty({\mathcal P},\vec{p}) \\[5mm]
\ds+\int_{G_2^{Leg}(V_2({\R}^4))\setminus\{\r=0\}} \lf[O(1)\,\chi\lf(\frac{\r(\vec{p})}{s}\rg)+O\lf(\frac{\r(\vec{p})}{s}\rg)\,\chi'\lf(\frac{\r(\vec{p})}{s}\rg)\rg]\ d{\mathbf v}_\infty({\mathcal P},\vec{p})\\[5mm]
\ds+\frac{1}{4}\,\int_{G_2^{Leg}(V_2({\R}^4))}\frac{\r^2(\vec{p})}{s^{2}}\,\chi''\lf(\frac{\r(\vec{p})}{s}\rg) \,\arctan\sigma(\vec{p})\ \frac{\nabla^{\mathcal P}\r}{\r}\cdot\frac{\nabla^{\mathcal P}\sigma}{1+\sigma^2}\ d{\mathbf v}_\infty({\mathcal P},\vec{p}) \\[5mm]
\ds-\frac{3}{4}\,\int_{G_2^{Leg}(V_2({\R}^4))}  \ \frac{\r(\vec{p})}{s}\,\chi'\lf(\frac{\r(\vec{p})}{s}\rg) \,\arctan\sigma(\vec{p})\,\frac{\nabla^{\mathcal P}\r}{\r}\cdot\frac{\nabla^{\mathcal P}\sigma}{1+\sigma^2(\vec{p})}\ \ d{\mathbf v}_\infty({\mathcal P},\vec{p}) \\[5mm]
\ds+\frac{1}{4}\,\int_{G_2^{Leg}(V_2({\R}^4))} \frac{\r(\vec{p})}{s}\chi'\lf(\frac{\r(\vec{p})}{s}\rg)\,\frac{|\ds\nabla^{\mathcal P}\sigma|^2}{(1+\sigma^2(\vec{p}))^2}\ \ d{\mathbf v}_\infty({\mathcal P},\vec{p}) \\[5mm]
\ds=4\,\int_{G_2^{Leg}(V_2({\R}^4))\setminus\{\r=0\}} \chi\lf(\frac{\r(\vec{p})}{s}\rg)\,\frac{|(\nabla^{\mathcal P}\r)^\perp|^2}{\r(\vec{p})^2}\ d{\mathbf v}_\infty({\mathcal P},\vec{p})
+\theta_0
\end{array}
\ee
Using thee fact that we have 
\[
|\nabla^{\mathcal P}\r|\le 1\quad\mbox{ and }\quad \lf|\frac{\nabla^{\mathcal P}\sigma}{1+\sigma^2(\vec{p})}\rg|\le \frac{2}{\r}
\]
we deduce
\[
\begin{array}{l}
\ds \theta_0+4\,\int_{G_2^{Leg}(V_2({\R}^4))\setminus\{\r=0\}} \chi\lf(\frac{\r(\vec{p})}{s}\rg)\,\frac{|(\nabla^{\mathcal P}\r)^\perp|^2}{\r(\vec{p})^2}\ d{\mathbf v}_\infty({\mathcal P},\vec{p})\le \frac{C}{s^2}\ \int_{s<\r<2s}\ d\mu_\infty+
\ \int_{0<\r<2s}\ d\mu_\infty
\end{array}
\]
Using
\[
\frac{|(\nabla^{\mathcal P}\r)^\perp|^2}{\r^2}=\frac{1}{4}\,\frac{|\nabla^{\mathcal P}\sigma|^2}{(1+\sigma^2)^2}+O(1)
\]
we deduce the first inequality in (\ref{mono}). For $2s<r$ we have
\be
\label{truncated-limit-3}
\begin{array}{l}
\ds-\int_{G_2^{Leg}(V_2({\R}^4))}\ \frac{\r(\vec{p})}{r}\,\chi'\lf(\frac{\r(\vec{p})}{r}\rg)\,\frac{1}{\r^2(\vec{p})}\ \lf[{|\nabla^{\mathcal P}\r|^2}+ \frac{\sigma(\vec{p})\,\arctan\sigma(\vec{p})}{\sqrt{1+\sigma^2(\vec{p})}} \rg]\ d{\mathbf v}_\infty({\mathcal P},\vec{p}) \\[5mm]
\ds+\int_{G_2^{Leg}(V_2({\R}^4))} \lf[O(1)\,\lf(\chi\lf(\frac{\r(\vec{p})}{r}\rg)-\chi\lf(\frac{\r(\vec{p})}{\eta}\rg)\rg)+O\lf(\frac{\r(\vec{p})}{r}\rg)\,\chi'\lf(\frac{\r(\vec{p})}{r}\rg)\rg]\ d{\mathbf v}_\infty({\mathcal P},\vec{p})\\[5mm]
\ds+\frac{1}{4}\,\int_{G_2^{Leg}(V_2({\R}^4))}\frac{\r^2(\vec{p})}{r^{2}}\,\chi''\lf(\frac{\r(\vec{p})}{r}\rg) \,\arctan\sigma(\vec{p})\ \frac{\nabla^{\mathcal P}\r}{\r}\cdot\frac{\nabla^{\mathcal P}\sigma}{1+\sigma^2}\ d{\mathbf v}_\infty({\mathcal P},\vec{p}) \\[5mm]
\ds-\frac{3}{4}\,\int_{G_2^{Leg}(V_2({\R}^4))}  \ \frac{\r(\vec{p})}{r}\,\chi'\lf(\frac{\r(\vec{p})}{r}\rg) \,\arctan\sigma(\vec{p})\,\frac{\nabla^{\mathcal P}\r}{\r}\cdot\frac{\nabla^{\mathcal P}\sigma}{1+\sigma^2(\vec{p})}\ \ d{\mathbf v}_\infty({\mathcal P},\vec{p}) \\[5mm]
\ds+\frac{1}{4}\,\int_{G_2^{Leg}(V_2({\R}^4))} \frac{\r(\vec{p})}{r}\chi'\lf(\frac{\r(\vec{p})}{r}\rg)\,\frac{|\ds\nabla^{\mathcal P}\sigma|^2}{(1+\sigma^2(\vec{p}))^2}\ \ d{\mathbf v}_\infty({\mathcal P},\vec{p}) \\[5mm]
\ds=4\,\int_{G_2^{Leg}(V_2({\R}^4))} \lf[\chi\lf(\frac{\r(\vec{p})}{r}\rg)-\chi\lf(\frac{\r(\vec{p})}{s}\rg)\rg]\,\frac{|(\nabla^{\mathcal P}\r)^\perp|^2}{\r(\vec{p})^2}\ d{\mathbf v}_\infty({\mathcal P},\vec{p})\\[5mm]
\ds-\int_{G_2^{Leg}(V_2({\R}^4))} \frac{\r(\vec{p})}{s}\,\chi'\lf(\frac{\r(\vec{p})}{s}\rg)\,\frac{1}{\r^2(\vec{p})}\,\lf[{|\nabla^{\mathcal P}\r|^2} +\frac{\sigma\,\arctan\sigma}{\sqrt{1+\sigma^2}}(\vec{p}) \rg]\ \ d{\mathbf v}_\infty({\mathcal P},\vec{p})\\[5mm]
\ds+\,4\,\int_{G_2^{Leg}(V_2({\R}^4))}\frac{\r(\vec{p})}{s}\, \chi'\lf(\frac{\r(\vec{p})}{s}\rg)\,\frac{|\nabla^{\mathcal P}\sigma|^2}{(1+\sigma^2(\vec{p}))^2}\ d{\mathbf v}_\infty({\mathcal P},\vec{p})\\[5mm]
\ds-\frac{3}{4}\int_{G_2^{Leg}(V_2({\R}^4))}\,\frac{\r(\vec{p})}{s}\,\chi'\lf(\frac{\r(\vec{p})}{s}\rg)\, \,\arctan\sigma(\vec{p})\ \frac{\nabla^{\mathcal P}\r}{\r(\vec{p})}\cdot\frac{\nabla^{\mathcal P}\sigma}{1+\sigma^2(\vec{p})}\ d{\mathbf v}_\infty({\mathcal P},\vec{p})  \\[5mm]
\ds+\frac{1}{4}\int_{G_2^{Leg}(V_2({\R}^4))}\frac{\r^2(\vec{p})}{s^{2}}\,\chi''\lf(\frac{\r(\vec{p})}{s}\rg)\, \,\arctan\sigma(\vec{p})\ \frac{\nabla^{\mathcal P}\r}{\r(\vec{p})}\cdot\frac{\nabla^{\mathcal P}\sigma}{1+\sigma^2(\vec{p})} \ d{\mathbf v}_\infty({\mathcal P},\vec{p})\\[5mm]
\ds+\int_{G_2^{Leg}(V_2({\R}^4))}\ \lf[O\lf(\frac{\r(\vec{p})}{s}\rg)\,\chi'\lf(\frac{\r(\vec{p})}{s}\rg)\rg]\ \ d{\mathbf v}_\infty({\mathcal P},\vec{p})
\end{array}
\ee
We have that $\chi'\le 0$ on ${\R}$ and supp$(\chi')\subset [1,2]$ hence we deduce
\[
\begin{array}{l}
\ds-\int_{G_2^{Leg}(V_2({\R}^4))} \frac{\r(\vec{p})}{s}\,\chi'\lf(\frac{\r(\vec{p})}{s}\rg)\,\frac{1}{\r^2(\vec{p})}\,\lf[{|\nabla^{\mathcal P}\r|^2} +\frac{\sigma\,\arctan\sigma}{\sqrt{1+\sigma^2}}(\vec{p}) \rg]\ \ d{\mathbf v}_\infty({\mathcal P},\vec{p})\\[5mm]
\ds\le \frac{C}{r^2}\ \int_{r<\r<2r}\ d\mu_\infty+
\ \int_{0<\r<2r}\ d\mu_\infty+  C\,\lf[   \frac{1}{s^2}\int_{s<\r<2s}\ d\mu_\infty\rg]^{1/2}\ \lf[ \int_{s<\r<2s}\ \frac{|\nabla^{\mathcal P}\sigma|^2}{(1+\sigma^2(\vec{p}))^2}\ d{\mathbf v}_\infty({\mathcal P},\vec{p})  \rg]^{1/2}
\end{array}
\]
We choose $\chi$ such that
\[
\chi'\le 0\quad\mbox{ on }{\R}\quad\mbox{and }\quad -\,\chi'> \frac{1}{2}\ \mbox{ on }[5/4, 7/4]\ .
\]
Since
\[
\frac{|\nabla^{\mathcal P}\r|^2}{\r^2}=\frac{|\nabla^H\r|^2}{\r^2}-\frac{1}{4}\,\frac{|\nabla^{\mathcal P}\sigma|^2}{(1+\sigma^2)^2}+O(1)=\frac{1}{\r^2}\,\frac{1}{\sqrt{1+\sigma^2}}-\frac{1}{4}\,\frac{|\nabla^{\mathcal P}\sigma|^2}{(1+\sigma^2)^2}+O(1)\ .
\]
Using remark~\ref{rm-II.1} we deduce
\[
\frac{1}{\r^2(\vec{p})}\,\lf[{|\nabla^{\mathcal P}\r|^2} +\frac{\sigma\,\arctan\sigma}{\sqrt{1+\sigma^2}}(\vec{p}) \rg]\ge \frac{1}{\r^2}-\frac{1}{4}\,\frac{|\nabla^{\mathcal P}\sigma|^2}{(1+\sigma^2)^2}+O(1)\ .
\]
We deduce that for $r<1$
\[
\begin{array}{l}
\ds \frac{1}{s^2}\int_{5s/4<\r<7s/4}\ d\mu_\infty\le \frac{C}{r^2}\ \int_{0<\r<2r}\ d\mu_\infty+  C\,\lf[   \frac{1}{s^2}\int_{s<\r<2s}\ d\mu_\infty\rg]^{1/2}\ \lf[ \int_{s<\r<2s}\ \frac{|\nabla^{\mathcal P}\sigma|^2}{(1+\sigma^2(\vec{p}))^2}\ d{\mathbf v}_\infty({\mathcal P},\vec{p})  \rg]^{1/2}\\[5mm]
 \ds +C\, \int_{s<\r<2s}\ \frac{|\nabla^{\mathcal P}\sigma|^2}{(1+\sigma^2(\vec{p}))^2}\ d{\mathbf v}_\infty({\mathcal P},\vec{p}) \ .
 \end{array}
\]
Let 
\[
A:=\sup_{2\,s<r} \frac{1}{s^2}\int_{5s/4<\r<7s/4}\ d\mu_\infty
\]
Using the first inequality in (\ref{mono}) we just proved, we deduce
\[
A\le  \frac{C}{r^2}\ \int_{0<\r<2r}\ d\mu_\infty+ C\  A^{1/2} \lf[ r^{-2}\int_{0<\r<2r}\ d\mu_\infty  \rg]^{1/2}
\]
Hence for any $s<r$ we have
\[
\frac{1}{s^2}\int_{5s/4<\r<7s/4}\ d\mu_\infty\le  \frac{C}{r^2}\ \int_{0<\r<2r}\ d\mu_\infty\ .
\]
This gives for any $s<r$
\[
\frac{1}{s^2}\int_{0<\r<s}\ d\mu_\infty=\frac{1}{s^2}\sum_{j=0}^\infty\int_{5^{j+1}s/7^{j+1}<\r<5^js/7^j}\ d\mu_\infty\le C\ \sum_{j=0}^\infty \lf(\frac{5}{7}\rg)^j\  \frac{1}{r^2}\ \int_{0<\r<2r}\ d\mu_\infty
\]
which implies the second inequality in (\ref{mono}). Lemma~\ref{lm-quasi-mono} is proved. \hfill$\Box$


\subsection{The finiteness and a global bound of the upper 2-density of $\mu_\infty$ at every point in $V_2({\R}^4)$.}

We denote by $\r_{\vec{p}_0}$ the Folland Kor\'anyi gauge with respect to $\vec{p}_0=(\vec{a}_0,\vec{b}_0)\in V_2({\R}^4)$. That is, for any $\vec{p}=(\vec{a},\vec{b})\in V_2({\R}^4)$
\[
\rho^2_{\vec{p}_0}(\vec{p})=|\vec{p}-\vec{p}_0|^2\quad,\quad \varphi_{\vec{p}_0}(\vec{p}):=\vec{a}\cdot \vec{b}_0-\vec{a}_0\cdot\vec{b}\quad\mbox{ and }\quad \r_{\vec{p}_0}^4:=\rho^4_{\vec{p}_0}+4\,\varphi_{\vec{p}_0}^2\ .
\]
For $\theta\in{\R}$ and $\vec{p}=(\vec{a},\vec{b})\in V_2({\R}^4)$ we denote
\[
\vec{p}^{\, \theta}:=(\cos\theta\,\vec{a}+\sin\theta\,\vec{b},-\sin\theta\,\vec{a}+\cos\theta\,\vec{b})\ .
\]
We shall now prove the following lemma
\begin{Lm}
\label{lm-finite-dens}
For any $\vec{p}_0\in V_2({\R}^4)$ the following holds
\be
\label{fin-dens}
\limsup_{r\rightarrow 0}\frac{1}{r^2}\int_{\r_{\vec{p}_0}<r} d\mu_\infty(\vec{p})<+\infty\ .
\ee
Hence there exists $C>1$ universal such that for any $\vec{p}_0\in V_2({\R}^4)$ and any $0<s<r$
\be
\label{mono-every}
C^{-1} \lf[\theta_0+\int_{0<\r_{\vec{p}_0}<s/2}\lf|\frac{\nabla^{\mathcal P}\sigma}{1+\sigma^2(\vec{p})}   \rg|^2
\ d{\mathbf v}_\infty({\mathcal P},\vec{p})\rg]\le \frac{1}{s^2}\ \int_{\r_{\vec{p}_0}<s}\ d\mu_\infty\le C\  \frac{1}{r^2}\ \int_{\r_{\vec{p}_0}<r}\ d\mu_\infty\ .
\ee
\end{Lm}
\noindent{\bf Proof of lemma~\ref{lm-finite-dens}.} Let $M>0$ and denote $\vec{p}$ such that
\[
E_{\infty}:=\lf\{\vec{p}\in V_2({\R}^4)\ ;\ \limsup_{r\rightarrow 0}\frac{1}{r^2}\int_{\r_{\vec{p}}<r} d\mu_\infty=+\infty\rg\}\ .
\]
Let $M>0$ and $\delta>0$. For any $\vec{p}\in E_\infty$ we choose $\delta>r_{\vec{p}}>0$ such that
\[
\frac{1}{r^2_{\vec{p}}}\int_{\r_{\vec{p}}<r_{\vec{p}}} d\mu_\infty>{M}\ .
\]
Let
\[
{\mathfrak M}_\infty:=\Pi^\ast\mu_\infty\ .
\]
where we recall that $\Pi$ is the tautological projection from $V_2({\R}^4)$ into $G_2({\R}^4)$ such that $\Pi(\vec{a},\vec{b})=\vec{a}\wedge\vec{b}$. Let $\vec{G}_{\vec{p}}:=\Pi(\vec{p})$.
We have for any pair $(\vec{a},\vec{b})$ and $(\vec{c},\vec{d})$ in $V_2({\R}^4)$
\[
|\vec{a}\wedge\vec{b}-\vec{c}\wedge\vec{d}|=|(\vec{a}-\vec{c})\wedge\vec{b}-\vec{c}\wedge(\vec{d}-\vec{b})|\le |\vec{a}-\vec{c}|+|\vec{d}-\vec{b}|\le \sqrt{2}\,\lf|(\vec{a},\vec{b})-(\vec{c},\vec{d})\rg|
\]
Let $\vec{G}_{\vec{p}}:=\Pi(\vec{p})$, because of the previous
\[
\Pi\lf(\{ \vec{q}\in V_2({\R}^4)\ ;\  \r_{\vec{p}}(\vec{q})<r_{\vec{p}}\}\rg)\subset B_{\sqrt{2} r_{\vec{p}}}(\vec{G}_{\vec{p}})
\]
Hence we have
\[
\frac{1}{r^2_{\vec{p}}}\int_{B_{\sqrt{2} r_{\vec{p}}}(\vec{G}_{\vec{p}})}d{\mathfrak M}_\infty>{M}
\]
Obviously $B_{\sqrt{2} r_{\vec{p}}}(\vec{G}_{\vec{p}})$ realizes an open cover of $\Pi(E_\infty)$. Extracting a Besicovitch covering $(B_{\sqrt{2} r_{i}}(\vec{G}_{i}))_{i\in I}$, we deduce that
\[
M\,\sum_{i\in I}r_i^2\le C\ \int_{G_2({\R}^4)}d{\mathfrak M}_\infty
\]
where $C$ is a constant depending only on $G_2({\R}^4)$ given by the Besicovitch covering theorem. Since $M>0$ and $\delta>0$ are arbitrary, we deduce
\[
{\mathcal H}^2(\Pi(E_{\infty}))=0\ .
\]
This implies in particular that $E_\infty$ has no interior. Recall that from lemma~\ref{lm-quasi-mono}, for any $\vec{p}\in V_2({\R}^4)\setminus E_\infty$, we have
\[
\forall\ r<1\quad\quad \frac{1}{r^2}\int_{0<\r_{\vec{p}}<r}d\mu_\infty\le C\, \mu_\infty(V_2({\R}^4))\ .
\]
where $C>0$. Moreover, since by assumption 
\[
\limsup_{r\rightarrow 0}\frac{1}{r^2}\int_{\r_{\vec{p}}<r} d\mu_\infty<+\infty\ ,
\]
we have
\[
\forall\,\vec{p}\in V_2({\R}^4)\setminus E_\infty\quad\quad \mu_\infty(\{\vec{p}\})=0\ .
\]
Hence 
\[
\forall\,\vec{p}\in V_2({\R}^4)\setminus E_\infty\quad\quad \forall\ r<1\quad\quad \frac{1}{r^2}\int_{\r_{\vec{p}}<r}d\mu_\infty\le C\, \mu_\infty(V_2({\R}^4))\ .
\]
Let now $\vec{p}_0=(\vec{a}_0,\vec{b}_0)\in E_\infty$ and $\vec{p}_k=(\vec{a}_k,\vec{b}_k)\rightarrow \vec{p}_\infty$ with $\vec{p}_k\in V_2({\R}^4)\setminus E_\infty$. Let $r>0$ and $k\in {\N}$ such that $\r_{\vec{p}_0}(\vec{p}_k)<2^{-1}\,C_0^{-1}\,r$ where $C_0>1$ is going to be fixed later on below. We have for any $\vec{p}\in V_2({\R}^4)$
\[
\r_{\vec{p}_0}(\vec{p})=\lf[|\vec{p}-\vec{p}_0|^4+4 \,\varphi^2_{\vec{p}_0}(\vec{p})\rg]^{1/4}=\lf[|\vec{p}-\vec{p}_0|^4+4 \,(\vec{a}\cdot\vec{b}_0-\vec{a}_0\cdot\vec{b})^2\rg]^{1/4}
\]
We write
\[
\begin{array}{l}
\ds\varphi_{\vec{p}_0}(\vec{p})=\vec{a}\cdot\vec{b}_0-\vec{b}\cdot\vec{a}_0=(\vec{a}-\vec{a}_0)\cdot\vec{b}_0-(\vec{b}-\vec{b}_0)\cdot\vec{a}_0\\[5mm]
\ds\quad=(\vec{a}-\vec{a}_0)\cdot(\vec{b}_0-\vec{b}_k)-(\vec{b}-\vec{b}_0)\cdot(\vec{a}_0-\vec{a}_k)
\ds+(\vec{a}-\vec{a}_0)\cdot\vec{b}_k-(\vec{b}-\vec{b}_0)\cdot\vec{a}_k\\[5mm]
\ds=(\vec{a}-\vec{a}_0)\cdot(\vec{b}_0-\vec{b}_k)-(\vec{b}-\vec{b}_0)\cdot(\vec{a}_0-\vec{a}_k)+\varphi_{\vec{p}_k}(\vec{p})+\varphi_{\vec{p}_0}(\vec{p}_k)
\end{array}
\]
Hence
\[
\begin{array}{l}
\ds|\varphi_{\vec{p}_0}(\vec{p})|\le \rho_{\vec{p}_0}(\vec{p})\ \rho_{\vec{p}_0}(\vec{p}_k)+|\varphi_{\vec{p}_k}(\vec{p})|+|\varphi_{\vec{p}_0}(\vec{p}_k)|\\[5mm]
\ds\le\rho_{\vec{p}_0}(\vec{p}_k)\ \rho_{\vec{p}_0}(\vec{p}_k)+\rho_{\vec{p}_k}(\vec{p})\ \rho_{\vec{p}_0}(\vec{p}_k)+|\varphi_{\vec{p}_k}(\vec{p})|+|\varphi_{\vec{p}_0}(\vec{p}_k)|
\end{array}
\]
Hence
\be
\label{quasi-dist}
\begin{array}{l}
\ds\r_{\vec{p}_0}(\vec{p})\le \rho_{\vec{p}_0}(\vec{p})+\sqrt{2\,|\varphi_{\vec{p}_0}(\vec{p})|}\\[5mm]
\ds \le \rho_{\vec{p}_0}(\vec{p}_k)+\rho_{\vec{p}_k}(\vec{p})+\sqrt{2}\,  \rho_{\vec{p}_0}(\vec{p}_k)+\sqrt{2\,\rho_{\vec{p}_k}(\vec{p})\ \rho_{\vec{p}_0}(\vec{p}_k)}+\sqrt{2\,|\varphi_{\vec{p}_k}(\vec{p})|}+\sqrt{2\,|\varphi_{\vec{p}_0}(\vec{p}_k)|}\\[5mm]
\ds\le C_0\ [\r_{\vec{p}_k}(\vec{p})+\r_{\vec{p}_0}(\vec{p}_k)]
\end{array}
\ee
where $C_0>1$ is universal. We choose $\vec{p}_k$ such that $\r_{\vec{p}_0}(\vec{p}_k)<2^{-1}\,C_0^{-1}\,r$. Hence
\[
\lf\{\vec{p}\in V_2({\R}^4)\ ;\ \r_{\vec{p}_0}(\vec{p})<r\rg\}\subset \lf\{\vec{p}\in V_2({\R}^4)\ ;\ \r_{\vec{p}_k}(\vec{p})<2^{-1}\,C_0^{-1}\,r\rg\}
\]
Since $\vec{p}_k\in V_2({\R}^4)\setminus E_\infty$ we have
\[
r^{-2}\int_{ \r_{\vec{p}_k}(\vec{p})<2^{-1}\,C_0^{-1}\,r}d\mu_\infty\le C\, \mu_\infty(V_2({\R}^4))
\]
where $C$ is universal. Hence we deduce
\[
r^{-2}\int_{ \r_{\vec{p}_0}(\vec{p})<\,r}d\mu_\infty\le C\, \mu_\infty(V_2({\R}^4))\ .
\]
This holds for any $0<r<1$ which contradicts the fact that $\vec{p}_0\in E_\infty$. Hence $E_\infty=\emptyset$ and lemma~\ref{lm-finite-dens} is proved. \hfill $\Box$.
\subsection{Bounding the limiting 2-density of $\mu_\infty$ from below}
We take again
\[
h_{r,\eta}:=[\chi(\r/r)-\chi(\r/\eta)]\, \arctan\sigma\ .
\]
we have
\[
\nabla^H\arctan\sigma=\frac{\rho^2}{\r^4}\,\nabla^H\varphi- \frac{1}{\r^4}\, \varphi\,\nabla^H\rho^2\ .
\]
Hence we have
\[
\begin{array}{l}
\ds\vec{X}_{h_{r,\eta}}:=J_H\nabla^Hh_{r,\eta}+\frac{h_{r,\eta}}{2}\, \vec{R}\\[5mm]
\ds\ =[\chi(\r/r)-\chi(\r/\eta)]\,\lf[  \frac{\rho^2}{\r^4}\,J_H\,\nabla^H\varphi- \frac{1}{\r^4}\, \varphi\,J_H\nabla^H\rho^2+ \arctan\sigma\ \frac{\vec{R}}{2} \rg]+[r^{-1}\,\chi'(\r/r)-\eta^{-1}\,\chi'(\r/\eta)]\, \arctan\sigma\,\nabla^H\r\ .
\end{array}
\]
Recall that
\[
\nabla^H\varphi=\sum_{i=1}^4v^i\,\nabla^Hw^i-w^i\,\nabla^Hv^i\ .
\]
Hence
\[
\lf\{
\begin{array}{l}
\ds|\vec{X}_{h_{r,\eta}}|\le C\, \sum_{j=0}^1 \r^{-1+j}\ \lf|r^{-j}\,\chi^{(j)}(\r/r)-\eta^{-j}\,\chi^{(j)}(\r/\eta)\rg|\\[5mm]
\ds |\nabla^H\vec{X}_{h_{r,\eta}}|\le C\, \sum_{j=0}^2 \r^{-2+j}\ \lf|r^{-j}\,\chi^{(j)}(\r/r)-\eta^{-j}\,\chi^{(j)}(\r/\eta)\rg| \\[5mm]
\ds |\nabla^H\nabla^H \vec{X}_{h_{r,\eta}}|\le C\, \sum_{j=0}^3 \r^{-3+j}\ \lf|r^{-j}\,\chi^{(j)}(\r/r)-\eta^{-j}\,\chi^{(j)}(\r/\eta)\rg|
\end{array}
\rg.
\]
We have also
\[
\nabla\vec{X}_h\cdot\vec{R}=\nabla\lf[J_H\nabla^Hh_{r,\eta}\rg]\cdot \vec{R}+\nabla h_{r,\eta}=-J_H\nabla^Hh_{r,\eta}\cdot \nabla^H\vec{R}+\nabla h_{r,\eta}
\]
which gives
\[
|\nabla\vec{X}_h\cdot\vec{R}|\le   C\,\r^{-1}\, [\chi(\r/r)-\chi(\r/\eta)]+C\,\lf|[r^{-1}\,\chi'(\r/r)-\eta^{-1}\,\chi'(\r/\eta)\rg|\ .
\]
Recall
\be
\label{forme-aire-bis}
\int_{\Sigma}<d\vec{w}_{r,\eta},d\vec{\La}_k>_g\ dvol_g=2\, \int_\Sigma\, <dh_{r,\eta},d\beta_k>_g\ dvol_g\ .
\ee
and 
\be
\label{euler-lagrange-bis}
\begin{array}{l}
\ds \int_{\Sigma}<d\vec{w}_{r,\eta}\cdot d\vec{\La}_k>_{g_k}\ dvol_{g_k}+e^{-\ep_k^{-2}}\, O(\|w\|_{\vec{\La}_k})\\[5mm]
\ds\quad=-4\,\ep^4_k\, \int_\Sigma (1+|d\vec{T}_k|^2_{g_k})\ \lf[-\,\lf<d\vec{T}_k\,\dot{\otimes}\,d\vec{T}_k, d\vec{w}\,\dot{\otimes}\,d\vec{\La}_k\rg>_g+\,d\vec{T}_k\cdot d\lf<d\vec{w}\wedge d\vec{\La}_k\rg>_{g_k}\rg]\ dvol_{g_k}\\[5mm]
\ds\quad-\ep^4_k\ \int_\Sigma\lf [(1+|d\vec{T}_k|^2_{g_k})^2- 4\, |d\vec{T}_k|^2_{g_k}\, (1+|d\vec{T}_k|^2_{g_k})\rg] <d\vec{w}\cdot d\vec{\La}_k>_{g_k}\ dvol_{g_k}\
\end{array}
\ee
and we recall  from   (\ref{pointwize-w})    and from (\ref{ptw-horiz-ref})
\be
\label{pointwize-w-ref}
|d\vec{\La}_k\,\dot{\otimes}\,d\vec{w}_{r,\eta}|_{g_k}\le \ |\nabla^H\vec{X}_{h_{r,\eta}}|
\ee
and
\be
\label{ptw-horiz-ref-bis}
\lf|d\vec{T}_k\cdot d\lf<d\vec{w}_{r,\eta}\wedge d\vec{\La}_k\rg>_g\rg|\le \lf|d\vec{T}_k\rg|_g^2 |\nabla^H\vec{X}_{h_{r,\eta}}|+\lf|d\vec{T}_k\rg|_g\ [|\nabla^H\nabla^H\vec{X}_{h_{r,\eta}}|+|\nabla\vec{X}_{h_{r,\eta}}\cdot\vec{R}|]
\ee
Hence
\[
\begin{array}{l}
\ds\lf|\int_\Sigma\, <dh_{r,\eta},d\beta_k>_g\ dvol_g\rg|\le C\, \ep^4_k\, \int_\Sigma (1+|d\vec{T}_k|^2_{g_k})^2\ |\nabla^H\vec{X}_{h_{r,\eta}}| \ dvol_{g_k}\\[5mm]
\ds\quad+C\, \ep^4_k\, \int_\Sigma (1+|d\vec{T}_k|^2_{g_k})^{3/2}\ [|\nabla^H\nabla^H\vec{X}_{h_{r,\eta}}|+|\nabla\vec{X}_{h_{r,\eta}}\cdot\vec{R}|] \ dvol_{g_k}+e^{-\ep_k^{-2}}\, O(\|w\|_{\vec{\La}_k})
\end{array}
\]
This gives
\be
\label{penalisation}
\begin{array}{l}
\ds\lf|\int_\Sigma\, <dh_{r,\eta},d\beta_k>_g\ dvol_g\rg|\le C\, \ep^4_k\, \int_{\Sigma\,\cap\, \{\eta<\r<2r\}} \r^{-2}\, (1+|d\vec{T}_k|^2_{g_k})^2\  \ dvol_{g_k}\\[5mm]
\ds\quad+C\, \ep^4_k\, \int_{\Sigma\,\cap\, \{\eta<\r<2r\}} \r^{-3}\,  (1+|d\vec{T}_k|^2_{g_k})^{3/2}\  \ dvol_{g_k}+e^{-\ep_k^{-2}}\, O(\|w\|_{\vec{\La}_k})\\[5mm]
\end{array}
\ee
We recall that the function $$d(\vec{p},\vec{q}):=\r_{\vec{p}}(\vec{q})=\lf[ |\vec{p}-\vec{q}|^4+4\,|\vec{a}\cdot\vec{d}-\vec{b}\cdot\vec{c}|^2    \rg]^{1/4}\quad\mbox{ where }\quad \vec{p}=(\vec{a},\vec{b}) \quad\mbox{ and }\vec{q}=(\vec{c},\vec{d}) $$ and we have seen in - see (\ref{quasi-dist}) - that it satisfies the quasi-distance condition
\be
\label{quasi-dist-bis}
\r_{\vec{p}}(\vec{q})=\r_{\vec{q}}(\vec{p})\quad\mbox{ and }\quad d(\vec{p}_1,\vec{p})\le C_0\,\lf[d(\vec{p}_1,\vec{q}) +d(\vec{p},\vec{q})\rg]
\ee
for some universal $C_0>1$. We are now proving the following lemma.
\begin{Lm}
\label{lm-good-points}
There exist two universal constants $C>0$ and  $c_0>0$ such that for any $\delta_k\rightarrow 0$ and for $k$ large enough and for any point $\vec{p}\in \vec{\La}_k(\Sigma)$ satisfying
\be
\label{good-1}
\forall\ r<1\quad\quad  {\ep^4_k}\, \int_{\Sigma\,\cap\, \{\r_{\vec{p}}< r\}}  (1+|d\vec{T}_k|^2_{g_k})^2\  \ dvol_{g_k}\le \delta_k\ \int_{\Sigma\,\cap\, \{\r_{\vec{p}}<5\,C_0^2\,r\}}  \ dvol_{g_k}\ ,
\ee
we have  for any $0<r<1/2 $
\be
\label{good-2}
\begin{array}{l}
\ds \int_{\Sigma\,\cap\, \{\r_{\vec{p}}<r\}}|\nabla^\Sigma\arctan\sigma_k|^2\ dvol_{g_k}\le \frac{C}{r^2}\int_{\Sigma\,\cap\, \{\r_{\vec{p}}<r\}}  \ dvol_{g_k}\ \\[5mm]
\ds \quad+C\, \lf[\delta\ \log\lf(\frac{r}{\ep}\rg)+\delta^{1/4} +\delta^{3/4}\ \log^{3/4}\lf(\frac{r}{\ep}\rg)\rg]\ \sup_{s<r}\ \frac{1}{s^2}\ \int_{\Sigma\,\cap\, \{\r_{\vec{p}}<s\}}  \ dvol_{g_k}\\[5mm]
\ds \quad+C\ e^{-\ep_k^{-2}}\ \ep_k^{-3}\ ,
\end{array}
\ee
and for any $0<s<r<1/2 $
\be
\label{good-2-b}
\begin{array}{l}
\ds \lf[c_0-\lf[\delta\ \log\lf(\frac{r}{\ep}\rg)+\delta^{1/4} +\delta^{3/4}\ \log^{3/4}\lf(\frac{r}{\ep}\rg)\rg]^{1/2}\rg]\ \frac{1}{s^2}\int_{\Sigma\,\cap\, \{\r_{\vec{p}}<s\}} dvol_{g_k}\\[5mm]
\ds\quad\quad\le \frac{C}{r^2}\int_{\Sigma\,\cap\, \{\r_{\vec{p}}<r\}}  \ dvol_{g_k}+C\ e^{-\ep_k^{-2}}\ \ep_k^{-3}\ .
\end{array}
\ee
\end{Lm}
\noindent{\bf Proof of lemma~\ref{lm-good-points}.}  We take $\vec{p}:=(\vec{\ep}_1,\vec{\ep}_2)$. Since the second fundamental form of $\vec{\La}_k$ is in $L^4$ the weak Immersion
is $C^1$. Hence $\vec{p}$ has finitely many pre-images by $\vec{\La}_k$ and we denote their number by $N_{\vec{p}}\in {\N}$. We choose $\vec{p}$ such that $N_{\vec{p}}\ge 1$.
Let $x_0\in \Sigma$ such that $\vec{\La}_k(x_0)=\vec{p}$. We take local conformal coordinates around $x_0$ such that $x_i(x_0)=0$ and we omit to write explicitly the subscript $k$. We can assume modulo rotation in the image and dilation of the coordinates that
$\p_{x_1}\vec{\La}(0)=(\vec{\ep}_3,0)$. Since both $J_H(\p_{x_1}\vec{\La}(0))=(0,\vec{\ep}_3)$ and $J_H(\p_{x_2}\vec{\La}(0))$ are orthogonal to $\vec{\La}_\ast T_{x_0}\Sigma$ we can assume modulo the action of rotations preserving $H$ and $J_H$ that $\p_{x_2}\vec{\La}(0)=(\vec{\ep}_4,0)$ and $J_H(\p_{x_2}\vec{\La}(0))=(0,\vec{\ep}_4)$. Hence, since $\vec{\La}$ is $C^1$ we have in the neighbourhood of $0$
\[
\vec{\La}(x_1,x_2)= \lf(\vec{\ep}_1+x_1\,\vec{\ep}_3+x_2\,\vec{\ep}_4+o(|x|),\vec{\ep}_2+o(|x|)\rg)
\]
Recall that for $\vec{a}=\vec{\ep}_1+\vec{v}$ and $\vec{b}:=\vec{\ep}_2+\vec{w}$ with $(\vec{a},\vec{b})\in V_2({\R}^4)$ recall that $\varphi=v^2-w^1$
\[
\nabla^H\varphi(\vec{a},\vec{b})=\sum_{i=1}^4v^i\,\nabla^H w^i-w^i\,\nabla^H v^i\ .
\]
Recall for  $\vec{p}=(\vec{a},\vec{b})=(\vec{\ep}_1+\vec{v},\vec{\ep}_2+\vec{w})$
\[
\begin{array}{l}
\ds H_{\vec{p}}=T_{\vec{p}}V_2({\R}^4)\cap (\vec{R})^\perp\\[5mm]
\ds\quad\quad=\lf\{(\vec{V},\vec{W})\in {\R}^8  \ ;\ V_1+\vec{V}\cdot\vec{v}=0\ ;\  V_2+\vec{V}\cdot\vec{w}=0\ ; \ W_1+\vec{W}\cdot\vec{v}=0\ ;\  W_2+\vec{W}\cdot\vec{w}=0\ \rg\}
\end{array}
\]
we have respectively 
\[
\nabla^H v^i=\pi_H(\vec{\ep}_i,0)\quad\mbox{ and }\quad\nabla^H w^i=\pi_H(0,\vec{\ep}_i)\ .
\]
First of all we have at the point $\vec{p}:=(\vec{a},\vec{b})$
\[
\nabla^{V_2({\R}^4)}v^i= (\vec{\ep}_i-a^i\,\vec{a}-b^i\,\vec{b},0)\quad\mbox{ and }\quad\nabla^{V_2({\R}^4)} w^i=(0,\vec{\ep}_i-a^i\,\vec{a}-b^i\,\vec{b})
\]
Hence
\be
\label{nabla-H}
\nabla^H v^i=(\vec{\ep}_i-a^i\,\vec{a}-b^i\,\vec{b},0)\quad\mbox{ and }\quad\nabla^{H} w^i=(0,\vec{\ep}_i-a^i\,\vec{a}-b^i\,\vec{b})
\ee
In particular at $\vec{p}=(\vec{a},\vec{b})=(\vec{\ep}_1+\vec{v},\vec{\ep}_2+\vec{w})$, thanks to (\ref{tan}) there holds
\be
\label{nabla-H-bis}
\lf\{
\begin{array}{l}
\nabla^Hv^1=(-\vec{v}-w^1\ (\vec{\ep}_2+ \vec{w})-v^1\,(\vec{\ep}_1+\vec{v}), 0)\\[3mm]
\nabla^H v^2=(-\vec{w}-v^2\,(\vec{\ep}_1+\vec{v})-w^2\,(\vec{\ep}_2+\vec{w}),0)\\[3mm]
\nabla^H v^3=(\vec{\ep}_3-v^3\, (\vec{\ep}_1+\vec{v})-w^3\, (\vec{\ep}_2+\vec{w}),0)\\[3mm]
\nabla^H v^4=(\vec{\ep}_4-v^4\, (\vec{\ep}_1+\vec{v})-w^4\, (\vec{\ep}_2+\vec{w}),0)\ ,
\end{array}
\rg.
\ee
 and 
\be
\label{nabla-H-ter}
\lf\{
\begin{array}{l}
\nabla^Hw^1=(0, -\vec{v}-w^1\ (\vec{\ep}_2+ \vec{w})-v^1\,(\vec{\ep}_1+\vec{v}))\\[3mm]
\nabla^H w^2=(0,-\vec{w}-v^2\,(\vec{\ep}_1+\vec{v})-w^2\,(\vec{\ep}_2+\vec{w}))\\[3mm]
\nabla^H w^3=(0, \vec{\ep}_3-v^3\, (\vec{\ep}_1+\vec{v})-w^3\, (\vec{\ep}_2+\vec{w}))\\[3mm]
\nabla^H w^4=(0,\vec{\ep}_4-v^4\, (\vec{\ep}_1+\vec{v})-w^4\, (\vec{\ep}_2+\vec{w}))\ .
\end{array}
\rg.
\ee
Observe that for $i=1...4$
\be
\label{J-H}
J_H\nabla^H v^i=\nabla^H w^i\ .
\ee
Thus in particular
\be
\label{horizont-grad}
\lf\{
\begin{array}{l}
\nabla^H v^1\circ\vec{\La}(x)=O(|x|)\ \\[3mm]\nabla^H v^2\circ\vec{\La}(x)=O(|x|)\ \\[3mm]\nabla^H v^3\circ\vec{\La}(x)=(\vec{\ep}_3,0)+O(|x|)\  \\[3mm]\nabla^H v^4\circ\vec{\La}(x)=(\vec{\ep}_4,0)+O(|x|)
\end{array}
\rg. 
\quad\mbox{ and }\quad\lf\{
\begin{array}{l}
\nabla^H w^1\circ\vec{\La}(x)=O(|x|)\ \\[3mm]\nabla^H w^2\circ\vec{\La}(x)=O(|x|)\ \\[3mm]\nabla^H w^3\circ\vec{\La}(x)=(0,\vec{\ep}_3)+O(|x|)\  \\[3mm]\nabla^H w^4\circ\vec{\La}(x)=(0,\vec{\ep}_4)+O(|x|)
\end{array}
\rg. 
\ee
Thus
\[
\lf\{
\begin{array}{l}
\nabla^\Sigma v^1\circ\vec{\La}(x)=O(|x|)\ \\[3mm]\nabla^\Sigma v^2\circ\vec{\La}(x)=O(|x|)\ \\[3mm]\nabla^\Sigma v^3\circ\vec{\La}(x)=(\vec{\ep}_3,0)+O(|x|)\  \\[3mm]\nabla^H v^4\circ\vec{\La}(x)=(\vec{\ep}_4,0)+O(|x|)
\end{array}
\rg. 
\quad\mbox{ and }\quad\lf\{
\begin{array}{l}
\nabla\Sigma w^1\circ\vec{\La}(x)=O(|x|)\ \\[3mm]\nabla^\Sigma w^2\circ\vec{\La}(x)=O(|x|)\ \\[3mm]\nabla^\Sigma w^3\circ\vec{\La}(x)=O(|x|)\  \\[3mm]\nabla^\Sigma w^4\circ\vec{\La}(x)=O(|x|)
\end{array}
\rg. 
\]
Hence we deduce in particular that
\[
w^3\circ\vec{\La}=O(|x|^2)\quad\mbox{ and }\quad w^4\circ\vec{\La}=O(|x|^2)
\]
Since
\[
\nabla^H\varphi\circ\vec{\La}=\sum_{i=3}^4v^i\circ\vec{\La}\,\nabla^H w^i\circ\vec{\La}-w^i\circ\vec{\La}\,\nabla^H v^i\circ\vec{\La}+O(|x|^2)\ .
\]
we deduce from the previous that
\[
\nabla^\Sigma\varphi= O(|x|^2) \quad\mbox{ and }\quad\varphi(x)=O(|x|^3)\ .
\]
Since $\rho^2=|x|^2\, (1+o(1))$ we deduce
\be
\label{arctan-bor}
|\nabla^\Sigma\arctan\sigma|=O(1)\quad\mbox{ and }\quad |\arctan\sigma|=O(|x|)\ .
\ee
Since  $\vec{\La}$ is a $C^1$ immersion we have 
\[
\lim_{s\rightarrow 0} \lf\| \frac{\rho}{\r}-1\rg\|_{L^\infty(\vec{\La}^{-1}B_\eta(\vec{p}))}=0
\]
and then, using again that $\vec{\La}$ is a $C^1$ immersion we have
\be
\label{density}
\lim_{\eta\rightarrow 0}\frac{1}{\eta^2}\,\int_{\r_{\vec{p}}<\eta}\ dvol_{g_k}=\lim_{\eta\rightarrow 0}\frac{1}{\eta^2}\,\int_{\rho_{\vec{p}}<\eta}\ dvol_{g_k}=\pi \, N_{\vec{p}}
\ee
We recall the identity (\ref{truncated})
\be
\label{truncated-bis}
\begin{array}{l}
\ds\int_\Sigma<d h_{r,\eta},d\beta>_g\ dvol_g-\int_\Sigma\frac{\r}{r}\,\chi'\lf(\frac{\r}{r}\rg)\,\frac{1}{\r^2}\,\lf[{|\nabla^\Sigma\r|^2} +\frac{\sigma\,\arctan\sigma}{\sqrt{1+\sigma^2}} \rg]\ dvol_g\\[5mm]
\ds+\int_\Sigma\ \lf[O(1)\chi\lf(\frac{\r}{r}\rg)+O\lf(\frac{\r}{r}\rg)\,\chi'\lf(\frac{\r}{r}\rg)\rg]\ dvol_g +\frac{1}{4}\,\int_\Sigma \frac{\r^2}{r^{2}}\,\chi''\lf(\frac{\r}{r}\rg) \,\arctan\sigma\,\frac{\nabla^\Sigma\r}{\r}\cdot\frac{\nabla^\Sigma\sigma}{1+\sigma^2}\ dvol_g \\[5mm]
\ds-\frac{3}{4}\,\int_\Sigma  \ \frac{\r}{r}\,\chi'\lf(\frac{\r}{r}\rg) \,\arctan\sigma\,\frac{\nabla^\Sigma\r}{\r}\cdot\frac{\nabla^\Sigma\sigma}{1+\sigma^2}\ dvol_g +\frac{1}{4}\,\int_\Sigma \frac{\r}{r}\chi'\lf(\frac{\r}{r}\rg)\,\frac{|\nabla^\Sigma\sigma|^2}{(1+\sigma^2)^2}\ dvol_g \\[5mm]
\ds=4\,\int_\Sigma[\chi(\r/r)-\chi(\r/\eta)]\,\frac{|(\nabla^\Sigma\r)^\perp|^2}{\r^2}\ dvol_g-\int_\Sigma\frac{\r}{\eta}\,\chi'\lf(\frac{\r}{\eta}\rg)\,\frac{1}{\r^2}\,\lf[{|\nabla^\Sigma\r|^2} +\frac{\sigma\,\arctan\sigma}{\sqrt{1+\sigma^2}} \rg]\ dvol_g\\[5mm]
\ds+4\,\int_\Sigma\frac{\r}{\eta}\, \chi'\lf(\frac{\r}{\eta}\rg)\,\frac{|\nabla^\Sigma\sigma|^2}{(1+\sigma^2)^2}\ dvol_g-\frac{3}{4}\int_\Sigma\,\frac{\r}{\eta}\,\chi'\lf(\frac{\r}{\eta}\rg)\, \,\arctan\sigma\,\frac{\nabla^\Sigma\r}{\r}\cdot\frac{\nabla^\Sigma\sigma}{1+\sigma^2}\ dvol_g  \\[5mm]
\ds+\frac{1}{4}\int_\Sigma\frac{\r^2}{\eta^{2}}\,\chi''\lf(\frac{\r}{\eta}\rg)\, \,\arctan\sigma\,\frac{\nabla^\Sigma\r}{\r}\cdot\frac{\nabla^\Sigma\sigma}{1+\sigma^2}\ dvol_g +\int_\Sigma\ \lf[O(1)\chi\lf(\frac{\r}{\eta}\rg)+O\lf(\frac{\r}{\eta}\rg)\,\chi'\lf(\frac{\r}{\eta}\rg)\rg]\ dvol_g
\end{array}
\ee
Using (\ref{arctan-bor}) we have
\[
\begin{array}{l}
\ds\int_\Sigma\frac{\r}{\eta}\,\chi'\lf(\frac{\r}{\eta}\rg)\,\frac{1}{\r^2}\,\lf[{|\nabla^\Sigma\r|^2} +\frac{\sigma\,\arctan\sigma}{\sqrt{1+\sigma^2}} \rg]\ dvol_g\\[5mm]
\ds=\frac{1}{\eta^2}\int_{\r<2\eta}\chi'\lf(\frac{\r}{\eta}\rg)\,\frac{\eta}{\r}\,|\nabla^\Sigma\r|^2 dvol_g+ O\lf(\frac{1}{\eta}\int_{\r<2\eta}dvol_g\rg)
\end{array}
\]
From (\ref{horizont-gradient}), (\ref{gg}) and  (\ref{arctan-bor})  we also have
\[
1+O(\r^2)=|\nabla^H\r|^2=|\nabla^\Sigma\r|^2+|(\nabla^\Sigma\r)^\perp|^2=|\nabla^\Sigma\r|^2+4^{-1}\,\r^2\,|\nabla^\Sigma(\arctan\sigma)|^2=|\nabla^\Sigma\r|^2+ O(\r^2)
\]
Combining the previous gives then
\be
\label{density-mooth}
\begin{array}{l}
\ds-\int_\Sigma\frac{\r}{\eta}\,\chi'\lf(\frac{\r}{\eta}\rg)\,\frac{1}{\r^2}\,\lf[{|\nabla^\Sigma\r|^2} +\frac{\sigma\,\arctan\sigma}{\sqrt{1+\sigma^2}} \rg]\ dvol_g\\[5mm]
\ds=-\frac{1}{\eta}\int_{\r<2\eta}\chi'\lf(\frac{\r}{\eta}\rg)\,\frac{\eta}{\r}\,\lf|\nabla^\Sigma\lf(\frac{\r}{\eta}\rg)\rg|\ dvol_g +o(1)=-\int_{0}^\infty\chi'(s)\ \frac{\ds{\mathcal H}^1(\{\r_{\vec{p}}=\eta\, s\})}{\eta\, s}\  ds+o(1)
\end{array}
\ee
Hence
\be
\label{sensity-smooth}
\lim_{\eta\rightarrow 0} -\int_\Sigma\frac{\r}{\eta}\,\chi'\lf(\frac{\r}{\eta}\rg)\,\frac{1}{\r^2}\,\lf[{|\nabla^\Sigma\r|^2} +\frac{\sigma\,\arctan\sigma}{\sqrt{1+\sigma^2}} \rg]\ dvol_g=2\,\pi \, N_{\vec{p}}
\ee
Thus
\be
\label{mon-2}
\begin{array}{l}
\ds\lim_{\eta\rightarrow 0}\,4\,\int_\Sigma[\chi(\r/r)-\chi(\r/\eta)]\,\frac{|(\nabla^\Sigma\r)^\perp|^2}{\r^2}\ dvol_g-\int_\Sigma\frac{\r}{\eta}\,\chi'\lf(\frac{\r}{\eta}\rg)\,\frac{1}{\r^2}\,\lf[{|\nabla^\Sigma\r|^2} +\frac{\sigma\,\arctan\sigma}{\sqrt{1+\sigma^2}} \rg]\ dvol_g\\[5mm]
\ds+4\,\int_\Sigma\frac{\r}{\eta}\, \chi'\lf(\frac{\r}{\eta}\rg)\,\frac{|\nabla^\Sigma\sigma|^2}{(1+\sigma^2)^2}\ dvol_g-\frac{3}{4}\int_\Sigma\,\frac{\r}{\eta}\,\chi'\lf(\frac{\r}{\eta}\rg)\, \,\arctan\sigma\,\frac{\nabla^\Sigma\r}{\r}\cdot\frac{\nabla^\Sigma\sigma}{1+\sigma^2}\ dvol_g  \\[5mm]
\ds+\frac{1}{4}\int_\Sigma\frac{\r^2}{\eta^{2}}\,\chi''\lf(\frac{\r}{\eta}\rg)\, \,\arctan\sigma\,\frac{\nabla^\Sigma\r}{\r}\cdot\frac{\nabla^\Sigma\sigma}{1+\sigma^2}\ dvol_g +\int_\Sigma\ \lf[O(1)\chi\lf(\frac{\r}{\eta}\rg)+O\lf(\frac{\r}{\eta}\rg)\,\chi'\lf(\frac{\r}{\eta}\rg)\rg]\ dvol_g\\[5mm]
\ds=4\,\int_\Sigma \chi\lf(\frac{\r}{r}\rg)\,\frac{|(\nabla^\Sigma\r)^\perp|^2}{\r^2}\ dvol_g+2\,\pi \, N_{\vec{p}}
\end{array}
\ee
We have also
\[
\begin{array}{l}
\ds\int_\Sigma<d h_{r,\eta},d\beta>_g\ dvol_g=\frac{1}{2}\int_\Sigma<d\vec{w}_{r,\eta},d\vec{\La}>_g\ dvol_g=\frac{1}{2}\int_\Sigma\vec{X}_{r,\eta}\circ\La_k\cdot\Delta_g\vec{\La}\ dvol_g\\[5mm]
\ds=-\int_\Sigma\nabla^H\lf[( \chi(\r/r)-\chi(\r/\eta) )  \arctan\sigma\rg]\cdot J_H\pi_{\vec{N}}\Delta_g\vec{\La} dvol_g\\[5mm]
\ds=-\int_\Sigma\nabla^\Sigma\lf[( \chi(\r/r)-\chi(\r/\eta) )  \arctan\sigma\rg]\cdot J_H\pi_{\vec{N}}\vec{H}_{\vec{\La}}\ dvol_g
\end{array}
\]
where we have used  the fact that $\vec{R}\circ\La\cdot\Delta\vec{\La}=0$ - see (\ref{leg-cond}) - and the fact that $J_H$ realizes an isomorphism between the horizontal projection of the normal bundle of the immersion and the tangent bundle to the immersion. We have then
\be
\label{mon-2a}
\begin{array}{l}
\ds\lf|\int_\Sigma<d h_{r,\eta},d\beta>_g\ dvol_g\rg|\le C\ \int_{\eta<\r<2\,r}|\nabla^\Sigma\arctan\sigma|\ \sqrt{1+|d\vec{T}|^2_{g}}\ dvol_g\\[5mm]
\ds+ C\ \int_{\eta<\r<2\,\eta}|x|^{-1}\,|\arctan\sigma|\ \sqrt{1+|d\vec{T}_k|^2_{g_k}}\ dvol_g+ C\ \int_{r<\r<2\,r}|x|^{-1}\,|\arctan\sigma|\ \sqrt{1+|d\vec{T}|^2_{g}}\ dvol_g\\[5mm]
\ds\le C\ \lf[\int_{\eta<\r<2\,r}|\nabla^\Sigma\arctan\sigma|^2\ dvol_g  \rg]^{1/2}\ \lf[ \int_{\eta<\r<2\,r}(1+|d\vec{T}|^2_{g})^2\ dvol_g  \rg]^{1/4}\ \lf[\int_{\eta<\r<2\,r}\ dvol_g \rg]^{1/4}\\[5mm]
\ds+ C\ \lf[\int_{\eta<\r<2\,\eta}|x|^{-2}\,|\arctan\sigma|^2\ dvol_g  \rg]^{1/2}\ \lf[ \int_{\eta<\r<2\,\eta}(1+|d\vec{T}|^2_{g})^2\ dvol_g  \rg]^{1/4}\ \lf[\int_{\eta<\r<2\,\eta}\ dvol_g \rg]^{1/4}\\[5mm]
\ds+ C\ \int_{r<\r<2\,r}|x|^{-1}\,|\arctan\sigma|\ \sqrt{1+|d\vec{T}|^2_{g}}\ dvol_g
\end{array}
\ee
Since $|\nabla^\Sigma\arctan\sigma|=O(1)$ and $|\arctan\sigma|=O(|x|)$ we deduce
\be
\label{mon-3}
\begin{array}{l}
\ds\limsup_{\eta\rightarrow 0}\lf|\int_\Sigma<d h_{r,\eta},d\beta>_g\ dvol_g\rg|\le  C\ \int_{r<\r<2\,r}|x|^{-1}\,|\arctan\sigma|\ \sqrt{1+|d\vec{T}|^2_{g}}\ dvol_g\\[5mm]
\ds +C\ \lf[\int_{\r<2\,r}|\nabla^\Sigma\arctan\sigma|^2\ dvol_g  \rg]^{1/2}\ \lf[ \int_{\r<2\,r}(1+|d\vec{T}|^2_{g})^2\ dvol_g  \rg]^{1/4}\ \lf[\int_{\r<2\,r}\ dvol_g \rg]^{1/4}
\end{array}
\ee
Combining (\ref{truncated-bis}), (\ref{mon-2}) and (\ref{mon-3}) gives
\be
\label{truncated-bis-re}
\begin{array}{l}
\ds 4\,\int_\Sigma \chi\lf(\frac{\r}{r}\rg)\,\frac{|(\nabla^\Sigma\r)^\perp|^2}{\r^2}\ dvol_g+2\,\pi \, N_{\vec{p}}\le 
-\int_\Sigma\frac{\r}{r}\,\chi'\lf(\frac{\r}{r}\rg)\,\frac{1}{\r^2}\,\lf[{|\nabla^\Sigma\r|^2} +\frac{\sigma\,\arctan\sigma}{\sqrt{1+\sigma^2}} \rg]\ dvol_g\\[5mm]
\ds+\int_\Sigma\ \lf[O(1)\chi\lf(\frac{\r}{r}\rg)+O\lf(\frac{\r}{r}\rg)\,\chi'\lf(\frac{\r}{r}\rg)\rg]\ dvol_g +\frac{1}{4}\,\int_\Sigma \frac{\r^2}{r^{2}}\,\chi''\lf(\frac{\r}{r}\rg) \,\arctan\sigma\,\frac{\nabla^\Sigma\r}{\r}\cdot\frac{\nabla^\Sigma\sigma}{1+\sigma^2}\ dvol_g \\[5mm]
\ds-\frac{3}{4}\,\int_\Sigma  \ \frac{\r}{r}\,\chi'\lf(\frac{\r}{r}\rg) \,\arctan\sigma\,\frac{\nabla^\Sigma\r}{\r}\cdot\frac{\nabla^\Sigma\sigma}{1+\sigma^2}\ dvol_g +\frac{1}{4}\,\int_\Sigma \frac{\r}{r}\chi'\lf(\frac{\r}{r}\rg)\,\frac{|\nabla^\Sigma\sigma|^2}{(1+\sigma^2)^2}\ dvol_g \\[5mm]
\ds +\, C\ \int_{r<\r<2\,r}|x|^{-1}\,|\arctan\sigma|\ \sqrt{1+|d\vec{T}|^2_{g}}\ dvol_g\\[5mm]
\ds +C\ \lf[\int_{\r<2\,r}|\nabla^\Sigma\arctan\sigma|^2\ dvol_g  \rg]^{1/2}\ \lf[ \int_{\r<2\,r}(1+|d\vec{T}|^2_{g})^2\ dvol_g  \rg]^{1/4}\ \lf[\int_{\r<2\,r}\ dvol_g \rg]^{1/4}

\end{array}
\ee
Using (\ref{gg}) we deduce for any $0<\al<1$
\[
\begin{array}{l}
\ds 4\,\int_{\r_{\vec{p}}<r} |\nabla^\Sigma\arctan\sigma|^2+O(1)\  dvol_g+2\,\pi \, N_{\vec{p}}\le C\ \r^{-2}\ \int_{r<\r_{\vec{p}}<2\,r} \  dvol_g\\[5mm]
\ds +C\ \r^{-1}\ \int_{r<\r<2\,r} \sqrt{1+|d\vec{T}|^2_{g}}\ dvol_g+C\, \al\, \int_{\r<2\,r}|\nabla^\Sigma\arctan\sigma|^2\ dvol_g \\[5mm]
\ds + C\ \al^{-1}\ \lf[ \int_{\r<2\,r}(1+|d\vec{T}|^2_{g})^2\ dvol_g  \rg]^{1/2}\ \lf[\int_{\r<2\,r}\ dvol_g \rg]^{1/2}\ .
\end{array}
\]
This gives
\[
\begin{array}{l}
\ds (4-C\al)\,\int_{\r_{\vec{p}}<r} |\nabla^\Sigma\arctan\sigma|^2 dvol_g+2\,\pi \, N_{\vec{p}}\le C\ \r^{-2}\ \int_{\r_{\vec{p}}<2\,r} \  dvol_g\\[5mm]
\ds +C\ \r^{-1}\ \int_{r<\r_{\vec{p}}<2\,r} \sqrt{1+|d\vec{T}|^2_{g}}\ dvol_g + C\ \al^{-2}\  \r^2\ \int_{\r_{\vec{p}}<2\,r}(1+|d\vec{T}|^2_{g})^2\ dvol_g   \\[5mm]
\ds\le C\ \r^{-2}\ \int_{\r_{\vec{p}}<2\,r} \  dvol_g+C\ \r^{-3/2}\  \lf[ \r^2 \int_{\r_{\vec{p}}<2\,r}(1+|d\vec{T}|^2_{g})^2\ dvol_g  \rg]^{1/4}\ \lf[\int_{\r_{\vec{p}}<2\,r}\ dvol_g \rg]^{3/4}\\[5mm]
\ds+ C\ \al^{-2}\  \r^2\ \int_{\r_{\vec{p}}<2\,r}(1+|d\vec{T}|^2_{g})^2\ dvol_g \ .
\end{array}
\]
We choose $C\al\le 2$ and we finally get using the assumption of the lemma
\[
 2\,\int_{\r_{\vec{p}}<r} |\nabla^\Sigma\arctan\sigma|^2 dvol_g+2\,\pi \, N_{\vec{p}}\le C\, \lf[1+\delta \frac{r^4}{\ep^4}\rg]\ \frac{1}{\r^2}\ \int_{\r_{\vec{p}}<10\,C_0^2\,r} \  dvol_g\ .
\]
Observe that $|\nabla^\Sigma\arctan\sigma|^2\le C\, \r^{-2}$ and this implies
\[
\int_{r<\r_{\vec{p}}<10\,C^2_0\,r} |\nabla^\Sigma\arctan\sigma|^2\ dvol_g\le \frac{C}{\r^2}\ \int_{r<\r_{\vec{p}}<10\,C_0^2\,r} \  dvol_g\ .
\]
Hence we get 
\be
\label{prem-conc}
\begin{array}{l}
\ds \forall\ r<C_0^{-2}\, \ep \,\delta^{-1/4}\quad\quad 2\,\int_{\r_{\vec{p}}<r} |\nabla^\Sigma\arctan\sigma|^2 dvol_g+2\,\pi \, N_{\vec{p}}\le C\ \frac{1}{\r^2}\ \int_{\r_{\vec{p}}<r} \  dvol_g\ .
\end{array}
\ee
Let $1>r>\ep \,\delta^{-1/4}$ and let $\ep\le s\le\ep \,\delta^{-1/4}$, from (\ref{mon-3})
\be
\label{mon-3-bis}
\begin{array}{l}
\ds\limsup_{\eta\rightarrow 0}\lf|\int_\Sigma<d h_{s,\eta},d\beta>_g\ dvol_g\rg|\le  C\ \int_{s<\r<2\,s}|x|^{-1}\,|\arctan\sigma|\ \sqrt{1+|d\vec{T}|^2_{g}}\ dvol_g\\[5mm]
\ds +C\ \lf[\int_{\r<2\,s}|\nabla^\Sigma\arctan\sigma|^2\ dvol_g  \rg]^{1/2}\ \lf[ \int_{\r<2\,s}(1+|d\vec{T}|^2_{g})^2\ dvol_g  \rg]^{1/4}\ \lf[\int_{\r<2\,s}\ dvol_g \rg]^{1/4}
\end{array}
\ee
We have obviously $h_{r,\eta}=h_{\ep,\eta}+h_{r,\ep}$. Combining (\ref{penalisation}) with (\ref{mon-3-bis}) gives then
\be
\label{penalisation-ter}
\begin{array}{l}
\ds\limsup_{\eta\rightarrow 0}\lf|\int_\Sigma\, <dh_{r,\eta},d\beta_k>_g\ dvol_g\rg|\le \limsup_{\eta\rightarrow 0}\lf|\int_\Sigma\, <dh_{\ep,\eta},d\beta_k>_g\ dvol_g\rg|\\[5mm]
\ds+\lf|\int_\Sigma\, <dh_{\ep,\eta},d\beta_k>_g\ dvol_g\rg|\le C\, \ep^4\, \int_{\Sigma\,\cap\, \{\ep<\r<2r\}} \r^{-2}\, (1+|d\vec{T}|^2_{g})^2\  \ dvol_{g}\\[5mm]
\ds+C\, \ep^4\, \int_{\Sigma\,\cap\, \{\ep<\r<2r\}} \r^{-3}\,  (1+|d\vec{T}_k|^2_{g})^{3/2}\  \ dvol_{g}+C\ \int_{\ep<\r<2\,\ep}|x|^{-1}\,|\arctan\sigma|\ \sqrt{1+|d\vec{T}|^2_{g}}\ dvol_g\\[5mm]
\ds +C\ \lf[\int_{\r<2\,\ep}|\nabla^\Sigma\arctan\sigma|^2\ dvol_g  \rg]^{1/2}\ \lf[ \int_{\r<2\,\ep}(1+|d\vec{T}|^2_{g})^2\ dvol_g  \rg]^{1/4}\ \lf[\int_{\r<2\,\ep}\ dvol_g \rg]^{1/4}\\[5mm]
\ds+e^{-\ep^{-2}}\ O(\|\vec{w}_{\ep,r}\|_{\vec{\La}_k})
\end{array}
\ee
where $\vec{w}_{\ep,r}:=\nabla\vec{X}_{h_{r,\ep}}$. We first estimate using (\ref{sss})
\be
\label{sss-bis}
\begin{array}{l}
\|\vec{w}_{\ep,r}\|_{\vec{\La}_k}\\[5mm]
\ds\le \lf[\int_\Sigma\lf[|\nabla^2 \vec{X}_{h_{r,\ep}}|^2\circ{\vec{\La}_k}+|\nabla \vec{X}_{h_{r,\ep}}|^2\circ{\vec{\La}_k}\ |1+|d\vec{T}_k||^2+|\vec{X}_{h_{r,\ep}}|^2\rg]^2\ dvol_{g}\rg]^{1/4}+\|\nabla \vec{X}_{h_{r,\ep}}\|_\infty\ .
\end{array}
\ee
This gives
\be
\label{sss-ter}
\|\vec{w}_{\ep,r}\|_{\vec{\La}_k}\le O(\ep^{-3})\ .
\ee
We have using the hypothesis of the lemma
\be
\label{err-1}
\begin{array}{l}
\ds  \ep^4\, \int_{\Sigma\,\cap\, \{\ep<\r_{\vec{p}}<2r\}} \r_{\vec{p}}^{-2}\, (1+|d\vec{T}|^2_{g})^2\  \ dvol_{g}=\ep^4\,\int_{\ep/2}^r \frac{ds}{s^2}\int_{\r=2\,s}\frac{(1+|d\vec{T}|^2_{g})^2}{|\nabla^\Sigma\r|} dl\\[5mm]
\ds=\ep^4 \frac{1}{4\,r^2}  \int_{\Sigma\,\cap\, \{\ep<\r_{\vec{p}}<2r\}}  (1+|d\vec{T}|^2_{g})^2\  \ dvol_{g}+2\,\ep^4\,\int_{\ep/2}^r \frac{ds}{s^3} \int_{\Sigma\,\cap\, \{\ep<\r_{\vec{p}}<2s\}} \, (1+|d\vec{T}|^2_{g})^2\  \ dvol_{g}\\[5mm]
\ds\le C\ \delta\ \log\lf(\frac{r}{\ep}\rg)\ \sup_{\ep<s<2\, C_0^2\,r}\ \frac{1}{s^2}\ \int_{\Sigma\,\cap\, \{\r_{\vec{p}}<s\}}  \ dvol_{g}
\end{array}
\ee
we have moreover
\be
\label{err-2}
\begin{array}{l}
 \ds\ep^4\, \int_{\Sigma\,\cap\, \{\ep<\r<2r\}} \r_{\vec{p}}^{-3}\,  (1+|d\vec{T}_k|^2_{g})^{3/2}\  \ dvol_{g}\\[5mm]
 \ds\le\ep^4\, \lf[\int_{\Sigma\,\cap\, \{\ep<\r<2r\}} \r_{\vec{p}}^{-6}\,    \ dvol_{g}\rg]^{1/4} \lf[\int_{\Sigma\,\cap\, \{\ep<\r<2r\}} \r_{\vec{p}}^{-2}\,  (1+|d\vec{T}_k|^2_{g})^{2}\  \ dvol_{g}\rg]^{3/4}
 \end{array}
\ee
We bound
\be
\label{err-3}
\begin{array}{l}
\ds\int_{\Sigma\,\cap\, \{\ep<\r<2r\}} \r_{\vec{p}}^{-6}\,    \ dvol_{g}=\int_{\ep/2}^r \frac{ds}{s^6}\int_{\r=2\,s}\frac{1}{|\nabla^\Sigma\r|} dl=(2r)^{-6}\ \int_{\Sigma\,\cap\, \{\ep<\r<2r\}} \,    \ dvol_{g}\\[5mm]
\ds +\,6\ \int_{\ep/2}^r \frac{ds}{s^7} \int_{\Sigma\,\cap\, \{\ep<\r_{\vec{p}}<2s\}} \,   \ dvol_{g}\le\ C\  \sup_{\ep<s<2r}\ \frac{1}{s^2}\ \int_{\Sigma\,\cap\, \{\r_{\vec{p}}<s\}}  \ dvol_{g}\ \int_{\ep/2}^r \frac{ds}{s^5}\ .
\end{array}
\ee
Combining (\ref{err-1}), (\ref{err-2}) and (\ref{err-3}) gives
\be
\label{err-4}
\begin{array}{l}
 \ds\ep^4\, \int_{\Sigma\,\cap\, \{\ep<\r<2r\}} \r_{\vec{p}}^{-3}\,  (1+|d\vec{T}_k|^2_{g})^{3/2}\  \ dvol_{g}\\[5mm]
\ds \le\,C\ \delta^{3/4}\ \log^{3/4}\lf(\frac{r}{\ep}\rg)\ \sup_{\ep<s<2\,C_0^2\,r}\ \frac{1}{s^2}\ \int_{\Sigma\,\cap\, \{\r_{\vec{p}}<s\}}  \ dvol_{g}\ .
\end{array}
\ee
We have also
\be
\label{err-5}
\begin{array}{l}
 \ds\int_{\ep<\r<2\,\ep}|x|^{-1}\,|\arctan\sigma|\ \sqrt{1+|d\vec{T}|^2_{g}}\ dvol_g\\[5mm]
 \ds\le C\ \ep^{-1}   \lf[\int_{\ep<\r<2\,\ep}\ dvol_g \rg]^{3/4}  \lf[ \int_{\ep<\r<2\,\ep}(1+|d\vec{T}|^2_{g})^2\ dvol_g  \rg]^{1/4}\le C\ \delta^{1/4}\ \ep^{-2}\ \int_{\ep<\r<2\,\ep}\ dvol_g\ . 
\end{array}
\ee
Using moreover (\ref{prem-conc}) for $r=2\ep$ we also get
\be
\label{err-6}
\begin{array}{l}
\ds\lf[\int_{\r<2\,\ep}|\nabla^\Sigma\arctan\sigma|^2\ dvol_g  \rg]^{1/2}\ \lf[ \int_{\r<2\,\ep}(1+|d\vec{T}|^2_{g})^2\ dvol_g  \rg]^{1/4}\ \lf[\int_{\r<2\,\ep}\ dvol_g \rg]^{1/4}\\[5mm]
\ds\quad\le C\ \delta^{1/4}\ \ep^{-2}\ \int_{\r<10\,C_0^2\,\ep}\ dvol_g
\end{array}
\ee
Combining (\ref{penalisation-ter})....(\ref{err-6}) we obtain
\be
\label{err-7}
\begin{array}{l}
\ds\limsup_{\eta\rightarrow 0}\lf|\int_\Sigma<d h_{s,\eta},d\beta>_g\ dvol_g\rg|\\[5mm]
\ds\le  C\ \lf[\delta\ \log\lf(\frac{r}{\ep}\rg)+\delta^{1/4} +\delta^{3/4}\ \log^{3/4}\lf(\frac{r}{\ep}\rg)\rg]\ \sup_{\ep<s<10\,C_0^2\,r}\ \frac{1}{s^2}\ \int_{\Sigma\,\cap\, \{\r_{\vec{p}}<s\}}  \ dvol_{g}+C\ e^{-\ep^{-2}}\ \ep^{-3}
\end{array}
\ee
Combining (\ref{truncated-bis}), (\ref{mon-2}) and (\ref{err-7}) we obtain
\be
\label{err-8}
\begin{array}{l}
\ds\int_{\r_{\vec{p}}<2\,r} |\nabla^\Sigma\arctan\sigma|^2 dvol_g\le C\  \frac{1}{{r}^2}\ \int_{\Sigma\,\cap\, \{\r_{\vec{p}}<2\,r\}}  \ dvol_{g}\\[5mm]
\ds\quad+C\, \lf[\delta\ \log\lf(\frac{r}{\ep}\rg)+\delta^{1/4} +\delta^{3/4}\ \log^{3/4}\lf(\frac{r}{\ep}\rg)\rg]\ \sup_{s<10\,C_0^2\, r}\ \frac{1}{s^2}\ \int_{\Sigma\,\cap\, \{\r_{\vec{p}}<s\}}  \ dvol_{g}\\[5mm]
\ds\quad+C\ e^{-\ep^{-2}}\ \ep^{-3} ,
\end{array}
\ee
where we have also used the obvious inequality
\[
\int_{r<\r_{\vec{p}}<2\,r} |\nabla^\Sigma\arctan\sigma|^2 dvol_g\le  \frac{1}{{r}^2}\ \int_{\Sigma\,\cap\, \{r<\r_{\vec{p}}<2\,r\}}  \ dvol_{g}\ .
\]
This implies (\ref{good-2}).

\medskip

Let $0<s_0<r$ such that
\be
\label{sup-1}
\frac{1}{s_0^2}\ \int_{\Sigma\,\cap\, \{\r_{\vec{p}}<s_0\}}  \ dvol_{g}:=\sup_{s< 10\, C_0^2\,r}\ \frac{1}{s^2}\ \int_{\Sigma\,\cap\, \{\r_{\vec{p}}<s\}}  \ dvol_{g}
\ee
Since $\chi$ has been chosen such that $\chi'\le 0$ and $\chi'(t)<-1/2$ on $[5/4,7/4]$ (\ref{truncated-bis}) implies for any $\eta>\ep$
\[
\begin{array}{l}
\ds 4\,\int_\Sigma[\chi(\r/r)-\chi(\r/\eta)]\, |\nabla^\Sigma\arctan\sigma|^2 \ dvol_g+ \eta^{-2}\, \int_{\Sigma\,\cap\{5\eta/4<\r<7\eta/4\}}\  dvol_g \\[5mm]
\ds\le\frac{C}{{r}^2}\ \int_{\Sigma\,\cap\, \{r<\r_{\vec{p}}<2\,r\}}  \ dvol_{g}+ C\, \lf[\int_{\Sigma\,\cap\{\eta<\r<2\eta\}} |\nabla^\Sigma\arctan\sigma|^2\ dvol_g\rg]^{1/2}\ \lf[ \eta^{-2}\, \int_{\Sigma\,\cap\{\eta<\r<2\eta\}}\  dvol_g  \rg]^{1/2}\\[5mm]
\ds+\ C\ \lf|\int_\Sigma<d h_{r,\eta},d\beta>_g\ dvol_g\rg|
\end{array}
\]
Combining (\ref{penalisation}), (\ref{err-1}) , (\ref{err-4}) and (\ref{sup-1}) for $\eta\ge\ep$ one obtains
\[
\lf|\int_\Sigma<d h_{r,\eta},d\beta>_g\ dvol_g\rg|\le C\ \lf[ \delta^{3/4}\ \log^{3/4}\lf(\frac{r}{\ep}\rg)+  \delta\ \log\lf(\frac{r}{\ep}\rg)\rg] \ \frac{1}{s_0^2}\ \int_{\Sigma\,\cap\, \{\r_{\vec{p}}<s_0\}}  \ dvol_{g}+C\ e^{-\ep^{-2}}\ \ep^{-3}
\]
While for $\eta<\ep$ combining (\ref{mon-2a}), (\ref{prem-conc}) and (\ref{err-8}) gives
\be
\label{mon-2a-bis}
\begin{array}{l}
\ds\lf|\int_\Sigma<d h_{\ep,\eta},d\beta>_g\ dvol_g\rg|\\[5mm]
\ds\le C\ \lf[\frac{1}{s_0^2}\ \int_{\Sigma\,\cap\, \{\r_{\vec{p}}<s_0\}}  \ dvol_{g}  \rg]^{1/2}\ \lf[ \ep^{-2}\ \int_{\eta<\r<2\,\ep}\ep^4\,(1+|d\vec{T}|^2_{g})^2\ dvol_g  \rg]^{1/4}\ \lf[ \ep^{-2}\,\int_{\eta<\r<2\,\ep} \ dvol_{g}
\rg]^{1/4}\\[5mm]
\ds+ C\ \ep^{-2}\,\int_{\ep<\r<2\,\ep}\sqrt{\ep^2\,(1+|d\vec{T}|^2_{g})}\ dvol_g\le C\  \delta^{1/4}\ \frac{1}{s_0^2}\ \int_{\Sigma\,\cap\, \{\r_{\vec{p}}<s_0\}}  \ dvol_{g}+C\ e^{-\ep^{-2}}\ \ep^{-3}
\end{array}
\ee
We assume first that $s_0<r$. For $\eta:={(7/5)}^{-j}\, s_0$ for $j$ going from $-1$ to $+\infty$ we  obtain, using again the fact that $h_{r,\eta}=h_{r,\ep}+h_{\ep,\eta}$,
\be
\label{sup-2}
\begin{array}{l}
\ds (7/5)^{2\,j}\, s_0^{-2}\, \int_{\Sigma\,\cap\{5 (7/5)^{-j}\, s_0/4<\r<7(7/5)^{-j}\, s_0/4\}}\  dvol_g \le\frac{C}{{r}^2}\ \int_{\Sigma\,\cap\, \{r<\r_{\vec{p}}<2\,r\}}  \ dvol_{g}\\[5mm]
\ds+ C\, \lf[\int_{\Sigma\,\cap\{\r<2r\}} |\nabla^\Sigma\arctan\sigma|^2\ dvol_g\rg]^{1/2}\ \lf[ s_0^{-2}\, \int_{\Sigma\,\cap\{\ep<\r<s_0\}}\  dvol_g  \rg]^{1/2}\\[5mm]
\ds+C\ \lf[ \delta^{1/4}+\delta^{3/4}\ \log^{3/4}\lf(\frac{r}{\ep}\rg)+  \delta\ \log\lf(\frac{r}{\ep}\rg)\rg] \ \frac{1}{s_0^2}\ \int_{\Sigma\,\cap\, \{\r_{\vec{p}}<s_0\}}  \ dvol_{g}\\[5mm]
\ds +C\ e^{-\ep^{-2}}\ \ep^{-3}
\end{array}
\ee
Multiplying by $(7/5)^{-2\,j}$ and summing from $j=-1$ up to $+\infty$, using the fact that $7(7/5)^{-j-1}=  5 (7/5)^{-j}$ we obtain
\be
\label{sup-3}
\begin{array}{l}
\ds \frac{1}{s_0^{2}}\, \int_{\Sigma\,\cap\{\r< s_0\}}\  dvol_g \le C\,\lf[\sum_{j=-1}^{+\infty}\lf(\frac{5}{7}\rg)^{2\,j}\rg]\ \frac{1}{{r}^2}\ \int_{\Sigma\,\cap\, \{r<\r_{\vec{p}}<2\,r\}}  \ dvol_{g}\\[5mm]
\ds+ C\,\lf[\sum_{j=-1}^{+\infty}\lf(\frac{5}{7}\rg)^{2\,j}\rg]\ \lf[\frac{1}{{r}^2}\ \int_{\Sigma\,\cap\, \{\r_{\vec{p}}<2\,r\}}  \ dvol_{g}\rg]^{1/2}\ \lf[ \frac{1}{s_0^{2}}\, \int_{\Sigma\,\cap\{\ep<\r<s_0\}}\  dvol_g  \rg]^{1/2}\\[5mm]
\ds+ C\,\lf[\sum_{j=-1}^{+\infty}\lf(\frac{5}{7}\rg)^{2\,j}\rg]\  \lf[\delta\ \log\lf(\frac{r}{\ep}\rg)+\delta^{1/4} +\delta^{3/4}\ \log^{3/4}\lf(\frac{r}{\ep}\rg)\rg]^{1/2}\  \frac{1}{s_0^{2}}\, \int_{\Sigma\,\cap\{\r<s_0\}}\  dvol_g \\[5mm]
\ds+C\ e^{-\ep^{-2}}\ \ep^{-3}
\end{array}
\ee
If $r<s_0<10\, C_0^2\, r$ we trivially write 
\[
 \frac{1}{s_0^{2}}\, \int_{\Sigma\,\cap\{\r< s_0\}}\  dvol_g \le \frac{1}{r^2}\int_{\Sigma\,\cap\{\r< 10\, C_0^2\ r\}}\  dvol_g\ .
\]
We then deduce (\ref{good-2-b}) and this concludes the proof of lemma~\ref{lm-good-points}.\hfill $\Box$

From the previous we shall deduce the following lemma
\begin{Lm}
\label{lower-dens-mu}
Under the previous notations there exists a universal constant $c_\ast>0$ such that
\be
\label{low-dens-mu}
\mbox{ for }\mu_\infty\,\mbox{ a. e. }\,\vec{p}\in V_2({\R}^4)\quad,\quad\quad\forall s>0\quad\mbox{then}\quad c_\ast<\frac{\ds\mu_\infty\lf(\{\r_{\vec{p}}<s\}\rg)}{s^2}\le C\ \mu_\infty(V_2({\R}^4))\  .
\ee
\end{Lm}
\noindent{\bf Proof of lemma~\ref{lower-dens-mu}.} Let $\vec{p}\in V_2({\R}^4)$ and $s>0$ such that $\mu_\infty\lf(\{\r_{\vec{p}}<s\}\rg)\ne 0$. Let $\delta_k\rightarrow 0$ given by
 \be
 \label{bord-10-bb}
\delta_k^2:=(\log\ep_k^{-1})^{-1}\, \ep^4_k\int_\Sigma(1+|d\vec{T}_k|^2)^2\ dvol_{g_k}
 \ee
 Observe that with this choice we garantee
 \be
 \label{condition-delta-bis}
 \lim_{k\rightarrow +\infty}\delta_k\,\log\ep_k^{-1}=0\quad\mbox{ and }\quad \lim_{k\rightarrow +\infty}\delta^{-1}_k\,{\ep^4_k\int_\Sigma(1+|d\vec{T}_k|^2)^2\ dvol_{g_k}}=0\ .
 \ee
 Let $E_k$ be the set of $\vec{p}\in V_2({\R}^4)$  such that there exists $s>0$ with
 \[
 \delta_k\, {\mu}_k(\{\r_{\vec{p}}<5\, C_0^2\,s\})\le \ep^4_k\int_{B_\rho(x_0)\cap\{\r_{\vec{p}}<s\}}(1+|d\vec{T}_k|^2)^2\ dvol_{g_k}
 \]
 Observe that for $\vec{p}\ne\vec{q}$
 \[
\sup_{\vec{o}\in V_2({\R}^4)}\ \frac{\r_{\vec{p}}(\vec{q})}{\max\{\r_{\vec{p}}(\vec{o}),\r_{\vec{q}}(\vec{o})\}}\le 2\, C_0
 \]
 We can then use the Vitali type Covering Lemma in Quasi-metric spaces given by Lemma 2.7 in \cite{AlMi} and cover $E_k$ by at most countably many balls $(B^\r_{5 C_0^2\,s_i}(\vec{p}_i))_{i\in I}$  for the quasi-distance $\r$
 such that the balls $B^\r_{5 C_0^2\,s_i}(\vec{p}_i)$ for $i\in I$ are disjoint to each other. Hence we have
 \[
 \begin{array}{l}
 \ds{\mu}_k(E_k)\le\sum_{i\in I}\ti{\mu}_k(\{\r_{\vec{p_i}}<5\, C_0^2\,s_i\})\le \delta_k^{-1}\ \sum_{i\in I}\ep^4_k\int_{B_\rho(x_0)\cap\{\r_{\vec{p}_i}<s_i\}}(1+|d\vec{T}_k|^2)^2\ dvol_{g_k}\\[5mm]
 \ds\quad\quad \quad\le \delta^{-1}_k\,{\ep^4_k\int_\Sigma(1+|d\vec{T}_k|^2)^2\ dvol_{g_k}}=o(1)\ .
 \end{array}
 \]
 Assume that there exists $k_0$ for any $k\ge k_0$  as large as we want 
 \[
 \vec{\La}_k(\Sigma)\cap B^\r_s(\vec{p})\subset E_k\cap B^\r_s(\vec{p})\ ,
  \]
  then we would have
 \[
 {\mu}_k\lf(\{\vec{q}\ ;\ \r_{\vec{p}}(\vec{q})< s\}\rg)\le {\mu}_k(E_k)\ \longrightarrow 0\ ,
 \]
 and this would imply
 \[
 \mu_\infty\lf(\{\r_{\vec{p}}<s\}\rg)= 0\ ,
 \]
which contradicts our assumption. Hence, we can find $k$ as large as we want such that
\[
\lf(\vec{\La}_k(\Sigma)\cap B^\r_s(\vec{p})\rg)\bigcap \lf(V_2({\R}^4)\setminus E_k\rg)\ne \emptyset
\]
For such a $k$ we consider then $\vec{p}_k\in \vec{\La}_k(\Sigma)\cap B^\r_s(\vec{p})$ and $\vec{p}_k\notin E_k$.  Applying lemma~\ref{lm-good-points} and more specifically (\ref{good-2-b}) to this point, for the choice of $\delta_k$ we are making we obtain that for any $\eta<s$ there holds
\[
\begin{array}{l}
\ds \lf[c_0-o(1)\rg]\ \frac{1}{\eta^2}\int_{\Sigma\,\cap\, \{\r_{\vec{p}_k}<\eta\}} dvol_{g_k}\\[5mm]
\ds\quad\quad\le \frac{C}{s^2}\int_{\Sigma\,\cap\, \{\r_{\vec{p}_k}<s\}}  \ dvol_{g_k}+C\ e^{-\ep_k^{-2}}\ \ep_k^{-3}\ .
\end{array}
\]
as explained in the first part of the proof of lemma~\ref{lm-good-points},  we have
\[
\lim_{\eta\rightarrow 0}\frac{1}{\eta^2}\int_{\Sigma\,\cap\, \{\r_{\vec{p}}<\eta\}} dvol_{g_k}=\pi N_{\vec{p}_k}\ge \pi
\]
Using one more time the fact that $\r$ is a quasi distance we deduce the existence of three universal  positive constants $c_0>0$ , $c_1>1$ and $C>0$ such that, for some subsequence $k'$ we have
\[
C\,\mu_{k'}(\{\r_{\vec{p}}< c_1\, s)\ge s^2\ \pi\, [c_0-o(1)]\ .
\]
Passing to the limit we have proved the following implication
\[
\mu_\infty\lf(\{\r_{\vec{p}}<s\}\rg)\ne 0\quad\Longrightarrow \quad \mu_\infty\lf(\{\r_{\vec{p}}\le c_1\,s\}\rg)\ge s^2\ \pi\, c_0\ .
\]
Assume there exists $G$, measurable, such that
\[
\mu_\infty(G)>0\quad\mbox{ and }\quad\forall\vec{p}\in G\ \quad\exists\, s_{\vec{p}}>0\quad\mbox{ s. t. } \mu_\infty\lf(\{\r_{\vec{p}}\le c_1\,s_{\vec{p}}\}\rg)< s^2_{\vec{p}}\ \pi\, c_0
\]
We then have
\[
\forall \ \vec{p}\in G\quad\exists\, s_{\vec{p}}>0\quad \mbox{s. t. }\quad\mu_\infty\lf(\{\r_{\vec{p}}<s_{\vec{p}}\}\rg)= 0
\]
It implies that $\mu_\infty(G)=0$ which is a contradiction and this concludes the proof of lemma~\ref{lower-dens-mu}.
  \hfill $\Box$
  
\section{The limit of the area density of $\vec{\La}_k$.}
\reset

In this section $\vec{\La}_k$ denotes an admissible sequence of almost critical points of $E_\ep$. We assume first that the conformal class of the induced metric $g_k:=\vec{\La}_k^\ast g_{V_2({\R}^4)}$ is constant (the general case will be treated in a subsequent section) and, thanks to the uniformization principle, we can then compose by a diffeomorphism such that there exists a fixed constant Gauss curvature $h$ on $\Sigma$ such that
\[
g_k=\vec{\La}_k^\ast \,g_{V_2({\R}^4)}=e^{2\la_k}\, h\ .
\]
We are interested in this section with the limiting behaviour of
\[
\nu_k:=\frac{1}{2} |d\vec{\La}_k|^2_h\ dvol_h=e^{2\la_k}\ dvol_h\ .
\]
Again we denote 
${\mathbf v}_{k}$ to be the  2-varifold associated to the immersion $\vec{\La}_k$
\[
\forall \ \Xi\in C^0(G_2^{Leg}(V_2({\R}^4)))\quad {\mathbf v}_{k}(\Xi):=\int_\Sigma\Xi((\vec{\La}_k)_\ast T_x\Sigma,\vec{\La}_k(x))\ dvol_{g_k}
\]
and by $\mu_k$ the Radon measure given by the integration along the fibers that is
\[
\forall \ A\subset V_2({\R}^4)\quad\mbox{ Borel }\quad\mu_k(A):={\mathbf v_k}(\pi^{-1}(A))\ .
\]
which is $\mu_k:=\pi_\ast {\mathbf v}_k$.  We clearly have for any $\Psi\in C^0(V_2({\R}^4))$
\[
<\mu_k,\Psi>=\int_\Sigma \Psi(\vec{\La}_k)\ dvol_{g_k}=<\nu_k,\Psi\circ\vec{\La}_k>\ ,
\]
or in other words
\[
(\vec{\La}_k)_\ast\nu_k=\mu_k\ .
\]
Modulo extraction of a subsequence we can assume that in the sense of Radon measure convergence we have the existence of $\nu_\infty$, $\mu_\infty$ and ${\mathbf v}_\infty$ such that
\[
\nu_k\rightharpoonup\nu_\infty\quad,\quad \mu_k\rightharpoonup\mu_\infty\quad\mbox{ and }\quad {\mathbf v}_k\rightharpoonup{\mathbf v}_\infty\ .
\]
Since the $W^{1,2}-$norm of $\vec{\La}_k$ is bounded on $(\Sigma,h)$ we can also assume that
\[
\vec{\La}_k\ \rightharpoonup\ \vec{\La}_\infty\mbox{ weakly in }W^{1,2}_h(\Sigma,V_2({\R}^4))\ .
\]
Thanks to the theorem of  Rellich Kondrachov we have that $\vec{\La}_k$ converges strongly to $\vec{\La}_\infty$ in $L^2(\Sigma)$ and for any smooth function $f$ on $V_2({\R}^4)$ (which is compact) we have that $f\circ\vec{\La}_k$ is also strongly converging in $L^2$ towards $f\circ\vec{\La}_\infty$. Hence in particular
\[
\vec{\La}_k^{\,\ast}\al\ \rightharpoonup\ \vec{\La}_\infty^{\,\ast}\al\quad\mbox{ in }{\mathcal D}'(\Sigma)\ .
\]
Since $0=\vec{\La}_k^{\,\ast}\al$ we deduce that $\vec{\La}_\infty^{\,\ast}\al=0$ and $\vec{\La}_\infty$ is then {\it weakly Legendrian}. There is however no reason a-priori for $\vec{\La}_\infty$ to be weakly conformal from $(\Sigma,h)$ into $V_2({\R}^4)$. We now prove the following lemma which is maybe the most important building-block of the present work after the almost monotonicity formula.


\begin{Lm}
\label{lm-energy-quant}
There exists a universal constant $c_Q>0$  such that, for any $x_0\in \Sigma$,  $\rho>0$  such that for some subsequence $k'$
\be
\label{hyp-1}
\limsup_{k'\rightarrow +\infty}\int_{\p B_\rho(x_0)}|d\vec{\La}_{k'}|_h\ dl_h<+\infty\ ,
\ee
and two positive constants $c_1, c_2>0$  depending only on 
\be
\label{const-c1-c2}
\Theta:=\limsup_{k'\rightarrow +\infty}\ \sup_{r<1}\,\sup_{ \vec{p}\in V_2({\R}^4)}\frac{\mu_k(\{\r_{\vec{p}}<r\})}{r^2}\ ,
\ee
such that if
\be
\label{hyp-2}
0<t:=\limsup_{k'\rightarrow +\infty}\sup_{\{\vec{p}\,,\,\vec{q}\,\in \ \vec{\La}_{k'}(\p B_\rho(x_0))\}}\r_{\vec{p}}(\vec{q}\,)\le c_1
\ee
then either
\be
\label{concl-1}
\nu_\infty(B_\rho(x_0))>c_Q\ ,
\ee
or
\be
\label{concl-2}
\limsup_{k'\rightarrow+\infty}\int_{ B_\rho(x_0)}|d\vec{\La}_{k'}|^2_h\ {\mathbf 1}_{\ds \lf\{x\in B_\rho(x_0)\ ;\ \r^\ast_{\vec{\La}_{k'}(x)}> c_2\, t\rg\}}\ dvol_h=0\ ,
\ee
where 
\[
\r^\ast_{\vec{p}}:=\inf\lf\{\r_{\vec{p}}(\vec{q})\quad;\quad \vec{q}\in \vec{\La}_{k'}(\p B_\rho(x_0)) \rg \}\ ,
\]
and $ {\mathbf 1}_{\lf\{x\ ;\ \r^\ast_{\vec{\La}_{k'}(x)}> c_2\, t\rg\}}$ is the characteristic function of the points $x\in B_\rho(x_0)$ such that $\r^\ast_{\vec{\La}_{k'}(x)}> c_2\, t$. In particular, if (\ref{concl-1}) does not hold, we have
\be
\label{concl-3}
\begin{array}{rl}
\ds\nu_\infty(B_\rho(x_0))&\ds\le\limsup_{k\rightarrow+\infty}\nu_k(B_\rho(x_0))=\limsup_{k\rightarrow+\infty}\lf[\lf(\lf.\vec{\La}_k\rg|_{B_\rho(x_0)}\rg)_\ast\nu_k\rg](V_2({\R}^4))\\[5mm]
\ds&\ds\le\limsup_{k\rightarrow +\infty}\mu_k\lf(\{\vec{p}\ ;\ \r^\ast_{\vec{p}}<c_2\, t\}\rg)\le C\, t^2\ \Theta\ .
\end{array}
\ee
\end{Lm}
\noindent{\bf Proof of Lemma~\ref{lm-energy-quant}.}
We denote by $k$ the subsequence given by the hypothesis of the lemma. Let $\vec{\gamma}(s)=(\vec{a}(s),\vec{b}(s))$ be a $C^1$ Legendrian path passing through $\vec{p}_0:=(\vec{\ep}_1,\vec{\ep}_2)$ at the time $s=0$. We denote $\vec{\gamma}(s):=(\vec{a}(s),\vec{b}(s))=(\vec{\ep}_1+\vec{v}(s),\vec{\ep}_2+\vec{w}(s))$. Recall 
\[
\varphi_{\vec{p}_0}(\vec{\gamma}(s))=v_2(s)-w_1(s)\ .
\]
Hence we have from (\ref{liouv})
\[
\frac{d\varphi_{\vec{p}_0}}{ds}=\frac{d(v_2-w_1)}{ds}=\sum_{i=1}^4 v_i(s)\,\frac{dw_i}{ds}-w_i(s)\,\frac{dv_i}{ds}\ .
\]
We deduce for any $s_0>0$
\[
\begin{array}{l}
\ds|\varphi_{\vec{p}_0}|(s_0)\le\,\int_0^{s_0}\lf|\frac{d\varphi_{\vec{p}_0}}{ds}\rg|\ ds\le \|\vec{v}\|_{L^\infty([0,t],{\R}^4)}\ \int_0^{s_0}\lf|\frac{d\vec{w}}{ds}\rg|\ ds+ \|\vec{w}\|_{L^\infty([0,s_0],{\R}^4)}\ \int_0^{s_0}\lf|\frac{d\vec{v}}{ds}\rg|\ ds\\[5mm]
\ds\le  2\ \int_0^{s_0}\lf|\frac{d\vec{v}}{ds}\rg|\ ds\int_0^{s_0}\lf|\frac{d\vec{w}}{ds}\rg|\ ds\le   \lf[\int_0^{s_0}\lf|\frac{d\vec{v}}{ds}\rg|\ ds\rg]^2+ \lf[\int_0^{s_0}\lf|\frac{d\vec{w}}{ds}\rg|\ ds\rg]^2\le 2\  \lf[\int_0^{s_0}\lf|\frac{d\vec{\gamma}}{ds}\rg|\ ds\rg]^2
\end{array}
\]
While we proved it for notations conveniences for $\vec{p}_0:=(\vec{\ep}_1,\vec{\ep}_2)$, it holds by homogeneity of $V_2({\R}^4)$ for any choice of $\vec{p}_0$.    This gives for $x(s):=\rho\,(\cos s,\sin s)$ and any $s\in[0,2\pi]$ and $\vec{p}_0:=\vec{\La}_{k'}(x(0)))$
\be
\label{res-00}
\begin{array}{rl}
\ds \r_{\vec{p}_0}\lf(\vec{\La}_{k'}(x(s))\rg)&\ds=\sqrt[4]{ \lf|\int_0^s\frac{d(\vec{\La}_{k'}\circ x)(s)}{ds}\ ds\rg|^4+4\,\varphi^2_{\vec{p}_0}\lf(\vec{\La}_{k'}(x(s))\rg)   } \\[5mm]
\ds\quad&\ds\le \sqrt[4]{17}\ \int_0^s\lf|\frac{d(\vec{\La}_{k'}\circ x)(s)}{ds}\rg|\ ds
\end{array}
\ee
We deduce from (\ref{hyp-1})
\be
\label{res-0000}
t:=\limsup_{k\rightarrow +\infty}\sup_{\{\vec{p}\,,\,\vec{q}\,\in \ \vec{\La}_{k'}(\p B_\rho(x_0))\}}\r_{\vec{p}}(\vec{q}\,)<+\infty\ .
\ee
Hence
\be
\label{res-1}
\sup_{\{\vec{p}\,,\,\vec{q}\,\in \ \vec{\La}_\infty(\p B_\rho(x_0))\}}\r_{\vec{p}}(\vec{q}\,)\le \limsup_{k\rightarrow +\infty}\sup_{\{\vec{p}\,,\,\vec{q}\,\in \ \vec{\La}_{k'}(\p B_\rho(x_0))\}}\r_{\vec{p}}(\vec{q}\,)= t\ .
\ee
We pick $x_1\in \p B_\rho(x_0)$ such that $\vec{\La}_{k}(x_1)\rightarrow \vec{\La}_\infty(x_1)=\vec{p}_1$ (since  $\vec{\La}_{k}$ is weakly converging to $\vec{\La}_\infty$ in $W^{1,2}(\Sigma)$ the trace on $\p B_\rho(x_0)$
is converging weakly in $H^{1/2}$ and hence almost everywhere thanks to Rellich Kondrachov).

We recall  - see (\ref{quasi-dist}) -  that the function $$d(\vec{p},\vec{q}):=\r_{\vec{p}}(\vec{q})=\lf[ |\vec{p}-\vec{q}|^4+4\,|\vec{a}\cdot\vec{d}-\vec{b}\cdot\vec{c}|^2    \rg]^{1/4}\quad\mbox{ where }\quad \vec{p}=(\vec{a},\vec{b}) \quad\mbox{ and }\vec{q}=(\vec{c},\vec{d}) $$ satisfies the quasi-distance condition
\be
\label{quasi-dist-rrr}
\r_{\vec{p}}(\vec{q})=\r_{\vec{q}}(\vec{p})\quad\mbox{ and }\quad d(\vec{p}_1,\vec{p})\le C_0\,\lf[d(\vec{p}_1,\vec{q}) +d(\vec{p},\vec{q})\rg]
\ee
for some universal $C_0>1$. In particular, for any $\vec{p}\in V_2({\R}^4)$ and any $\vec{q}\in \vec{\La}_k(\p B_\rho(x_0))$ there holds
\[
\r_{\vec{p}_1}(\vec{p}\,)\le C_0\ \r_{\vec{p}_1}(\vec{q})+C_0\ \r_{\vec{p}}(\vec{q}\,)\ ,
\]
which implies for $k$ large enough
\[
\r_{\vec{p}_1}(\vec{p})> 8\,C_0 \ t\quad\Longrightarrow\quad\inf_{\vec{q}\in \vec{\La}_{k}(\p B_\rho(x_0))}\r_{\vec{p}}(\vec{q})\ge 4\ t\ .
\]
We denote by $\ti{\mu}_k:=(\vec{\La}_k| B_\rho(x_0))_\ast\nu_k$ and we denote by $\ti{\mathbf v}_k$ the corresponding varifold and $\ti{\mu}_\infty$ is the limit of $\ti{\mu}_k$.

We introduce for the same function $\chi$ as in the previous section which is equal to 1 in $[0,1]$ equal to zero on $[2,+\infty)$
and satisfying $\chi'\le 0$ and consider
\[
\zeta(r):=1-\chi\lf(\frac{r}{8\, C_0\,s}\rg)\ .
\]
where $s>t$ will be chosen later. The function $\zeta$ has been constructed in such a way that, for $k'$ large enough, $\zeta\circ \r_{{\vec{p}}_1}$ is supported away from an open neighbourhood of $\vec{\La}_{k}(\p B_\rho(x_0))$. More precisely
\be
\label{zeta-1}
\zeta(\r_{{\vec{p}}_1}(\vec{p}))\ne 0\quad\Longrightarrow\quad \inf_{\vec{q}\in \vec{\La}_{k}(\p B_\rho(x_0))}\r_{\vec{p}}(\vec{q})\ge 2\ t\ .
\ee
We shall now be considering the infinitesimal variation for $E_{\ep_k}$ at the restriction of $\vec{\La}_k$ to $B_\rho(x_0)$ for the Hamiltonian given by 
\[
\ti{h}_{r,\eta}:=\zeta\, h_{r,\eta}=\zeta(\r_{\vec{p}_1})\ [ \chi(\r_{\vec{p}}/r)- \chi(\r_{\vec{p}}/\eta)]\ \arctan\sigma_{\vec{p}}
\]
 for some $\vec{p}$, $r>0$ and $\eta>0$ to be fixed later such that 
 \[
 \lf\{
 \begin{array}{l}
\ds \zeta\circ\r_{\vec{p}_1}(\vec{q})\equiv1\quad\mbox{ for }\quad \r_{\vec{p}}(\vec{q})<2\,\eta\ ,\\[5mm]
\ds \zeta\circ\r_{\vec{p}_1}(\vec{q})\equiv1\quad\mbox{ for }\quad \r_{\vec{p}}(\vec{q})>r\ .
\end{array}
\rg.
 \]
  We denote
 \[
 \r^{\,\ast}_{\vec{p}}:=\inf\lf\{\r_{\vec{p}}(\vec{q})\ ;\ \vec{q}\in \vec{\La}_{k'}(\p B_\rho(x_0))\rg\}\ .
 \]
 Because of (\ref{zeta-1}) we have $\r^{\,\ast}_{\vec{p}}>2\,t$. We consider respectively
 \[
 \vec{X}_{\ti{h}_{r,\eta}}:= J_H\nabla^H\ti{h}_{r,\eta}+\frac{\ti{h}_{r,\eta}}{2}\ \vec{R}\quad\mbox{ and }\quad\vec{w}_{{\ti{h}_{r,\eta}}}:= \vec{X}_{\ti{h}_{r,\eta}}\circ \vec{\La}_k\ .
 \]
 To simplify the notations, in the rest of the proof of the lemma we write
 \[
 \r:=\r_{\vec{p}}\quad ,\quad \sigma:=\sigma_{\vec{p}}\quad ,\quad \r_1:=\r_{\vec{p}_1}\quad ,\quad \sigma_1:=\sigma_{\vec{p}_1}\quad\mbox{ and }\quad \r^\ast:=\r^\ast_{\vec{p}}\ .
 \]
 We multiply (\ref{identite-monotone}) by $\zeta(\r_1)\ [\chi(\r/r)-\chi(\r/\eta)]$, (\ref{autre-identite}) by $\r^{-3}\,\arctan\sigma\, \zeta(\r_1)\ [r^{-1}\chi'(\r/r)-\eta^{-1}\chi'(\r/\eta)]$ and (\ref{autre-identite}) for $\r_1$ and $\sigma_1$ instead of $\r$ and $\sigma$    by $\r^{-3}_1\ (8\, C_0\, s)^{-1}\,\zeta'(\r_1)\ [\chi(\r/r)-\chi(\r/\eta)] \ \arctan\sigma$ this gives
\[
\begin{array}{l}
\ds<d \ti{h}_{r,\eta},d\beta>_g+\zeta(\r_1)\ [\chi(\r/r)-\chi(\r/\eta)]\,\Delta_g\log\r    =   4\,\zeta(\r_1)\ [\chi(\r/r)-\chi(\r/\eta)]\,\frac{|(\nabla^\Sigma\r)^\perp|^2}{\r^2}\\[5mm]
\ds+\,O(1)\,\zeta(\r_1)\ [\chi(\r/r)-\chi(\r/\eta)]+O(\r)\, \zeta(\r_1)\ [r^{-1}\chi'(\r/r)-\eta^{-1}\chi'(\r/\eta)]\\[5mm]
\ds+\,\zeta(\r_1)\ [r^{-1}\chi'(\r/r)-\eta^{-1}\chi'(\r/\eta)]\,\r^{-1}\,\arctan\sigma\,\frac{\sigma}{\sqrt{1+\sigma^2}}\\[5mm]
\ds+\,4^{-1}\,\r^{-3} \ \zeta(\r_1)\  [r^{-1}\chi'(\r/r)-\eta^{-1}\chi'(\r/\eta)]\ \arctan\sigma\,\mbox{div}^\Sigma\lf(\r^4\frac{\nabla^\Sigma\sigma}{1+\sigma^2}\rg)\\[5mm]
\ds +\,(8\, C_0\,s)^{-1}\,\zeta'(\r_1)\ [\chi(\r/r)-\chi(\r/\eta)]\ \,\arctan\sigma\, O(\r_1)\\[5mm]
\ds+\ (8\, C_0\,s)^{-1}\,\zeta'(\r_1)\ [\chi(\r/r)-\chi(\r/\eta)]\ \,\arctan\sigma\ \r_1^{-1}\,\frac{\sigma_1}{\sqrt{1+\sigma^2_1}}\\[5mm]
\ds+\,4^{-1}\,\ (8\, C_0\, s)^{-1}\, \r_1^{-3}\, \zeta'(\r_1)\ [\chi(\r/r)-\chi(\r/\eta)]\ \,\arctan\sigma\,\,\mbox{div}^\Sigma\lf(\r^4_1\frac{\nabla^\Sigma\sigma_1}{1+\sigma^2_1}\rg)
\end{array}
\]
Integrating over $B_\rho(x_0)$ gives
\be
\label{mono-bord}
\begin{array}{l}
\ds\int_{B_\rho(x_0)}<d \ti{h}_{r,\eta},d\beta>_g\ dvol_g-\int_{B_\rho(x_0)}\,\zeta(\r_1)\ \lf[\frac{\r}{r}\,\chi'\lf(\frac{\r}{r}\rg)-\frac{\r}{\eta}\,\chi'\lf(\frac{\r}{\eta}\rg)\rg]\frac{|\nabla^\Sigma\r|^2}{\r^2}\ dvol_g\\[5mm]
\ds+\int_{B_\rho(x_0)}\ O(1)\,\zeta(\r_1)\ \ [\chi(\r/r)-\chi(\r/\eta)]+\, \zeta(\r_1)\,[O(\r/r)\,\chi'(\r/r)-O(\r/\eta)\chi'(\r/\eta)]\ dvol_g \\[5mm]
\ds+4^{-1}\,\int_{B_\rho(x_0)} \zeta(\r_1)\,\nabla^\Sigma\lf[\r^{-3} \ [r^{-1}\chi'(\r/r)-\eta^{-1}\chi'(\r/\eta)]\ \arctan\sigma\rg]\,\r^4\frac{\nabla^\Sigma\sigma}{1+\sigma^2}\ dvol_g \\[5mm]
\ds-\int_{B_\rho(x_0)} \zeta(\r_1)\ \lf[\frac{\r}{r}\,\chi'\lf(\frac{\r}{r}\rg)-\frac{\r}{\eta}\,\chi'\lf(\frac{\r}{\eta}\rg)\rg]\,\frac{1}{\r^{2}}\,\frac{\sigma\,\arctan\sigma}{\sqrt{1+\sigma^2}}\ dvol_g \\[5mm]
\ds-4\,\int_{B_\rho(x_0)}\, \zeta(\r_1)\ [\chi(\r/r)-\chi(\r/\eta)]\,\frac{|(\nabla^\Sigma\r)^\perp|^2}{\r^2}\ dvol_g\\[5mm]
\ds=-\,4^{-1}\,\int_{B_\rho(x_0)}\,  \frac{\zeta'(\r_1)}{8\, C_0\,s}\  \lf[\frac{\r}{r}\chi'\lf(\frac{\r}{r}\rg)-\frac{\r}{\eta}\ \chi'\lf(\frac{\r}{\eta}\rg)\rg]\ \arctan\sigma\,\frac{\nabla^\Sigma\r_1\cdot\nabla^\Sigma\sigma}{1+\sigma^2}\ dvol_g\\[5mm]
\ds+\int_{B_\rho(x_0)}\,   \lf[\frac{3}{4}\,\frac{\zeta'(\r_1)}{8\, C_0\,s}- \frac{\r_1\,\zeta''(\r_1)}{(8\, C_0\,s)^2}  \rg]\ [\chi(\r/r)-\chi(\r/\eta)]\ \frac{\nabla^\Sigma\r_1\cdot\nabla^\Sigma\sigma_1}{1+\sigma^2_1}\ dvol_g\\[5mm]
\ds+\int_{B_\rho(x_0)}\, \,\frac{\zeta'(\r_1)}{8\, C_0\,s}\ [\chi(\r/r)-\chi(\r/\eta)]\ \lf[  O(\r_1)+\r^{-1}\,\nabla^\Sigma\r_1\cdot\nabla^\Sigma\r + \arctan\sigma\ \r_1^{-1}\,\frac{\sigma_1}{\sqrt{1+\sigma^2_1}}  \rg]\ dvol_g
\end{array}
\ee
We have the following bounds (using (\ref{mono-every}) in lemma~\ref{lm-finite-dens}), since $\zeta'\circ\r_1$ is supported in the domain $8\,C_0\, s< \r_1<16\,C_0\,s$
\be
\label{bord-1}
\begin{array}{l}
\ds\limsup_{k\rightarrow +\infty}\lf|\int_{B_\rho(x_0)}\, \,\frac{\zeta'(\r_1)}{8\, C_0\,s}\ [\chi(\r/r)-\chi(\r/\eta)]\ \lf[  O(\r_1)+\r^{-1}\,\nabla^\Sigma\r_1\cdot\nabla^\Sigma\r \rg]\ \ dvol_g\rg|\\[5mm]
\ds\le C\, s^{-1}\, (1+\r^{-1}_\ast)\,\ti{\mu}_\infty(\r_1< 16\,C_0\,s)\le C\, s\,(1+\r^{-1}_\ast)\, {\mu}_\infty(\r_1^{-1}[0,1])\ .
\end{array}
\ee
We are choosing $\eta$ and $r$ such that $\zeta'(\r_1)\,\chi'\lf(\frac{\r}{r}\rg)\equiv 0$ and  $\zeta'(\r_1)\,\chi'\lf(\frac{\r}{\eta}\rg)\equiv 0$
hence
\be
\label{bord-2}
\begin{array}{l}
\ds\lf|4^{-1}\,\int_{B_\rho(x_0)}\,  \frac{\zeta'(\r_1)}{8\, C_0\,s}\  \lf[\frac{\r}{r}\chi'\lf(\frac{\r}{r}\rg)-\frac{\r}{\eta}\ \chi'\lf(\frac{\r}{\eta}\rg)\rg]\ \arctan\sigma\,\frac{\nabla^\Sigma\r_1\cdot\nabla^\Sigma\sigma}{1+\sigma^2}\ dvol_g\rg|=0\ .
\end{array}
\ee
We have also 
\be
\label{bord-3}
\begin{array}{l}
\ds\limsup_{k\rightarrow +\infty}\lf|\int_{B_\rho(x_0)}\,   \lf[\frac{3}{4}\,\frac{\zeta'(\r_1)}{8\, C_0\,s}- \frac{\r_1\,\zeta''(\r_1)}{(8\, C_0\,s)^2}  \rg]\ [\chi(\r/r)-\chi(\r/\eta)]\ \frac{\nabla^\Sigma\r_1\cdot\nabla^\Sigma\sigma_1}{1+\sigma^2_1}\ dvol_g\rg|\\[5mm]
\ds\le C\, (C_0\,s)^{-1}\ [\ti{\mu}_\infty(8\,C_0\,s<\r_1< 16\,C_0\,s)]^{1/2}\ \lf[\int_{8\,C_0\,s<\r_1(\vec{q})< 16\,C_0\,s}\frac{|\nabla^{\mathcal P}\sigma_1|^2}{(1+\sigma_1^2)^2}\ d\ti{\mathbf v}_\infty({\mathcal P},\vec{q})\rg]^{1/2}
\end{array}
\ee
 We recall from (\ref{id-sigma}) that away from $\r=0$
\be
\label{id-sigma-bis}
|\nabla^H\arctan\sigma|=\lf|2\,\frac{\rho^2}{\r^4}\,\nabla\varphi- \frac{2}{\r^4}\, \varphi\,\nabla\rho^2\rg|\le\frac{C}{\r} \ .
\ee
We have
\[
\lf\| \nabla^H\lf[  [\chi(\r/r)-\chi(\r/\eta)]\  \arctan\sigma\rg]\rg\|_{L^\infty(\r_1<16\,C_0\ s)}=\lf\| \nabla^H  \arctan\sigma\rg\|_{L^\infty(\r_1<16\,C_0\ s)}\ .
\]
Hence we deduce from (\ref{id-sigma-bis})
\be
\label{bd}
\lf\| \nabla^H\lf[  [\chi(\r/r)-\chi(\r/\eta)]\  \arctan\sigma\rg]\rg\|_{L^\infty(\r_1<16\,C_0\ s)}\le \frac{C}{\r^\ast}\ .
\ee
Hence
\be
\label{bord-3-a}
\begin{array}{l}
\ds\zeta'(\r_1(\vec{q}))\ne 0\ \Longrightarrow\ \\[5mm]
\ds \lf|[\chi(\r(\vec{q})/r)-\chi(\r(\vec{q}/\eta)]\  \arctan\sigma(\vec{q})-[\chi(\r(\vec{p}_1)/r)-\chi(\r(\vec{p}_1)/\eta)]\  \arctan\sigma(\vec{p}_1)\rg|\\[5mm]
\ds\quad=\lf|\arctan\sigma(\vec{q})- \arctan\sigma(\vec{p}_1)\rg|\le C\ C_0\, \frac{s}{\r^\ast}\ .
\end{array}
\ee
We consider (\ref{autre-identite}) with $\r_1$ and $\sigma_1$ instead of $\r$ and $\sigma$ 
\[
\r_1^3\, <d\r_1,d\beta>=\r_1^2\,\frac{\sigma_1}{\sqrt{1+\sigma_1^2}}+O(\r_1^4)+4^{-1}\,\mbox{div}^\Sigma\lf(\r_1^4\frac{\nabla^\Sigma\sigma_1}{1+\sigma_1^2}\rg)
\]
and we multiply this identity by $\zeta'(\r_1)\ \r_1^{-3}$ and integrate over $B_\rho(x_0)$ this gives
\be
\label{bord-4}
\begin{array}{l}
\ds\int_{B_\rho(x_0)}<d(\zeta\circ\r_1), d\beta>\ dvol_g=\int_{B_\rho(x_0)}\zeta'(\r_1)\, \lf[\r_1^{-1}\,\frac{\sigma_1}{\sqrt{1+\sigma_1^2}}+O(\r_1)\rg]\ dvol_g\\[5mm]
\ds +\frac{1}{4}\, \int_{B_\rho(x_0)}\zeta'(\r_1)\ \r_1^{-3}\,\mbox{div}^\Sigma\lf(\r_1^4\frac{\nabla^\Sigma\sigma_1}{1+\sigma_1^2}\rg)\ dvol_g
\end{array}
\ee
Since $\zeta\circ\r_1$ is a smooth function supported away from $\vec{\La}_k(\p B_\rho(x_0))$ lemma~\ref{lm-passage} is implying
\be
\label{bord-5}
\lim_{k\rightarrow 0}\int_{B_\rho(x_0)}<d(\zeta\circ\r_1), d\beta>\ dvol_g=0
\ee
Hence, after integrating by parts the last integral in (\ref{bord-4}) one gets
\be
\label{bord-6}
\begin{array}{l}
\ds\limsup_{k\rightarrow 0}\lf|\int_{B_\rho(x_0)}\zeta'(\r_1)\, \r_1^{-1}\,\frac{\sigma_1}{\sqrt{1+\sigma_1^2}}\ dvol_g\rg|\\[5mm]
\ds\quad\le [\ti{\mu}_\infty(8\,C_0\,s<\r_1< 16\,C_0\,s)]^{1/2}\ \lf[\int_{8\,C_0\,s<\r_1(\vec{q})< 16\,C_0\,s}\frac{|\nabla^{\mathcal P}\sigma_1|^2}{(1+\sigma_1^2)^2}\ d\ti{\mathbf v}_\infty({\mathcal P},\vec{q})\rg]^{1/2}
\end{array}
\ee
Combining (\ref{bord-3-a}) and (\ref{bord-6}) gives
\be
\label{bord-7}
\begin{array}{l}
\ds\limsup_{k\rightarrow 0}\lf|\int_{B_\rho(x_0)}\, \,\frac{\zeta'(\r_1)}{8\, C_0\,s}\ [\chi(\r/r)-\chi(\r/\eta)]\  \arctan\sigma\ \r_1^{-1}\,\frac{\sigma_1}{\sqrt{1+\sigma^2_1}} \ dvol_g\rg|\\[5mm]
\ds\le \frac{|[\chi(\r(\vec{p}_1)/r)-\chi(\r(\vec{p}_1)/\eta)]\  \arctan\sigma(\vec{p}_1)|}{8\, C_0\,s}\, \limsup_{k\rightarrow 0}\lf|\int_{B_\rho(x_0)}\zeta'(\r_1)\, \r_1^{-1}\,\frac{\sigma_1}{\sqrt{1+\sigma_1^2}}\ dvol_g\rg|\\[5mm]
\ds +C\ C_0\, s\ (\r^\ast)^{-1} \limsup_{k\rightarrow 0}\lf|\int_{B_\rho(x_0)}\,  \lf|\frac{\zeta'(\r_1)}{8\, C_0\,s}\rg|\  \r_1^{-1}\,\lf|\frac{\sigma_1}{\sqrt{1+\sigma^2_1}}\rg|   \ dvol_g \rg|\\[5mm]
\ds\le C\ (C_0\, s)^{-1}\  [\mu_\infty(C_0\,s<\r_1< 2\,C_0\,s)]^{1/2}\ \lf[\int_{C_0\,s<\r_1(\vec{q})< 2\,C_0\,s}\frac{|\nabla^{\mathcal P}\sigma_1|^2}{(1+\sigma_1^2)^2}\ d\ti{\mathbf v}_\infty({\mathcal P},\vec{q})\rg]^{1/2}\\[5mm]
\ds\ +C\ (C_0\, s)^{-1}\ (\r^\ast)^{-1} \ \ti{\mu}_\infty(C_0\,s<\r_1< 2\,C_0\,s)\le C \ C_0\, s\ (\r^\ast)^{-1} \ {\mu}_\infty(\{\r_1<1\})
\end{array}
\ee
Let $N\in {\N}$ such that $\r^\ast >2^{N+1}t$, because of lemma~\ref{lm-finite-dens} there exists $j_k\in\{1\cdots N\}$ such that 
\be
\label{bord-8}
\int_{C_0\,2^{j_k}\, t<\r_1(\vec{q})< 2^{j_k+1}\,C_0\,t}\frac{|\nabla^{\mathcal P}\sigma_1|^2}{(1+\sigma_1^2)^2}\ d\ti{\mathbf v}_\infty({\mathcal P},\vec{q})\le \frac{C}{N}\ {\mu}_\infty(\{\r_1<1\})
\ee
We take $s:=2^{j_k}\,t$ and we shall fix $N$ later on. Combining (\ref{bord-1})...(\ref{bord-6}) we obtain
\be
\label{bord-9}
\begin{array}{l}
\ds\limsup_{k\rightarrow+\infty}\lf|-\,4^{-1}\,\int_{B_\rho(x_0)}\,  \frac{\zeta'(\r_1)}{8\, C_0\,s}\  \lf[\frac{\r}{r}\chi'\lf(\frac{\r}{r}\rg)-\frac{\r}{\eta}\ \chi'\lf(\frac{\r}{\eta}\rg)\rg]\ \arctan\sigma\,\frac{\nabla^\Sigma\r_1\cdot\nabla^\Sigma\sigma}{1+\sigma^2}\ dvol_g\rg.\\[5mm]
\ds+\int_{B_\rho(x_0)}\,   \lf[\frac{3}{4}\,\frac{\zeta'(\r_1)}{8\, C_0\,s}- \frac{\r_1\,\zeta''(\r_1)}{(8\, C_0\,s)^2}  \rg]\ [\chi(\r/r)-\chi(\r/\eta)]\ \frac{\nabla^\Sigma\r_1\cdot\nabla^\Sigma\sigma_1}{1+\sigma^2_1}\ dvol_g\\[5mm]
\ds\lf.+\int_{B_\rho(x_0)}\, \,\frac{\zeta'(\r_1)}{8\, C_0\,s}\ [\chi(\r/r)-\chi(\r/\eta)]\ \lf[  O(\r_1)+\nabla^\Sigma\r_1\cdot\nabla^\Sigma\r + \arctan\sigma\ \r_1^{-1}\,\frac{\sigma_1}{\sqrt{1+\sigma^2_1}}  \rg]\ dvol_g\rg|\\[5mm]
\ds\le C\, 2^{N}\, t\ {\mu}_\infty(\{\r_1<1\})+C \ C_0\, 2^{N}\, t\ (\r^\ast)^{-1} \ {\mu}_\infty(\{\r_1<1\})+ C \ \frac{1}{\sqrt{N}}\,{\mu}_\infty(\{\r_1<1\})\ .
\end{array}
\ee
 Let $\delta_k\rightarrow 0$ given by
 \be
 \label{bord-10}
\delta_k^2:=(\log\ep_k^{-1})^{-1}\, \ep^4_k\int_\Sigma(1+|d\vec{T}_k|^2)^2\ dvol_{g_k}
 \ee
 Observe that with this choice we garantee
 \be
 \label{condition-delta}
 \lim_{k\rightarrow +\infty}\delta_k\,\log\ep_k^{-1}=0\quad\mbox{ and }\quad \lim_{k\rightarrow +\infty}\delta^{-1}_k\,{\ep^4_k\int_\Sigma(1+|d\vec{T}_k|^2)^2\ dvol_{g_k}}=0\ .
 \ee
 Let $F_k$ be the set of $\vec{p}\in V_2({\R}^4)$ with $\r_{\vec{p}_1}(\vec{p})> N\, 2^{N+1}t$ such that there exists $s>0$ with
 \[
 \delta_k\, \ti{\mu}_k(\{\r_{\vec{p}}<5\, C_0^2\,s\})\le \ep^4_k\int_{B_\rho(x_0)\cap\{\r_{\vec{p}}<s\}}(1+|d\vec{T}_k|^2)^2\ dvol_{g_k}\ .
 \]
 Observe that for $\vec{p}\ne\vec{q}$
 \[
\sup_{\vec{o}\in V_2({\R}^4)}\ \frac{\r_{\vec{p}}(\vec{q})}{\max\{\r_{\vec{p}}(\vec{o}),\r_{\vec{q}}(\vec{o})\}}\le 2\, C_0\ .
 \]
 We can then use the Vitali type Covering Lemma in Quasi-metric spaces given by Lemma 2.7 in \cite{AlMi} and cover $F_k$ by at most countably many balls $(B^\r_{5 C_0^2\,s_i}(\vec{p}_i))_{i\in I}$  for the quasi-distance $\r$
 such that the balls $B^\r_{5 C_0^2\,s_i}(\vec{p}_i)$ for $i\in I$ are disjoint to each other. Hence we have
 \[
 \begin{array}{l}
 \ds\ti{\mu}_k(F_k)\le\sum_{i\in I}\ti{\mu}_k(\{\r_{\vec{p_i}}<5\, C_0^2\,s_i\})\le \delta_k^{-1}\ \sum_{i\in I}\ep^4_k\int_{B_\rho(x_0)\cap\{\r_{\vec{p}_i}<s_i\}}(1+|d\vec{T}_k|^2)^2\ dvol_{g_k}\\[5mm]
 \ds\quad\quad \quad\le \delta^{-1}_k\,{\ep^4_k\int_\Sigma(1+|d\vec{T}_k|^2)^2\ dvol_{g_k}}=o(1)\ .
 \end{array}
 \]
 Assume that for every $k$ and every $\vec{p}$ satisfying  $\r_1(\vec{p})>2^N\,N\, t$ either  $\vec{p}\in F_k$ or $\vec{p}\in V_2({\R}^4)\setminus \vec{\La}_k(B_\rho(x_0))$ then
 \[
 \ti{\mu}_k\lf(\{\vec{p}\ ;\ \r_1(\vec{p})>2^N\,N\, t\}\rg)=\ti{\mu}_k(F_k)\ \longrightarrow 0\ ,
 \]
 from which we deduce
 \be
 \label{loc}
 \ti{\mu}_\infty(\{\vec{p}\ ;\ \r_1(\vec{p})>2^N\, t\})=0\ .
 \ee
 Alternatively, assume $ \vec{\La}_k(B_\rho(x_0)) \cap \{\vec{p}\ ;\ \r_1(\vec{p})>2^N\,N\, t\}\cap (V_2({\R}^4)\setminus F_k)\ne\emptyset$, for $\vec{p}_k$ in this intersection we can apply lemma~\ref{lm-good-points} for $r=t$  and any $s<t$ which gives
 \be
\label{good-2-bb}
\begin{array}{l}
\ds \lf[c_0-\lf[\delta\ \log\lf(\frac{t}{\ep}\rg)+\delta^{1/4} +\delta^{3/4}\ \log^{3/4}\lf(\frac{t}{\ep}\rg)\rg]^{1/2}\rg]\ \frac{1}{s^2}\int_{B_\rho(x_0)\,\cap\, \{\r_{\vec{p}}<s\}} dvol_{g_k}\\[5mm]
\ds\quad\quad\le \frac{C}{t^2}\int_{B_\rho(x_0)\,\cap\, \{\r_{\vec{p}}<t\}}  \ dvol_{g_k}+C\ e^{-\ep_k^{-2}}\ \ep_k^{-3}\ \ .
\end{array}
\ee
 Since $\vec{p}_k\in\vec{\La}_k(B_t(x_0))\setminus  \vec{\La}_k(\p B_t(x_0))$ we have using (\ref{density})
 \[
\lim_{s\rightarrow 0} \frac{1}{s^2}\int_{B_\rho(x_0)\,\cap\, \{\r_{\vec{p}_k}<s\}} dvol_{g_k}= \pi\,N_{\vec{p}_k}\ge \pi\ .
 \]
 This gives
 \be
 \label{good-3-bb}
 \begin{array}{l}
\ds \pi\,\lf[c_0-\lf[\delta\ \log\lf(\frac{t}{\ep}\rg)+\delta^{1/4} +\delta^{3/4}\ \log^{3/4}\lf(\frac{t}{\ep}\rg)\rg]^{1/2}\rg]\\[5mm]
\ds\quad\quad\quad\le \frac{C}{t^2}\int_{B_\rho(x_0)\,\cap\, \{\r_{\vec{p}_k}<t\}}  \ dvol_{g_k}+C\ e^{-\ep_k^{-2}}\ \ep_k^{-3}\ .
\end{array}
 \ee
 We deduce
 \be
 \label{good-4-bb}
  \pi\,c_0\le\liminf_{k\rightarrow +\infty} \frac{C}{t^2}\int_{B_\rho(x_0)\,\cap\, \{\r_{\vec{p}_k}<t\}}  \ dvol_{g_k}\ .
 \ee
 Assume that for $k$ large enough there exists $ \vec{p}_k\in\vec{\La}_k(B_\rho(x_0)) \cap \{\vec{p}\ ;\ \r_1(\vec{p})>2^N\, t\}\cap (V_2({\R}^4)\setminus F_k)$. If $r^\ast_{\vec{p}_k}>1$ then, since in the domain $0<r_{\vec{p}_k}<1/2$  we have $\zeta\circ \r_1\equiv 1$ we can apply (\ref{good-2-b}) as such and we get
\be
\label{good-2-b-bis}
\begin{array}{l}
\ds\pi \lf[c_0-\lf[\delta\ \log\lf(\frac{r}{\ep}\rg)+\delta^{1/4} +\delta^{3/4}\ \log^{3/4}\lf(\frac{r}{\ep}\rg)\rg]^{1/2}\rg]\ \\[5mm]
\ds\quad\quad\le{C_\ast}\int_{B_{\rho}(x_0)\,\cap\, \{\r_{\vec{p}_k}<1/2\}}  \ dvol_{g_k}+C\ e^{-\ep_k^{-2}}\ \ep_k^{-3}\ .
\end{array}
\ee
and passing to the limit we get
\be
\label{fff}
\nu_\infty(B_{\rho}(x_0))\ge \pi\,\frac{c_0}{C_\ast}
\ee
 and by choosing $0<c_Q<\frac{c_0}{C_\ast}$ we have proved the lemma in this case. 
 
 \medskip
 
 We now fix $N\in {\N}^\ast$ such that
  \be
  \label{N}
 C \ \frac{1}{\sqrt{N}}\,\,{\mu}_\infty(\{\r_1<1)\le \frac{\pi}{4}\,c_0\ ,
 \ee
 and assume now $0<2^N\,t<r^\ast_{\vec{p}_k}<1$. We extract a subsequence that we keep denoting $\vec{p}_k$ such that $\vec{p}_k\rightarrow \vec{p}_\infty$. We apply  (\ref{bord-9}) for $\vec{p}=\vec{p}_\infty$ that we combine
 with the the proof of lemma~\ref{lm-good-points} but with an error term coming from $\p \vec{\La}_k(B_t(x_0))$ in the r.h.s. of (\ref{mono-bord}) which is controlled by (\ref{bord-9}) we get
 \be
 \label{good-5-bb}
 \begin{array}{l}
 \ds\limsup_{k\rightarrow +\infty} \frac{C}{t^2}\int_{B_\rho(x_0)\,\cap\, \{\r_{\vec{p}_\infty}<2\,t\}}  \ dvol_{g_k}
 \ds- C\, 2^{N}\, t\ {\mu}_\infty(\{\r_1<1\})-C \ C_0\, 2^{N}\, t\ (\r^\ast_{\vec{p}_\infty})^{-1} \ {\mu}_\infty(\{\r_1<1\})\\[5mm]
 \ds- C \ \frac{1}{\sqrt{N}}\,{\mu}_\infty(\{\r_1<1)\le 
 C_\ast\, {\mu}_\infty(\r< 1)\ .
 \end{array}
 \ee
From (\ref{good-4-bb}) we have
\be
\label{good-5-xxx}
 \liminf_{k\rightarrow +\infty}\frac{C}{t^2}\int_{B_\rho(x_0)\,\cap\, \{\r_{\vec{p}_\infty}<2\,t\}}  \ dvol_{g_k}\ge \pi\,c_0
\ee
Combining (\ref{good-5-bb}) and (\ref{good-5-xxx}) is giving finally
\be
\label{good-6-xxx}
\begin{array}{l}
\ds\pi\, c_0- C\, 2^{N}\, t\ {\mu}_\infty(\{\r_1<1\})-C \ C_0\, 2^{N}\, t\ (\r^\ast_{\vec{p}_\infty})^{-1} \ {\mu}_\infty(\{\r_1<1\})\\[5mm]
 \ds- C \ \frac{1}{\sqrt{N}}\,{\mu}_\infty(\{\r_1<1)\le 
 C_\ast\, \ti{\mu}_\infty(V_2({\R}^4)=C_\ast  \nu_\infty(B_\rho(x_0))\ .
\end{array}
\ee
 We restrict to $0<t$ such that
 \be
 \label{N-1}
C\, 2^{N}\, t\ {\mu}_\infty(\{\r_1<1\})\le \frac{\pi}{4}\,c_0\ ,
 \ee
and we choose $\r^\ast_{\vec{p}_\infty}$ large enough compare to $t$ so that
\be
\label{N-2}
C \ C_0\, 2^{N}\, t\ (\r^\ast_{\vec{p}_\infty})^{-1}\ {\mu}_\infty(\{\r_1<1\})\le  \frac{\pi}{4}\,c_0\ .
\ee
Combining (\ref{N})  assuming $t:= \limsup_{k'\rightarrow +\infty}\sup_{\{\vec{p}\,,\,\vec{q}\,\in \ \vec{\La}_{k'}(\p B_\rho(x_0))\}}\r_{\vec{p}}(\vec{q}\,)$ satisfies (\ref{N-1}), assume one can find a subsequence of points 
$\vec{p}_k\in\vec{\La}_k(B_\rho(x_0)) \cap \{\vec{p}\ ;\ \r_1(\vec{p})>2^N\, t\}\cap (V_2({\R}^4)\setminus F_k)$ converging to a point $\vec{p}_\infty$ such that (\ref{N-2}) is satisfied, then we obtain
\be
\label{N-3}
 \frac{\pi}{4}\,c_0\le C_\ast  \nu_\infty(B_\rho(x_0))\ .
\ee
This concludes the proof of  lemma~\ref{lm-energy-quant}.\hfill $\Box$

\medskip
\begin{Lm}
\label{lm-abs-cont}
The limiting measure $\nu_\infty$ decomposes as follows
\be
\label{dec-1}
\nu_\infty= f\ dvol_h+\sum_{i=1}^Q c_i\ \delta_{x_i}\ ,
\ee
where $f\in L^1(\Sigma,{\R}_+)$ and $c_i\ge c_Q$ (where $c_Q$ is the constant given in lemma~\ref{lm-energy-quant}). Moreover the weak $W^{1,2}-$limit $\vec{\La}_\infty$ of the sequence $\vec{\La}_k$ is in $C^0_{loc}(\Sigma\setminus \{x_1\cdots x_Q\},V_2({\R}^4))$. Finally for any open set $U$ in $\Sigma$  with $x_i\notin\ov{U}$ for any $i=1\cdots Q$ there holds 
\be
\label{limit-meas-1}
\lim_{k\rightarrow +\infty} (\vec{\La}_k)_\ast(\nu_k\res U)\ =(\vec{\La}_\infty)_\ast(\nu_\infty\res U)\ .
\ee
\end{Lm}
\noindent{\bf Proof of lemma~\ref{lm-abs-cont}} Let $x\in\Sigma$ and let $r>0$. Using Fubini combined with the mean value argument, there exists $\rho\in(r/2,r)$ and a subsequence that we keep denoting $k$ such that
\[
\int_{\p B_\rho(x)}|d\vec{\La}_k|_h^2\ dl_h\le \frac{C}{r}\,\int_{B_r(x)}|d\vec{\La}_k|_h^2\ dvol_h\le 2\,\frac{C}{r}\,\mu_k(V_2({\R}^4))\ .
 \]
 and simultaneously
 \[
 \int_{\p B_\rho(x)}|d\vec{\La}_\infty|_h^2\ dl_h\le \frac{C}{r}\,\int_{B_r(x)}|d\vec{\La}_\infty|_h^2\ dvol_h\le 2\,\frac{C}{r}\,\mu_\infty(V_2({\R}^4))\ .
 \]
 Hence the restriction of $\vec{\La}_k$ to $\p B_\rho(x)$ is sequentially weakly pre-compact in $W^{1,2}(\p B_\rho(x), V_2({\R}^4))$ as well as in $C^{0,1/2}(\p B_\rho(x),V_2({\R}^4))$. By Arzela-Ascoli we deduce that $\vec{\La}_k$ is strongly pre-compact in $C^{0,\al}(\p B_\rho(x),V_2({\R}^4))$ for any $\al<1/2$. Since $\vec{\La}_k$ converges weakly towards $\vec{\La}_\infty$  in $W^{1,2}(B_\rho(x),V_2({\R}^4))$, by continuity of the trace operation into $W^{1/2,2}(\p B_\rho,V_2({\R}^4))$ we deduce
\be
\label{holder-cv}
\forall\al<1/2\quad\quad\vec{\La}_k\ \longrightarrow\ \vec{\La}_\infty\quad\mbox{ strongly in }C^{0,\al}(\p B_\rho(x),V_2({\R}^4))\ .
\ee 
and moreover
\be
\label{nu-bound-1}
\begin{array}{l}
\ds\int_{\p B_\rho(x)}|d\vec{\La}_\infty|_h^2\ dl_h\le \liminf_{k\rightarrow +\infty}\int_{\p B_\rho(x)}|d\vec{\La}_k|_h^2\ dl_h\le  \frac{2}{r}\,\liminf_{k\rightarrow +\infty}\int_{B_r(x)}|d\vec{\La}_k|_h^2\ dvol_h\\[5mm]
\ds\quad\le\frac{4}{r}\,\nu_\infty(\ov{B_r(x)})\ .
\end{array}
\ee
Cauchy Schwartz inequality gives
\be
\label{nu-bound-2}
\int_{\p B_\rho(x)}|d\vec{\La}_k|_h\ dl_h\le C\, \rho^{1/2}\, \lf[ \int_{\p B_\rho(x)}|d\vec{\La}_k|_h^2\ dl_h  \rg]^{1/2}\le C\, \sqrt{\frac{\rho}{r}}\ \sqrt{\int_{B_r(x)}|d\vec{\La}_k|_h^2\ dvol_h}\ .
\ee
Hence
\[
\limsup_{k\rightarrow +\infty}\int_{\p B_\rho(x)}|d\vec{\La}_k|_h\ dl_h<+\infty\ .
\]
Let
\[
t:=\lim_{k\rightarrow +\infty}\sup_{\{\vec{p}\,,\,\vec{q}\,\in \ \vec{\La}_{k'}(\p B_\rho(x_0))\}}\r_{\vec{p}}(\vec{q}\,)=\sup_{\{\vec{p}\,,\,\vec{q}\,\in \ \vec{\La}_{\infty}(\p B_\rho(x_0))\}}\r_{\vec{p}}(\vec{q}\,)
\]
Applying lemma~\ref{lm-energy-quant} we deduce that
\be
\label{altern}
\mbox{either }\quad\nu_\infty(B_\rho(x))>c_Q\quad\mbox{ or }\quad\nu_\infty(B_\rho(x))\le C\, t^2\, V_2({\R}^4)\ .
\ee 
Assuming
\[
\nu_\infty(\{x\})>0
\]
the second alternative cannot happen because, taking $r$ small enough we can make $t$ as small as we want. Indeed
\be
\label{schwar}
t\le C\, \int_{\p B_\rho(x)}|d\vec{\La}_\infty|_h\ dl_h\le C\, \sqrt{\frac{\rho}{r}}\ \sqrt{\int_{B_r(x)}|d\vec{\La}_\infty|_h^2\ dvol_h}\quad\longrightarrow\ 0\quad\mbox{ as }r\rightarrow 0\ .
\ee
Hence we have proved the following quantization phenomenon
\be
\label{quanti}
\nu_\infty(\{x\})>0\quad\Longrightarrow\quad \nu_\infty(\{x\})>c_Q\ ,
\ee
and this can only happen at at most finitely many points. We denote by $x_1\cdots x_Q$ the set of atoms (when it is not empty). Let $K$ be a compact set containing no atom for $\nu_\infty$ we claim that there exists $C$ independent of $K$ such that
\be
\label{lebesgue}
\nu_\infty(K)\le C\,\int_K|d\vec{\La}_\infty|^2_h\ dvol_h\ .
\ee
This claim is obviously implying the decomposition (\ref{dec-1}) of the measure $\nu_\infty$.
\medskip

\noindent{\bf Proof of (\ref{lebesgue}).} Since there are at most finitely many atoms, there exists $r_0>0$ such that $K\cap \ov{B_{r_0}(x_i)}=\emptyset$ for any $i=1\cdots Q$. Let $U$ be the open set
given by
\[
U:=\Sigma\setminus \bigcup_{i=1}^Q \ov{B_{r_0}(x_i)}\ ,
\]
and we choose $V$ open, arbitrary, such that $K\subset V$ and $\ov{V}\subset U$. Let $r>0$ such that $B_{5r}(y)\subset V$ for any $y\in K$ and 
\[
\sup_{y\in K}\nu_\infty(B_{5r}(y))<c_Q\ .
\]
We take a maximal finite family of balls $(B_{r}(y_j))_{j\in J}$ for $y_j\in K$ and the distance between two centers $y_j\ne y_j'$ is at least $r$ . Since it is a maximal family we have
\[
K\subset \bigcup_{j\in J} B_{r}(y_j)\quad\mbox{ and }\quad\sum_{j\in J}{\mathbf 1}_{B_{5r}(y_j)}\le C
\] 
where $C>0$ is universal and ${\mathbf 1}_{B_{5r}(y_j)}$ denotes the characteristic function of the ball $B_{5r}(y_j)$. For each $j\in J$ we chose $\rho_j\in [2\,r,4\,r]$ and we extract a subsequence that we keep denoting $k$ (recall that the cardinal of $J$ is finite) such that $\forall j\in J$ there holds
\[
\int_{\p B_{\rho_j}(y_j)}|d\vec{\La}_k|_h^2\ dl_h\le \frac{C}{r}\,\int_{B_{4\,r}(y_j)}|d\vec{\La}_k|_h^2\ dvol_h\le 2\,\frac{C}{r}\,\mu_k(V_2({\R}^4))\ ,
 \]
and simultaneously
\[
 \int_{\p B_{\rho_j}(y_j)}|d\vec{\La}_\infty|_h^2\ dl_h\le \frac{C}{r}\,\int_{B_{4\,r}(y_j)}|d\vec{\La}_\infty|_h^2\ dvol_h\le 2\,\frac{C}{r}\,\mu_\infty(V_2({\R}^4))\ .
\]
Hence, as for the proof of (\ref{holder-cv}), we have 
\[
\vec{\La}_k\longrightarrow \vec{\La}_\infty\quad\mbox{ in }C^0\lf(\bigcup_{j\in J}\p  B_{\rho_j}(y_j), V_2({\R}^4)\rg)
\]
and because of (\ref{schwar}) we have
\[
\begin{array}{rl}
\ds t_j&\ds:= \lim_{k\rightarrow +\infty}\sup_{\{\vec{p}\,,\,\vec{q}\,\in \ \vec{\La}_{k}(\p B_{\rho_j}(y_j))\}}\r_{\vec{p}}(\vec{q}\,)=\sup_{\{\vec{p}\,,\,\vec{q}\,\in \ \vec{\La}_{\infty}(\p B_{\rho_j}y_j))\}}\r_{\vec{p}}(\vec{q}\,)\\[5mm]
\ds&\ds\le C\ \sqrt{\int_{B_{4r}(y_j)}|d\vec{\La}_\infty|_h^2\ dvol_h}
\end{array}
\]
Because of lemma~\ref{lm-energy-quant}, since for every $j\in J$ we have $\nu_\infty(B_{\rho_j}(y_j))<c_Q$,  using more precisely (\ref{concl-3}) we deduce
\[
\begin{array}{l}
\ds\nu_\infty(K)\le \sum_{j\in J}\nu_\infty(B_{\rho_j}(y_j))\le C\sum_{j\in J}t_j^2\ \mu_\infty(V_2({\R}^4))\,\le C\,\mu_\infty(V_2({\R}^4))\,\sum_{j\in J}\int_{B_{4r}(y_j)}|d\vec{\La}_\infty|_h^2\ dvol_h\\[5mm]
\ds\le C\ \mu_\infty(V_2({\R}^4))\,\int_\Sigma \sum_{j\in J}{\mathbf 1}_{B_{4r}(y_j)}\ |d\vec{\La}_\infty|_h^2\ dvol_h\le \ C\,\mu_\infty(V_2({\R}^4))\,\int_V |d\vec{\La}_\infty|_h^2\ dvol_h
\end{array}
\]
Since this holds for any open set $V$ containing $K$ we deduce (\ref{lebesgue}) and the claim is proved. We deduce in particular that 
\[
\int_K |d\vec{\La}_\infty|_h^2\ dvol_h=0\quad \Longrightarrow\quad \nu_\infty(K)=0\ .
\]
and we deduce that away from the atoms the measure $\nu_\infty$ is absolutely continuous with respect to $dvol_h$. Again from lemma~\ref{lm-energy-quant} we have
\be
\label{concl-22}
\limsup_{k\rightarrow+\infty}\int_{ B_{\rho_j}(y_j)\cap \r^j_{\vec{\La}_{k}(x)}> c_2\, t_j}|d\vec{\La}_{k}|^2_h\ dvol_h=0\ ,
\ee
where 
\[
\r^j_{\vec{p}}:=\inf\lf\{\r_{\vec{p}}(\vec{q})\quad;\quad \vec{q}\in \vec{\La}_{k}(\p B_{\rho_j}(y_j)) \rg \}\ .
\]
Let $\chi_j(\vec{p}):=(\r^j_{\vec{p}}-c_2\,t_j)^+$ we have that the restriction of $\chi_j\circ \vec{\La}_k$ to $B_{\rho_j}(y_j)$ is strongly converging to zero in $W^{1,2}$ and we deduce that for all $j\in J$
\[
\forall \,\vec{p}\in \vec{\La}_\infty(B_{\rho_j}(y_j))\quad\mbox{and }\quad\forall \,\vec{q}\in \vec{\La}_{\infty}(\p B_{\rho_j}(y_j))\quad \r_{\vec{p}}(\vec{q})< C\, (c_2+1)\, t_j
\]
where $C$ denotes the constant in the triangular inequality of the quasi distance $\r$. Thus for $x$ and $y$ in $K$ such that $|x-y|<r$ we have that $x$ and $y$ belong to two balls $B_{2\,r}(y_j)$ which intersect each-other
and using again the fact that $\r$ is a quasi-distance we deduce that
\[
\r_{\vec{\La}_\infty(x)}(\vec{\La}_\infty(y))< C\, \max_{j\in J}\{t_j\}\le  C\ \max_{j\in J}\sqrt{\int_{B_{4r}(y_j)}|d\vec{\La}_\infty|_h^2\ dvol_h}
\]
Let now $\ep>0$, we choose $r>0$ such that
\[
\max_{y\in K} \sqrt{\int_{B_{4r}(y)}|d\vec{\La}_\infty|_h^2\ dvol_h}<\ep
\]
we have that
\[
\forall\, x,y\in K\quad|x-y|<r\ \Longrightarrow \ \r_{\vec{\La}_\infty(x)}(\vec{\La}_\infty(y))<\ep\ .
\]
This implies that $\vec{\La}_\infty$ is continuous on $K$. 

\medskip

Finally we have
\[
\begin{array}{l}
\ds\int_K|\vec{\La}_k-\vec{\La}_\infty|\ d\nu_k\le \sum_{j\in J}\int_{ B_{\rho_j}(y_j)}|\vec{\La}_k-\vec{\La}_\infty|\ d\nu_k\le C\ \sum_{j\in J} t_j\ \int_{ B_{\rho_j}(y_j)\cap \r^j_{\vec{\La}_{k}(x)}\le c_2\, t_j}|d\vec{\La}_{k}|^2_h\ dvol_h\\[5mm]
\ds\ +\sum_{j\in J} \int_{ B_{\rho_j}(y_j)\cap \r^j_{\vec{\La}_{k}(x)}>c_2\, t_j}|d\vec{\La}_{k}|^2_h\ dvol_h
\end{array}
\]
From (\ref{concl-22}) we deduce
\[
\label{concl-222-bis}
\begin{array}{l}
\ds\limsup_{k\rightarrow+\infty}\int_K|\vec{\La}_k-\vec{\La}_\infty|\ d\nu_k\le C\, \max_{j\in J}t_j\ \mu_\infty(V_2({\R}^4))\\[5mm]
\ds\le   C\ \mu_\infty(V_2({\R}^4))\ \max_{j\in J}\sqrt{\int_{B_{4r}(y_j)}|d\vec{\La}_\infty|_h^2\ dvol_h}\le   C\ \mu_\infty(V_2({\R}^4))\ \max_{y\in K}\sqrt{\int_{B_{4r}(y)}|d\vec{\La}_\infty|_h^2\ dvol_h}
\end{array}
\]
By making $r$ go to zero implies
\be
\label{concl-222}
\lim_{k\rightarrow+\infty}\int_K|\vec{\La}_k-\vec{\La}_\infty|\ d\nu_k=0\ ,
\ee
For any $\Psi\in C^0(V_2({\R}^4))$, because of (\ref{concl-222}) and since $\vec{\La}_\infty$ is continuous on $K$, the weak convergence in Radon measure is implying
\be
\label{concl223}
\lim_{k\rightarrow+\infty}\int_K\Psi(\vec{\La}_k)\ d\nu_k=\lim_{k\rightarrow+\infty}\int_K\Psi(\vec{\La}_\infty)\ d\nu_k=\int_K\Psi(\vec{\La}_\infty)\ d\nu_\infty
\ee
This holds for any compact not containing the atoms of $\nu_\infty$ and this is closing the proof of lemma~\ref{lm-abs-cont} .\hfill $\Box$

\medskip
\begin{Lm}
\label{lm-integer-mult}
There exists
$N\in L^\infty(\Sigma,{\N}^\ast)$ such that
\[
\nu_\infty= \, N\,  \frac{|d\vec{\La}_\infty\,\dot{\wedge} \,d\vec{\La}_\infty|_h}{2}\ dvol_h+\sum_{i=1}^Q c_i\ \delta_{x_i}\ .
\]
\end{Lm}
\noindent{\bf Proof of lemma~\ref{lm-integer-mult}.}
We take local conformal coordinates $(x_1,x_2)$ for $h$ around a point which is $0=(0,0)$ in these coordinates and we can write in these coordinates 
\[
h=e^{2\la}\ [dx_1^2+dx_2^2]\ .
\]
We can assume again that $\vec{\La}_\infty(0)=(\vec{\ep}_1,\vec{\ep}_2)$. We assume that $0$ is not one of the atoms of $\nu_\infty$, we assume moreover that $0$ is a Lebesgue point for $f$ in such a way that
\be
\label{exist-dens}
\lim_{r\rightarrow 0}\frac{\nu_\infty(B_r(0))}{\pi\,r^2}=f(0)\, e^{2\la(0)}\ .
\ee
This is satisfied ${\mathcal H}^2-a.e.$ on $\Sigma$. The goal in this lemma is to prove that there exists a non zero integer $N$ such that
\be
\label{but}
f(0)=N\, \lf|\p_{x_1}\vec{\La}_\infty(0)\wedge\p_{x_2}\vec{\La}_\infty(0)\rg|\ e^{-2\la(0)}\ .
\ee
 where we recall the notation $\nu_\infty=f\ dvol_h+\sum_{i=1}^Qc_i\ \delta_{x_i}$.  We can also assume that
 \be
 \label{pt-cv}
 \lim_{k\rightarrow+\infty}\vec{\La}_k(0)=\vec{\La}_\infty(0)\ .
\ee
 We can assume without loss of generality that $0$ is a Lebesgue point for $d\vec{\La}_\infty$
and that $\vec{\La}_\infty$  is approximate differentiable at $0$. This holds for ${\mathcal H}^2-$a.e. points in $\Sigma$ as well since $\vec{\La}_\infty\in W^{1,2}$ (see \cite{EG}  section 6.1). Hence we have
\be
\label{lebesgue-pt}
\lim_{r\rightarrow 0}\int_{B_r(0)}\frac{|d\vec{\La}_\infty(x)-d\vec{\La}_\infty(0)|^2}{r^2}\ dx^2=0\ .
\ee
as well as
\be
\label{app-diff}
\lim_{r\rightarrow 0}\int_{B_r(0)}\frac{\lf|\vec{\La}_\infty(x)-\vec{\La}_\infty(0)-d\vec{\La}_\infty(0)\cdot x\rg|}{r^3}\ dx^2=0\ ,
\ee
and
\be
\label{app-diff-bis}
\lim_{r\rightarrow 0}\int_{B_r(0)}\frac{\lf|\vec{\La}_\infty(x)-\vec{\La}_\infty(0)-d\vec{\La}_\infty(0)\cdot x\rg|^2}{r^4}\ dx^2=0\ .
\ee
Since $0$ is chosen to be a point of approximate differentiability of $\vec{\La}_\infty$
Hence in particular for any $t>0$ there exists $r\in [t/2,t]$ and $\theta_r\in[0,2\pi]$ such that
\[
\lim_{r\rightarrow 0}\int_0^{2\pi} \lf|d\vec{\La}_\infty(r,\theta)-d\vec{\La}_\infty(0)\rg|\ d\theta=0
\]
and
\[
\lim_{r\rightarrow 0}r^{-1}\ \int_0^{2\pi} \lf|\vec{\La}_\infty(r,\theta)-\vec{\La}_\infty(0)-r\, \cos\theta\,\p_{x_1}\vec{\La}_\infty(0)-r\, \sin\theta\,\p_{x_2}\vec{\La}_\infty(0)\rg|\ d\theta=0\ .
\]
By the mean value theorem there exists $\theta_r\in[0,2\pi]$ such that
\[
\lf|\vec{\La}_{\infty}(r,\theta_r)-\vec{\La}_\infty(0)-r\, \cos\theta_r\,\p_{x_1}\vec{\La}_\infty(0)-r\, \sin\theta_r\,\p_{x_2}\vec{\La}_\infty(0)\rg|=o(r)\ .
\]
We have
\[
\begin{array}{rl}
\ds\forall\,\theta\in [0,2\pi]\ \quad\vec{\La}_\infty(r,\theta)-\vec{\La}_\infty(r,\theta_r)&\ds=r\,\int_{\theta_r}^\theta -\sin\,{\phi}\ \partial_{x_1}\vec{\La}_\infty(r,\phi)+\cos\,{\phi}\ \partial_{x_2}\vec{\La}_\infty(r,\phi)\ d\phi\\[5mm]
\ds\quad&\ds=r\,\int_{\theta_r}^\theta -\sin\,{\phi}\ \partial_{x_1}\vec{\La}_\infty(0)+\cos\,{\phi}\ \partial_{x_2}\vec{\La}_\infty(0)\ d\phi +o(r)\ .
\end{array}
\]
Combining the previous we obtain
\be
\label{bord-cont}
\lf\|\vec{\La}_\infty(r,\theta)-\vec{\La}_\infty(0)-r\, \cos\theta\,\p_{x_1}\vec{\La}_\infty(0)-r\, \sin\theta\,\p_{x_2}\vec{\La}_\infty(0)\rg\|_{L^\infty([0,2\pi])}=o(r)\ .
\ee
As in the proof of the previous lemma we can also choose $r$ such that  there exists a subsequence (we keep denoting $k$) satisfying
\be
\label{bord-conv}
\int_{\p B_r(0)}|d\vec{\La}_k|^2\ dl\le r^{-1}\ \int_{B_{4\,r}(0)\setminus B_{r/4}(0)}|\nabla \vec{\La}_k|^2\ dx^2
\ee
The radii such that (\ref{bord-cont}) and (\ref{bord-conv}) hold are called good radii.
\[
\|\vec{\La}_\infty(y)-\vec{\La}_\infty(x)\|_{L^\infty(\p B_r(0))\times \p B_r(0))}\le \pi\, |d\vec{\La}_\infty(0)|\ r +o(r)
\]
we deduce from (\ref{loc}) that there exists a constant $C>0$ independent of $r$ such that
 \be
 \label{loc-bis}
 \ti{\mu}_\infty\lf(\lf\{\vec{p}\ ;\ \inf\lf\{\r_{\vec{q}}(\vec{p})\ ; \vec{q}\in \p\vec{\La}_\infty(\p B_r(0))\rg\}>C |d\vec{\La}_\infty(0)|\ r\rg\}\rg)=0\ .
 \ee
 where we adopt again the notation
 \[
 \ti{\mu}_\infty:=\lim_{k\rightarrow +\infty}\ \lf(\lf.\vec{\La}_k\rg|B_r(0)\rg)_\ast\nu_k\ .
 \]
 Since $0$ is a Lebesgue point for $f$ we have
 \[
 \begin{array}{l}
\ds \nu_\infty(B_{4^{-4\,l}}(0))=\int_{B_{4^{-4\,l}}(0)}f(x)\ e^{2\la}\ dx^2\\[5mm]
\ds\quad=\int_{B_{4^{-4\,l}}(0)}f(0)\ e^{2\la}\ dx^2+o(4^{-4\,l})
\end{array}
 \]
 (recall $\nu_\infty(\{0\})=0$) we can choose $t_j=4^{-4\,j}\rightarrow 0$ and from the previous fact we have
 \be
 \label{choix-tj}
 \nu_\infty\lf(B_{4t_j}(0)\rg)= O(t_j^2)
 \ee
 Finally Taking  $\delta_k\rightarrow 0$ given by
 \be
 \label{re-bord-10}
\delta_k^2:=(\log\ep_k^{-1})^{-1}\, \ep^4_k\int_\Sigma(1+|d\vec{T}_k|^2)^2\ dvol_{g_k}
 \ee
Recall that this choice is made in such a way that
 \be
 \label{recondition-delta}
 \lim_{k\rightarrow +\infty}\delta_k\,\log\ep_k^{-1}=0\quad\mbox{ and }\quad \lim_{k\rightarrow +\infty}\delta^{-1}_k\,{\ep^4_k\int_\Sigma(1+|d\vec{T}_k|^2)^2\ dvol_{g_k}}=0\ .
 \ee
 Let $F_k$ be the set of $\vec{p}\in V_2({\R}^4)$ with  such that there exists $s>0$ with
 \[
 \delta_k\, {\mu}_k(\{\r_{\vec{p}}<5\, C_0^2\,s\})\le \ep^4_k\int_{B_\rho(x_0)\cap\{\r_{\vec{p}}<s\}}(1+|d\vec{T}_k|^2)^2\ dvol_{g_k}
 \]
 Observe that for $\vec{p}\ne\vec{q}$
 \[
\sup_{\vec{o}\in V_2({\R}^4)}\ \frac{\r_{\vec{p}}(\vec{q})}{\max\{\r_{\vec{p}}(\vec{o}),\r_{\vec{q}}(\vec{o})\}}\le 2\, C_0
 \]
 We can then use the Vitali type Covering Lemma in Quasi-metric spaces given by Lemma 2.7 in \cite{AlMi} and cover $F_k$ by at most countably many balls $(B^\r_{5 C_0^2\,s_i}(\vec{p}_i))_{i\in I}$  for the quasi-distance $\r$
 such that the balls $B^\r_{5 C_0^2\,s_i}(\vec{p}_i)$ for $i\in I$ are disjoint to each other. Hence we have
 \[
 \begin{array}{l}
 \ds{\mu}_k(F_k)\le\sum_{i\in I}\ti{\mu}_k(\{\r_{\vec{p_i}}<5\, C_0^2\,s_i\})\le \delta_k^{-1}\ \sum_{i\in I}\ep^4_k\int_{B_\rho(x_0)\cap\{\r_{\vec{p}_i}<s_i\}}(1+|d\vec{T}_k|^2)^2\ dvol_{g_k}\\[5mm]
 \ds\quad\quad \quad\le \delta^{-1}_k\,{\ep^4_k\int_\Sigma(1+|d\vec{T}_k|^2)^2\ dvol_{g_k}}=o(1)\ .
 \end{array}
 \]
  We construct a sequence of ``good radii'' $r_j\in [t_j/2,t_j]$ a sub sequence $\phi_j(k)$ and an increasing sequence index $k_j$ as follows. Assuming $(r_l)_{l<j}$ as well as $(\phi_l(k))_{l<j}$ and $(k_l)_{l<j}$ have been constructed, we choose $r_j$ as explained above such that
 \be
\label{bord-cont-2}
\lf\|\vec{\La}_\infty(r_j,\theta)-\vec{\La}_\infty(0)-r_j\, \cos\theta\,\p_{x_1}\vec{\La}_\infty(0)-r_j\, \sin\theta\,\p_{x_2}\vec{\La}_\infty(0)\rg\|_{L^\infty([0,2\pi])}=o(r_j)\ .
\ee
and such that there exists a subsequence of $\phi_{j-1}(k)$  that we denote $\phi_{j}(k)$ such that
\be
\label{w12-bord-control}
\int_{\p B_{r_j}(0)}|d\vec{\La}_{\phi_{j}(k)}|^2\ dl\le r_j^{-1}\ \int_{B_{4\,r_j}(0)\setminus B_{r_j/4}(0)}|\nabla \vec{\La}_{\phi_{j}(k)}|^2\ dx^2\ .
\ee
 Combining ({concl-222}), (\ref{pt-cv}), (\ref{choix-tj}) with (\ref{w12-bord-control}), we can ``wait enough'' and pick $k_j$ large enough so that
 \begin{itemize}
 \item
 \be
 \label{item-0}
 |\vec{\La}_{\phi_{j}(k_j)}(0)-\vec{\La}_\infty(0)|=o(r_j)\ ,
 \ee
 \item
 \be
 \label{item-1}
\lf| \nu_{\phi_{j}(k_j)}(B_{r_j}(0))- \nu_{\infty}(B_{r_j}(0))\rg|=o(r_j^2)\ ,
 \ee
 \item
 \be
 \label{item-2}
 \int_{\p B_{r_j}(0)}|d\vec{\La}_{\phi_{j}(k_j)}|^2\ dl\le  r_j^{-1}\  O(t^2_j) =O(r_j)\ ,
 \ee
 \item
 \be
 \label{item-2b}
 \int_{B_{r_j}}|\vec{\La}_{\phi_{j}(k_j)}-\vec{\La}_\infty|^2\ dx^2=o(r_j^4)\ ,
 \ee
 \item
 \be
 \label{item-3}
 \int_{{B_{r_j}(0)}}|\vec{\La}_{\phi_{j}(k_j)}-\vec{\La}_\infty|\ d\nu_{\phi_{j}(k_j)}=o(r^3_j)\ ,
 \ee
  \item
 \be
 \label{item-3-bis}
 \int_{{B_{r_j}(0)}}|\vec{\La}_{\phi_{j}(k_j)}-\vec{\La}_\infty|^2\ d\nu_{\phi_{j}(k_j)}=o(r^4_j)\ ,
 \ee
 \item
 \be
 \label{item-4}
 r_{j}^{-1}\,\sqrt{\ep_{\phi_{j}(k_j)}}\longrightarrow 0\ ,
 \ee
 \item
 \be
 \label{item-5}
\frac{\ep^4_{\phi_{j}(k_j)}}{r_j^2}\int_{B_{r_j}(0)}(1+|d\vec{T}_{\phi_{j}(k_j)}|^2_{g_{\phi_{j}(k_j)}})^2 \ dvol_{g_{\phi_{j}(k_j)}}=o\lf(\frac{1}{\log \ep^{-1}_{\phi_{j}(k_j)}}\rg)\ ,
 \ee
 \item
 \be
 \label{item-6}
 {\mu_{{\phi_{j}(k_j)}}\lf(F_{\phi_{j}(k_j)}\rg)}=o(r_j^2)\ .
 \ee
 \end{itemize}
Denoting $\ti{\mathbf v}_{k,r}$ the varifold generated by the restriction to $B_r(0)$ of $\vec{\La}_k$ we introduce the dilated varifold
 \[
 \hat{\mathbf v}_{k,r}:= (D_{r^{-1}})_\ast \ti{\mathbf v}_{k,r}
 \]
 Where 
 \[
 \begin{array}{rcl} 
\ds  D_{r^{-1}}\ : \quad\quad V_2({\R}^4)\ &\ \longrightarrow \ &\ r^{-1}(V_2({\R}^4)-\vec{\La}_\infty(0))\\[3mm]
 \ds (\vec{\ep}_1+\vec{v},\vec{\ep}_2+\vec{w})\ &\ \longrightarrow \ &\ (r^{-1}\vec{v},r^{-1}\vec{w})
 \end{array}
 \]
 Let $\vec{v}_r:=r^{-1}\vec{v}$, $\vec{w}_r:=r^{-1}\vec{w}$ and introduce $\varphi_r:=r^{-2}\, \varphi=r^{-2}(v^2-w^1)$. The inverse of $D_r$ is $(D_r)^{-1}$ such that 
 \[
 (D_{r^{-1}})^{-1}(\vec{v}_r,\vec{w}_r):=r\,(\vec{v}_r,\vec{w}_r)+(\vec{\ep}_1,\vec{\ep}_2)\ .
 \]
  we have that
  \[
  \begin{array}{l}
 \ds \lf((D_{r^{-1}})^{-1}\rg)^\ast \al= \lf((D_{r^{-1}})^{-1}\rg)^\ast\lf[(\vec{\ep}_1+\vec{v})\cdot d\vec{w}-(\vec{\ep}_2+\vec{w})\cdot d\vec{v}\rg]\\[5mm]
 \ds = r\ (dw^1_r-dv^2_r)+r^2 \,\lf[\vec{v}_r\cdot d\vec{w}_r-\vec{w}_r\cdot d\vec{v}_r\rg]= r^{2}\, \al_r
  \end{array}
  \]
  where
 \[
 \al_r:=-d\varphi_r+\vec{v}_r\cdot d\vec{w}_r-\vec{w}_r\cdot d\vec{v}_r
 \]
Observe that {\bf $ \hat{\mathbf v}_{k,r}$ is Legendrian for $\al_r$}. We denote moreover
 \[
 \r_r:=\lf[( \lf(|\vec{v}_r|^2+|\vec{w}_r|^2\rg)^2+4\,\varphi_r^2\rg)^{1/4}=r^{-1}\,\r
 \]
For the above constructed subsequence $\phi_{j}(k_j)$ we introduce respectively
\[
\hat{\La}_j(y):=r_j^{-1}\,\, (\vec{\La}_{\phi_{j}(k_j)}(r_j\,y)-\vec{\La}_{\infty}(0))\quad,\quad \hat{\La}_{\infty,j}(y):=r_j^{-1}\,\, (\vec{\La}_{\infty}(r_j\,y)-\vec{\La}_{\infty}(0))
\]
and
\[
\hat{\nu}_j:= \frac{1}{2}\,|\nabla_y \hat{\La}_j|^2\ dy^2=\frac{1}{2\, r_j^2}\,|\nabla_x {\La}_{\phi_{j}(k_j)}|^2\ dx^2\quad,\quad\hat{\mathbf v}_j:= \hat{\mathbf v}_{\phi_{j}(k_j),r_j}\quad\mbox{ and }\quad\al_j:=\al_{r_j}\ .
\]
We write
\[
\vec{v}_j:=\vec{v}_{r_j}=r_j^{-1}\,\vec{v}\quad,\quad   \vec{w}_j:=\vec{w}_{r_j}=r_j^{-1}\,\vec{w} \quad ,\quad\varphi_j:=\varphi_{r_j}:=r_j^{-2}\, \varphi\quad\mbox{ and }\quad \r_j:=r_j^{-1}\,\r \ .
\]
We see $\al_j$ as a one form in ${\R}^9$ given by
\[
\al_j:=-d\varphi_j+\vec{v}_j\cdot d\vec{w}_j-\vec{w}_j\cdot d\vec{v}_j\ .
\]
We also denote
\[
\hat{v}_j:=\vec{v}\circ \hat{\La}_j=r_j^{-1}\,\vec{v}\circ\vec{\La}_j\quad\mbox{ and }\quad\hat{w}_j:=\vec{w}\circ \hat{\La}_j=r_j^{-1}\,\vec{w}\circ\vec{\La}_j
\]
and naurally
\[
\hat{v}_{l,j}:=r_j^{-1}\,v_l\circ\vec{\La}_j\quad,\quad\hat{w}_{l,j}:=r_j^{-1}\,w_l\circ\vec{\La}_j\quad,\quad\hat{\varphi}_j:=r_j^{-2}\, \varphi\circ\vec{\La}_j\quad\mbox{ and }\quad \hat{\r}_j:=r_j^{-1}\,\r\circ\vec{\La}_j\ .
\]
We have then in particular $\hat{\La}_j=(\hat{v}_j,\hat{w}_j)$.
Because of (\ref{lebesgue-pt}), (\ref{app-diff}), (\ref{bord-cont-2}), (\ref{item-0})...(\ref{item-5}) we have respectively that
\be
\label{dila-conv}
\int_{B_1(0)} |\nabla\hat{\La}_{\infty,j}(y)-\nabla\vec{\La}_\infty(0)|^2\ dy^2\longrightarrow 0\quad,\quad\int_{B_1(0)}|\hat{\La}_{\infty,j}(y)-\nabla\vec{\La}_\infty(0)\cdot y|\ dy^2\longrightarrow 0\ .
\ee
and
\be
\label{conv-bord-ellipse}
\lf\|\hat{\La}_{\infty,j}(1,\theta)- \cos\theta\,\p_{x_1}\vec{\La}_\infty(0)- \sin\theta\,\p_{x_2}\vec{\La}_\infty(0)\rg\|_{L^\infty([0,2\pi])}=o(1)
\ee
and
\be
\label{blow-1}
\lim_{j\rightarrow+\infty}\int_{B_1(0)}|\hat{\La}_j-\hat{\La}_{j,\infty}|^2\ dy^2=0
\ee
as well as
\be
\label{blow-1-bis}
\lim_{j\rightarrow+\infty}\int_{B_1(0)}|\hat{\La}_j-\hat{\La}_{j,\infty}|^2\ d\hat{\nu}_j=0
\ee
Moreover the map $\hat{\La}_j$ is a critical point of 
\[
\hat{E}_j(\hat{\La}):=\int_{B_1(0)}\ dvol_{g_{\hat{\La}}}+\hat{\ep}^{\,4}_j\int_{B_1(0)}(r_j^2+|d\hat{T}|^2_{g_{\hat{\La}}})^2\ dvol_{g_{\hat{\La}}}
\]
for perturbations within $r_j^{-1}(V_2({\R}^4)-\vec{\La}_\infty(0))$ preserving the Legendrian condition relative to $\al_j$ supported away from $\hat{\La}_j(\p B_1(0))$ and satisfying
\be
\label{ener-asym}
\limsup_{j\rightarrow+\infty}\hat{E}_j(\hat{\La}_j)<+\infty\quad,\quad\hat{\ep}^{\,4}_j\,\int_{B_1(0)}(r_j^2+|d\hat{T}_j|^2_{g_{\hat{\La}_j}})^2\ dvol_{g_{\hat{\La}}}=o\lf(\frac{1}{\log \hat{\ep}^{-1}_{j}}\rg)\ 
\ee
where
\[
\hat{T}_j:=\frac{\p_{x_1}\hat{\La}_j\wedge\p_{x_2}\hat{\La}_j }{|\p_{x_1}\hat{\La}_j\wedge\p_{x_2}\hat{\La}_j|}\quad\mbox{ and }\quad      \hat{\ep}_j:=\frac{\ep_{\phi_{j}(k_j)}}{\sqrt{r_j}}\ .
\]
Observe that $r_j^{-1}(V_2({\R}^4)-\vec{\La}_\infty(0))$ equipped with $\al_j$ in the local chart given by $(\varphi_{j},v^3_{j},w^3_{j},v^4_{j},w^4_{j})$ is converging in any $C^l$ norm locally towards the Heisenberg group : ${\R}^5:=\{(\varphi_\infty,v^3_\infty,w^3_\infty,v^4_\infty,w^4_\infty)\}$
equipped with $\al_\infty:=-d\varphi_\infty+v^3_\infty\,dw^3_\infty-w^3_\infty\, dv^3_\infty+v^4_\infty\,dw^4_\infty-w^4_\infty\, dv^4_\infty$. We can extract a subsequence that we keep denoting $j$ such that
\[
\hat{\nu}_j\rightharpoonup\hat{\nu}_\infty \mbox{ weakly as Radon measure on $B_1(0)$}\quad\mbox{ and }\quad\hat{\mathbf v}_j\rightharpoonup \hat{\mathbf v}_\infty\mbox{ as varifold in ${\mathcal M}(G_2({\R}^8))$ }
\]
The previous lemma apply and we can then deduce that there exists $\hat{f}$ supported in $B_1(0)$ and finitely many points $y_j$ such that
\[
\hat{\nu}_\infty\res B_1(0)= \hat{f}\ dy^2+\sum_{j=1}^{\hat{Q}}c_j\ \delta_{y_j} 
\]
where the masses $c_j$ are bounded from below by a universal constant $c_Q>0$. By possibly having dilated a bit less in the image (i.e. instead of taking $\hat{\La}_j:=r_j^{-1}(\vec{\La}_k(r_j\,y)-\vec{\La}_\infty(0)$ one would have taken
$\hat{\La}_j:=\beta\,r_j^{-1}(\vec{\La}_k(r_j\,y)-\vec{\La}_\infty(0)$ for some $\beta<1$ we can assume that $\hat{\nu}_\infty(B_1(0))<c_Q$ and assume that there is no mass concentration. Hence $\hat{\nu}_\infty$ is absolutely continuous with respect to Legesgue. Let  $\hat{\La}_{\infty,\infty}(y)=\nabla\vec{\La}(0)\cdot y$. Because of (\ref{app-diff-bis}) and (\ref{blow-1}) there holds respectively 
\[
\int_{B_1(0)} |\hat{\La}_{j}(y)-\hat{\La}_{\infty,\infty}(y)|\ dy^2=o(1)\ .
\]
Hence we have
\be
\label{weak-conv-1}
\hat{\La}_{j}\ \rightharpoonup \ \hat{\La}_{\infty,\infty}\quad\quad\mbox{ weakly in }W^{1,2}(B_1(0))\ .
\ee
We can apply\footnote{Lemma~\ref{lm-abs-cont} is proved in a compact ambiant space $V_2({\R}^4$, while we are now considering an ``expanding space'' converging to the Heisenberg group ${\mathbb H}^2$. Nevertheless the proof of Lemma~\ref{lm-abs-cont} is the same except that one has to consider compactly supported test functions exclusively.} Lemma~\ref{lm-abs-cont} and in particular (\ref{concl223}) for compactly supported $\Psi$ to deduce
 that
\be
\label{conv-meas}
\lim_{j\rightarrow +\infty} (\hat{\La}_j)_\ast(\hat{\nu}_j\res B_1(0))=(\hat{\La}_{\infty,\infty})_\ast(\hat{\nu}_\infty\res B_1(0))
\ee
This implies that $|\hat{\mathbf v}_\infty|$ is supported in the (possibly degenerated) ellipse ${\mathcal E}:=\{\vec{p}:=\nabla\vec{\La}(0)\cdot y\ ;\ |y|<1\}$. This ellipse is itself contained in the Legendrian hyperplane given by $\varphi_\infty=0$.

\medskip

For every $\vec{p}\in r_j^{-1}([V_2({\R}^4)\setminus F_k]-\vec{\La}_\infty(0))$
and any $s>0$ we have from (\ref{good-2-b}) 
\be
\label{regood-2-b}
\begin{array}{l}
\ds \lf[c_0-\lf[\delta_k\ \log\lf(\frac{s\,r_j}{\ep_k}\rg)+\delta_k^{1/4} +\delta_k^{3/4}\ \log^{3/4}\lf(\frac{s\,r_j}{\ep_k}\rg)\rg]^{1/2}\rg]\ \frac{1}{s^2}\int_{\Sigma\,\cap\, \{\r_{\vec{p}}<s\, r_j\}} dvol_{g_k}\\[5mm]
\ds\quad\quad\le C \int_{\Sigma\,\cap\, \{\r_{\vec{p}}<1\}}  \ dvol_{g_k}+C\ e^{-\ep_k^{-2}}\ \ep_k^{-3}\le  C\, E_{\ep_k}\ .
\end{array}
\ee
where $\delta_k$ is given by (\ref{re-bord-10}) (we have simply written $k$ for $\phi_{j}(k_j)$) in the two lines above ). Because of (\ref{item-6}) we deduce from (\ref{regood-2-b}) that for $|\hat{\mathbf v}_\infty|$ a.e. $\vec{p}$ and any $s>0$ there holds
\[
|\hat{\mathbf v}_\infty|(B^\r_s(\vec{p}))\le C\ s^2\ .
\]
Assume first that the rank of $\nabla\vec{\La}_0$ is less or equal than one that is, the ellipse ${\mathcal E}$ is degenerated and is either a segment or a point. We cover the support of $|\hat{\mathbf v}_\infty|$ by at most $s^{-1}$ balls and we deduce
\[
|\hat{\mathbf v}_\infty|({\R}^5)\le C\, s
\]
This  holds for any $s$ we deduce
\[
(\hat{\La}_{\infty,\infty})_\ast(\hat{\nu}_\infty\res B_1(0))\le|\hat{\mathbf v}_\infty|({\R}^5)=0
\]
Hence $\hat{\nu}_\infty\res B_1(0)=0$ in that case which implies $f=0$ (thanks to (\ref{exist-dens}) and (\ref{item-1}) ) as stated in the lemma. 

\medskip

We work now under the assumption that  $\mbox{rank}(\nabla\vec{\La}(0))=2$ without loss of generality one can assume that 
\[
{\mathcal P}_\infty:=\mbox{Span}(\p_{x_1}\vec{\La}_\infty(0),\p_{x_2}\vec{\La}_\infty(0))=\mbox{Span}\lf((\vec{\ep}_3,0),(\vec{\ep}_4,0)\rg)
\]
where we recall that we are using the following notations for the coordinates in ${\R}^8$ in the neighbourhood of $(\vec{\ep}_1,\vec{\ep}_2)$ : 
\[
\vec{p}:= ((1+v^1)\ \vec{\ep}_1+v^2\ \vec{\ep}_2+v^3\ \vec{\ep}_3+v^4\ \vec{\ep}_4, w^1\ \vec{\ep}_1+(1+w^2)\ \vec{\ep}_2+w^3 \vec{\ep}_3+w^4\ \vec{\ep}_4)
\]
We claim that
\be
\label{rec-1}
\hat{\mathbf v}_\infty=\delta_{{\mathcal P}_\infty}\otimes|\hat{\mathbf v}_\infty|\ .
\ee
This is equivalent to
\be
\label{rec-2}
\begin{array}{l}
\ds\int_{B_1(0)}|\nabla \hat{v}_{j}^1|^2+|\nabla \hat{w}_{j}^1|^2+|\nabla \hat{v}_{j}^2|^2+|\nabla \hat{w}_{j}^2|^2+|\nabla \hat{w}_{j}^{3}|^2+|\nabla \hat{w}_{j}^{4}|^2+|\nabla\hat{\varphi}_j|^2 \ dy^2=o(1)\ .
\end{array}
\ee
We know from (\ref{weak-conv-1}) that $\hat{\La}_j$ is weakly converging towards  $\hat{\La}_{\infty,\infty}$ in $W^{1,2}(B_1(0))$. Because of (\ref{item-2}) we have
\[
\limsup_{j\rightarrow +\infty}\int_{\p B_1(0)}|d\hat{\La}_j|^2\ dl<+\infty
\]
Hence $\hat{\La}_j$ is strongly pre-compact in $C^{0,\mu}(\p B_1(0))$ for any $\mu<1/2$\be
\label{cv-bord}
\lim_{j\rightarrow +\infty}\|\hat{\La}_j-\hat{\La}_{\infty,\infty}\|_{L^\infty(\p B_1(0))}=0\ .
\ee
We have $\hat{\nu}_\infty(B_1(0))<c_Q$. Let $t:=\mbox{diam}(\hat{\La}_{\infty,\infty}(\p B_1(0)))$. Thanks to (\ref{concl-2}), there exists $c_2$
\be
\label{concl-2-reu}
\lim_{j\rightarrow+\infty}\int_{ B_1(0)}|\nabla \hat{\La}_{j}|^2\ {\mathbf 1}_{\om_j(t)}\ dy^2=0\ ,
\ee
where $ {\mathbf 1}_{\om_j(t)}$ is the characteristic function of $$\om_j(t):={\lf\{x\in B_1(0)\ ;\ \r^\ast_j({\hat{\La}_{j}(x)})> c_2\, t\rg\}}\ , $$  and where 
\[
\r^\ast_{j}(\hat{p}):=\inf\lf\{\r_{j,\hat{p}}(\hat{q})\quad;\quad \hat{q}\in \hat{\La}_{j}(\p B_1(0)) \rg \}\ .
\]
For any $\hat{p},\hat{q}\in r_j^{-1}\,(V_2({\R}^4)-\vec{\La}_\infty(0))$ we are using the notation
\[
\r_{j,\hat{p}}(\hat{q}):=r_j^{-1}\ \r_{\vec{p}_j}(\vec{q}_j)\quad\mbox{and }\quad \vec{p}_j:=r_j\,\hat{p}+\vec{\La}_\infty(0)\quad\mbox{ and }\quad \vec{q}_j:=r_j\,\hat{q}+\vec{\La}_\infty(0)\ .
\]
Recall that in a neighbourhood of $\vec{\La}_j(0)=(\vec{\ep}_1,\vec{\ep}_2)$ there holds
\[
\nabla^Hv^1=O(\rho)\quad,\quad\nabla^H w^2=O(\rho)\quad\mbox{ and }\quad\nabla^H(w^1+v^2)=O(\rho)
\]
This implies
\be
\label{negl-1}
\lf\{
\begin{array}{l}
\ds|\nabla\hat{v}^{1}_{j}|=|\nabla^Hv^1\circ\vec{\La}_j\cdot\nabla\hat{\La}_j|=O(\rho\circ\vec{\La}_j)\ |\nabla\hat{\La}_j|\\[5mm]
\ds|\nabla\hat{w}^2_{j}|=|\nabla^Hw^2\circ\vec{\La}_j\cdot\nabla\hat{\La}_j|=O(\rho\circ\vec{\La}_j)\ |\nabla\hat{\La}_j|\\[5mm]
\ds|\nabla(\hat{w}_{j}^1+\hat{v}_{j}^2)|=|\nabla^H(w^1+v^2)\circ\vec{\La}_j\cdot\nabla\hat{\La}_j|=O(\rho\circ\vec{\La}_j)\ |\nabla\hat{\La}_j|
\end{array}
\rg.
\ee
Hence
\be
\label{first-coo}
\ds\int_{B_1(0)}|\nabla \hat{v}_{j}^1|^2+|\nabla \hat{w}_{j}^2|^2+|\nabla (\hat{w}_{j}^1+\hat{v}^2_j)|^2 \ dy^2=o(1)\ .
\ee
We consider  $\xi(z_1,z_3)$ to be an arbitrary smooth non negative cut-off function on ${\mathcal P}_\infty$ such that Supp$(\xi)\subset {\mathcal E}_\infty$ . We introduce the following Hamiltonian function
\be
\label{VI.83-d}
\begin{array}{l}
\ds h_j:=\chi(\r^\ast_j/s)\  \xi(v^3_{j},v_{j}^4)\ \lf[ v^3_j\,w_j^3+v_j^4\,w^4_j-\varphi_j\rg]\\[5mm]
\ds\quad=\chi(\r^\ast/r_j\,s)\ r_j^{-2}\  \xi(r_j^{-1}\,v^3,r_j^{-1}\,v^4)\ \lf[ v^3\,w^3+v^4\,w^4-v^2+w^1\rg]\ ,
\end{array}
\ee
where $\chi$ is a cut-off function equal to one on $[0,1]$ and supported in $[2,+\infty)$. We choose $s>c_2\, t$ large enough and fixed such that  $\chi(\r^\ast/s)$ is identically equal to 1 on ${\mathcal E}_\infty$. Because of 
(\ref{concl-2-reu}) there holds for any smooth function $f$
\be
\label{VI.83-e}
\lim_{j\rightarrow+\infty}\int_{ B_1(0)}|\nabla \hat{\La}_{j}|^2\ f(|\hat{\La}_j|)\ [|\chi'(\r^\ast/s)|+|\chi''(\r^\ast/s)|]\circ\hat{\La}_j\ dy^2=0\ .
\ee
For $j$ large enough, because of (\ref{cv-bord}),  $\hat{\La}_j(\p B_1(0))\cap \mbox{Supp}(\xi)=\emptyset$. We have
\[
\begin{array}{l}
\ds\vec{X}_{h_j}:=r_j^{-1}\,s^{-1}\ J_H \nabla^H\r_\ast\ \chi'(\r^\ast_j/s)\ \xi(v^3_{j},v_{j}^4)\ \lf[ v^3_j\,w_j^3+v_j^4\,w^4_j-\varphi_j\rg]\\[5mm]
\ds\quad+\chi(\r^\ast/r_j\,s)\ r_j^{-2}\  J_H\,\nabla^H\lf[\xi(r_j^{-1}\,v^3,r_j^{-1}\,v^4)\ \lf[ v^3\,w^3+v^4\,w^4-v^2+w^1\rg]\rg]+\frac{h_j}{2}\ (\vec{\ep}_2+\vec{w},-\vec{\ep}_1-\vec{v})\ .
\end{array}
\] 
which gives in particular
\be
\label{VI.83-ee}
\|\nabla^H\vec{X}_{h_j}\|_{L^\infty (B_{R\,r_j}(\vec{\ep}_1,\vec{\ep}_2))}\le C\, (R\,r_j)^{-2}\quad\mbox{ and }\quad \|\nabla^H(\nabla^H\vec{X}_{h_j})\|_{L^\infty (B_{R\,r_j}(\vec{\ep}_1,\vec{\ep}_2))}\le C\, (R\,r_j)^{-3}
\ee
Let
\[
\hat{X}_{h_j}(\hat{p}):=r_j\, \vec{X}_{h_j}(r_j\,\hat{p}+\vec{\La}_\infty(0))
\]
Recall
\[
\al_j:=-d\varphi_j+\vec{v}_j\cdot d\vec{w}_j-\vec{w}_j\cdot d\vec{v}_j
\]
and denote $$H_j:=\mbox{Ker}(\al_j)\quad\mbox{  and }\quad J_{H_j}:= (D_{r_j^{-1}})_\ast\circ J_H\circ  (D_{r_j^{-1}})^{-1}_\ast\ .$$ We have in particular $\hat{X}_{h_j}\circ\hat{\La}_j(y)=r_j\, \vec{X}_{h_j}(\vec{\La}_j)(r_j\,y)$ and with the previously introduced notation there holds
\[
\begin{array}{l}
\ds\hat{X}_{h_j}(\hat{p})= r_j\,J_H \nabla^Hh_j(r_j\,\hat{p}+\vec{\La}_\infty(0))+r_j\,\frac{h_j}{2}\ (\vec{\ep}_2+\vec{w}( r_j\,\hat{p}+\vec{\La}_\infty(0)),-\vec{\ep}_1-\vec{v}(r_j\,\hat{p}+\vec{\La}_\infty(0) ))\\[5mm]
\ds\quad= J_{H_j}\nabla^{H_j} (r_j\,\hat{p}+\vec{\La}_\infty(0))+r_j\, \frac{\hat{h}_j}{2}\ (\vec{\ep}_2+\vec{w}( \vec{p}),-\vec{\ep}_1-\vec{v}(\vec{p} ))
\end{array}
\]
where we are using the induced metric from ${\R}^8$ on $r_j^{-1}\, (V_2({\R}^4)-\vec{\La}_\infty(0)) $.  $\hat{h}_j:=h_j( r_j\,\hat{\La}_j+\vec{\La}_\infty(0) )$. Observe that
\[
\nabla^{H_j}\hat{X}_{h_j}= r_j^2\,(\nabla^H \vec{X}_{h_j})(r_j\,\hat{p}+\vec{\La}_\infty(0))\quad\mbox{ and }\quad\nabla^{H_j}(\nabla^{H_j}\hat{X}_{h_j})= r_j^3\,\nabla^H(\nabla^H \vec{X}_{h_j})(r_j\,\hat{p}+\vec{\La}_\infty(0))
\]
which gives thanks to (\ref{VI.83-ee})
\be
\label{VI.83-f}
\|\nabla^{H_j}\hat{X}_{h_j}\|_{L^\infty (B_{R}(0))}+\|\nabla^{H_j}(\nabla^{H_j}\hat{X}_{h_j})\|_{L^\infty (B_{R}(0))}\le C\,
\ee
Thanks to the entropy condition (\ref{ener-asym}) and the computations  (\ref{smoother}), the criticality of $\hat{\La}_j$ of $\hat{E}_j$ for the Hamiltonian deformations generated by $\hat{X}_{h_j}$  is then giving
\be
\label{VI.83-g}
\int_{B_1(0)}<d\hat{\La}_j,d(\hat{X}_{h_j}\circ\hat{\La}_j)>\ dy^2=o(1)\ .
\ee
Observe that
\be
\label{VI.83-h}
\begin{array}{l}
\ds\lf<d(\hat{v}_j,\hat{w}_j),d\lf(   \hat{h}_j \, (\vec{\ep}_2+r_j\,\hat{w}_j, -\vec{\ep}_1+r_j\,\hat{v}_j) \rg)\rg>\\[5mm]
\ds\quad= r_j\,\hat{h}_j \ \lf<d(\hat{v}_j,\hat{w}_j),d(\hat{w}_j,-\hat{v}_j)\rg>- (d\hat{v}_j^2-d\hat{w}_j^1+r_j\,\hat{w}_j\cdot d\hat{v}_j-r_j\,\hat{v}_j\cdot d\hat{w}_j)\cdot d \hat{h}_j
\end{array}
\ee
Recall $\varphi_j=r_j^{-2}\,\varphi=r_j^{-2}\, (v^2-w^1)=r_j^{-1}\,(v^2_j-w^1_j)$. Since $\hat{\La}_j^\ast  \al_j=0$  there holds
\be
\label{VI.83-i}
\begin{array}{l}
r_j^{-1}\,d(\hat{v}^2_j-\hat{w}^1_j)=-\hat{w}_j\cdot d\hat{v}_j+\hat{v}_j\cdot d\hat{w}_j\ .
\end{array}
\ee
Combining (\ref{VI.83-e}), (\ref{VI.83-g}), (\ref{VI.83-h}) and (\ref{VI.83-i}) we deduce
\be
\label{VI.83-j}
\int_{B_1(0)}\chi\lf(\frac{\r^\ast}{r_j\,s}\rg)\ \frac{1}{r_j}\ \,\lf<d\hat{\La}_j,d \lf( J_H\,\nabla^H\lf[\xi(r_j^{-1}\,v^3,r_j^{-1}\,v^4)\ \lf[ v^3\,w^3+v^4\,w^4-\varphi\rg]\rg]\circ\vec{\La}_j\rg)\rg>\ dy^2=o(1)\ .
\ee
We have
\be
\label{VI.84-1}
\begin{array}{l}
\ds\int_{B_1(0)}\chi\lf(\frac{\r^\ast}{r_j\,s}\rg)\ \frac{1}{r_j}\ \,\lf<d\hat{\La}_j,d \lf( J_H\,\nabla^H\lf[\xi(r_j^{-1}\,v^3,r_j^{-1}\,v^4)\rg]\ \lf[ v^3\,w^3+v^4\,w^4-\varphi\rg]\circ\vec{\La}_j\rg)\rg>\ dy^2\\[5mm]
\ds=\int_{B_1(0)}\chi\lf(\frac{\r^\ast}{r_j\,s}\rg)\,\lf<d\hat{\La}_j,d \lf( \lf[\nabla^H w^3(\vec{\La}_j)\,\p_{v_j^3}\xi(\hat{v})+\nabla^H w^4(\vec{\La}_j)\,\p_{v_j^4}\xi(\hat{v})\rg]\, \lf[ \hat{v}_j^3\,\hat{w}_j^3+\hat{v}_j^4\,\hat{w}_j^4-\hat{\varphi}_j\rg]\rg)\rg>\ dy^2
\end{array}
\ee
From (\ref{conv-meas}), since $\hat{\La}_{\infty,\infty}(B_1)\subset \mbox{Span}\{(\vec{\ep}_3,0),(\vec{\ep}_4,0)\}$, we deduce for $q>0$
\be
\label{VI.83-l}
\int_{B_1(0)}\chi\lf(\frac{\r^\ast}{r_j\,s}\rg)\ \lf[|\hat{v}_j^1|^q+|\hat{v}_j^2|^q+|\hat{w}_j^1|^q+|\hat{w}_j^2|^q+|\hat{w}_j^3|^q+|\hat{w}_j^4|^q+|\hat{\varphi}_j|^q\rg]\ |d\hat{\La}_j|^2\ dy^2=o(1)\ .
\ee
Recall also from (\ref{nabla-H-bis}) and (\ref{nabla-H-ter}) that $\nabla^H w^i\circ(\hat{\La}_j)$ are polynomials of $\hat{v}^l_j$ and $\hat{w}^l_j$ (of grade at most 3). Hence we deduce 
\be
\label{VI.84-2}
\begin{array}{l}
\ds\int_{B_1(0)}\chi\lf(\frac{\r^\ast}{r_j\,s}\rg)\ \frac{1}{r_j}\ \,\lf<d\hat{\La}_j,d \lf( J_H\,\nabla^H\lf[\xi(r_j^{-1}\,v^3,r_j^{-1}\,v^4)\rg]\ \lf[ v^3\,w^3+v^4\,w^4-\varphi\rg]\circ\vec{\La}_j\rg)\rg>\ dy^2\\[5mm]
\ds=\int_{B_1(0)}\chi\lf(\frac{\r^\ast}{r_j\,s}\rg)\,\lf<d\hat{\La}_j,\lf[\nabla^H w^3(\vec{\La}_j)\,\p_{v_j^3}\xi(\hat{v})+\nabla^H w^4(\vec{\La}_j)\,\p_{v_j^4}\xi(\hat{v})\rg]\, \lf[ \hat{v}_j^3\,d\hat{w}_j^3+\hat{v}_j^4\,d\hat{w}_j^4-d\hat{\varphi}_j\rg]\rg>\ dy^2\\[5mm]
\ds\quad+\ o(1)\ .
\end{array}
\ee
We  write 
\[
\hat{v}_j^3\,d\hat{w}_j^3+\hat{v}_j^4\,d\hat{w}_j^4-d\hat{\varphi}_j=\hat{w}_j^3\,d\hat{v}_j^3+\hat{w}_j^4\,d\hat{v}_j^4-\hat{v}_j^1\,d\hat{w}_j^1+\hat{w}_j^1\,d\hat{v}_j^1-\hat{v}_j^2\,d\hat{w}_j^2+\hat{w}_j^2\,d\hat{v}_j^2
\]
Substituting this expression in (\ref{VI.84-2}) and using (\ref{VI.83-l}) we obtain finally
\be
\label{VI.84-2-bis}
\begin{array}{l}
\ds\int_{B_1(0)}\chi\lf(\frac{\r^\ast}{r_j\,s}\rg)\ \frac{1}{r_j}\ \,\lf<d\hat{\La}_j,d \lf( J_H\,\nabla^H\lf[\xi(r_j^{-1}\,v^3,r_j^{-1}\,v^4)\rg]\ \lf[ v^3\,w^3+v^4\,w^4-\varphi\rg]\circ\vec{\La}_j\rg)\rg>\ dy^2\\[5mm]
\ds=o(1)\ .
\end{array}
\ee
Combining (\ref{VI.83-j})  and (\ref{VI.84-2-bis}) is giving
\be
\label{VI.84-3}
\int_{B_1(0)}\chi\lf(\frac{\r^\ast}{r_j\,s}\rg)\ \frac{1}{r_j}\ \,\lf<d\hat{\La}_j,d \lf( \lf[\xi(r_j^{-1}\,v^3,r_j^{-1}\,v^4)\ J_H\,\nabla^H\lf[ v^3\,w^3+v^4\,w^4-\varphi\rg]\rg]\circ\vec{\La}_j\rg)\rg>\ dy^2=o(1)\ .
\ee
Recall from (\ref{nabla-H-bis}) and (\ref{nabla-H-ter})
\be
\label{VI.84-4}
\begin{array}{l}
\ds\nabla^H\varphi-\sum_{i=3}^4 v^i\,\nabla^Hw^i-w^i\,\nabla^H v^i=v^1\,\nabla^H w^1-w^1\,\nabla^H v^1+v^2\,\nabla w^2-w^2\,\nabla^H v^2\\[5mm]
\ds\quad=(w^1\,\vec{v}+(w^1)^2\ (\vec{\ep}_2+ \vec{w})+w^1\,v^1\,(\vec{\ep}_1+\vec{v}), -v^1\,\vec{v}-v^1\,w^1\ (\vec{\ep}_2+\vec{w})-(v^1)^2\,(\vec{\ep}_1+\vec{v}))\\[5mm]
\ds\quad+(w^2\,\vec{w}+w^2\,v^2\,(\vec{\ep}_1+\vec{v})+(w^2)^2\,(\vec{\ep}_2+\vec{w}),-v^2\,\vec{w}-(v^2)^2\,(\vec{\ep}_1+\vec{v})-v^2\,w^2\,(\vec{\ep}_2+\vec{w}))
\end{array}
\ee
Hence
\be
\label{VI.84-5}
\begin{array}{l}
\ds   \nabla^H\lf[ v^3\,w^3+v^4\,w^4-\varphi\rg]-2\,w_3\,\nabla^H v^3-2\,w_4\,\nabla^H v^4\\[5mm]
\ds\quad=(w^1\,\vec{v}+(w^1)^2\ (\vec{\ep}_2+ \vec{w})+w^1\,v^1\,(\vec{\ep}_1+\vec{v}), -v^1\,\vec{v}-v^1\,w^1\ (\vec{\ep}_2+\vec{w})-(v^1)^2\,(\vec{\ep}_1+\vec{v}))\\[5mm]
\ds\quad+(w^2\,\vec{w}+w^2\,v^2\,(\vec{\ep}_1+\vec{v})+(w^2)^2\,(\vec{\ep}_2+\vec{w}),-v^2\,\vec{w}-(v^2)^2\,(\vec{\ep}_1+\vec{v})-v^2\,w^2\,(\vec{\ep}_2+\vec{w}))
\end{array}
\ee
Hence
\be
\label{VI.84-6}
\int_{B_1(0)}\chi\lf(\frac{\r^\ast}{r_j\,s}\rg)\ \frac{1}{r_j}\ \,\lf<d\hat{\La}_j,d \lf( \lf[\xi(r_j^{-1}\,v^3,r_j^{-1}\,v^4)\ \lf[ w^3\,\nabla^Hw^3+w^4\,\nabla^Hw^4\rg]\rg]\circ\vec{\La}_j\rg)\rg>\ dy^2=o(1)\ .
\ee
Using one more time  (\ref{VI.83-l}) we finally obtain
\be
\label{VI.84-7}
\int_{B_1(0)}\chi\lf(\frac{\r^\ast}{r_j\,s}\rg)\ \,\lf[|d\hat{w}^3_j|^2+|d\hat{w}^4_j|^2 \rg]\ \xi(\hat{v}_j^3,\hat{v}_j^4)\ \ dy^2=o(1)\ .
\ee
Combining (\ref{first-coo}) and (\ref{VI.84-7}) is giving finally (\ref{rec-2}) which itself is implying (\ref{rec-1}).

Any infinitesimal variation $\phi_t$ such that supp$((\phi_t)_\ast |\hat{\mathbf v}_\infty|)\subset {\mathcal E}_\infty$ is preserving the legendrian property of the varifold hence
$\hat{\mathbf v}_\infty$ is stationary in a classical sense within the Legendrian two plane spanned by $\p_{x_1}\vec{\La}(0)$ and $\p_{x_2}\vec{\La}(0)$ relative to its boundary $\partial{\mathcal E}_\infty$. We deduce using the constancy theorem ( see \cite{Sim} chapter 8 theorem 41.1) that there exists a constant $\theta_0>0$ such that
\[
|\hat{\mathbf v}_\infty|=\theta_0\ d{\mathcal H}^2\res {\mathcal E}\ .
\]
It remains to prove that $m$ is a positive integer. 

\medskip

For any $0<\beta<1$ we denote 
\[
{\mathcal E}_\beta:=\lf\{ \vec{p}:=\nabla\vec{\La}_\infty(0)\cdot y\quad; \quad |y|<\beta\rg\}
\]
Let $a(v^3_\infty,v^4_\infty)\in C^\infty_0({\mathcal E}_\beta(0))$ and $b(v^3_\infty,v^4_\infty)\in C^\infty_0({\mathcal E}_\beta)$ arbitrary. Hence, because of (\ref{cv-bord}), for $j$ large enough, there holds
\be
\label{compatib}
a\circ \hat{\La}_j\equiv 0\ \quad\mbox{ and }\quad b\circ\hat{\La}_j\equiv 0\quad\quad\mbox{ on }\p B_1(0)\quad.
\ee
We introduce now the following Hamiltonian 
\[
h_j(\vec{v}_j,\vec{w}_j):= \chi\lf(\frac{\r^\ast}{r_j\,s}\rg)\, \lf[a(v^3_j,v^4_j)\ w^3_j+b(v^3_j,v^4_j)\ w^4_j\rg]\ .
\]
Arguing exactly as for the previous Hamiltonian gives then
\be
\label{divzero}
\lim_{j\rightarrow 0}\int_{B_1(0)}\chi\lf(\frac{\r^\ast}{r_j\,s}\rg)\, \nabla\lf[ a(\hat{v}^3_j,\hat{v}^4_j) \,\vec{\ep}_3\rg]\,\nabla\ti{\La}_j+\nabla\lf[ b(\hat{v}^3_j,\hat{v}^4_j)\, \vec{\ep}_4\rg]\,\nabla\ti{\La}_j\ dx^2=0\ \\[5mm]
\ee
where we have introduced the notation $\ti{\La}_j:=(\hat{v}^3_j,\hat{v}^4_j)$ for the projection on the two Legendrian plane ${\mathcal L}$ generated by $(0,\vec{\ep}_3,0)$ and $(0,\vec{\ep}_4,0)$. 
Recall that $\hat{\La}_j$ is conformal on $B_1(0)$. Denote
\[
\hat{\La}_j^\ast g_{{\R}^8}=e^{2\hat{\la}_j}\ [dx_1^2+dx_2^2]\ .
\]
Because of (\ref{rec-2}) we have in particular for $l=3,4$
\[
\lim_{j\rightarrow +\infty}\int_{B_1(0)}|e^{2\hat{\la}_j}\,(\vec{\ep}_l,0)-(\vec{\ep}_l,0)\cdot\nabla\hat{\La}_j\ \nabla\hat{\La}_j|\ dx^2=0\ ,
\]
which implies for $l=3,4$
\be
\label{sau-1}
\lim_{j\rightarrow +\infty}\int_{B_1(0)}|e^{2\hat{\la}_j}-|\nabla \hat{v}_j^l|^2|\ dx^2=0\ .
\ee
Observe that we have also
\[
\lim_{j\rightarrow +\infty}\int_{B_1(0)}\lf|e^{\hat{\la}_j}\,(\vec{\ep}_l,0)-(\vec{\ep}_l,0)\cdot\nabla\hat{\La}_j\ \nabla\hat{\La}_j\ e^{-\hat{\la}_j}\rg|^2\ dx^2=0\ .
\]
Using the fact that $(\vec{\ep}_3,0)\cdot(\vec{\ep}_4,0)=0$ we deduce from the previous line
\be
\label{sau-2}
\lim_{j\rightarrow +\infty}\int_{B_1(0)}|\nabla\hat{v}_j^4\cdot\nabla\hat{v}_j^3|\ dx^2=0\ .
\ee
We have
\[
\begin{array}{l}
\ds\p_{x_1}\lf[ a(\hat{v}^3_j,\hat{v}^4_j) \rg]\,\p_{x_1}\hat{v}_j^3+\p_{x_1}\lf[ b(\hat{v}^3_j,\hat{v}^4_j) \rg]\,\p_{x_1}\hat{v}_j^4+\p_{x_2}\lf[ a(\hat{v}^3_j,\hat{v}^4_j) \rg]\,\p_{x_2}\hat{v}_j^3+\p_{x_2}\lf[ b(\hat{v}^3_j,\hat{v}^4_j) \rg]\,\p_{x_2}\hat{v}_j^4\\[5mm]
=\p_{\hat{v}_j^3}a\, (\p_{x_1}\hat{v}_j^3)^2+\p_{\hat{v}_j^4}a\, \p_{x_1}\hat{v}_j^4\,\p_{x_1}\hat{v}_j^3
+\p_{\hat{v}_j^3}a\, (\p_{x_2}\hat{v}_j^3)^2+\p_{\hat{v}_j^4}a\, \p_{x_2}\hat{v}_j^4\,\p_{x_2}\hat{v}_j^3\\[5mm]
\ds+\p_{\hat{v}_j^3}b\, \p_{x_1}\hat{v}_j^3\,\p_{x_1}\hat{v}_j^4+\p_{\hat{v}_j^4}b\, (\p_{x_1}\hat{v}_j^4)^2
+\p_{\hat{v}_j^3}b\, \p_{x_2}\hat{v}_j^3\,\p_{x_2}\hat{v}_j^4+\p_{\hat{v}_j^4}b\, (\p_{x_2}\hat{v}_j^4)^2\\[5mm]
\ds=e^{2\hat{\la}_j}\,\lf[\p_{\hat{v}_j^3}a\, +\p_{\hat{v}_j^4}b\rg]+\p_{\hat{v}_j^3}a\ [|\nabla \hat{v}_j^3|^2-e^{2\hat{\la}_j}]+\p_{\hat{v}_j^4}b\ [|\nabla \hat{v}_j^4|^2-e^{2\hat{\la}_j}]+[\p_{\hat{v}_j^4}a+\p_{\hat{v}_j^3}b]\, \nabla\hat{v}_j^4\cdot\nabla\hat{v}_j^3\ .
\end{array}
\]
Combining this identity with (\ref{divzero}), (\ref{sau-1}) and (\ref{sau-2}) gives then
\[
\lim_{j\rightarrow +\infty}\int_{B_1(0)} \lf[\p_{\hat{v}_j^3}a\, +\p_{\hat{v}_j^4}b\rg]\circ\hat{\La}_j\ |\p_{x_1}\hat{\La}_j\wedge\p_{x_1}\hat{\La}_j|\ dx_1\wedge dx_2=0\ .
\]
The area formula gives then
\be
\label{divz}
\lim_{j\rightarrow +\infty}\int_{\mathcal L}\ \lf[\p_{y_1}a\, +\p_{y_2}b\rg](y)\ N_j(y)\ dy^2=0
\ee
where $N_j(y)$ is the number of $(x_1,x_2)$ in $B_1(0)$ such that $(\hat{v}_3(x_1,x_2),\hat{v}_4(x_1,x_2))=(y_1,y_2)$\ . At this stage we have used a uniform $L^\infty$ control of $(\nabla_y^la,\nabla_y^lb)$ for $l\le 3$
and we have proved
\be
\label{divz-a}
\lim_{j\rightarrow +\infty}\sup_{\sum_{l=0}^3\|\nabla_y^l(a,b)\|_\infty\le 1}\int_{\mathcal L}\ \lf[\p_{y_1}a\, +\p_{y_2}b\rg](y)\ N_j(y)\ dy^2=0
\ee
Let's take now a small parameter $\tau_j\rightarrow 0$ and consider the convolution $(a_j,b_j)$ of $a$ and $b$ respectively with $\chi_j(y_1,y_2):=\tau_j^{-1}\chi((y_1/\tau_j,y_2/\tau_j))$ where $\chi\in C^\infty_0({\R}^2)$ with $\int_{{\R}^2} \chi(y)\ dy^2=1$. We have obviously
\[
\|\nabla_y[(a,b)-(a_j,b_j)]\|_{L^\infty({\R}^2)}=o(1)\ \|\nabla_y(a,b)\|_{L^\infty({\R}^2)}
\]
which implies in particular
\[
\lim_{j\rightarrow +\infty}\int_{\mathcal L}\ \lf[\p_{y_1}(a-a_j)\, +\p_{y_2}(b-b_j)\rg](y)\ N_j(y)\ dy^2=o(1)\  \|\nabla_y(a,b)\|_{L^\infty({\R}^2)}\ .
\]
We have also for $l\le 3$
\[
\|\nabla^l_y(a_j,b_j)\|_{L^\infty({\R}^2)}\le \tau_j^{-l+1}\  \|\nabla_y(a,b)\|_{L^\infty({\R}^2)}
\]
We choose $\tau_j$ in such a way that $\tau_j^2$ goes infinitesimaly slower to zero than all the terms in the estimates above which go to zero. In particular we require
\[
\tau_j^{-2}\hat{\ep}^{\,4}_j\,\int_{B_1(0)}(r_j^2+|d\hat{T}_j|^2_{g_{\hat{\La}_j}})^2\ dvol_{g_{\hat{\La}}}=o(1)\  \|\nabla_y(a,b)\|_{L^\infty({\R}^2)}\ .
\]
In this way we deduce foor any $\beta<1$
\be
\label{divz-b}
\lim_{j\rightarrow +\infty}\sup_{\|\nabla_y(a,b)\|_\infty\le 1\ ;\ \mbox{Supp}(a,b)\subset {\mathcal E}_\beta}\int_{\mathcal L}\ \lf[\p_{y_1}a\, +\p_{y_2}b\rg](y)\ N_j(y)\ dy^2=0\ .
\ee
Using Lemma A.7 of \cite{PiRi1} we obtain
\be
\label{weakl1}
\lim_{j\rightarrow +\infty}\lf| N_j- |{\mathcal E}_{1/2}|^{-1}\int_{{\mathcal E}_{1/2}} N_j(y)\ dy^2\rg|_{L^{1,\infty}({\mathcal E}_{1/2})}=0
\ee
The area formula gives also
\[
N_j(y)\ dy^2\rightharpoonup |\hat{\mathbf v}_\infty|=\theta_0 \ d{\mathcal H}^2\res{\mathcal E}\ \mbox{ in }(C^0({\mathcal E}_{1/2}))^\ast
\]
This implies in particular
\[
\lim_{j\rightarrow +\infty}|{\mathcal E}_{1/2}|^{-1}\int_{{\mathcal E}_{1/2}} N_j(y)\ dy^2=\theta_0
\]
and from (\ref{weakl1}), we obtain
\[
\lim_{j\rightarrow +\infty}\mbox{dist}\lf(|{\mathcal E}_{1/2}|^{-1}\int_{{\mathcal E}_{1/2}} N_j(y)\ dy^2,{\N}\rg)=0
\]
hence
\[
\theta_0\in {\N}^\ast\ .
\]
Because of (\ref{item-1}) we have
\[
\pi\,\theta_0\ |\p_{x_1}\vec{\La}_\infty(0)\wedge\p_{x_2}\vec{\La}_\infty(0)|= \theta_0\,|{\mathcal E}|=\lim_{j\rightarrow +\infty}\hat{\nu}_j(B_1(0))=\lim_{r\rightarrow 0}\frac{\nu_\infty(B_r(x))}{r^2}=\pi\, f(0)\, e^{2\la(0)}
\]
 This implies (\ref{but}). Since $\mu_k=(\vec{\La}_k)_\ast \nu_k\rightharpoonup (\vec{\La}_\infty)_\ast \nu_\infty$, we have that $\mu_\infty=(\vec{\La}_\infty)_\ast \nu_\infty$. Thus for $|d\vec{\La}_\infty|^2_h\ dvol_h$ almost every $x\in \Sigma$
\[
\pi\,N_x\le \limsup_{t\rightarrow 0}\frac{\mu_\infty(B^\r_t(\vec{\La}_\infty(x)))}{t^2}\le C\,\mu_\infty(V_2({\R}^4))
\] 
Hence $N$ is in $L^\infty(\Sigma)$. This concludes the proof of lemma~\ref{lm-integer-mult}. \hfill $\Box$

We now identify the limiting varifold ${\mathbf v}_\infty$ as being the varifold associated to $\vec{\La}_\infty$ with the weight $N$. Precisely we have.
\begin{Lm}
\label{lm-vari-lim}
We denote by ${\mathbf v}_{\infty}^\eta$ the limit of the varifold induced by $\vec{\La}_k$ on $\Sigma\setminus \cup_{i=1}^QB_\eta(x_i)$ where $\eta>0$ is arbitrary
\[
\forall \ \Xi\in C^0(G_2^{Leg}(V_2({\R}^4)))\quad {\mathbf v}^\eta_{\infty}(\Xi):=\int_{\Sigma\setminus \cup_{i=1}^QB_\eta(x_i)} N\ \Xi((\vec{\La}_\infty)_\ast T_x\Sigma,\vec{\La}_\infty(x))\ dvol_{g_\infty}\ .
\]
\end{Lm}
\noindent{\bf Proof of lemma~\ref{lm-vari-lim}} For almost every $x_0\in \Sigma$ one has
\be
\label{lebesgue-pt-rep}
\lim_{r\rightarrow 0}\int_{B_r(x_0)}\frac{|d\vec{\La}_\infty(x)-d\vec{\La}_\infty(0)|^2}{r^2}\ dx^2=0\ .
\ee
as well as
\be
\label{app-diff-rep}
\lim_{r\rightarrow 0}\int_{B_r(x_0)}\frac{\lf|\vec{\La}_\infty(x)-\vec{\La}_\infty(0)-d\vec{\La}_\infty(0)\cdot x\rg|}{r^3}\ dx^2=0
\ee
We take such a point $x_0$ and we consider conformal coordinates for $h$ such that $x_0=(0,0)$. As in the proof of the previous lemma, we have a sequence of radii such that
\be
\label{bord-cont-rep}
\lf\|\vec{\La}_\infty(r_j,\theta)-\vec{\La}_\infty(x_0)-r_j\, \cos\theta\,\p_{x_1}\vec{\La}_\infty(0)-r_j, \sin\theta\,\p_{x_2}\vec{\La}_\infty(0)\rg\|_{L^\infty([0,2\pi])}=o(r_j)\ .
\ee
Let $\delta>0$. We claim first that for any $\Xi\in C^0(G_2^{Leg}(V_2({\R}^4)))$, for $j$ large enough
\[
\limsup_{k\rightarrow+\infty}\int_{B_{r_j}(0)}\lf|\Xi((\vec{\La}_k)_\ast T_x\Sigma,\vec{\La}_k(x))-\Xi((\vec{\La}_\infty)_\ast T_x\Sigma,\vec{\La}_\infty(x))\rg|\ dvol_{g_k}\le \delta\ r_j^2
\]
This claim is proved as in the proof of (5.10) in \cite{Piga} and the lemma follows for the same reasons as the one given in the proof of proposition 5.7 of \cite{Piga}.\hfill $\Box$
\section{Proof of theorem~\ref{th-1} for $N^5=V_2({\R}^4)$.}
Let ${\mathcal A}$ be an admissible family in ${\mathfrak M}:={\mathfrak E}_{\Sigma,Leg}^{2,4}$ such that
\[
\beta:=\inf_{A\in {\mathcal A}}\ \sup_{\vec{\La}\in A}\int_\Sigma dvol_{g_{\vec{\La}}}>0
\]
We apply proposition~\ref{pr-stmono} for $f(t):=e^{-1/t^2}$ and we obtain a sequence $\ep_k\rightarrow 0$ as well as $\vec{\La}_k\in {\mathfrak E}_{\Sigma,Leg}^{2,4}$ satisfying
\be
\label{energybound-bis}
\limsup_{k\rightarrow +\infty}E_{\ep_k}(\vec{\La}_k)=\beta
\ee
moreover
\be
\label{alm-crit-bis}
\|dE_{\ep_k}(\vec{\La}_k)\|_{\vec{\La}_k}\le e^{-\ep_k^{-2}}\ ,
\ee
and satisfying the entropy condition
\be
\label{entropy-bis}
\ep_k^4\,\int_\Sigma (1+|d\vec{T}_k|^2_{g_k})^2\ dvol_{g_k}=o\lf(\frac{1}{\log \ep_k^{-1}}\rg)\ .
\ee
The section III and IV apply. At the concentration points $x_i$ we proceed to a blow-up analysis exactly as in the proof of lemma III.6 in \cite{Riv-1}, using (\ref{concl-3}) of lemma~\ref{lm-energy-quant} in the same way as
the {\it global energy quantization lemma} III.3 is used in \cite{Riv-1}, to prove that there is no area dissipation in the neck regions (see also 6.3 in \cite{Piga}). Finally the case when the Riemann structure induced by
 $\vec{\La}_k$ is degenerating is treated as in section IV of \cite{Riv-1} or section 6 in \cite{Piga}. This is ending the proof of theorem~\ref{th-1}.\hfill $\Box$
 
 \section{The proof of theorem~\ref{th-1} in the general case.}
 \reset
 We recall the definition
 \begin{Dfi}
 \label{df-sasakian}
 A compact Riemannian manifold $(N^5,g)$ is {\bf Sasakian} if its metric cone $(C(S):={\R}_+^\ast\times N^5,\ov{g}=dt^2+t^2\,g)$ is K\"ahler.
 \end{Dfi}
 As a consequence of the definition $C(S)$ is equipped with an integrable complex structure ${J}$ and a K\"ahler $2-$form. On $\{1\}\times N^5\simeq N^5$ we define the Reeb vector-field
 \[
 \vec{R}:=\sqrt{2}\,{J}(\p_t)\ ,
 \]
 as well as the $1-$form
 \[
 \al:= \sqrt{2}\, J\, dt\ ,
 \]
 so that $\al(\vec{R})=2$. In $N^5$ we define $H:=\mbox{ker}(\al)$ and $J^H$ to be the restriction to $H$ of $J$ (bearing in mind that $J$ realizes an isometry of $H$). Observe that, since $\ov{g}$ is compatible with $J$
 the explicit expression of $\ov{g}$ gives
 \[
 \forall \vec{X}\in H\quad \quad 0=\ov{g}(\p_t,J\vec{X})=-g(J\p_t,\vec{X})
 \]
Hence $H\mbox{ker}(\al)$ coincides with the orthogonal to the Reeb vector-field $\vec{R}$. The symplectic form $\Om$  on $(C(S):={\R}_+^\ast\times N^5,\ov{g}=dt^2+t^2\,g)$ is given by
 \[
 \Om :=dt\wedge J\,dt+\frac{t^2}{2\,\sqrt{2}}\,d\al=\frac{i}{2}\,\ov{\p}\p t^2\ .
 \]
 Since $\Om$ is non degenerate we deduce that $\al$ is contact on $N^5$ that is
 \[
 \al\wedge d\al\wedge d\al\ne 0\ .
 \]
 A 5-Sasakian manifold is then often seen as the collection of $(N^5, g,\al,\vec{R},J^H)$ satisfying the above conditions.   In \cite{GKN} the existence of a ``Sasakian potential'' associated to every Sasakian structure  has been discovered (see also a presentation in \cite{spar} section 1.2).  At the neighbourhood of any point $\vec{p}$ there exists on $N^5$ a local chart $(\psi,z_1=x_1+iy_1,z_2=x_2+iy_2 )$ 
 where $\psi(\vec{p})=0$, $z_j(\vec{p})=0$ and a real valued function $h$ depending only on $(x_1,y_1,x_2,y_2)$ such that
 \be
 \label{VII.1}
 \lf\{
 \begin{array}{l}
 \ds \vec{R}=\sqrt{2}\,\p_{\psi}\\[3mm]
 \ds H^{0,1}=(H\otimes {\C})^{(1,0)}\quad\mbox{ is generated by } \p_{z_j}-\sqrt{2}^{-1}\,\al(\p_{z_j})\ \p_\psi\\[3mm]
 \ds J^H(\p_{z_j}-\sqrt{2}^{-1}\,\al(\p_{z_j})\ \p_\psi)= i\, \lf(\p_{z_j}-\sqrt{2}^{-1}\,\al(\p_{z_j})\ \p_\psi\rg)\\[5mm]
 \ds \al= \sqrt{2}\, d\psi+i\, \sqrt{2}\,\lf[\sum_{j=1}^2\ \p_{z_j}h\ dz_j-\p_{\ov{z}_j}h\ d\ov{z}_j\rg]
 \end{array}
 \rg.
 \ee
 Moreover $h$ is a K\"ahler potential for the {\it transversally K\"ahler foliation} $\psi=cte$ and we can also assume\footnote{This change of coordinate is called  {\it special foliated coordinates} in \cite{VeZe}} (modulo an holomorphic change of coordinates $z\rightarrow w$ on the leafs)  that 
 \be
 \label{VII.2}
 \frac{\p^2 h}{\p{z_j}\p{\ov{z}_k}}(0)=\delta_{jk}
 \ee
 so that the transverse K\"ahler metric on the leafs $\psi=cte$ satisfies 
 \[
 g^T(X,Y)=\frac{1}{2\,\sqrt{2}} d\al(X,J Y)=g^T_{j\ov{k}}\ dz_j\otimes d\ov{z}_k\quad\mbox{ and }\quad g^T_{j\ov{k}}(0)=\delta_{jk}
 \]
 Now we introduce 
 \[
 \varphi:=\psi +i\, \sum_{j=1}^2\p_{z_j}h(0)\, z_j-i\, \sum_{j=1}^2\p_{\ov{z}_j}h(0)\, \ov{z}_j+ \frac{i}{2}\, \sum_{j,k=1}^2\p^2_{z_jz_k}h(0)\ z_j\,z_k- \frac{i}{2}\, \sum_{j,k=1}^2\p^2_{\ov{z}_j\ov{z}_k}h(0)\ \ov{z}_j\,\ov{z}_k\ .
 \]
 We have
 \be
 \label{VII.3}
 \begin{array}{l}
\ds\sqrt{2}^{-1}\, \al=  d\varphi +i\,\sum_{j=1}^2\ \lf[\p_{z_j}h-\p_{{z}_j}h(0)-\sum_{k=1}^2\ \frac{\p^2 h}{\p{z_j}\p{{z}_k}}(0){z}_k\rg] \ dz_j\\[5mm]
\ds- i\ \sum_{j=1}^2\ \lf[\p_{\ov{z}_j}h-\p_{{\ov{z}}_j}h(0)-\sum_{k=1}^2\p^2_{\ov{z}_j\ov{z}_k}h(0)\ \ov{z}_k\rg] \   d\ov{z}_j\\[5mm]
\ds =d\varphi+i\,\sum_{j=1}^2\ (\ov{z}_j+O(|z|^2))\ dz_j-i\,\sum_{j=1}^2\ ({z}_j+O(|z|^2))\ d\ov{z}_j
\end{array}
 \ee
 Hence by dilating at ${\vec{p}}$ we are asymptotically converging to the Heisenberg group ${\mathbb H}^2$ as in the special case $V_2({\R}^4)$ above and the proof of theorem~\ref{th-1} given for $V_2({\R}^4)$ in sections II...V can  be adapted word by word with the same error terms (which were not taking into account the particular structure of $V_2({\R}^4)$ but only the convergence to ${\mathbb H}^2$ in the blow-up).

\bigskip

\noindent{\bf Conflict of interest} The author declares no conflict of interest.

\end{document}